\documentclass[12pt]{amsart}

\usepackage{amsmath,amssymb}

\input{melanie.sty}

\usepackage{fullpage,pdfsync}

\usepackage[capitalise]{cleveref}
\usepackage[]{mdframed}

\newtheorem{definition}[equation]{Definition}

\newcommand{\details}[1]{}

\newcommand{\category}{C}
\newcommand{\cL}{\mathcal{L}}
\newcommand{\cW}{\mathcal{W}}
\newcommand{\W}{\mathcal{W}}
\newcommand{\Img}{\operatorname{Im}}
\newcommand{\chr}{\operatorname{char}}
\newcommand{\Ex}{\mathscr{E}_\cL}
\newcommand{\Surj}{\operatorname{Surj}}
\newcommand{\Sp}{\operatorname{Sp}}
\newcommand{\ASp}{\operatorname{ASp}}
\newcommand{\Ps}{\operatorname{Ps}}
\newcommand{\Br}{\operatorname{Br}}
\newcommand{\Cat}{\operatorname{C}_{\Gamma,n}}
\newcommand{\sspan}{\operatorname{span}}
\newcommand{\AV}{\operatorname{AV}}
\newcommand{\av}{\operatorname{av}}
\newcommand{\Hur}{\operatorname{Hur}}

\newcommand{\bG}{\mathbf G}
 
\newcommand{\bbC}{\mathbb C} 
\newcommand{\bH}{\mathbf H}
\newcommand{\bX}{\mathbf X}

\newcommand{\rDisc}{\operatorname{rDisc}}
\newcommand{\Ext}{\operatorname{Ext}}
\newcommand{\oftau}{} 
\DeclareMathOperator{\bGal}{\mathbf{Gal}}
\newcommand{\Hf}{H_{\operatorname{fppf}}}
\newcommand{\Hfc}{H_{\operatorname{fppf},c}}
\newcommand{\tA}{\operatorname{tr}_{\operatorname{AV}}}
\newcommand{\Cl}{\mathrm{Cl}}

\renewcommand\C{\mathcal{C}}

\usepackage{mathrsfs}
\usepackage{physics}
\usepackage{tikz-cd}

\title{Distributions of unramified extensions of global fields}

\author{Will Sawin}
\address{Department of Mathematics\\
Princeton University \\
Fine Hall, Washington Road \\
Princeton, NJ 08540 USA}  
\email{wsawin@math.princeton.edu}

\author{Melanie Matchett Wood}
\address{Department of Mathematics\\
Harvard University\\
Science Center Room 325\\
1 Oxford Street\\
Cambridge, MA 02138 USA}  
\email{mmwood@math.harvard.edu}

\begin{document}

\begin{abstract}
Given a finite group $\Gamma$,  we prove results on the distribution of the prime-to-$q|\Gamma|$ part of fundamental groups of $\Gamma$-covers of the projective line $\P^1_{\F_q}$ over a finite field $\F_q$ as $q\ra\infty$.    Equivalently, this is a result on the distribution of the Galois groups of maximal unramified extensions of $\Gamma$-extensions of $\F_q(t)$, and thereby motivates a new conjecture on the distribution of Galois groups of maximal unramified extensions of $\Gamma$-extensions of a number field.   In particular, this allows us to see and predict the effect of roots of unity in the base field on such distributions.
We introduce the idea to study these groups along with the class in their 3rd homology group that arises from Artin-Verdier Duality. This invariant refines the lifting invariant that, in the function field setting, corresponds to stable components of Hurwitz space.
One major input into our function field results is an application of our recently developed methods to determine a distribution of groups (or more general algebraic structures) from its moments.   We prove non-existence results in the number field case that support our conjectures in the case where our conjectures predict certain kinds of groups occur with probability zero.
\end{abstract}

\maketitle


\section{Introduction}

This paper is devoted to understanding the distribution of the Galois group of the maximal unramified extension of a random number field, and specifically to understanding the influence of the roots of unity on the distribution. We make conjectural predictions, which are motivated by function field evidence in the $q\to\infty$ limit, and also checked against some numerical evidence.

It is useful to first consider the abelianization of the Galois group of the maximal unramified extension. By class field theory, this is the class group. Since work of Achter \cite{Achter2006} and Malle \cite{Malle2008}, it has been understood that the roots of unity of a field affect the distribution of the class groups of random extensions of a global field. It was not immediately clear how to give a formula for the probability of obtaining a given class group that correctly accounts for roots of unity, but a series of works have made predictions for this in special cases compatible with numerical and function-field evidence. Recently, the authors~\cite{Sawin2023} gave a prediction for the distribution of the Sylow $p$-subgroup of the class group of Galois extensions with Galois group $\Gamma$ of a fixed number field, incorporating roots of unity, as long as $p$ is prime to the order of $\Gamma$.

Predictions avoiding the influence of roots of unity for the distribution of the Galois group of the maximal unramified extension of a random extension of $\mathbb Q$ were recently made by Liu, Wood, and Zurieck-Brown~\cite{Liu2024} in the case when the extension is split at $\infty$ and Liu and Willyard~\cite{LiuWillyard} in general. Incorporating roots of unity requires introducing several new ideas, which we will discuss shortly.

We now give an example of a consequence of our main theorem in the function field case, which concerns the maximal $2$-group quotient of the Galois group of the maximal unramified extension, or, equivalently, the Galois group of the composition of all unramified Galois extensions whose degree is a power of $2$. 

 Let $q$ be a prime power, not divisible by $2$ or $3$, and let $\mathbb F_q$ be the finite field of $q$ elements. Let $\infty$ be the unique place of $\mathbb F_q(t)$ which does not arise from a prime ideal of $\mathbb F_q[t]$, i.e. the place corresponding to the valuation $v(f) = -\deg f$. For $K$ an extension of $\mathbb F_q(t)$, we let $\rDisc K$ be the product of the norms of the places of $\mathbb F_q(t)$ that ramify in $K$ (the norm of a place $v$ of $\mathbb F_q(t)$ being $q^{\deg v}$).  Let $E_{\mathbb Z/3} ( q^m, \mathbb F_q(t))$ be the set of Galois extensions of degree $3$ of $\mathbb F_q(t)$ (which necessarily have Galois group $\mathbb Z/3$), split at $\infty$, such that $\rDisc K=q^m$. 

For $K$ an extension of $\mathbb F_q(t)$, let  $K^{ \operatorname{un},2}$ be the composition of all finite Galois extensions of $K$ which are everywhere unramified, split at all places lying above $\infty$, and have degree a power of $2$.  Thus $\Gal(K^{ \operatorname{un},2}/K)$ is either a finite $2$-group or an infinite pro-$2$ group.  The following describes the probability of obtaining certain finite  $2$-groups, and is a sample of a kind of corollary one can obtain from our main result.

\begin{theorem}\label{intro-bb} Let $H$ be a finite $2$-group. Then there exists $p_H\in [0,1]$ such that
\[ \lim_{ \substack{ q\to\infty \\ q \equiv 3 \bmod 4 \\  3\nmid q }}  \limsup_{b\to\infty} \frac{\sum_{m\leq b}  \abs {\{ K \in E_{\mathbb Z/3}( q^m , \mathbb F_q(t))  \mid  \Gal( K^{ \operatorname{un},2} /K) \cong H  \} }} { \sum_{m\leq b}\abs{ E_{\mathbb Z/3}( q^m , \mathbb F_q(t))}} =p_H \] and

\[ \lim_{ \substack{ q\to\infty \\ q \equiv 3 \bmod 4 \\  3\nmid q }} \liminf_{b \to\infty} \frac{\sum_{m\leq b}  \abs {\{ K \in E_{\mathbb Z/3}( q^m , \mathbb F_q(t))  \mid  \Gal( K^{ \operatorname{un},2} /K) \cong H  \} }} { \sum_{m\leq b}\abs{ E_{\mathbb Z/3}( q^m , \mathbb F_q(t))}} =p_H \] and if $|H| \leq 8$ then
\[ p_H =  \prod_{k=0}^\infty (1+ 4^{-k-\frac{3}{2}})^{-1}  \cdot \begin{cases} 1 & \textrm{if } H \textrm{ is trivial} \\ \frac{1}{12} & \textrm{if } H \textrm{ is the Klein four group}  \\ \frac{7}{2^5 \cdot 3} & \textrm{if } H \textrm{ is the eight-element quaternion group}  \\ 0 & \textrm{ for all other } H \textrm{ with } |H| \leq 8 \end{cases} .\]
\end{theorem}
Based on this, we conjecture that the probability that for $K$ a random cyclic cubic extension of $\mathbb Q$ (necessarily split at $\infty$), the Galois group of the composition of all finite Galois extension of $K$ which are everywhere unramified (and thus split at all places lying above $\infty$), and have degree a power of $2$, is given the same formula $p_H$. 

For $H$ the trivial group, having the Galois group of $K^{  \operatorname{un},2}$ isomorphic to $H$ is equivalent to the $2$-rank of the class group of $K$ vanishing, and so in that case the conjecture agrees with an earlier conjecture of Malle~\cite[Equation (1) on p. 2827]{Malle2008}.

For $H$ the Klein four group or eight-element quaternion group, this agrees with an earlier conjecture of Boston and Bush~\cite[Table 7]{Boston2021a}. We discuss the exact relationship to their work in more detail later. The conjecture of Boston and Bush was based on matching numerical data (after guessing based on Malle's conjecture~\cite{Malle2010} that the probabilities are likely rational numbers with small denominators  times  $\frac{1}{6}\prod_{k=0}^\infty (1+ 4^{-k-\frac{3}{2}})^{-1}$) rather than any theoretical model. Hence our conjecture also matches Boston and Bush's data, as well as providing a conceptual explanation for their conjecture.

The condition $3\nmid q$ in Theorem \ref{intro-bb} is necessary to use comparison of fundamental groups between zero and positive characteristic.
The condition $q\equiv 3 \bmod 4$ ensures that $\mathbb F_q(t)$ has exactly $2$ roots of unity of order a power of $2$, the same number that $\mathbb Q$ has, which is the case that models the influence of the roots of unity of $\mathbb Q$ on the statistics. 
However, our main result in the function field case allows for arbitrary groups of roots of unity in the base field.
The restriction that $|H|\leq 8$ is only there to give some specific examples.  Our results give a method to enable the computation of $p_H$ for any $2$-group $H$.

An interesting consequence of Theorem \ref{intro-bb} is that for $H$ the trivial group, the Klein group, or the quaternion group, as long as $q$ is sufficiently large, coprime to $3$, and congruent to $3$ mod $4$, there exist infinitely many cyclic cubic extensions $K/\mathbb F_q(t)$ with $ \Gal(K^{ \operatorname{un},2} /K) \cong H$. This result appears to be new. Many other existence results for function fields with Galois groups of maximal unramified extensions of a specific form follow from our main results. 

The extensive prior work on the distribution of class groups and their non-abelian analogues, e.g. \cite{Cohen1984, Cohen1990, Malle2008, Malle2010,Garton2015, Adam2015, Boston2017a,Lipnowski2019, Lipnowski2020,Boston2021,Liu2022,Liu2024,LiuWillyard}, 
has usually constructed probability distributions either directly by guessing a formula for the probability of a given group or in terms of explicit constructions of random groups by generators and random relations, such as the cokernel of a random matrix. In the case of nonabelian groups, with the influence of roots of unity, it is not at all clear how to guess the probabilities or find a suitable construction of random relations to accurately model the Galois group. Liu~\cite{Liu2022} has given a Galois group model in one case, but it is not clear how to use this model to calculate the probability of obtaining a given group.

Instead, we use the moment method introduced by the authors in ~\cite{Sawin2022}. A classical problem in probability theory is to give criteria for when a distribution is uniquely determined by its moments. Subsequent work defined a notion of moments for probability distributions on the set of isomorphism classes of finite abelian groups, and gave criteria for when these probability distributions are determined by their moments \cite{Ellenberg2016, Wood2017}. 
More recently, similar results were obtained for nonabelian groups~\cite{Boston2017, Sawin2020}. The paper \cite{Sawin2022} generalizes those results and also strengthens them by giving criteria for a distribution with a given moments to exist and, more crucially, a formula for the probability that distribution assigns to a given group.

Thus, to construct a probability distribution that conjecturally models the distribution of the Galois group, we need only conjecture a formula for the moments of the probability distribution, and apply the results of \cite{Sawin2022} to calculate the probabilities. The moments of the Galois group, in this sense, may be calculated in the function field setting in the $q\to\infty$ limit, and the first part of the paper is devoted to computing these moments, which eventually gives a theorem instead of a conjecture in that case.

However, using the formulas of \cite{Sawin2022} to compute the probabilities from the moments in this situation is highly nontrivial. 
It turns out to be helpful in our arguments to consider the Galois group with an extra invariant, the \emph{Artin-Verdier fundamental class}. The statement of Theorem \ref{intro-bb} above does not incorporate this extra structure, but we calculate the probabilities $p_H$ by first calculating probabilities of obtaining a group together with an Artin-Verdier fundamental class, and then summing over possible choices of Artin-Verdier fundamental class.  A sum corresponding to possible values of this invariant naturally arises when one tries to compute the probabilities, and it does not seem possible to give a formula for the probabilities in general without that sum, so it is natural to give a conceptual interpretation of this sum by incorporating the Artin-Verdier fundamental class. The construction and consideration of this class along with the unramified Galois groups is one of the key new ideas of this paper. As we will discuss later, this invariant refines several invariants that have been previously considered.

Let $K$ be a number field, let $H$ be the Galois group of an unramified extension  $L/K$ and let $n$ be a positive integer.  The extension $L/K$ gives a map $H^3 ( H, \mathbb Z/n) \to H^3( \mathcal O_K, \mathbb Z/n)$ of \'{e}tale cohomology. If $L/K$ is split at each infinite place then it is easy to see that the composition $H^3 ( H, \mathbb Z/n) \to H^3( \mathcal O_K, \mathbb Z/n)\to H_3(K_v, \mathbb Z/n)$ vanishes for each infinite place $v$ of $K$. With slightly more work, one can see (in Lemma \ref{L:compactsupports}) that in the split at $\infty$ case, the map $H^3 ( H, \mathbb Z/n) \to H^3( \mathcal O_K, \mathbb Z/n)$ factors through a natural map  $H^3 ( H, \mathbb Z/n) \to H^3_c( \mathcal O_K, \mathbb Z/n)$.

If $K$ contains the $n$th roots of unity then, after fixing a generator for the $n$th roots of unity of $K$, we have a composed map \[ H^3 ( H, \mathbb Z/n) \to H^3_c ( \mathcal O_K, \mathbb Z/n)\to H^3_c(\mathcal O_K, \mu_n) \to \mathbb Z/n \] where the first arrow is from $L/K$ as above, the second arrow arises from the map of sheaves $\mathbb Z/n\to\mu_n$ sending $1$ to our fixed generator, and the last arrow comes from Artin-Verdier duality.
We call this map the \emph{Artin-Verdier trace} for $L/K$,
 and the corresponding element in $H_3( H, \mathbb Z/n)$ the \emph{Artin-Verdier fundamental class}  for $L/K$.  More details are given in Section~\ref{S:AV}.

\subsection{The main theorems and conjectures}

 Our main conjecture, combined with our main theorem, predict the distribution of Galois groups together with two types of extra structure -- first, the natural action of $\Gamma$, and second, the Artin-Verdier fundamental class. To express this extra structure, we work with the category of finite $n$-oriented $\Gamma$-groups.
 
 We fix a positive integer $n$ and a finite group $\Gamma$. 
We define a  \emph{$\Gamma$-group} to be a group $H$ together with an action of $\Gamma$ on $H$. 
 We define a \emph{finite $n$-oriented $\Gamma$-group} to be a finite $\Gamma$-group $H$ of order prime to $\abs{\Gamma}$ along with a $\Gamma$-invariant class $s_H \in H_3 ( H, \mathbb Z/n)$. A morphism of $n$-oriented $\Gamma$-groups $(G, s_{G} ) \to (H, s_{H})$ is a $\Gamma$-equivariant group homomorphism $f\colon G\to H$ such that $f_* s_{G} = s_H$.

We define a profinite $n$-oriented $\Gamma$-group to be a profinite group $X$, the inverse limit of finite groups of order prime to $|\Gamma|$, with an action of $\Gamma$ by automorphisms, and a $\Gamma$-invariant class $s_X \in H_3(X, \mathbb Z/n)$. We define morphisms the same as in the finite case. 

We write $n$-oriented $\Gamma$-groups with boldface letters such as $\mathbf H$, and write $\Sur(\mathbf X, \mathbf H)$ for the number of surjections from $\mathbf X$ to $\mathbf H$ in the category of profinite $n$-oriented $\Gamma$-groups, and $\Aut(\bH)$ for the group of automorphisms in the same category.
We use the corresponding non-bolded letter to denote the underlying group with a $\Gamma$-action, e.g. $H$ for $\mathbf H$.

 The conjecture is motivated by an analogous function field theorem, which we now describe the notation for and state. 
  Let $q$ be a prime power, which will always be prime to $\abs{\Gamma}$, and let $\mathbb F_q$ be the finite field with  $q$ elements. An extension of $\mathbb F_q(t)$, given as a Galois subfield of the separable closure of $\mathbb F_q(t)$, together with a choice of isomorphism of the Galois group to $\Gamma$, is called a $\Gamma$-extension of $\mathbb F_q(t)$. For $K$ a $\Gamma$-extension of $\mathbb F_q(t)$, we let $\rDisc K$ be the product of the norms of the places of $\mathbb F_q(t)$ that ramify in $K$ (the norm of a place $v$ of $\mathbb F_q(t)$ being $q^{\deg v}$).  Let $E_\Gamma( q^m, \mathbb F_q(t) )$ be the set of $\Gamma$-extensions
  $K$ of $\mathbb F_q(t)$, split at $\infty$ (the unique place of $\mathbb F_q(t)$ that does not arise from a prime ideal of $\mathbb F_q[t]$), such that  $\rDisc K = q^m$.

For $K$ a $\Gamma$-extension of $\mathbb F_q(T)$, let $K^{ \operatorname{un}, \abs{\Gamma}'} $ be the composition of all finite Galois extensions of $K$ with degree prime to $|\Gamma|$ that are everywhere unramified and split over $\infty$.

If $\mathbb F_q(t)$ contains the $n$th roots of unity, i.e. $n \mid q-1$, then since $K^{ \operatorname{un}, \abs{\Gamma}'}/K$ is an extension with Galois group $\Gal( K^{ \operatorname{un}, \abs{\Gamma}'} /K) $, it has an associated Artin-Verdier fundamental class $s\in H_3(\Gal( K^{ \operatorname{un}, \abs{\Gamma}'} /K)^\Gamma  ,\Z/n)^\Gamma$.
This makes $\bGal( K^{ \operatorname{un}, \abs{\Gamma}'} /K)$ naturally an $n$-oriented $\Gamma$-group.
(The action of $\Gamma$ is by using a homomorphic section of $\Gal( K^{ \operatorname{un}, \abs{\Gamma}'} /\F_q(t))\ra\Gal(K/\F_q(t))=\Gamma$
and conjugating. By the Schur--Zassenhaus theorem \cite[Theorem 2.3.15]{Ribes2010}, this section is well-defined up to conjugation by elements of $\Gal( K^{ \operatorname{un}, \abs{\Gamma}'} /K)$, and thus the resulting  $n$-oriented $\Gamma$-group is well-defined up to isomorphism.)

\begin{theorem}[Theorem \ref{ff-theorem}]\label{intro-ff-theorem} Let $\Gamma$ be a finite group, $n$ a positive integer coprime to $\abs{\Gamma}$, and $\mathbf H$ a finite $n$-oriented $\Gamma$-group with $H_\Gamma=1$.  Let $q$ be a prime power. As long as $q$ satisfies the congruence conditions $ (q, |\Gamma| |H|)=1$, $ q \equiv 1 \bmod n$,  and $(q-1, |H|)= (n,|H|)$ and $q$ is sufficiently large depending on $\Gamma, H$ we have 
\[ \lim_{b \to\infty}\frac{\sum_{m\leq b}   \sum_{K \in E_\Gamma( q^m , \mathbb F_q(t))}  \Sur ( \bGal ( K^{ \operatorname{un}, \abs{\Gamma}'} / K), \mathbf H)  }{ \sum_{m\leq b}\abs{ E_\Gamma( q^m , \mathbb F_q(t))}} =  \frac{ \abs{H^\Gamma} \abs{ H^2( H\rtimes \Gamma , \mathbb Z/n) } }{  \abs{H} \abs{ H^3( H\rtimes \Gamma , \mathbb Z/n) }  }.\] \end{theorem}

We are now ready to state the main conjecture in the number field case. Fix a number field $k$ and a finite group $\Gamma$. We will assume that $k$ lacks a nontrivial everywhere unramified extension, which is satisfied for instance if $k=\mathbb Q$. Let $n$ be the number of roots of unity of order prime to $|\Gamma|$ contained in $k$. A $\Gamma$-extension of $k$ is a Galois subfield $K$ of the algebraic closure of $k$ together with an isomorphism $\Gal(K/k) \cong \Gamma$.  For each $\Gamma$-extension $K$ of $k$, let $K^{ \operatorname{un}, \abs{\Gamma}'}$ be the maximal extension of $K$, unramified at all places including the infinite ones, that is a limit of finite Galois extensions of degree prime to $|\Gamma|$. Then $\bGal ( K^{ \operatorname{un}, \abs{\Gamma}'} / K)$ is a profinite $n$-oriented $\Gamma$-group.

We write $\operatorname{Spec}(k \otimes\mathbb R)$ for the set of infinite places of $k$. 
For $H$ a finite $n$-oriented $\Gamma$-group, let \[ b_{\mathbf H}=  \frac{ \abs{H^\Gamma} \abs{ H^2( H\rtimes \Gamma , \mathbb Z/n) } }{ \abs{ H^3( H\rtimes \Gamma , \mathbb Z/n) }  \prod_{ v \in \operatorname{Spec}(k\otimes \mathbb R)}  \abs{H^{ \gamma_v}} }\] if $ H_\Gamma$ is trivial and  $b_{\mathbf H}=0$ if $H_\Gamma$ is nontrivial. Note that $b_{\mathbf H}$ only depends on the underlying $\Gamma$-group $H$, not on the $n$-orientation.

\begin{conjecture}\label{intro-moments-conjecture} Let $k$ be a number field that lacks any nontrivial everywhere unramified field extension, and let $\Gamma$ be a finite group.
For each infinite place $v$ of $k$, fix a conjugacy class $\gamma_v $ of $\Gamma$ of order at most $2$ if $v$ is real or order $1$ if $v$ is complex.
 Let $n$ be the number of roots of unity of order prime to $\abs{\Gamma}$ contained in $k$. Let $\mathbf H$ be a finite $n$-oriented $\Gamma$-group. 

Among $\Gamma$-extensions $K/k$  such that for each infinite place $v$ of $k$ the conjugacy class of complex conjugation at the place $v$ is in $\gamma_v$, counted in a natural way in which the proportion of fields containing any fixed subfield larger than $k$ is $0$,
$$\E(\Sur ( \bGal ( K^{ \operatorname{un}, \abs{\Gamma}'} / K), \mathbf H)) =b_{\mathbf H}.$$
\end{conjecture}

We are vague about the method of counting fields so as to separate the question of the distribution of 
$\bGal ( K^{ \operatorname{un}, \abs{\Gamma}'} / K)$, from the separate interesting question 
 of under which methods of counting field extensions does one avoid ever having a fixed non-trivial subextension 
 occurring in a positive proportion of fields.  Only under such counting methods do we expect to obtain the distributional results of 
 $\bGal ( K^{ \operatorname{un}, \abs{\Gamma}'} / K)$ conjectured above.  See \cite[Section 6]{Bartel2020} for an example of how a positive proportion of fields containing a fixed subfield can violate even the Cohen-Lenstra-Martinet heuristics for class group distributions,  and \cite{Koymans2023a,Wang2025} for more examples of the fixed subfield issue.  See \cite[Section 1.4]{Wang2025} and \cite{Gundlach2022} for suggestions of general ways to count fields to avoid this issue, though for any specific $\Gamma$ there may be many such possible ways to count $\Gamma$-extensions.  

Now we explain how the moments in Conjecture~\ref{intro-moments-conjecture}
determine a distribution. We first define the space our distributions will be on. We say a profinite $n$-oriented $\Gamma$-group is \emph{small} if it has finitely many open subgroups of each index, i.e if the underlying profinite group is small in the classical sense.

A \emph{level} $\cL$ (of the category of finite $\Gamma$-groups) is a finitely generated formation of finite $\Gamma$-groups (a \emph{formation} is a set of isomorphism classes of $\Gamma$-groups closed under taking quotients and finite subdirect products).  For example, the group $\Z/p\Z$ (with trivial action) generates the level of 
elementary abelian $p$-groups with trivial $\Gamma$-action.

We let $\mathcal P_{\Gamma,n}$ be the set of isomorphism classes of small profinite $n$-oriented $\Gamma$-groups.  A group $X\in \mathcal{P}_{\Gamma,n}$ has a pro-$\cL$ completion $X^\cL$ which is the inverse limit of the (continuous) quotients of $X$ in $\cL$, and this completion is a finite $n$-oriented $\Gamma$-group since it is finite by \cite[Lemma 5.6]{Sawin2022} and inherits a homology class by pushforward from $X$. 

We consider a topology on $\mathcal{P}_{\Gamma,n}$ generated by $\{X \,|\, X^\cL\isom F\}$ as $\cL$ varies through levels and $F$ varies through finite $n$-oriented $\Gamma$-groups whose $\Gamma$-group is in $\cL$. We also consider the associated Borel $\sigma$-algebra.

In the function field case, we have the following probabilistic result.

\begin{theorem}\label{T:intro-measure-ff} 
There exists a unique measure $\nu_{\Gamma, n, \{1\}}$ on $\mathcal P_{\Gamma,n}$ such that for every finite $n$-oriented $\Gamma$-group $\bH$
$$
\int_{X\in \mathcal{P}_{\Gamma,n} } \Sur(X,\bH) d\nu_{\Gamma,n,\{1\}} =\begin{cases} \frac{ \abs{H^\Gamma} \abs{ H^2( H\rtimes \Gamma , \mathbb Z/n) } }{  \abs{H} \abs{ H^3( H\rtimes \Gamma , \mathbb Z/n) }}& \textrm{if } H_\Gamma = 1 \\ 0 & \textrm{if } H_\Gamma \neq 1 \end{cases} 
$$

For any level $\cL$ of the category of finite $\Gamma$-groups, 
and finite $n$-oriented $\Gamma$-group $\bH$ whose underlying $\Gamma$-group is in $\cL$, the quantity $\nu_{\Gamma, n, \{1\}}(\{X|X^\cL\isom  \bH\})$ is described by an explicit formula given later in \S\ref{ss-intro-formulas}. Furthermore, in the same setting, letting $M$ be the product of all primes dividing the orders of elements of $\cL$, we have

\[ \lim_{ \substack{ q\to\infty \\ (q, |\Gamma| M )=1 \\   (q-1, nM) = n }} \limsup_{b \to\infty} \frac{\sum_{m\leq b}  \abs{ \{ K \in  E_\Gamma( q^m , \mathbb F_q(t))  \mid  \bGal ( K^{ \operatorname{un}, \abs{\Gamma}'} / K)^\cL \isom \bH   \}} }{ \sum_{m\leq b}\abs{ E_\Gamma( q^m , \mathbb F_q(t))}} =  \nu_{\Gamma, n, \{1\}}(\{X|X^\cL\isom  \bH\})\] and

\[ \lim_{ \substack{ q\to\infty \\ (q, |\Gamma| M )=1 \\   (q-1, nM) = n }} \liminf_{b \to\infty} \frac{\sum_{m\leq b} \abs{ \{ K \in  E_\Gamma( q^m , \mathbb F_q(t))  \mid  \bGal ( K^{ \operatorname{un}, \abs{\Gamma}'} / K)^\cL \isom \bH   \}} }{ \sum_{m\leq b}\abs{ E_\Gamma( q^m , \mathbb F_q(t))}} =  \nu_{\Gamma, n, \{1\}}(\{X|X^\cL\isom  \bH\})\]\end{theorem}

The analogous result in the number field case is as follows.
\begin{theorem}\label{T:intro-measure-nf}
Let $n$ be a positive integer, $\Gamma$ a finite group of order prime to $n$,  and $\overline{\gamma}$ a tuple of conjugacy classes of $\Gamma$ of order at most $2$, with at least one trivial element unless $n=1$.

Then there exists a unique measure $\nu_{\Gamma, n, \overline{\gamma}} $  on $\mathcal P_{\Gamma,n}$ such that for every finite $n$-oriented $\Gamma$-group $\bH$ we have $$
\int_{X\in \mathcal{P}_{\Gamma,n} } \Sur(X,\bH) d\nu_{\Gamma, n, \overline{\gamma}} =b_{\mathbf H} .$$ 

For any level $\cL$ of the category of finite $\Gamma$-groups, 
and $\bH$ a finite $n$-oriented $\Gamma$-groups whose underlying $\Gamma$-group is in $\cL$, the quantity $\nu_{\Gamma, n, \overline{\gamma}} (\{X|X^\cL\isom  \bH\})$ is described by an explicit formula given later in \S\ref{ss-intro-formulas}.

Furthermore, in the same setting, assuming Conjecture \ref{intro-moments-conjecture} we have, in the same ordering of fields as in the conjecture
\[ \Prob( \bGal ( K^{ \operatorname{un}, \abs{\Gamma}'} / K)^{\cL} \isom \bH )= \nu_{\Gamma, n, \overline{\gamma}} (\{X|X^\cL\isom  \bH\}). \]
 \end{theorem}

Suppose $k$ is a number field that has exactly $n$ roots of unity of order prime to $\abs{\Gamma}$, and such that the elements of $\overline{\gamma}$ can be placed in bijection with the infinite places of $k$ with all the non-trivial conjugacy classes corresponding to real places of $k$. 
Then either $n\geq 3$, and $k$ has all complex places, so all $\gamma_v$ are trivial, or $|\Gamma|$ is odd, so all $\gamma_v$ are trivial, or $n=1$.  So the hypothesis on $\overline{\gamma}$ holds in all cases of interest.

 These theorems have three parts: The first states that there exists a unique measure with the relevant moments, the second gives a formula for that measure (the full details of which we delay to \S\ref{ss-intro-formulas}), and the third says that this measure gives the correct probabilities of obtaining a given group, in a double limit sense in the function field case and conditionally on Conjecture \ref{intro-moments-conjecture} in the number field case.

We point out that, even though the explicit formulas in \S\ref{ss-intro-formulas} require splitting into many different cases, in each case, it is easy to determine when the probability is zero and when it is positive. Whether the probabilities are positive or zero is particularly important because, when the probability is positive, from Theorem \ref{T:intro-measure-ff} we obtain an existence result for function fields with $ \bGal ( K^{ \operatorname{un}, \abs{\Gamma}'} / K)^\cL \isom \bH $, and, when the probability is zero, Theorem \ref{T:zero} gives a corresponding nonexistence result.

\subsection{Relation with prior work}\label{ss-relations}

A number of prior works have considered the distribution of the Galois group of the maximal unramified extension of a random number field, and even more have considered the distribution of its abelianization, i.e., the class group. We review here the relationship to that prior work, with more details given later in \S\ref{s-prior-work}. Crucially, we study the Galois group together with some extra data, the Artin-Verdier trace, which generalizes some extra data considered in prior works.

In 1984,  Cohen and Lenstra \cite{Cohen1984} gave conjectures for the distribution of the odd parts of class groups of imaginary and real quadratic fields,
as well as for any finite abelian group $A$,  the prime-to-$|A|$ parts of class groups of totally real $A$-extensions of $\Q$.  Cohen and Martinet \cite{Cohen1990} generalized these conjectures
to the situation of an arbitrary number field $K_0$ as a base field, and arbitrary group $\Gamma$, 
giving conjectures for distributions of parts of class groups of order relatively prime to $|\Gamma|$ among $\Gamma$-extensions of a fixed $K_0$ with any fixed behavior at the infinite places of $K_0$.    However, it was noted that these conjectures appeared to be wrong 
for the parts of class groups at primes dividing the number of roots of unity in the base field, in empirical number field work of Malle \cite{Malle2008, Malle2010}
and 
theoretical function field work of Achter  \cite{Achter2006}
and Garton \cite{Garton2015}.
There have since been many papers aimed at correcting these class group distribution conjectures in the presence of roots of unity, including work of   
Malle \cite{Malle2008,Malle2010} in his original papers, Garton \cite{Garton2015}, Adam and Malle \cite{Adam2015}, 
Lipnowski, the first author, and Tsimerman \cite{Lipnowski2020}, and most recently of the two current authors \cite{Sawin2023}.  See \cite[Section 1.2]{Sawin2023} for more discussion of the relationship between these papers.

Our  work~\cite{Sawin2023} made conjectures for the $|\Gamma|$-prime part of the class group of Galois $\Gamma$-extensions of any fixed number field, without any extra data beyond the $\Gamma$ action. The predictions of this paper can be compared with \cite{Sawin2023} by projecting the measure $\nu_{\Gamma, n, \overline{\gamma}}$ onto abelian groups and then forgetting the Artin-Verdier trace. This produces a distribution whose moments agree with the moments conjectured in \cite{Sawin2023}, and since it is shown in \cite{Sawin2023} that this distribution is determined by its moments, it follows that the two distributions agree. The paper \cite{Sawin2023} itself built on the prior work described above giving conjectures for distributions of class groups, agreeing with most of it but disagreeing with some, and it follows that the current paper has the same relationships with prior conjectures, as described in detail in \cite[Section 12]{Sawin2023}.

Lipnowski, Sawin, and Tsimerman~\cite{Lipnowski2020} considered the class group with two invariants, the $\omega$ and $\psi$ invariants, each describable in a different way as bilinear forms on the dual of the class group, with a relationship between them. We show that both of these invariants can be calculated from the Artin-Verdier invariant -- roughly, the Artin-Verdier invariant gives a linear form on the third cohomology, and we can use this to obtain a bilinear form on the first cohomology in two different ways, either cupping two classes and then taking a Bockstein homomorphism or cupping one class with a Bockstein homomorphism applied to the other. These two constructions turn out to give the two bilinear forms. This differs from the approach in~\cite{Lipnowski2020} to defining these invariants in that work also considered Bockstein homomorphisms in flat cohomology arising from exact sequences of group schemes and here we use only Bockstein homomorphisms in group cohomology arising from exact sequences of finite groups. Thus, the Artin-Verdier trace generalizes the $\omega$ and $\psi$ invariants and gives additional data even in the case of abelian groups (since from the Artin-Verdier trace we also obtain a trilinear form arising from cup product of three classes in first cohomology), which is not computable from only the $\omega$ and $\psi$ invariants, though this does not occur in the $\Gamma = \mathbb Z/2$ setting considered in~\cite{Lipnowski2020}.

Morgan and Smith~\cite{MorganSmith} obtained the $\psi$ invariant as part of a general formalism of Cassels-Tate pairings. They also defined~\cite[Notation 6.3 and Proposition 6.4]{MorganSmith} an invariant that satisfies properties identical to the key properties of the $\omega$ invariant. They observed that it seems likely that this invariant is equal to the $\omega$ invariant, but did not prove this. Since we prove that a pairing defined via Artin-Verdier duality agrees with the $\omega$ invariant, it also seems likely that this pairing agrees with the pairing $\operatorname{CTP}_{D(m)}$ of \cite[Proposition 6.4]{MorganSmith}, but we also do not prove this.

There have also been a number of papers on conjectures for the distribution of the non-abelian generalization of the class group, i.e. the Galois group of the maximal unramified extension of a random number field.
Boston, Bush, and Hajir \cite{Boston2017a,Boston2021} gave conjectures for the distribution of the $p$-class tower group (i.e. the Galois group of the maximal unramified Galois extension of degree a power of $p$) of quadratic fields, for $p$ odd.  
Liu, the second author,
and Zurieck-Brown~\cite{Liu2024} considered the Galois group of the maximal unramified extension of $\Gamma$-extensions of $\Q$ for an arbitrary finite group $\Gamma$, making conjectures similar to our own, except that they only considered extensions of degree relatively prime to the number of roots of unity in $k$ (i.e., in the case $k=\mathbb Q$, relatively prime to $2$). In that paper, as in this one, the conjectures were justified by moment theorems in a function field model. Since the moments conjectured in this paper agree with the moments conjectured in \cite{Liu2024} in the relevant special case, and we check the distribution is uniquely determined by its moments, it follows that the conjectural distributions from the two papers agree in this special case. A key difference is that \cite{Boston2017a,Boston2021,Liu2024} all construct their distributions from a model of random profinite groups based on generators and explicit random relations, whereas our distribution is constructed using a general formalism to find a distribution with given moments. In this special case of extensions of degree relatively prime to the number of roots of unity in the base field, the Artin-Verdier trace is trivial. Correspondingly, \cite{Boston2017a,Boston2021,Liu2024} did not consider any extra data on their profinite group beyond the $\Gamma$ action.

Boston and Bush~\cite{Boston2021a} considered the 2-class tower group  of cyclic cubic fields, finding empirical data and making conjectures based on that. They consider multiple variants -- the wide class group tower where the extension is assumed split at infinite places, the narrow class group tower where the extension is not assumed split at infinite places, and the joint distribution of the wide and narrow tower as a pair of groups. Since we consider extensions that are split at $\infty$, our conjectures apply to the wide case. As mentioned above, we check that the measure $\nu_{\mathbb Z/3\mathbb Z, \mathbb Q, (1)}$ arising from our Conjecture \ref{intro-moments-conjecture} agrees with Boston and Bush's numerical data for the smallest two rank two $2$-groups that can appear. The Artin-Verdier trace in this setting is nontrivial, so applying our conjecture requires summing over possible values of the Artin-Verdier trace.

In the nonabelian setting, the study of extra invariants started with unpublished work of Ellenberg, Venkatesh, and Westerland. Building on work of Serre \cite{Serre1990} and Fried \cite{Fried1995a}, they defined a lifting invariant which separates components of a Hurwitz space parameterizing covers of the projective line. This lifting invariant was defined in more detail and studied by the second author in \cite[Section 3]{Wood2019} and \cite{Wood2021}.  In \cite{Wood2019} the second author gives conjectures on the moments of the Galois group of the maximal unramified extension of a random quadratic field, including extensions of even degree, and thus in the case where roots of unity are relevant.  There are refined conjectures in \cite[Section 5]{Wood2019} that describe the moments with a fixed lifting invariant.
These conjectures are based on function field theorems. 
Liu~\cite{Liu2022} defined a $\omega$-invariant of an unramified Galois extension $L/K$ of a $\Gamma$-extension $K / \mathbb F_q(t)$, for a general finite group $\Gamma$, and showed it agrees with the prime-to-the-order-of-$\Gamma$ part of the lifting invariant of $L/\mathbb F_q(t)$ as defined in \cite{Wood2021}.
This allowed Liu to prove function field theorems on counting unramified extensions with a fixed $\omega$-invariant and motivated 
 refined conjectures in the number field case.
 We show that the Artin-Verdier trace may be used to calculate Liu's $\omega$-invariant. In the abelian group case, all these invariants specialize to the $\omega$-invariant. As we mentioned above, the Artin-Verdier trace gives additional data not computable from these invariants.
 Combined with Liu's result, the Artin-Verdier trace may be used to compute the prime-to-the-order-of-$\Gamma$ part of the lifting invariant of \cite{Wood2021}, and thus the Artin-Verdier trace  refines and generalizes the previously studied invariants on Galois groups of maximal unramified extensions. 

\subsection{Evidence}

In addition to the function field evidence discussed above, we find evidence for Conjecture \ref{intro-moments-conjecture} by calculating the probabilities $\nu_{\Gamma, n, \overline{\gamma}} (\{X|X^\cL\isom  \bH\})$, which that conjecture predicts is the probability that $\bH$ appears as the maximal quotient in $\cW$ of the Galois groups of maximal unramified extensions, and checking these predictions. When the predicted probability is zero, we can check the prediction by proving that $\bH$ does not appear as the maximal quotient in $\cW$ of a Galois group of a maximal unramified extension. We do this for every group whose predicted probability is zero. When the predicted probability is nonzero, we must instead compare against empirical data. We do not generate new data in this paper, but rather compare against some data produced in prior work.

More precisely, in Proposition \ref{nzp-criterion} below, we find necessary and sufficient conditions to have $\nu_{\Gamma, n, \overline{\gamma}} (\{X|X^\cL\isom  \bH\})=0$. In Proposition \ref{nzp-converse}, we check that groups violating these conditions in fact do not appear as  $\Gal ( K^{ \operatorname{un}, \abs{\Gamma}'} / K)^{\cL}$ (except if $K$ has additional roots of unity not contained in $k$, which should happen with density zero regardless). To do this, we need the following, purely Galois-theoretic, statement.

\begin{theorem}[Theorems \ref{ne-nd-prof}, \ref{ne-d-prof}, and \ref{rou-descent}]\label{T:zero}  Let $k$ be a number field, $\Gamma$ a finite group, $K$ an extension of $k$ with Galois group $\Gamma$. Let $p$ be a prime not dividing $\abs{\Gamma}$ and let $V$ be an irreducible representation of  $\Gal ( K^{ \operatorname{un}, \abs{\Gamma}'} / K) \rtimes \Gamma$ over $\mathbb F_p$. Assume $k$ admits no nontrivial unramified extensions.

If $k$ contains the $p$th roots of unity, then
\[  |H^1 (  \Gal ( K^{ \operatorname{un}, \abs{\Gamma}'} / K) \rtimes \Gamma , V^\vee)| |H^2 (  \Gal ( K^{ \operatorname{un}, \abs{\Gamma}'} / K) \rtimes \Gamma , V)^{s_K}| \] \[\leq |H^1 (  \Gal ( K^{ \operatorname{un}, \abs{\Gamma}'} / K) \rtimes \Gamma , V)|   \frac{  \prod_{ v \textrm{ archimedean place of } k} | V ^{ \Gal(k_v)}| }{|V^{ \Gal ( K^{ \operatorname{un}, \abs{\Gamma}'} / K) \rtimes \Gamma}|}  .\]
 
If $k$ does not contain the $p$th roots of unity,  and  $K$ does not contain the $p$th roots of unity, we have 
\[ |H^2( \Gal ( K^{ \operatorname{un}, \abs{\Gamma}'} / K) \rtimes \Gamma , V)| \leq |H^1( \Gal ( K^{ \operatorname{un}, \abs{\Gamma}'} / K) \rtimes \Gamma , V)|    \frac{ \prod_{ v \textrm{ archimedean place of } k} | V ^{ \Gal(k_v)}|}{ |V^{ \Gal ( K^{ \operatorname{un}, \abs{\Gamma}'} / K ) \rtimes \Gamma}|} .\]
\end{theorem}   

This verifies a prediction of Conjecture \ref{intro-moments-conjecture} that is spelled out in Theorem~\ref{T:intro-measure-nf}. 

The proof of Theorem~\ref{T:zero} 
is relatively short, and relies on comparing the Galois cohomology groups appearing in its statement to flat cohomology groups of number fields so that Artin-Verdier duality and Euler characteristic formulas for flat cohomology may be applied. The significance of Theorem~\ref{T:zero} is that, according to Conjecture \ref{intro-moments-conjecture} and Theorem~\ref{T:intro-measure-nf}, it should be the \emph{only} restriction on $ \Gal ( K^{ \operatorname{un}, \abs{\Gamma}'} / K)$ which is detectable by examining finite quotients and holds for an arbitrary number field.

We now turn to nonzero probabilities. 

First observe that $\Gal ( K^{ \operatorname{un}, \abs{\Gamma}'} / K)$ is the prime-to-$\abs{\Gamma}$ part of the class group of $K$, so Conjecture \ref{intro-moments-conjecture} implies a conjectural distribution for the prime-to-$\abs{\Gamma}$ part of the class group. As mentioned in \S\ref{ss-relations}, this agrees with the conjectural distribution of class groups in \cite{Sawin2023}. Malle~\cite{Malle2008,Malle2010} has previously produced a large amount of data on the distribution of the class groups of number fields containing roots of unity, and \cite{Sawin2023} checks that this data agrees with the conjectures of \cite{Sawin2023}, and hence also with the conjectures of this paper. For the distribution with extra invariants, note that Lipnowski, Sawin, and Tsimerman~\cite{Lipnowski2020} produced data on the distribution of the $\omega$ and $\psi$ invariants on class groups of quadratic extensions of $\mathbb Q(\mu_3)$, checking that this data agrees with the conjectures of \cite{Lipnowski2020}, and hence also with the conjectures of this paper.

Genuinely nonabelian data was produced by Boston and Bush~\cite{Boston2021a}, specifically on the Galois group of the maximal unramified Galois extension with degree a power of $2$ of cyclic cubic fields. As mentioned above, we check that our predicted probabilities agree with Boston and Bush's data for the two smallest $2$-groups considered in \cite{Boston2021a}. Checking agreement for further groups would require more intricate calculations in group cohomology.

\subsection{Ideas of the proof}

Our computation of the function field moment builds on the strategy of prior works such as \cite{Liu2024}, \cite{Liu2022}, and \cite{LandesmanLevy}, suitably modified to keep track of the Artin-Verdier trace. These works showed that the moment may be interpreted as a ratio of counts of $\mathbb F_q$-points on certain Hurwitz spaces. In the large $q$ limit, the number of points of the Hurwitz spaces may be estimated using a Lang-Weil-like argument if one can count the geometrically irreducible connected components of Hurwitz spaces. This is done using topological theorems to count the connected components of Hurwitz spaces over $\mathbb C$, deducing counts for the connected components of Hurwitz spaces over $\overline{\mathbb F}_q$, and then understanding the action of Frobenius on these components to count the fixed points, since Frobenius-fixed points over $\overline{\mathbb F}_q$ correspond to geometrically irreducible components over $\mathbb F_q$. For fixed $q$, Lang-Weil does not suffice and control of the low-degree cohomology groups of individual components is needed in addition to the count of the geometrically irreducible components. For Hurwitz spaces of arbitrary groups, this was handled in the breakthrough work of \cite{LandesmanLevy}, which enabled a more refined counting theorem.

We check that the moments we are interested in can be calculated as counts of $\mathbb F_q$-points on certain covers of Hurwitz spaces.
In particular, we require the use of new spaces in order to track the Artin-Verdier trace.
 To count the geometrically irreducible connected components of these covers of Hurwitz spaces, we calculate the monodromy group of this covering. This may be done by reduction to $\mathbb C$, which requires studying a purely topological question.

Fix a finite group $H$ with an action of a finite group $\Gamma$ and assume that $H_\Gamma=1$. We define an invariant, the braid fundamental class, associated to a braid together with a homomorphism from the fundamental group of the complement of that braid to $H \rtimes \Gamma$, that, for each strand of the braid, sends a small loop around that strand to an element of a conjugate of $\Gamma$. Using this data, we construct a $3$-manifold as covering of $S^3$, branched along the braid, with monodromy group contained in $\Gamma$. This 3-manifold admits a natural homomorphism from its fundamental group to $H$, or in other words, a map to the classifying space $BH$. The braid fundamental class is the fundamental class of this 3-manifold inside the homology group $H_3(H,\mathbb Z)$. 

The braid fundamental class is $\Gamma$-invariant, and our first key topological theorem shows that every $\Gamma$-invariant class in $H_3(H,\mathbb Z)$ is the braid fundamental class of some braid. To prove this, we first show that every $\Gamma$-invariant class in $H_3(H,\mathbb Z)$ is the fundamental class of an unbranched covering with monodromy group $\Gamma$ of a 3-manifold. We replace the base 3-manifold with $S^3$ by Dehn surgery, which requires introducing curves where the covering branches. This produces a link in $S^3$, which can be expressed as the closure of the braid.

Our second key topological theorem computes the stable low-degree homology of certain covers of Hurwitz space. These covers may also be described using the braid fundamental class: The fundamental group of a configuration space of points in the plane is a braid group, Since the $H\times \Gamma$-Hurwitz space is a covering space of configuration space, its fundamental group is a subgroup of the braid group. There is a natural homomorphism from the fundamental group of the complement of a braid in this subgroup to $H \rtimes \Gamma$. The braids in this subgroup with trivial braid fundamental class form a further subgroup, which corresponds to a finite covering of Hurwitz space. It is this covering we compute the cohomology of.

We use the results of \cite{LandesmanLevy} to compute the homology of this cover of Hurwitz space. It is not immediately obvious how to do this, since the rational homology of a space only gives a lower bound for the homology of its finite covering spaces. However, we show that this covering of Hurwitz space is itself covered by a different Hurwitz space of a larger group, generalizing an argument of \cite{Lipnowski2020}. This gives both upper bounds and lower bounds for the homology, which match each other, giving a computation of the homology of the cover.

Since our covering has abelian monodromy group, its rational cohomology splits as a sum of the cohomology of certain rank one local systems on Hurwitz space. Our homology computation is equivalent to a vanishing result for the low-degree cohomology of the nontrivial rank one local systems appearing in this decomposition. Very recently, similar results for low-degree cohomology of certain rank one local systems arising from 2-cocycles on racks were proven by Ellenberg and Shusterman~\cite{Ellenberg2026}, as a special case of vanishing results for low-degree cohomology of braided vector spaces.  However, these results do not seem to apply in our case, as the local systems we consider do not arise from 2-cocycles on racks. Thus we require a different approach, which the covering by a Hurwitz space of a larger group provides.

Our proof of the probabilistic theorem requires evaluating formulas from \cite{Sawin2022} that express the probability of obtaining a group $H$ as a weighted sum of moments of various extensions $G$ of $H$ by other groups $F$. Our formula for the moment of $G$ involves the size of the group cohomology group $H^2(G \rtimes \Gamma, \mathbb Z/n)$. It is therefore natural to apply the Lyndon-Hochschild-Serre spectral sequence, which computes the group cohomology of the extension $G \rtimes \Gamma$ of $H \rtimes \Gamma$ by $F$ in terms of the cohomology of $H$ and $F$: More precisely, it is a spectral sequence converging to $H^{p+q} ( G \rtimes \Gamma, \mathbb Z/n)$ whose second page is $H^p ( H \rtimes \Gamma, H^q( F, \mathbb Z/n))$.

Applying the spectral sequence requires computing $H^p ( H \rtimes \Gamma, H^q( F, \mathbb Z/n))$ for $p,q$ small and computing several differentials. Computing these differentials, and then using the differentials to calculate the moments, is the most technically difficult part of the paper. We use a description of the differentials due to Huebschmann \cite{Huebschmann1981}, but need to make this description more concrete and simplify it, taking advantage of the fact that only certain groups may appear as $F$ in our setup.

The reader interested in the proof may benefit from first reading our previous paper \cite{Sawin2024} on the distribution of the fundamental groups of 3-manifolds. This paper included a similar Lyndon-Hochschild-Serre spectral sequence argument, but it was considerably simpler as several difficulties that occur for Galois groups of number fields do not appear for fundamental groups of 3-manifolds. First, we need to keep track of the action of $\Gamma$, which did not exist in the 3-manifold setting we studied (though could, of course, if one studied $3$-manifolds with actions of $\Gamma$). Second, there are multiple cases in this work depending on the number of roots of unity, but only one case needs to be considered for 3-manifolds. (Roughly, an oriented 3-manifold is analogous to a number field containing all roots of unity that are possibly relevant.) Third, in \cite{Sawin2024} we used Poincar\'e duality to obtain strong restrictions on the fundamental group of a 3-manifold. We only needed to consider groups $H$ satisfying these restriction, which forced certain extension groups to vanish. In this paper, several calculations involve these extension groups, which need not vanish. (Artin-Verdier duality can be used to prove restrictions on the Galois group of the maximal unramified extension of a number field, as we do in Lemma \ref{ne-d-prof}, but owing to the influence of Archimedean places, these restrictions are weaker than in the 3-manifold case, and do not show that these extension groups vanish.) For all these reasons, the proof of \cite{Sawin2024} can serve as a simpler model for the proofs in this paper.

\subsection{Notation}\label{s-notation}
Before giving the formulas for the probabilities in our main theorems, we must introduce some notation.  In order to have things in one place, we collect here notation that is used throughout the paper.

Throughout the paper $\Gamma$ is a finite group and $n$ is a positive integer relatively prime to $|\Gamma|$.
We write $\Cat$ for the category of finite $n$-oriented $\Gamma$-groups.
The letter $q$ always denotes a prime power and $\F_q$ the finite field of order $q$.
For any field $K$ that contains the $n$th roots of unity, we fix a generator $\xi$ for the $n$th roots of unity in $K$.

{\bf Asymptotic notation:} We write 
$f(X) = O(g(X))$ to mean there exists a constant $C$ such that $|f(X)| \leq Cg(X)$ for all values of $X$, and
$f(P,X) = O_P(g(P,X))$ to mean that there is a function $C(P)$ of parameters $P$ such that $|f(P,X)| \leq C(P)g(P,X)$ for all values of $X$.
Sometimes we write $O$ for $O_P$ when we specify in words what the constant may depend on.

{\bf Group theory:} Given a group $\Gamma$, a \emph{$\Gamma$-group} $G$ is a group $G$ with an action of $\Gamma$ by automorphisms.  
A \emph{morphism of $\Gamma$-groups} is a $\Gamma$-equivariant morphism of the underlying groups, and we write $\Aut_\Gamma(G)$ for the automorphisms of a $\Gamma$-group $G$.  When $\bG$ is an $n$-oriented $\Gamma$-group, we write $\Aut(\bG)$ for the group of automorphisms in that category,  without any subscript since the bolding serves as a reminder.
A \emph{$\Gamma$-subgroup} $H$ of $G$ is a subgroup of $G$ closed under the $\Gamma$-action, and a \emph{$\Gamma$-group quotient} is the image of a surjective $\Gamma$-group morphism.  
A $\Gamma$-group 
is a \emph{simple $\Gamma$-group} if it is nontrivial and contains no proper, nontrivial normal $\Gamma$-subgroup.  
A $\Gamma$-group is \emph{semisimple} if it is a finite direct product of simple $\Gamma$-groups.

A \emph{$[\Gamma]$-group} is a group $G$ together with a homomorphism $\rho: \Gamma \ra \Out(G)$, though we often leave $\rho$ implicit.
An \emph{isomorphism of $[\Gamma]$-groups} from $(G,\rho)$ to $(G,\rho')$ is given by a group isomorphism $f:G\ra G'$ such that the induced map $f_*:\Out(G)\isom \Out(G')$ satisfies $f_* \circ \rho=\rho'$.
We write $\Aut_{[\Gamma]}(G)$ for the automorphisms of a $[\Gamma]$-group $G$.
We have that $\Gamma $ acts on the set of normal subgroups of an $[\Gamma ]$-group $G$, and we say a nontrivial $[\Gamma]$-group is \emph{simple} if it has no nontrivial proper fixed points for this action.  
 A \emph{semisimple} $[\Gamma]$-group is a finite direct product of simple $[\Gamma]$-groups.  
If we have an exact sequence of groups $1\ra N \ra G \ra  \Gamma \ra 1$, then $N$ is naturally a $[\Gamma ]$-group.  A normal subgroup $N'$ of $N$ is fixed by the $\Gamma$-action if and only if it is normal in $G$. 
Hence $N$ is a simple (semisimple) $G$-group (via conjugation) if and only if $N$ is a simple (semisimple) $[\Gamma]$-group.

Simple finite $\Gamma$-groups or $[\Gamma]$-groups are characteristically simple as groups, and hence their underlying group is product of isomorphic finite simple groups.

We write $\Out(G)$ for the group of outer automorphisms of a group, $\Gamma$-group, or $[\Gamma]$-group $G$, where $\Out(G)$ only depends on the underlying group structure of $G$.

We often work with profinite groups. 
Whenever we refer to morphisms of profinite groups, we always mean continuous homomorphisms.
When considering the homology and cohomology of profinite groups, we always mean the continuous cohomology. When discussing quotients of profinite groups, we always mean quotients by closed normal subgroups, which are themselves profinite groups. 

For an abelian group $G$, we write $G^\vee$ for $\Hom(G,\Q/\Z)$.

{\bf Invariants and coinvariants:}
For a $\Gamma$-group $G$, we write $G^\Gamma$ for the invariants, and for $\gamma\in\Gamma$, we write $G^\gamma$ for the elements fixed by $\gamma$.  We write $G_\Gamma$ for the coinvariants, i.e. the maximal $\Gamma$-invariant quotient of $G$, which is the quotient of $G$ by the normal subgroup generated by the elements $g^{-1}\gamma(g)$ for each $g\in G$ and $\gamma\in \Gamma$.

When $G$ is abelian and $(|G|,|\Gamma|)=1$, we have that the natural map $G^\Gamma \ra G_\Gamma$ is an isomorphism, and 
$(G_\Gamma)^\vee\isom (G^\vee)_\Gamma$.

{\bf Products:}
By convention, a product $\prod_{k=a}^{b}$ where $a>b$ is $1$.

{\bf $q$-series:}
 We define an adjusted $q$-Pochhammer symbol $(q)_n:=\prod_{i=1}^n (1-q^{-i})$ for any positive integer $n$.  By convention $(q)_0=1$. 

{\bf Tensor powers, $\Delta$:}
If $A$ is a finitely generated abelian group, there is an $\Aut(A)$-equivariant map $\Delta: A\tensor A \ra A\tensor A$ given by $a\tensor b \mapsto a\tensor b - b\tensor a$, which descends to an injection $\wedge^2 A \ra A^{\tensor 2}$.
We define $\wedge_2 A$ to be the image of $\Delta$, and $\Delta$ induces an isomorphism $\wedge^2 A \isom \wedge_2 A$.

{\bf Cohomology of schemes:} When taking cohomology groups of a scheme, we always mean \'{e}tale cohomology groups, except in Sections \ref{s-prior-work} and \ref{S:NE}  where we will consider flat cohomology groups, denoted with the subscript $\operatorname{fppf}$, as well.

{\bf Affine Symplectic Group:} 
We follow Gurevich and Hadani \cite{Gurevich2012}.
Given a finite dimensional $\F_2$-vector space $V$, and a full rank form $\omega\in \wedge_2 V^\vee$, we let $\Sp(V)$ be the set of automorphisms of $V$ preserving $\omega$ (and so this definition will always require an implicit $\omega$, which should be clear from context).  If $[f]\in H^2(V,\Z/4)$ is the unique class represented by a bilinear form $f$ such that $\Delta f=\omega$ (which exists by Lemma~\ref{L:H2struc}), then we define $\ASp(V)$ to be the 
automorphisms of the extension $E$ of $V$ by $\Z/4$ corresponding to $f$ that fix $\Z/4$ pointwise and act on $V$ through $\Sp(V)$.
We have an exact sequence
$$
1\ra \Hom(V,\F_2) \ra \ASp(V) \ra \Sp(V) \ra 1.
$$

For more about the affine symplectic group, including the proof that this exact sequence exists, see \cite[\S2.1]{Sawin2024}.
This exact sequence represents a class $\Phi \in H^2 ( \Sp(V), \Hom(V, \mathbb F_2))$ which combined with $\omega^{-1} \colon \Hom (V, \F_2) \to V$ gives a class $\omega^{-1}_*(\Phi )\in H^2(\Sp(V), V)$.

{\bf Representations over finite fields:}
All representations of groups considered in this paper are finite-dimensional.

Let $p$ be a prime and $V$ an irreducible representation of a finite group $\Pi$ over $\F_p$, with $\kappa=\End_\Pi(V)$.
(In most of this paper, $\Pi$ will be $H\rtimes \Gamma$, where $\Gamma$ is a finite group and $H$ a finite $\Gamma$-group.)
We write $V^\vee$ for the dual representation $\Hom(V,\F_p)$ (which under the inclusion $\F_p\ra\Q/\Z$ agrees with our notation above), 
and say $V$ is \emph{self-dual} if $V\isom V^\vee$ as $\Pi$-representations.
Note that $V$ is self-dual if and only if $(V\tensor V)^\Pi\ne 0$.
The trace map $\kappa\ra\F_p$ gives an $\Pi$-equivariant isomorphism $\Hom_\kappa(W,\kappa)\isom W^\vee$.

We call $V$ \emph{unitary}  if $(V\tensor V)^\Pi\ne 0$ but $(V\tensor_\kappa V)^\Pi= 0$ (i.e. $V$ is self-dual over $\F_p$ but not over $\kappa$).
We call $V$ \emph{symmetric}  if $(V\tensor_\kappa V)^\Pi\ne 0$, but $(\wedge_\kappa ^2V)^\Pi=0$.
We call $V$ \emph{symplectic}  if $(\wedge_\kappa ^2 V)^\Pi\ne0$. 

For self-dual irreducible representations $V$, we define a modified Frobenius-Schur indicator $\epsilon_V$.
In characteristic $2$, the value of $\epsilon_V$ will actually also depend on whether $4\mid n$, and will not be defined for every representation.

When $p$ is odd,  we call all  self-dual irreducible representations  \emph{non-anomalous}, and we define $\epsilon_V$ as follows:
\begin{itemize}
\item $\epsilon_V=0$ if $V$ is unitary,
\item $\epsilon_V=1$ if $V$ is symmetric, and
\item $\epsilon_V=-1$ if $V$ is symplectic.
\end{itemize}

Now we consider the case $p=2$.
We say $V$ is \emph{A-symplectic} if $V$ is symplectic and a map $\Pi\ra\Sp(V)$
for an invariant form $\phi\in \wedge_\kappa^2 V$ factors through $\Pi\ra\ASp(V)$.
We say $V$ is \emph{intermediate} if $V$ is symplectic but not $A$-symplectic.  
We say $V$ is \emph{$\F_2$-orthogonal} 
if there is a non-zero invariant quadratic form in $\Sym^2 V$.
We say $V$ is \emph{anomalous}
\begin{itemize}
\item if $4\mid n$ and $V$ is intermediate or,
\item  if $2\mid n$ but $4\nmid n$, and $V$ is intermediate and $\F_2$-orthogonal.
\end{itemize}
Any other $V$ is called \emph{non-anomalous}.

When $p=2$ and $4\mid n$, we define $\epsilon_V$ for non-anomalous self-dual $V$  
as follows:
\begin{itemize}
\item $\epsilon_V=0$ if $V$ is unitary,
\item $\epsilon_V=1$ if $V$ is symmetric (which is the same as trivial),  and
\item $\epsilon_V=-1$ if $V$ is A-symplectic.
\end{itemize}

When $p=2$ and $2\mid n$ but $4\nmid n$, we define $\epsilon_V$ 
for non-anomalous self-dual $V$
as follows:
\begin{itemize}
\item $\epsilon_V=0$ if $V$ is unitary and $\F_2$-orthogonal,
\item $\epsilon_V=1$ if $V$ is symmetric (which is the same as trivial)
\item $\epsilon_V=1$ if $V$ is non-$\F_2$-orthogonal, and
\item $\epsilon_V=-1$ if $V$ is A-symplectic and $\F_2$-orthogonal.
\end{itemize}

For anomalous $V$, $\epsilon_V$ is not defined.

\subsection{Explicit formulas for the probabilities}\label{ss-intro-formulas}

We now describe the terminology needed to give the formulas for the measures $\nu_{\Gamma,n,\{1\}}$ and $\nu_{ \Gamma, n,\overline{\gamma}}$. We fix a finite group $\Gamma$, a positive integer $n$ prime to $\abs{\Gamma}$, and a multiset $U$ of conjugacy classes of $\gamma$, which will specialize to $U= \{1\}$ in the function field case or $U= \overline{\gamma}$ in the number field case. Associated to the multiset $U$ we have a function on $\Gamma$-groups
$$G^{\cdot U} :=\frac{\prod_{\gamma\in U} |G^\gamma|}{|G^\Gamma|}.$$

Let $\cL$ be a level of the category of finite $\Gamma$-groups.
Let  $\mathbf H$ be a finite $n$-oriented $\Gamma$-group whose underlying $\Gamma$-group is in $\cL$.

A finite simple abelian $H \rtimes \Gamma$-group $V$ is a group of the form $\mathbb F_\ell^{ d}$, for some prime $\ell$, with an irreducible action of $H \rtimes \Gamma$. We say $V$ is \emph{admissible} if $V \rtimes H$, viewed as a $\Gamma$-group by its natural embedding in $V \rtimes (H\rtimes \Gamma) = (V \rtimes H) \rtimes \Gamma$, is in $\cL$.

Let $V_1, \dots V_r$ be {representatives of the isomorphism classes of} admissible finite simple abelian $ H \rtimes \Gamma$-groups.

We say a finite simple nonabelian $[H\rtimes\Gamma ]$-group $N$, of order prime to $\abs{\Gamma}$ is \emph{admissible} if for one, equivalently, every, lift of $\Gamma$ to $H\rtimes \Gamma \times_{\Out(N) } \operatorname{Aut}( N)$, the complementary normal subgroup $H \times_{\Out(N)} \operatorname{Aut}(N)$ is in $\cL$.  (By Schur-Zassenhaus, different lifts are conjugate and thus the complements are isomorphic as $\Gamma$-groups.)

Let $N_1,\dots, N_s$ be {representatives of the isomorphism classes of} admissible finite simple nonabelian $[H\rtimes\Gamma ]$-groups.
Note our definitions of admissible implicitly depend on $\cL$ and $H$.

Our formula for $\nu ( \{ \mathbf X \mid \mathbf X^{\cL} \cong \bH \}) $ will be a product of local factors corresponding to $V_1,\dots, V_r$ and $N_1,\dots, N_s$.

To define these local factors, we have the following invariants. For $V$ an admissible finite simple abelian $H \rtimes \Gamma$-group, let $\kappa_V$ be the ($H \rtimes \Gamma$-equivariant) endomorphisms of $V$, necessarily a finite field. Let $q_V= \abs{\kappa_V}$. Let $H^2( H \rtimes \Gamma, V)^{\cL}$ be the set of cohomology classes whose associated extension of $H$ by $V$ as a $\Gamma$-group lies in $\cL$ (the association is described explicitly in Lemma~\ref{L:diffext}).  Note that $H^2( H \rtimes \Gamma, V)^{\cL}$ is a $\kappa_V$-vector space (Lemma~\ref{L:basicH2}). 

The orientation $s_H \in H_3( H, \mathbb Z/n)$ induces a map $H^3( H, \mathbb Z/n) \to \mathbb Z/n$, which we also denote by $s_H$ (see Lemma~\ref{L:homdual}).  
We also write $s_H$ for the composite of this map with the pullback isomorphism $H^3(H\rtimes\Gamma,\Z/n)\ra H^3( H, \mathbb Z/n)^\Gamma$ (Lemma~\ref{L:GammaH}), giving a map $s_H: H^3(H\rtimes\Gamma,\Z/n) \ra \Z/n$.

If $\kappa_V$ has characteristic dividing $n$ and $V$ is a nontrivial representation, let $ H^2( H \rtimes \Gamma, V)^{\cL, s_H} $ be the set of cohomology classes $\alpha\in H^2( H \rtimes \Gamma, V)^{\cL}$ such that $s_H (\alpha \cup \beta)=0$ for all $\beta$ in $H^1( H \rtimes \Gamma, \Hom(V,\Z/n))\isom H^1( H \rtimes \Gamma, V^\vee)$.

If $V$ is a trivial representation $\mathbb F_\ell$ with characteristic dividing $n$, let $ H^2( H \rtimes \Gamma, V)^{\cL, s_H} $ be the set of cohomology classes $\alpha\in H^2( H \rtimes \Gamma, V)^{\cL}$ such that $s_H (\alpha \cup \beta)=0$ for all $\beta$ in $H^1( H \rtimes \Gamma, \Hom(V,\Z/n))\isom H^1( H \rtimes \Gamma, V^\vee)$ and $s_H (\mathcal B (\alpha))=0$ for $\mathcal B \colon H^2( H\rtimes \Gamma, \mathbb F_\ell) \to H^3( H \rtimes \Gamma, \mathbb Z/n)$ the connecting homomorphism associated to the exact sequence $0\to \mathbb Z/n \to \mathbb Z/n\ell \to \mathbb Z/\ell \to 0$.

If $\kappa_V$ has characteristic not dividing $n$, let  $ H^2( H \rtimes \Gamma, V)^{\cL, s_H}= H^2( H \rtimes \Gamma, V)^{\cL}$.

Let\[z_V =\dim_{\kappa_V}  H^2( H \rtimes \Gamma, V)^{\cL, s_H}\] and \[h_V = \dim_{\kappa_V} H^1( H\rtimes \Gamma, V^\vee) - \dim_{\kappa_V} H^1(H\rtimes \Gamma, V).\]Let $u_V$ be such that \[q_V^{u_V} =\frac{ \abs{V^\Gamma}}{\abs{V ^{H\rtimes \Gamma}} } V^{ \cdot U}.\]

For $V$ nontrivial, of characteristic prime to $n$, let \[ w_{V} = \prod_{j=1}^\infty \left(1-q_V^{-j-u_V} \frac{ | H^2( H \rtimes \Gamma, V)^{\cL} |}{  |H^1( H \rtimes \Gamma, V)|   } 
 \right).    \]

For $V$ trivial, let \[ w_V = \prod_{k=0}^{z_V-1} (1-q_V^{k-u_V}).\]

If the characteristic of $V$ divides $n$ and $V$ is admissible but $V^\vee$ is not, let \[ w_{V} = \prod_{j=1}^\infty \left(1-q_V^{h_V +z_V -j - u_V } \right). \] 

If the characteristic of $V$ divides $n$ and $V^\vee$ is admissible but not isomorphic to $V$,  then let  \[ w_V = \sqrt{ \frac{ (q_V)_\infty (q_V)_{2u_V}}{ (q_V)_{u_V-h_V-z_V} (q_V)_{u_V+h_V - z_{V^\vee}}}}\] if $-u_V + z_{V^\vee} \leq h_V \leq u_V -z_V$ and $0$ otherwise. (Using that $V$ is a non-trivial $H\rtimes \Gamma$-representation of characteristic not dividing $|\Gamma|$, it is straightforward to check that $h_{V^\vee} = -h_V$ and $u_V=u_{V^\vee}$,
 so that $w_{V}=w_{V^\vee}$.)

If the characteristic of $V$ divides $n$ and $V$ is self-dual, non-trivial, and non-anomalous, let
 \[ w_{V}   =
 \prod _{k=0}^{\infty }(1+q_V^{-k-\frac{\epsilon_{V}+1}{2} -u_V})^{-1}
 \prod _{k=0}^{z_V-1 }{(1-q_V^{k -u_V})}.
 \] 
 
 If the characteristic of $V$ divides $n$ and $V$ is anomalous, there is a particular class $\omega_*^{-1}(\Phi)\in H^2(\Sp(V) ,V)$ defined in \S\ref{s-notation}, which pulls back to a class in $H^2( H \rtimes \Gamma, V)$. If $\omega_*^{-1}(\Phi)\in H^2(H\rtimes \Gamma,V)^{\cL,s_H}$, 
$$
w_{V}=
\prod _{k=0}^{\infty }(1+q_V^{-k-u_V})^{-1}
\prod_{k=0}^{z_V-1}(1-q_V^{k-u_V}) 
$$
If $\omega_*^{-1}(\Phi)\not\in H^2(H\rtimes \Gamma,V)^{\cL,s_H}$,
 \[ w_{V}   =
 \prod _{k=0}^{\infty }(1+q_V^{-k-1 -u_V})^{-1}
 \prod _{k=0}^{z_V-1 }{(1-q_V^{k -u_V})}.
 \] 
 One can observe that we could define $\epsilon_V$ to be $-1$ in case $\omega_*^{-1}(\Phi)\in H^2(H\rtimes \Gamma,V)^{\cL,s_H}$, and $1$ in case $\omega_*^{-1}(\Phi)\not\in H^2(H\rtimes \Gamma,V)^{\cL,s_H}$, which would make the same formula valid in the anomalous and non-anomalous self-dual cases, but we have not made this definition as that would make $\epsilon_V$ depend on $\cL$ and $s_H$.

For $N$ an admissible nonabelian group,  let $L_N$ be the number of lifts of $H \rtimes \Gamma \to \Out(N)$ to $H \rtimes \Gamma \to \Aut(N)$. Given a choice of lift, $N$ has the structure of a $\Gamma$-group, and so we can define $N^{\cdot U}$. However, since $|\Gamma|$ is relatively prime to $|N|$, any two lifts of $\Gamma$ from $\Out(N)$ to $\Aut(N)$ are conjugates by an inner automorphism (which follows from the Schur-Zassenhaus Theorem applied to $\Aut(N)\times_{\Out(N)} \Gamma$).
Thus the $\Gamma$-groups produced by the two lifts are isomorphic, and in fact $N^{\cdot U}$ depends only on the structure of $N$ as an $[H \rtimes \Gamma]$-group.  Thus, we will use the notation $N^{\cdot U}$ without choosing a lift.

If the natural map \[\delta\colon H^2 ( N, \mathbb Z/n)^{H \rtimes \Gamma} \to H^3 ( H \rtimes \Gamma, \mathbb Z/n) \] from the spectral sequence computing $H^3( \Aut(N)  \times_{ \Out(N)} (H \rtimes \Gamma) , \mathbb Z/n)$, composed with the pullback
$H^3 ( H \rtimes \Gamma, \mathbb Z/n)\ra H^3 ( H, \mathbb Z/n)$ and then 
 $s_H$,
 is nontrivial, then we define $w_{N} =1$. If the map $H\rtimes \Gamma \to \operatorname{Out}(N)$ is trivial, then we define $w_N =1$. In any other case, we define  \[w_{N}  = e^{ - \frac{  L_{N} \left| H^2(N,\mathbb Z/n)^{H \rtimes \Gamma} \right|  } { |N|   |Z_{ \Out(N)  } ( H \rtimes \Gamma)|  N^{\cdot U} }} .\] 

\begin{example} The map $\delta$ is trivial for most small finite simple groups, but one can check it can be nontrivial for $N= A_6$ or $N= PSU_3( \mathbb F_8)$. The first example was pointed out to us by Pham Huu Tiep, and we found the second using the Atlas of Finite Simple Groups. They can be recognized in the Atlas from the fact that some central extension of $N$ does not extend to a central extension of some particular $S$ with $N \subset S \subseteq \operatorname{Aut}(N)$. This reflects the fact that the map $H^{2} ( S, \mathbb Z/n) \to H^2(N, \mathbb Z/n)$ is not surjective for some $n$, which implies the relevant class in $H^2( N, \mathbb Z/n)$ must not lie in the kernel of every differential in the spectral sequence.  \end{example}

With these definitions we have the following.

\begin{theorem}\label{T:intro-measure-general} Fix a finite group $\Gamma$, a positive integer $n$ prime to $\abs{\Gamma}$, and a multiset $U$ of conjugacy classes of $\Gamma$.  Assume that each nonzero representation of $\Gamma$ of characteristic dividing $n$ contains a nonzero vector fixed by some element of $U$.

For a finite $n$-oriented $\Gamma$-group $\bG$, let \[M_{\bG}= \begin{cases} \frac{ |H^2(G\rtimes \Gamma, \mathbb Z/n)| }{ |H^3(G\rtimes \Gamma,\mathbb Z/n) G^{\cdot U}} & \textrm{if } G_\Gamma = 1\\ 0 & \textrm{if }G_\Gamma \neq 1\end{cases} .\]

There exists a unique measure $\nu_{\Gamma, n,U} $  on $\mathcal P_{\Gamma,n}$ such that for every finite $n$-oriented $\Gamma$-group $\bG$ we have $$
\int_{X\in \mathcal{P}_{\Gamma,n} } \Sur(X,\bG) d\nu_{\Gamma, n, U} =M_{\bG}. $$ 

For any level $\cL$ of the category of finite $\Gamma$-groups, finite $n$-oriented $\Gamma$-group $\bH$ whose underlying $\Gamma$-group is in $\cL$, we have

\[\nu_{\Gamma, n, U} (\{X|X^\cL\isom  \bH\}) =  \frac{ M_{\bH}}{\abs{\Aut(\bH)}} \prod_{i=1}^r w_{V_i} \prod_{i=1}^s w_{N_i}\]
where $V_1,\dots, V_r$ are representatives of the isomorphism classes of admissible finite simple abelian $H \rtimes \Gamma$-groups and $N_1,\dots, N_s$ are representatives of the isomorphism classes of admissible  finite simple non-abelian $[H\rtimes \Gamma]$-groups.  \end{theorem}

\subsection{Outline of the paper}
In Section~\ref{S:AV}, we define precisely the Artin-Verdier fundamental class for an extension $L/K$, and prove basic properties about this class.
In Section~\ref{s-prior-work}, we show how previously defined additional structure on Galois groups of unramified extensions,  including the bilinearly enhanced structure of Lipnowski-Sawin-Tsimerman \cite{Lipnowski2020} and the lifting invariant of  \cite{Liu2022},  is determined by the Artin-Verdier fundamental class.
We also work out explicit cases of our predictions for the $2$-class tower group of cyclic degree $3$ extensions, and compare to the data, results, and predictions of Boston and Bush \cite{Boston2021a}.
In Section~\ref{S:evidence}, we prove Theorem~\ref{intro-ff-theorem}, our main result on function field moments.  In Section~\ref{S:GroupTheory}, we collect some standard results in group theory that we use repeatedly in our arguments.  In Section~\ref{S:LHS}, we analyze the $E^3_{0,2}$ term in the Lyndon-Hochschild-Serre spectral sequence, especially for an extension of a group $H$ by a semisimple abelian $H$-group.  That is a key ingredient in our argument for Theorems~\ref{intro-bb}, \ref{T:intro-measure-ff}, \ref{T:intro-measure-nf}, and \ref{T:intro-measure-general} in Section~\ref{S:momprob}, where we work out explicitly the distributions of oriented groups determined by the moments we have found.  
In Section~\ref{S:NE}, we prove Theorem~\ref{T:zero},  giving number field evidence for our conjecture in the form of non-existence proofs when our conjecture predicts probability zero.

\subsection*{Acknowledgements} 
We would like to thank Johan de Jong and Yuan Liu for helpful conversations, Aaron Landesman, for helpful comments on an earlier version of this manuscript, and Jiuya Wang for help with corecting Conjecture  \ref{intro-moments-conjecture}.
The first author was supported by a Clay Research Fellowship, NSF DMS-2101491 and DMS-2502029 and a Sloan Fellowship. The second author was partially supported by a Packard Fellowship for Science and Engineering,  NSF  DMS-2052036 and DMS-2140043,  the Radcliffe Institute for Advanced Study at Harvard University, and a MacArthur Fellowship.

\tableofcontents

\section{The Artin-Verdier trace and fundamental class}\label{S:AV}

In this section we establish some of the basic properties of the components of the Artin-Verdier trace map and fundamental class.

Let $K$ be a global field, of characteristic not dividing a positive integer $n$.  If $K$ is a number field, let $X_K=\Spec \O_K$.  If $K$ is a function field with field of constants $k$, let $X_K$ be the unique connected smooth complete curve over $k$ having $K$ as its function field.  
There is a natural isomorphism on the compactly supported \'etale cohomology \cite[Chapter II, Section 3]{Milne2006}
\begin{equation}\label{E:basic-AV}
H^3_c(X_K,\G_m)\isom \Q/\Z,
\end{equation}
that is in particular invariant under automorphisms of $K$.  
The map of $n$th roots of unity $\mu_n \ra \G_m$ gives a map $H^3_c(X_K,\mu_n)\ra H^3_c(X_K,\G_m) $, and taking the composite with the map above we obtain a map $D\colon H^3_c(X_K,\mu_n)\ra \frac{1}{n}\Z/\Z\stackrel{\times n}{\ra} \Z/n$.  

If $K$ contains the $n$th roots of unity then,  we fix throughout the paper a generator $\xi$ for the $n$th roots of unity of $K$, which in particular gives a morphism $H^3_c ( X_K, \mathbb Z/n)\to H^3_c(X_K,\mu_n)$ that we will use.

Let $L/K$ be a Galois extension that is unramified everywhere, with $H=\Gal(L/K)$.
We will show in Remark~\ref{R:defineor} that there is a natural map $H^3(H,\Z/n) \ra H^3(X_K,\Z/n)$. 
Moreover, if  $L/K$ is split completely over all real places of $K$, then we show  in Lemma~\ref{L:compactsupports} that we can refine this to a map to cohomology with compact supports
$$
H^3(H,\Z/n)\ra H^3_c(X_K,\Z/n). 
$$

Then, associated to this extension $L/K$,
we have a composed map \[ H^3 ( H, \mathbb Z/n) \to H^3_c ( X_K, \mathbb Z/n)\to H^3_c(X_K ,\mu_n) \stackrel{D}{ \to} \mathbb Z/n \] where each map is  discussed above.
We call this map
$\AV_{L/K} :H^3 ( H, \mathbb Z/n)\ra \mathbb Z/n$
 the \emph{Artin-Verdier trace} for $L$ over $K$, 
 and the corresponding element $\av_{L/K}$ (via the duality in Lemma~\ref{L:homdual}) in $H_3( H, \mathbb Z/n)$ the \emph{Artin-Verdier fundamental class} for $L$ over $K$.

If $K$ is a Galois extension of some other global field $F$ with $\Gamma=\Gal(K/F)$, such that $L/F$ is Galois, we will show in Lemma~\ref{L:AVinvariant} that the Artin-Verdier fundamental class for $L/K$ is in $H_3(H,\Z/n)^{\Gal(K/F)}$.

\begin{lemma}\label{L:BHpullback}
Let $X$ be a connected scheme, quasi-projective over an affine scheme.
Let $Y/X$ be a finite  \'etale morphism with a group of automorphisms $H$ of $Y$ over $X$ that acts simply transitively on each geometric fiber.   Let $p$ be a non-negative integer. 
Let $A$ be an $H$-module, viewed as a representation of $\pi_1(X)$ and thus an \'etale sheaf on $X$ via the covering $Y/X$ (and a choice of a basepoint in $Y$ over a basepoint in $X$, which gives a map $\pi_1(X)\ra H$).
Then there is a homomorphism
$$\phi_{Y/X}=\phi_{Y/X,A}: H^p(H,A) \ra H^p(X,A).$$ 
If we have $f:Z\ra X$,  such that $Z$ is a connected scheme, quasi-projective over an affine scheme,  then 
$f^* \phi_{Y/X}=\phi_{Z\times_X Y/Z}$ (with compatible choices of basepoints).
If $A'$ is an $H$-module and $f:A\ra A'$ is an $H$-module homomorphism, then we have $f_* \phi_{Y/X,A}= \phi_{Y/X,A'}f_*$.
If $A''$ is another $H$-module, and we let $B$ denote the connecting homomorphism for an exact sequence of $H$-modules
$$
0\ra A \ra A' \ra A''\ra 0 
$$
in either the cohomology of $H$ or $X$, then we have $\phi_{Y/X,A}B=B\phi_{Y/X,A''}$.
If $a\in H^*(H,A)$ and $b\in H^*(H,A')$, then $\phi_{Y/X}(a \cup b)= \phi_{Y/X}(a) \cup \phi_{Y/X}(b)\in H^*(X,A\tesnor A')$.

If $Y'/X$ is another finite  \'etale morphism with a group of automorphisms $H'$ of $Y'$ over $X$ that acts simply transitively on each geometric fiber, with the map $Y'\ra X$ factoring through $Y\ra X$ (and $Y'\ra Y$ basepoint preserving), and with associated quotient $q: H'\ra H$ inducing $q^* :H^p(H,A)\ra H^p(H',A)$, then $\phi_{Y'/X} q^*=\phi_{Y/X}$. 
\end{lemma}
\begin{proof}
We interpret  $H^*(H, A)$ as the cohomology of the standard complex $C^*(H,A)$ of homogeneous cocycles, i.e.  $H$-equivariant functions from $H^{p+1}$ to $A$, where the action of $H$ on $H^p$ is given by left multiplication on each coordinate. 
We will map the complex $C^*(H,A)$ to the  (\'etale) \v{C}ech complex  $C^*(Y/X,A)$ for the cover $Y/X$.  

Choose a basepoint $x\in X$, and let $F_Y$ be the fiber over $x$ in $Y$.  Pick a basepoint $y\in F_Y$.  We then have a bijection
$F_Y \ra H$ taking $h^{-1}y\mapsto h$ for $h\in H$. The choice of basepoint $y\in F_Y$ also gives an associated group homomorphism
$\pi_1(X)\ra H$ such that for $g\in \pi_1(X)$, if $gy=h_g^{-1}y$ for $h_g\in H$, the $g\mapsto h_g$.  

By the fact that $\pi_1(X)$ is an automorphism of the fiber functor, we have for $g\in \pi_1(X)$, that
$g(h^{-1}y)= (h_g h)^{-1} y$ for all $h\in H$.  
The sections of $C^*(Y/X,A)$ are given by maps of $\pi_1(X)$-sets from $F_Y^{p+1} \ra A$, which under the identification of $F_Y$ with $H$ above are precisely maps of $\pi_1(X)$-sets $H^{p+1}\ra A$, with $\pi_1(X)$ acting through the map $\pi_1(X)\ra H$ and then left multiplication on each component.  This gives a map $C^*(H,A)\ra C^*(Y/X,A)$ taking a map $H^{p+1}\ra A$ to the same map $H^{p+1}\ra A$.
Moreover this map manifestly respects the boundary maps in these complexes.
Thus we obtain a homomorphism from
$H^p(H,A)$ to the cohomology of the complex $C^p(Y/X,A)$, which maps to the \v{C}ech cohomology, and thus by the assumption on $X$, the \'etale cohomology
$H^p(X,A)$.

The compatibility on pullbacks can be checked directly from the construction.
The compatibilities with a change in coefficients and the connecting homomorphisms are immediate from the construction.
Also immediate from the construction is the compatibility of the cup product on group cohomology and the cup product on \v{C}ech cochains, which by
\cite[Corollary 3.10]{Swan1999} agrees with the cup product in \'{e}tale cohomology.
The compatibility between different covers is also immediate from the construction.

(In the language of stacks, the covering $Y\ra X$ defines a map $ X \to BH$, and the group cohomology of $H$ is the cohomology of $BH$, with the standard group cohomology cochain complex corresponding to the Cech complex for the covering $\operatorname{pt}\to BH$.  In this language $\phi_{Y/X}$ is just the pullback along $X \to BH$.)
\end{proof}

\begin{remark}\label{R:defineor}
If $L/K$ is a Galois extension that is unramified everywhere with $H=\Gal(L/K)$, then since $H^3(H,\Z/n)$ is the direct limit of $H^3(H_i,\Z/n)$
over finite quotients $H_i$ of $H$,  from Lemma~\ref{L:BHpullback} and the compatibility between different quotients, we get a map
$$
H^3(H,\Z/n) \ra H^3(X_K,\Z/n).
$$
More generally, if $V$ is a finite $H$-module, then $H^3(H,V)$ is the direct limit of $H^3(H_i,V)$ over finite quotients $H_i$ of $H$ such that the action of $H$ on $V$ factors through $H_i$, and so we have
$$
H^3(H,V) \ra H^3(X_K,V).
$$
\end{remark}

When $X_K=\Spec \O_K$, for a number field $K$, 
and the \'etale morphsim is from an extension of fields that is split completely at all real places,
we can refine the map above to a map to cohomology with compact supports.

\begin{lemma}\label{L:compactsupports} Let $L/K$ be a  Galois unramified extension of a number field $K$, split completely over all real places of  $K$.
Let $X_K=\Spec \O_K$.
Let $A$ be an $\Gal(L/K)$-module.
A choice of algebraic closure $\bar{K}$ and a choice of an embedding $L\ra \bar{K}$, gives a homomorphism $\pi_1(X)\ra \Gal(L/K)$,
and thus $A$ is also a representation of $\pi_1(X_K)$ and thus an \'etale sheaf on $X_K$.

For all $p\geq2$, there is a homomorphism
$$
\phi_{L/K}^c: H^p(\Gal(L/K), A) \ra H^p_c(X_K,A).
$$ whose composition with $H^p_c(X_K,A)\to H^p(X_K,A)$ is the usual pullback map $\phi_{\Spec \O_L/\Spec \O_K}$ of Lemma~\ref{L:BHpullback} (or Remark~\ref{R:defineor}). 

This map has the following compatibilities: For an inclusion of fields $K\to K'$ such that $L' = L \otimes_K K'$ is a field, this map forms a commutative square 
\[\begin{tikzcd}  H^p(\Gal(L/K), A)
\arrow[r,"\phi_{L/K}^c"] &H^p_c(X_K,A)
\\ H^p(\Gal(L'/K'), A)  \arrow[u] \arrow[r,"\phi_{L'/K'}^c"] &
 H^p_c(X_{K'},A)  \arrow[u]\end{tikzcd} \]
with the natural pull-backs.
For $L'$ a subfield of $L$ that is Galois over $K$, this map forms a commutative triangle with the usual pull-back map $H^p(\Gal(L'/K), A) \to H^p(\Gal(L/K), A) $
\[\begin{tikzcd}  
H^p(\Gal(L'/K), A)  \arrow[d] \arrow[r,"\phi_{L'/K}^c"] &
 H^p_c(X_{K},A) 
\\ 
 H^p(\Gal(L/K), A).
\arrow[ur,"\phi_{L/K}^c"] 
 \end{tikzcd} \]
\end{lemma}

\begin{proof}
As in Remark~\ref{R:defineor}, we can immediately reduce to showing the case when $L/K$ is finite, as the final statement of the lemma gives the compatibility in towers $L/L'/K$ that we need to pass to the direct limit.

For $p\geq 1$, the compactly supported cohomology $ H^p_c(X_K,A)$ is defined \cite[Section II:2]{Milne2006} as the cohomology of the translated mapping cone of the \v{C}ech complex $C^*(X_K,A)$ representing the cohomology of $X_K$ to the sum over infinite places $v$ of $K$ 
of the
 \v{C}ech complex $C^*(\Spec  K_v,A)$.
 
Let $f$ be the map $\coprod_{v\mid\infty} \Spec K_v \ra X_K$ and let $L_v=K_v\tesnor_K L$.
A class in the translated mapping cone described above is a pair $(a,b)\in C^*(X_K,A) \times \prod_{v\mid\infty} C^{*-1}(\Spec K_v,A)$, and it is a cocycle in particular if $a$ is a cocycle
in $C^*(X_K,A)$ such that $f^*a=-d(b),$ where $d$ is the differential.

Let $A^*:=C^*(X_L/X_K,A)$ and $B^*=\coprod_{v\mid\infty} C^{*}(\Spec L_v/\Spec K_v,A)$.
First,  we will find a homomorphism from the subgroup of cocycles in $C^*(\Gal(L/K),A)$ to the subgroup of  cocycles in the mapping cone of
the \v{C}ech complexes for coverings
 $f^*:  A^*
\ra B^*$.

From the proof of Lemma~\ref{L:BHpullback}, we have a homomorphism of complexes  $\Phi: C^*(\Gal(L/K),A) \ra A^*$.
Since $L/K$ is split completely at $v$,  the complex $B^p$ is exact at $p\geq 1$.
From this it follows, for $p\geq 2$, that given a cocycle $f^*a\in A^p,$ there is a  $b\in B^{p-1}$, unique up to coboundaries,  such that
$f^*a=-d(b)$.  So for $p\geq 2$,  the map $\alpha \mapsto  (\Phi \alpha,b)$,  such that $f^*\Phi \alpha=-d(b)$,
gives us a homorphism from the subgroup of 
cocycles in  $C^p(\Gal(L/K),A)$ to the quotient of the group of elements $(a,b)\in A^p\times B^{p-1}$ such that $f^*a=-d(b)$
by the group $0\times dB^{p-2}$.  
We can check that a coboundary $d\alpha$ is mapped to
the coboundary $d(-\Phi \alpha,0)=(d\Phi \alpha, -f^* \Phi \alpha)$,  using the definition of the boundary map of a mapping cone.

The maps from \v{C}ech complexes for coverings to  \v{C}ech complexes then give us 
a map of complexes from the mapping cone for 
$f^*:  A^*\ra B^*$ to the mapping cone for 
$f^*:  C^*(X_K,A)
\ra \prod_{v\mid\infty} C^{*}(\Spec K_v,A)$.
Combining with the above, we obtain a homomorphism  $$H^p(\Gal(L/K), A) \ra H^p_c(X_K,A).$$
It is clear from the construction that the composition with $H^p_c(X_K,A)\to H^p(X_K,A)$ gives the map of Lemma~\ref{L:BHpullback}.
The compatibilities both follow in a straightforward way from the compatibility of $\Phi$ with pullback.
\end{proof}

\begin{lemma}\label{L:AVinvariant}
Let $n$ be a positive integer.
Let $L/F$ be a finite Galois extension of global fields, with $K$ an intermediate field such that $K/F$ is Galois and $K$ contains $n$th roots of unity.  
The Artin-Verdier trace and the Artin-Verdier fundamental class for $L/K$ are $\Gal(K/F)$-invariant.
\end{lemma}
\begin{proof}
Let $\Gamma=\Gal(K/F)$.
For any lift $\tilde{\gamma}\in \Gal(L/F)$ of a $\gamma\in\Gamma,$ we have $\tilde{\gamma}:L\ra L$, and $\tilde{\gamma}:K \ra K$, and thus by Lemma~\ref{L:compactsupports} we have that the map $H^3 ( H, \mathbb Z/n) \to H^3_c ( X_K, \mathbb Z/n)$ 
associated to $L/K$
is a $\Gamma$-equivariant map
for the natural actions of $\Gamma$ on the left and right (where $\Gamma$ acts on $H^3 ( H, \mathbb Z/n)$ by lift to $\Gal(L/F)$ and conjugation on $H$, not depending on the lift since $H$'s self-conjugation acts trivially on $H^3 ( H, \mathbb Z/n)$).
 The map $H^3_c ( X_K, \mathbb Z/n)\to H^3_c(X_K,\mu_n)$ is $\Gamma$-equivariant, and $H^3_c(X_K ,\mu_n) \stackrel{D}{ \to} \mathbb Z/n$ is $\Gamma$-invariant.  We conclude that the Artin-Verdier trace $H^3(\Gal(L/K),\Z/n)\ra \Z/n$ and the Artin-Verdier fundamental class  for $L/K$ are $\Gamma$-invariant.
\end{proof}

\section{Relations with prior invariants}\label{s-prior-work}

This section will explain in detail how our work relates to  prior papers. 
We will explain in Section~\ref{SS:LST} how our Artin-Verdier fundamental class refines the bilinearly enhanced group structure that Lipnowski,  Tsimerman, and the first author \cite{Lipnowski2020} defined in the case that $|\Gamma|=2$.  In Section~\ref{SS:Liu}, we will explain how our Artin-Verdier fundamental class refines lifting invariants defined by Liu \cite{Liu2022} and the second author \cite{Wood2019,Wood2021}.  In Section~\ref{ss-bb}, we will explain how our conjectures relate to some conjectures made by Boston and Bush \cite{Boston2021a} in the case that $\Gamma$ is the cyclic group of order $3$.
The relationship with \cite{Liu2022} will be used in subsequent proofs, as we will cite some results from \cite{Liu2022}, while the others are used only to provide context and numerical evidence.

\subsection{Bilinearly enhanced groups, after Lipnowski-Sawin-Tsimerman}\label{SS:LST} In the paper \cite{Lipnowski2020},  two invariants are defined on the class group of a number field,  called the $\omega$ and $\psi$-invariants, which can be calculated using the Artin-Verdier trace.  In this subsection, we first review the definitions of the $\omega$ and $\psi$ invariants of \cite{Lipnowski2020}, then explain how to derive these invariants from the Artin-Verdier trace, and finally prove some auxiliary results about the information captured by the Artin-Verdier trace but not the $\omega$ and $\psi$ invariants.

Let $K$ be a number field and $n$ be the order of the group of roots of unity in $K$. Let $\ell$ be an odd prime and let $v\geq 1$ be the $\ell$-adic valuation of $n$, so that $K$ contains the $\ell^v$th roots of unity. 
In particular, this implies that $K$ has no real places.
The paper \cite{Lipnowski2020} defines two invariants, the $\omega$ and $\psi$ invariant, on the $\ell$-part of the class group of $K$. We will use in this subsection the same notation used in that paper, except that their $n$ is our $v$.

In this section, we will consider flat cohomology $\Hf^i(\Spec \O_K,F)$ along with \'etale cohomology, as the definitions of \cite{Lipnowski2020} use flat cohomology.  Recall that when $F$ is a smooth, quasi-projective, commutative group scheme, the flat and \'etale cohomology groups coincide \cite[III, Theorem 3.9]{Milne1980}.  We will use this in particular when $F=\Z/k$ or $\G_m$. Since $K$ has no real places,  cohomology coincides with compactly support cohomology in both the flat and \'etale cases.
We write $\tA$ for the map (isomorphism) $H^3(\Spec \O_K,\G_m) \ra \Q/\Z$ \cite[Chapter II, Section 3]{Milne2006}, and with a slight abuse of notation, also
use $\tA$ to denote the maps $H^3(\Spec \O_K,\mu_k) \ra \Q/\Z$ and $\Hf^3(\Spec \O_K,\mu_k) \ra \Q/\Z$ obtained by first mapping to cohomology with coefficients in $\G_m$ and then applying the above map.

The $\omega$-invariant on the class group $\mathrm{Cl}(K)$ of a number field $K$ containing the $\ell^v$'th roots of unity is an element of $(\wedge^2 \mathrm{Cl}(K))[\ell^v]$. To define this invariant \cite[Definition 4.5]{Lipnowski2020}, one uses \cite[Lemma 4.2 and Corollary 4.3]{Lipnowski2020}, which show that such an element is equivalent to a tuple of, for each $m$, a symplectic bilinear form $\omega_{m,K} \colon \mathrm{Cl}(K)^\vee [\ell^m] \times \mathrm{Cl}(K)^\vee [\ell^m] \to \mathbb Q/\mathbb Z [\ell^v]$, such that for all $a \in \mathrm{Cl}(K)^\vee[\ell^m], b\in \mathrm{Cl}(K)^\vee [\ell^{m+1}]$ we have $\omega_{m,K}(a,\ell b) = \omega_{m+1,K}(a,b)$. 
Let $\tilde{K}$ be the Hilbert class field of $K$, and we use class field theory to write $\Cl(K)=\Gal(\tilde{K}/K)$.  
For any integer $k$, we write $\phi=\phi_{\tilde{K}/K}: H^i(\Cl(K),\Z/k )
\ra H^i(\Spec \O_K,\Z/k)$ for the map from
Lemma~\ref{L:BHpullback}.
Recall our fixed generator $\xi$ of the $n$th roots of unity in $K$.
Let $\zeta_m \in \Hf^1( \operatorname{Spec} \mathcal O_K, \mu_{\ell^m})$ correspond to the $\mu_{\ell^m}$ torsor consisting of the $\ell^m$th roots of  the generator $\xi^{n/\ell^v}$ of $\mu_{\ell^v}$ (as in \cite[III, Proposition 4.6]{Milne1980}). 
One constructs the pairing $\omega_{m,K}$ by the formula
\[ \omega_{m,K}(a,b) = -\frac{1}{2} \tA (\zeta_m \cup \phi(a) \cup \phi(b)) \] 
for $a,b \in \mathrm{Cl}(K)^\vee [\ell^m]=H^1(\Cl(K),\Z/\ell^m) $.
 Since $\zeta_m\in \Hf^1$,  the cup products are in flat cohomology.
  In \cite[Lemma 4.4]{Lipnowski2020}, these are checked to satisfy the conditions required to define an element of $(\wedge^2 \mathrm{Cl}(K))[\ell^v]$.

The $\psi$-invariant is a homomorphism $\psi\colon \mathrm{Cl}(K)^\vee [\ell^v] \to \mathrm{Cl}(K) [\ell^v] $. This data is equivalent by duality to a pairing $\langle , \rangle \colon \mathrm{Cl}(K)^\vee [\ell^v] \times \mathrm{Cl}(K)^\vee  \to \mathbb Q/\mathbb Z$ given by the formula $\langle a,b\rangle = b ( \psi(a))$. 
Again, we use class field theory to write $\Cl(K)=\Gal(\tilde{K}/K)$.  
We define the $\psi$ map  \cite[Definition 4.1]{Lipnowski2020} as the composition of $\phi:\mathrm{Cl}(K)^\vee [\ell^v]=H^1(\Cl(K),\Z/\ell^v)\ra
  H^1( \operatorname{Spec} \mathcal O_K, \mathbb Z/\ell^v)$, the cup product $H^1( \operatorname{Spec} \mathcal O_K, \mathbb Z/\ell^v)\to H^1( \operatorname{Spec} \mathcal O_K, \mu_{\ell^v})$ with the generator $\xi^{n/\ell^v}$ of $\mu_{\ell^v}$, 
and the map $H^1( \operatorname{Spec} \mathcal O_K, \mu _{\ell^v})\to H^1(\operatorname{Spec} \mathcal O_K, \mathbb G_m ) = \mathrm{Cl}(K)$ from the map $\mu_{\ell^v}\ra\G_m$. 

Next, we will give formulas for the bilinear forms $\omega_{m,K}$ and $\langle, \rangle$ in terms of the Artin-Verdier trace and the Bockstein homomorphism $B \colon H^i (\mathrm{Cl}(K) , \mathbb Z/\ell^m) \to H^{i+1} (\mathrm{Cl}(K), \mathbb Z/\ell^v)$ arising from the short exact sequence $\mathbb Z/\ell^v\to \mathbb Z/\ell^{v+m} \to \mathbb Z/\ell^m $. 
 Let $s_{\mathrm{Cl}(K)} \in H_3( \mathrm{Cl}(K), \mathbb Z/n)$ be the Artin-Verdier fundamental class for $\tilde{K}/K$.
Let $s_{\mathrm{Cl}(K)} : H^3(\mathrm{Cl}(K),\Z/n)\ra \Z/n$ be the usual induced map, and let
$s^\ell_{\mathrm{Cl}(K)} : H^3(\mathrm{Cl}(K),\Z/\ell^v)\ra \Q/\Z[\ell^v]$ be the composition of the
 map 
$H^3(\mathrm{Cl}(K),\Z/\ell^v)\ra H^3(\mathrm{Cl}(K),\Z/n)$ by multiplication by $n/\ell^v$, the 
   map $s_{\mathrm{Cl}(K)}$, and $\Z/n \ra \Q/\Z$ given by division by $n$.

\begin{theorem}\label{lst-tau} Let $\ell$ be an odd prime and $K$ a number field containing the $\ell^v$th roots of unity, and not the $\ell^{v+1}$th roots of unity, for some $v\geq 1$. 
For all positive integers $m$ and  $a,b \in H^1 (\mathrm{Cl}(K), \mathbb Z/\ell^m)$ we have
\[ \omega_{m,K} (a,b) =-\frac{1}{2} s^\ell_{\mathrm{Cl}(K)}  ( B(a\cup b))\]
and for $a \in  H^1 (\mathrm{Cl}(K), \mathbb Z/\ell^v)$ and $b\in  H^1 (\mathrm{Cl}(K), \mathbb Z/\ell^m)$ we have
\[ \langle a,b\rangle  = s^\ell_{\mathrm{Cl}(K)}  ( a \cup B(b)).\] 
\end{theorem}

Thus, the $\omega$ and $\psi$ invariants both can be calculated by applying $s_{\mathrm{Cl}(K)} $ to classes in $H^3 ( \mathrm{Cl}(K), \mathbb Z/\ell^v)$ arising from the cup product and Bockstein homomorphism. 

Motivated by Theorem \ref{lst-tau}, we define $\omega$ and $\psi$ invariants for an oriented $\Gamma$-group $\mathbf H =(H,s_H)$ by the same formulas:
\[ \omega_{m} (a,b) =-\frac{1}{2} s_{H}^\ell  ( B(a\cup b)),  \quad \textrm{ for $a,b\in H^1(H,\Z/\ell^m)$ }\]
\[ \langle a,b\rangle = b(\psi(a))= s_{H}^\ell  ( a \cup B(b)), \quad \textrm{ for $a\in H^1(H,\Z/\ell^v)$ and $b\in H^1(H,\Z/\ell^m)$ } \]
in the first case assuming that $\ell$ is an odd prime. 
Note that the $\Gamma$ action plays no role in the definitions, but the fact that $s_H$ is $\Gamma$-invariant does imply that the pairings $\omega_{m} $ and $\langle, \rangle$ are $\Gamma$-invariant. 

To prove Theorem \ref{lst-tau}, the key lemma is the following, which generalizes \cite[Lemma 6.18]{Lipnowski2020}.

\begin{lemma}\label{exact-sequence-duality} Let $X$ be a scheme and let $0\to G_1 \to G_2 \to G_3 \to 0$ be a short exact sequence of finite flat commutative
 group schemes over $X$. Let $0 \to G_3^\vee \to G_2^\vee \to G_1^\vee \to 0$ be the Cartier dual sequence, which is exact (e.g. by \cite[Proposition 1.2.2]{OhAbelian}).
 
Let $B\colon \Hf^i (X, G_3) \to \Hf^{i+1} ( X, G_1) $ be the connecting homomorphism and similarly for $B^\vee \colon \Hf^i ( X, G_1^\vee) \to \Hf^{i+1} ( X, G_3^\vee ) $.

For $i, j$ nonnegative integers, $\alpha \in \Hf^i ( X, G_3)$ ,and $\beta \in  \Hf^{j}  ( X, G_1^\vee)$ we have \[ \alpha \cup  B^\vee \beta= (-1)^{i+1} B \alpha \cup \beta .\]

For $X$ the spectrum of a ring of integers of the number field,  if $i\geq 2$, then the same claim is true for $\alpha \in \Hfc^i( X, G_3)$  and $\beta \in  \Hf^{j}  ( X, G_1^\vee)$, where the cup product of a compactly-supported cohomology class and an ordinary cohomology class is a compactly-supported cohomology class.
\end{lemma} 

\begin{proof} The crux of the argument is to consider Cech cohomology and then use the product rule for cup products. 

Choose cocycles in the Cech complex of a suitable hypercovering representing $\alpha$ and $\beta$ \cite[01GU]{stacks}. 
Lift $\alpha$ to a cochain $\tilde{\alpha}$ in $C^i ( X, G_2)$ and lift $\beta$ to a cochain $\tilde{\beta}$ in $C^{j} ( X, G_2^\vee)$.  Then $\tilde{\alpha} \cup \tilde{\beta} $ is a cochain in $C^{i+j} (  X, \mathbb G_m)$ so $d (\tilde{\alpha} \cup \tilde{\beta} ) $ is a coboundary and thus represents the zero class in cohomology. Hence
\[  d (\tilde{\alpha} \cup \tilde{\beta} )= d\tilde{\alpha} \cup \tilde{\beta} + (-1)^i  \tilde{\alpha} \cup d \tilde{\beta} \] represents the zero class in cohomology \cite[01FP]{stacks}.

Now by the definition of the Bockstein homomorphism, $d \tilde{\alpha}$ is the pushforward of a cocycle representing $B \tilde{\alpha}$ along $G_1 \to G_2$. The multiplication map $G_2 \times G_2^\vee \to \mathbb G_m$, restricted to $G_1 \times G_2^\vee$, factors through the multiplication map $G_1 \times G_1^\vee \to \mathbb G_m$, so the multiplication map $C^{i+1} (X, G_2)\times C^{j} ( X, G_2^\vee)
 \to C^{i+j+1} ( X, \mathbb G_m)$, restricted to $C^{i+1} (X, G_1)\times C^{j} ( X, G_2^\vee)$, factors through the multiplication map  $C^{i+1} (X, G_1)\times C^{j} ( X, G_1^\vee)
 \to C^{i+j+1}( X, \mathbb G_m)$ and the projection $C^{j} ( X, G_2^\vee)\to C^{j} ( X, G_1^\vee)$ applied to $\tilde{\beta}$ recovers $\beta$. This implies $d\tilde{\alpha} \cup \tilde{\beta}$ represents the cohomology class $B\alpha \cup \beta$. 
 
 Symmetrical reasoning shows that $\tilde{\alpha} \cup d \tilde{\beta}$ represents the cohomology class $\alpha \cup B^\vee \beta$. Combining these gives the statement.
 
We now consider the case where $X$ is $\Spec \O_L$,  and $\alpha$ lies in the compactly supported cohomology, and $i\geq 2$. Compactly supported cohomology of the ring of integers of a number field is defined as the mapping cone of the natural map from the usual cohomology complex to a complex representing Tate cohomology of the Galois groups of the places at infinity. Thus cochains in degree $i$ for the compactly supported cohomology consist of cochains in degree $i$ for the usual cohomology together with cochains in degree $i-1$ for the Tate cohomology. Cochains in degree $i-1\geq 1$ for the Tate cohomology are the same as cochains in degree $i-1$ for the ordinary group cohomology which may be computed by Cech cohomology of the same hypercovering. So we may choose a cocycle representing the cohomology class $\alpha$ that is a pair consisting of a cochain in $C^{i}(\Spec \mathcal O_L; G_2)$ and a cochain in $C^{i-1}$ of the Cech complex at the infinite places. For lifting to $\tilde{\alpha}$, applying the cup product, and applying the differential, we may proceed as above.\end{proof}

\begin{proof}[Proof of Theorem \ref{lst-tau}] 

We apply Lemma \ref{exact-sequence-duality} to $X=\operatorname{Spec} \mathcal O_K$ and $G_1 = \mathbb Z/\ell^v$ and $G_2 = \mathbb Z/\ell^{m+v} $ and $G_3 = \mathbb Z/\ell^{m}$. Let $B$ be the connecting homomorphism of this short exact sequence and let $B^\vee$ be the connecting homomorphism of the dual short exact sequence $0 \to \mu_{\ell^m} \to \mu_{\ell^{v+m}} \to \mu_{\ell^v}\to 0$. For $i,j$ satisfying $i+j=2$ and classes $\alpha \in \Hf^i(X, \mathbb Z/\ell^m)$ and $\beta \in \Hf^j(X, \mu_{\ell^v})$, this gives the identity in $\Hf^3(\operatorname{Spec} \mathcal O_K, \mathbb G_m)$ 
\begin{equation}\label{AV-eq}  \alpha \cup  B^\vee \beta= (-1)^{i+1} B \alpha \cup \beta. 
\end{equation}
We apply \eqref{AV-eq} to establish both parts of Theorem \ref{lst-tau}.

The class called $\zeta_m$ in the definition of $\omega_{m,K}$ is the connecting homomorphism $B^\vee$ applied to the generator $\xi^{-n/\ell^v}$ 
of $\Hf^0( \operatorname{Spec} \mathcal O_K, \mu_{\ell^v})$. Applying \eqref{AV-eq} with $i=2$ we see that
\begin{align*} \omega_{m,K} (a,b) =- \frac{1}{2}\tA (\zeta_m \cup \phi(a) \cup \phi(b))&= \frac{1}{2}\tA (B^\vee \xi^{n/\ell^v} \cup \phi(a) \cup \phi(b))\\
&= \frac{1}{2}\tA (  \phi(a) \cup \phi(b) \cup B^\vee \xi^{n/\ell^v})\\
&= -\frac{1}{2}\tA (  B(\phi(a) \cup \phi(b)) \cup \xi^{n/\ell^v})\\
&= -\frac{1}{2}\tA (  \phi(B(a \cup b)) \cup \xi^{n/\ell^v}).
\end{align*}
The  last line follows from Lemma~\ref{L:BHpullback}, which gives the compatibility of $\phi$ with cup products and the Bockstein homomorphism.
Finally, we have a commutative diagram relating the definition of $s^\ell_{\Cl(K)}$ to the definition of $s_{\Cl(K)}$ (using the functoriality of $\phi$ in the coefficients from Lemma~\ref{L:BHpullback})
\[
\begin{tikzcd}
H^3(\Cl(K),\Z/n) \arrow[r,"\phi"] & H^3(\Spec \O_K,\Z/n) \arrow[r,"1\mapsto \xi"] & H^3(\Spec \O_K,\mu_n) \arrow[r,"\tA"] & \Q/\Z\arrow[r,"\times n"] & \Z/n \arrow[dl,"\times \frac{1}{n}"] \\
H^3(\Cl(K),\Z/\ell^v) \arrow[r,"\phi"] \arrow[u,"\times \frac{n}{\ell^v}"] & H^3(\Spec \O_K,\Z/\ell^v) \arrow[r,"1\mapsto \xi^{n/\ell^v}"] \arrow[u,"\times \frac{n}{\ell^v}"] & H^3(\Spec \O_K,\mu_{\ell^v}) \arrow[r,"\tA"] \arrow[u] &\Q/\Z. \arrow[u, equal] & 
\end{tikzcd}
\]
The composite of the top line above is $s_{\Cl(K)}$, and then going from the bottom left to the bottom right is  $s^\ell_{\Cl(K)}$, which is also the map $x\mapsto \tA( \phi(x)\cup\xi^{n/\ell^v})$.  Thus we conclude that
$$
 \omega_{m,K} (a,b) = -\frac{1}{2} s^\ell_{\Cl(K)}(  B(a \cup b)).
$$

Now we consider the claim about the pairing $\langle, \rangle$. First \cite[Proposition 6.3]{Lipnowski2020} implies that $\langle a, b\rangle $ for 
$a \in \Cl_K^\vee[\ell^v]$ and $b   \in \Cl_K^\vee[\ell^m]$ 
can be calculated by first cupping $\phi(a)$ with the fixed generator $\zeta$ to obtain a class in $\Hf^1( \operatorname{Spec} \mathcal O_K, \mu_{\ell^v})$, then applying the map $\Hf^1( \operatorname{Spec} \mathcal O_K, \mu _{\ell^v})\to \Hf^1(\operatorname{Spec} \mathcal O_K, \mathbb G_m ) $, then applying the Kummer map $\Hf^1(\operatorname{Spec} \mathcal O_K, \mathbb G_m ) \to \Hf^2(\operatorname{Spec} \mathcal O_K, \mu_{\ell^m} )  $, and then taking the cup product with $\phi(b)$, and then taking Artin-Verdier trace $\tA$.  Note that the two  maps
$\Hf^1( \operatorname{Spec} \mathcal O_K, \mu _{\ell^v})\to \Hf^1(\operatorname{Spec} \mathcal O_K, \mathbb G_m ) $
and $\Hf^1(\operatorname{Spec} \mathcal O_K, \mathbb G_m ) \to \Hf^2(\operatorname{Spec} \mathcal O_K, \mu_{\ell^m} ) $
are 
 both are part of the Kummer sequence but, even if $v=m$, these are not adjacent arrows in the Kummer sequence and so their composition is not automatically zero. Instead, one arises from a map of group schemes $\mu_{\ell^v} \to \mathbb G_m$ and the next from a short exact sequence $\mu_{\ell^m} \to \mathbb G_m \to \mathbb G_m$. Their composition therefore arises from the pullback of the short exact sequence along the map of group schemes, which is $\mu_{\ell^m} \to \mu_{\ell^{v+m}} \to \mu_{\ell^v}$. In other words, the composition of the two maps is $B^\vee$. Using this, \eqref{AV-eq}, and the compatibility of $\phi$ with cup products and the Bockstein homomorphism, we have
\begin{align*} \langle a,b\rangle = \tA ( B^\vee(\phi(a) \cup \xi^{n/\ell^v})\cup \phi(b))
&= \tA ( \phi(a) \cup \xi^{n/\ell^v}\cup B(\phi(b)))\\
&= \tA ( \phi(a \cup B(b)) \cup \xi^{n/\ell^v})\\
&=s^\ell_{\Cl(K)}(a \cup B(b)). \qedhere \end{align*}\end{proof}

The paper \cite[Lemma 6.20]{Lipnowski2020} 
proved a certain compatibility between the $\omega$ and $\psi$ invariants of class groups. We will show this compatibility for general oriented groups
follows straightforwardly from our definition of the general version of these invariants. The proof follows the same idea as that in \cite{Lipnowski2020}, but is a bit more transparent as certain non-split group schemes appearing in \cite{Lipnowski2020} are replaced with cyclic groups.  The upshot is that the data of a bilinearly-enhanced abelian group from \cite{Lipnowski2020} all comes from the data of an oriented  group.

\begin{proposition}\label{P:psiomega} Let $H$ be a finite group, $n$ a positive integer, $s_H \in H_3(H, \mathbb Z/n)$ a class. Let $\ell$ be an odd prime and $v$ the $\ell$-adic valuation of $n$.  Let $r$ be a non-negative integer and let $m=v+r$.
For $\omega_m$ and $\langle, \rangle$ the pairings defined in terms of $s_H^\ell$, we have, for $a,b \in H^1(H, \mathbb Z/\ell^m)$
and $\bar{a},\bar{b}\in H^1(H, \mathbb Z/\ell^v)$ their images under the map $1\mapsto 1$ on coefficients,
 \[   \langle \bar{a}, b \rangle- \langle \bar{b}, a\rangle = 2 \omega_m (a,b).\]
\end{proposition}
\begin{proof} Choose cocycles representing $a,b$ (and call these $a,b$ by a slight abuse of notation) and lift these to cochains $\tilde{a},\tilde{b} \in C^1(H, \mathbb Z/\ell^{v+m})$. Then $\tilde{a} \cup \tilde{b} \in C^2(H, \mathbb Z/\ell^{v+m})$ lifts $a \cup b$.

Hence, by the definition of Bockstein homomorphism, $d( \tilde{a} \cup \tilde{b})$ is the pushforward from $C^3 (H , \mathbb Z/\ell^v) $ to $C^3(H, \mathbb Z/\ell^{v+m})$ of a cycle representing $B(a \cup b)$. The Leibnitz formula for cup product gives 
\[ d ( \tilde{a} \cup \tilde{b}) = d (\tilde{a}) \cup \tilde{b}- a \cup d (\tilde{b} ).\]

By the definition of Bockstein homomorphism, $d(\tilde{\alpha})$ is the pushforward from $C^2(H, \mathbb Z/\ell^v)$ of a cocycle representing $B(a)$. The multiplication map $\mathbb Z/ \ell^{v+m} \times \mathbb Z/\ell^{v+m} \to \mathbb Z/\ell^{v+m} $ composed with the multiplication by $\ell^m$ map $t: \mathbb Z/\ell^v \to \mathbb Z/\ell^{v+m}$ gives a bilinear form $\mathbb Z/\ell^v \times \mathbb Z/\ell^{v+m} \to \mathbb Z/\ell^{v+m}$. This can be expressed as the composition of first the mod $\ell^v$ map $r: \mathbb Z/\ell^{v+m} \to \mathbb Z/\ell^v$, then the multiplication map $\mathbb Z/\ell^v \times \mathbb Z/\ell^v \to \mathbb Z/\ell^v$, followed by the multiplication by $\ell^m$ map $t\colon \mathbb Z/\ell^v \to \mathbb Z/\ell^{v+m}$.  

Applying this identity to cocycles, the cup product map $C^2(H, \mathbb Z/\ell^{v+m}) \times C^1(H, \mathbb Z/\ell^{v+m}) \to C^3(H, \mathbb Z/\ell^{v+m})$ composed with  $t_*: C^2(H, \mathbb Z/\ell^v) \to C^2(H, \mathbb Z/\ell^{v+m})$ is the composition of $r_*:C^1(H, \mathbb Z/\ell^{v+m})\to C^1(H, \mathbb Z/\ell^v)$, the cup product map $C^2(H, \mathbb Z/\ell^{v}) \times C^1(H, \mathbb Z/\ell^{v}) \to C^3(H, \mathbb Z/\ell^{v})$, and  $t_* :C^3(H, \mathbb Z/\ell^v) \to C^3(H, \mathbb Z/\ell^{v+m})$.  In other words, for $x\in C^2(H, \mathbb Z/\ell^v)$ and $y\in C^1(H, \mathbb Z/\ell^{v+m})$, we have $t_*(x)\cup y=t_*(x \cup r_*(y))$.  

Applying this to $(t^{-1})_*d(\tilde{a})$ (which is a cochain representing $B(a)$) and $\tilde{b}$,  we have $d(\tilde{a})\cup \tilde{b}=t_*((t^{-1})_*d(\tilde{a}) \cup r_*(\tilde{b}))$. 
We observe that since $m =v+r\geq v$, the projection $r_*:C^1(H, \mathbb Z/\ell^{v+m})\to C^1(H, \mathbb Z/\ell^v)$ factors through $C^1(H, \mathbb Z/\ell^m)$.  
Thus $r_*(\tilde{b})=\bar{b}$.
So we have $d(\tilde{a})\cup \tilde{b}=t_*((t^{-1})_*d(\tilde{a}) \cup \bar{b})$ and similarly $\tilde{a}\cup d(\tilde{b})=t_*(\bar{a} \cup (t^{-1})_*d(\tilde{b})).$  Applying $t_*^{-1}$ to the Leibnitz formula then gives
$$
t_*^{-1}(d ( \tilde{a} \cup \tilde{b}))= (t^{-1})_*d(\tilde{a}) \cup \bar{b} - \bar{a} \cup (t^{-1})_*d(\tilde{b}),
$$
which in cohomology gives
$$
B(a\cup b)=B(a)\cup \bar{b} -\bar{a}\cup B(b).
$$
Taking $s_H^\ell$ of both sides we obtain
\[ s_H^\ell( B( a\cup b) ) = s_H^\ell (  \bar{b} \cup B(a)  ) - s_H^\ell( \bar{a} \cup B(b))\] which by definition gives.
\[ -2 \omega_m (a,b)= \langle \bar{b}, a\rangle - \langle \bar{a}, b \rangle. \qedhere\] \end{proof}

\begin{remark}We remark on the case $\ell=2$, not considered in \cite{Lipnowski2020}. The factor $\frac{1}{2}$ in the formula for $\omega_m$ giving in Theorem \ref{lst-tau} is unreasonable to include in this case, and can be removed (it exists only for consistency with \cite{Lipnowski2020}, who included it for consistency with a certain random matrix model.) 

The pairing $a ,b \mapsto  s_{H}  ( B(a\cup b))$ is clearly antisymmetric. For pairings on $2$-power torsion groups, being antisymmetric is a strictly weaker condition than being alternating, i.e., the pairing of any element with itself being $0$. It is straightforward to check the pairing $a ,b \mapsto  s_{H}  ( B(a\cup b))$ is alternating. Any class $a \in H^1(H, \mathbb Z/2^m)$ is the pullback of a universal class in $H^1( \mathbb Z/2^m, \mathbb Z/2^m)$ along a homomorphism $H \to \mathbb Z/2^m$, and so $B( a\cup a)$ is the pullback of a class in the image of the Bockstein map $H^2(\mathbb Z/2^m, \mathbb Z/2^m) \to H^3( \mathbb Z/2^m, \mathbb Z/2^v)$, but that map is zero, as may be calculated using the standard periodic resolution for group cohomology of cyclic groups.

If this pairing were checked to be equivalent to the pairing constructed in \cite[Proposition 6.4]{MorganSmith}, then the fact that this pairing is alternating would be equivalent to \cite[Proposition 6.6]{MorganSmith}, and the argument above would give a new proof of that proposition. \end{remark}

It is straightforward to check that for an abelian group $H$, the $\omega$ and $\psi$ invariants typically do not contain all the information included in $s_H$.

\begin{example} Let $\ell$ be an odd prime, $n=\ell$, and $H =( \mathbb Z/\ell)^d$. Then the K\"unneth formula and the fact that $H^i (\mathbb Z/\ell, \mathbb Z/\ell) \cong \mathbb Z/\ell$ for all $i \geq 0$ gives $H^3( H, \mathbb Z/\ell) \cong (\mathbb Z/\ell)^{\frac{ d(d+1)(d+2)}{6}}$.  Thus there are $\ell^{\frac{ d(d+1)(d+2)}{6}}$ choices of orientation on $H$.
 On the other hand, the $\psi$ invariant is a linear form $(\mathbb Z/\ell)^d \to \mathbb (Z/\ell)^d$ and the $\omega$ invariant is determined by $\omega_1$ and thus by $\psi$ using Proposition~\ref{P:psiomega}.  So the $\omega,\psi$ invariants are parametrized by $\Hom ((\mathbb Z/\ell)^d, (\mathbb Z/\ell)^d)\cong (\mathbb Z/\ell)^{d^2}$. Thus the number of choices of orientation $s_H$ for each $\omega, \psi$ is at least $\ell^{ \frac{ d(d+1)(d+2)}{6} - d^2} = \ell^{ \frac{ d(d-1)(d-2)}{6}}$. \end{example}

In fact, for a finite  abelian group $H$, an orientation $s_H$ is determined by $\omega$ and $\psi$ together with certain trilinear forms on $H^1(H, \mathbb Z/\ell^m)$.

\begin{lemma}\label{h3-abelian-classification} Let $H$ be a finite abelian group and $n$ a nonnegative integer. The group $H^3 (H, \mathbb Z/n)$ is generated by classes of the following three types, over all primes $\ell\mid n$ with $v$ the  $\ell$-adic valuation $v$ of $n$, under the $\times n/\ell^v$ map on coefficients
$H^3 (H, \mathbb Z/\ell^v)\ra H^3 (H, \mathbb Z/n)$: 
\begin{enumerate}
\item For $m \leq v$, the image under $H^3(H, \mathbb Z/\ell^m) \stackrel{\times \ell^{v-m}}{\to} H^3(H, \mathbb Z/\ell^v)$ of the cup product of three classes in $H^1(H, \mathbb Z/\ell^m)$.
\item Classes of the form $B ( a \cup b)$ for $a,b\in H^1(H, \mathbb Z/\ell^m)$.
\item Classes of the form $a \cup B(b)$ for $a \in H^1(H, \mathbb Z/\ell^v)$ and $b\in H^1(H, \mathbb Z/\ell^m)$.
\end{enumerate}

In particular, a class $s_H \in H_3(H, \mathbb Z/n)$ is uniquely determined by its pairings with classes of these three types. \end{lemma}

\begin{proof} Since $H^3( H, \mathbb Z/n)$ is generated by $H^3(H, \mathbb Z/\ell^v)$ for $\ell$ dividing $n$, we may as well assume $n=\ell^v$. Since $\ell$-power cohomology classes of finite abelian groups are always pullbacks from a finite abelian $\ell$-group, we may assume $H$ is an $\ell$-group.

The short exact sequence $0 \to \mathbb Z \to \mathbb Z \to \mathbb Z/n \to 0$ induces a long exact sequence \[\dots \to H^3( H, \mathbb Z)\to H^3(H, \mathbb Z/n) \to H^4(H,\mathbb Z) \to H^4(H,\mathbb Z) \to \dots \] and thus in particular an exact sequence
\[ H^3( H, \mathbb Z)/ n H^3( H, \mathbb Z)    \stackrel{t_0} \to H^3(H, \mathbb Z/n) \stackrel{B_v}{\to} H^4(H,\mathbb Z)[n] \to 0.\] We will check that classes of type (2) span the image of $t_0$, and that $B_v$ applied to classes of type (1) and (3) generate $H^4(H,\mathbb Z)[n]$. It immediately follows that classes of types (1), (2), and (3) together generate $H^3(H, \mathbb Z/n)$.

To do this, we use a description of the integral cohomology ring of a finitely generated abelian group due to Huebschmann~\cite{Huebschmann91}. We first review this description, though we must adjust Huebschmann's notation slightly to harmonize it with our own. Write $H = \prod_{i=1}^r \mathbb Z/d_i$ where $d_i \mid d_{i+1}$. Let $A(H)$ be the graded algebra generated by variables $x_1,\dots, x_r$ in degree $1$ and $\zeta_1,\dots, \zeta_r$ in degree $2$, subject to the relations that the $\zeta_i$ commute with everything, the $x_i$ anticommute with each other, $x_i^2 = \frac{ d_i (d_i-1)}{2} \zeta_i$, and $d_i \zeta_i=0$. Also for $i_1,\dots,i_t$ a tuple of integers satisfying $1\leq i_1<\dots<i_t \leq r$, let $Z_{{i_1} ,{i_2}, \dots ,{i_t} }=B_{d_{i_1}}(X_{i_1} \cup \dots \cup X_{i_t})\in H^{t+1}(H, \mathbb Z)$, for certain
 the classes $X_i\in H^1(H,\Z/d_i)$. It is not necessary for our argument to know what classes these are, but presumably $X_i$ is dual to the standard generator of $\mathbb Z/d_i$, although it is not easy to check this from \cite{Huebschmann91}.

Then \cite[Theorem B]{Huebschmann91} states that the integral cohomology ring $H^*(H, \mathbb Z)$ is isomorphic to the subring of $A(H)$ generated by the expressions
\[ \tilde{\zeta}_{x_{i_1} x_{i_2} \dots x_{i_t} } = \sum_{s=1}^t (-1)^{s-1} \frac{d_s}{d_1} \zeta_{i_s} x_{i_1} \dots x_{i_{s-1}} x_{i_{s+1}}\dots x_{i_t} \] for $i_1,\dots,i_t$ a tuple of integers satisfying $1\leq i_1<\dots<i_t \leq r$, and furthermore states that this isomorphism sends $Z_{{i_1},{i_2}, \dots ,{i_t} }$ to $\tilde{\zeta}_{x_{i_1} x_{i_2} \dots x_{i_t} }$. In particular, the cohomology ring is generated by the classes $Z_{{i_1},{i_2}, \dots ,{i_t}}$ in degree $t+1$.

Since we have generators in degree $2,3,$ and $4$, but not $1$, the only way to obtain a class in degree $3$ is from a generator $Z_{{i} ,{j}}$ in degree $3$, and the only ways to obtain a class in degree $4$ is from a generator $Z_{{i} ,j ,k}$ in degree $4$ or the product $Z_{i} Z_{j}$ of two generators in degree $2$. In the last case we do not assume $i<j$ but without loss of generality assume $i\leq j$. To calculate the $n$-torsion in $H^4( H,\mathbb Z)$, we need to understand the relations between these classes. Note that $A(H)$ is generated as an abelian group by monomials in the classes $x_1,\dots, x_r,\zeta_1,\dots, \zeta_r$ where no $x_i$ appears twice, subject to the relation that a monomial containing $\zeta_i$ is $d_i$-torsion. Thus a polynomial is $n$-torsion if and only if each monomial appearing is $n$-torsion. Since no two of the classes $\tilde{\zeta}_{x_{i} x_j x_k}$ or  $\tilde{\zeta}_{x_i} \tilde{\zeta}_{x_j}$ share a monomial, a sum of them is $n$-torsion if and only if each summand is $n$-torsion. So the $n$-torsion in $H^4(H,\mathbb Z)$ is generated by the classes $\frac{d_i}{\gcd(d_i,n)} Z_{{i} ,j ,k} $ and $\frac{d_i}{\gcd(d_i,n)}Z_{i} Z_{j} $.

We will now show that the images of classes $Z_{i,j}$ under $t_0$ are classes of type (2) and the classes $\frac{d_i}{\gcd(d_i,n)} Z_{i,j,k} $ and $\frac{d_i}{\gcd(d_i,n)}Z_{i} Z_{j} $ arise from applying $B_v$ to classes of type (1) and (3).
This  requires using the compatibility of many different Bockstein maps.  Let $B_{a,b}$ denote the connecting Bockstein homomorphisms for the short exact sequence $\Z/\ell^b\stackrel{\times \ell^a}{\ra}\Z/\ell^{a+b} \ra \Z/\ell^a$, and $B_a=B_{a,\infty}$ analogously for the sequence $\Z\stackrel{\times \ell^a}{\ra}\Z\ra \Z/\ell^a$.
Our previous map $B$ is $B_{m,v}$ in this notation (though $m$ can vary).
We let $t_k:H^*(H,\Z/\ell^a) \ra H^*(H,\Z/\ell^b)$ denote the map induced by the map of coefficients $\Z/\ell^a \ra \Z/\ell^b$ (or $\Z\ra \Z/\ell^b$) sending $1$ to $\ell^k$ (for $a+k\leq b$), slightly overloading the notation.  
It is straightforward to check from the definitions that $t_k B_{a,b}=B_{c,d} t_{c+k-a}$ for any $a,b,c,d,k\geq 0$ with $c+k-a\geq 0$.   (Roughly, $B_{a,b}$ comes about by lifting, taking $d$, and dividing by $\ell^a$.)  Note that the identity map is a special case of $t_0$.

We first consider a class $Z_{i,j}$. Let $m$ be such that $\ell^m = d_i$. 
Let $\bar{X_j}=t_0(X_j)\in H^1(H,\Z/\ell^m)$.
We have $Z_{i,j}=B_m(X_{i} \cup X_j)=B_m(X_{i} \cup \bar{X}_j)$.
Thus in $H^3(H, \mathbb Z/\ell^v)$ we have $t_0(Z_{i,j})=
t_0B_m(X_{i} \cup \bar{X}_j)=B_{m,v}(X_{i} \cup \bar{X}_j)$, a class of type $(2)$, as desired.


We now consider a generator $\frac{d_i}{\gcd(d_i,n)} Z_{i,j,k} $ of $H^4(H,\mathbb Z)[n]$. Let $m$ be such that $\ell^m=d_i$. If $m \leq v$ then $\gcd(d_i, n) =\ell^m$ and  $\frac{d_i}{\gcd(d_i,n)} Z_{i,j,k}=Z_{i,j,k}$.   We write $\bar{X}_j$ and $\bar{X}_k$ as above.
We have $Z_{i,j,k}=B_m(X_{i} \cup X_j \cup X_k)=B_m(X_i \cup \bar{X}_j \cup \bar{X}_k)=B_{v}( t_{v-m} (X_i \cup \bar{X}_j \cup \bar{X}_k) )$
is the image of a class of type (1) under $B_v$.
If $m>v$ then our generator of $H^4(H,\mathbb Z)[n]$ is $\frac{d_i}{\gcd(d_i,n)} Z_{i,j,k}=\ell^{m-v} Z_{i ,j ,k}$. 
We write $\bar{X}_i$ for $t_0(X_i)\in  H^3(H, \mathbb Z/\ell^v)$, and similarly for $\bar{X}_j, \bar{X}_k$.
We have $\ell^{m-v} Z_{i ,j ,k}=t_{m-v} B_m(X_{i} \cup X_j \cup X_k)=B_v(t_0(X_i \cup {X}_j \cup {X}_k ))
=B_v(\bar{X}_i \cup \bar{X}_j \cup \bar{X}_k )$.  In particular $\ell^{m-v} Z_{i ,j ,k}$ is $B_v$ applied to a class of type (1) with $m=v$.

Finally, consider a class $\frac{d_i}{\gcd(d_i,n)}Z_iZ_j$. We first choose a class $a \in H^1(H, \mathbb Z/\ell^v)$ so that $\frac{d_i}{\gcd(d_i,n)}Z_i=B_v(a).$
If $d_i =\ell^c$ for $c\geq v$ then we can take $a=t_0(X_i)$, and   $\frac{d_i}{\gcd(d_i,n)}Z_i=\ell^{c-v}B_c(X_i)=B_v(a)$.
If $d_i =\ell^c$ for $c<v$ then we can take $a=t_{v-c}(X_i)$, and $\frac{d_i}{\gcd(d_i,n)}Z_i=B_c(X_i)=B_v(a)$.
We let $m$ be such that $\ell^m = d_j$ and set $b = X_j \in H^1(H, \mathbb Z/\ell^m)$, 
and in $H^2(H, \mathbb Z/\ell^v)$ we have that $t_0(Z_j)=B_{m,v}(b)=B(b).$
Then to confirm that $\frac{d_i}{\gcd(d_i,n)}Z_iZ_j$ arises as the image under $B_v$ of a class of type (3), it suffices to check that $B_v( a \cup B(b))=\frac{d_i}{\gcd(d_i,n)}Z_i \cup Z_j$. To do this, choose a cocycle in $C^1( H, \mathbb Z/\ell^v)$ representing $a$ and a lift to a cochain $\tilde{a} \in C^1(H, \mathbb Z)$. Also choose a cocycle 
$\in C^2(H, \mathbb Z)$
representing $Z_j $, which we will refer to for simplicity as $Z_j$. Then $\tilde{a} \cup Z_j$ is a lift to $C^3(H, \mathbb Z)$ of a cochain representing $a \cup B(b)\in H^3(H,\Z/\ell^v)$.  Thus  $d (\tilde{a} \cup Z_j)$ is the image under the multiplication-by-$\ell^v$ map $C^4(H,\mathbb Z)\to C^4(H,\mathbb Z)$ of a cochain representing $B_v(a \cup B(b))\in H^4 (H, \mathbb Z)$. We have
\[ d( \tilde{a} \cup Z_j) =  d(\tilde{a}) \cup Z_j -  \tilde{a} \cup d(Z_j) = d(\tilde{a}) \cup Z_j -  \tilde{a} \cup 0 = d(\tilde{a})  \cup Z_j\] since $Z_j$ is a cocycle. Now $d(\tilde{a})$ is the image under the multiplication by $\ell^v$ map of a cocycle representing $\frac{d_i}{\gcd(d_i,n)}Z_i$, so that  $d (\tilde{a} ) \cup Z_j$ is the image under multiplication by $\ell^v$ of a cochain representing $\frac{d_i}{\gcd(d_i,n)}Z_i \cup Z_j$. 
So  $d (\tilde{a} \cup Z_j)$ is the image under the multiplication-by-$\ell^v$ map $C^4(H,\mathbb Z)\to C^4(H,\mathbb Z)$ of both 
a cochain representing $B_v(a \cup B(b))$ and a cochain representing $\frac{d_i}{\gcd(d_i,n)}Z_i \cup Z_j$.
Since multiplication by $\ell^v$ is injective on cochains, this gives that $B_v( a \cup B(b))=\frac{d_i}{\gcd(d_i,n)}Z_i \cup Z_j$. 
Thus,  the classes $\frac{d_i}{\gcd(d_i,n)}Z_iZ_j$ are in the image of classes of type (3) under $B_v$, completing the proof.
\end{proof}

However, when we consider a group $ H \rtimes \Gamma$ with $H$ abelian, if $ \Gamma \cong \mathbb Z/2$ then the trilinear forms on $H^1(H, \mathbb Z/\ell^m)$ vanish so the $\omega$ and $\psi$ invariants do determine $\tau$.

\begin{lemma}\label{lst-case} Assume $\Gamma\cong \mathbb Z/2$ and let $H$ be a finite abelian $\Gamma$-group of odd order such that $H_\Gamma$ is trivial. The group $H^3 (H, \mathbb Z/n)^\Gamma$ is generated by classes of types (2) and (3) of Lemma \ref{h3-abelian-classification}. In particular, a class $s_H \in H_3(H, \mathbb Z/n)^\Gamma$ is uniquely determined by its $\omega$ and $\psi$ invariants. \end{lemma}

\begin{proof} Since $H$ has odd order, $H^3( H, \mathbb Z/n)$ has odd order, and hence splits into summands where $\Gamma$ acts with eigenvalue $+1$ and $-1$ respectively. Since $\Gamma$ acts on $H$ with eigenvalue $-1$, it acts on $H^1(H, \mathbb Z/\ell^m)$ with eigenvalue $-1$, and hence on cup products of three classes in $H^1(H, \mathbb Z/\ell^m)$ with eigenvalue $(-1)^3=-1$. Similarly, classes of types (2) and (3) have eigenvalue $(-1)^2=1$. So any invariant class is a sum of classes of type (2) and (3). The last claim follows because, again for odd order reasons $H_3(H, \mathbb Z/n)^\Gamma$ and $H^3(H, \mathbb Z/n)^\Gamma$ are dual. \end{proof}

Furthermore, for $\mathcal L$ a level consisting only of abelian $\Gamma$-groups, our formula for the probability depends only on $\psi$ and the number of automorphisms of the oriented group, as shown in the following lemma.  One result of this is that we see that different $s_H$ with the same $\omega,\psi$ do occur, and so we have found an invariant that carries more information than the billinearly enhanced group structure of \cite{Lipnowski2020}, even on class groups.
On the other hand, when for a given $\Gamma$-group $H$, the probabilities for all possibly orientations $s_H$
with a given $\psi$ are summed, the following lemma gives that the resulting formula is still a nice product.  Thus for class groups the entire refinement to $s_H$, or even to the pair $\omega, \psi$, is not necessary to obtain product formula probabilities.

\begin{lemma} Let $\mathcal L$ be a level in the category of finite $\Gamma$-groups that consists only of abelian $\Gamma$-groups (or, equivalently, a level in the category of finite abelian $\Gamma$-groups). Let $\mathcal W$ be the set of finite $n$-oriented $\Gamma$-groups whose underlying $\Gamma$-group is in $\mathcal L$. Let $H \in \mathcal L$ be of order relatively prime to $|\Gamma|$. 
For any $n$-orientation $s_H$ on $H$, let $\nu_{\Gamma, n, U} (\{X|X^\cL\isom  (H,s_H)\})$  be the probability from Theorem~\ref{T:intro-measure-general}, which is also called $v_{\cW,(H,s_H)}$ in Theorem~\ref{T:measure}.
Then $  | \Aut(H, s_H)| v_{\cW,(H,s_H)}$ 
depends on $s_H$ only through $\psi \colon H^\vee [\ell^v] \to H[\ell^v]$ for primes $\ell$ dividing $n$.

\end{lemma}

\begin{proof} Let $\mathbf H = (H, s_H)$.  Theorem~\ref{T:intro-measure-general} expresses $v_{\mathcal W, \mathbf H}$ as $\frac{M_{\mathbf H} }{|\Aut(\mathbf H)|}$ times a product of local factors $w_{V_i}$ and $w_{N_i}$. We have cancelled the $|\Aut(\mathbf H)|$ factor and the formula for $M_{\mathbf H}$ 
manifestly involves only $H$ and not $s_H$, so it suffices to check the $w_{V_i}$ and $w_{N_i}$ factors depend on $s_H$ only through $\psi$. Also if $H_\Gamma$ is nontrivial then $M_{\mathbf H}=0$ by definition and thus $v_{\mathcal W, \mathbf H}=0$ independently of $s_H$, so we may assume $H_\Gamma$ is trivial.

Since $\mathcal L$ consists only of abelian groups, there are certainly no admissible non-abelian groups $N_i$, and no $N_i$ terms. Furthermore, for $V_i$ a representation of $H\rtimes \Gamma$, if $V_i$ is admissible then $V_i \rtimes H \in \mathcal L$ is abelian so the action of $H$ on $V_i$ is trivial, i.e., $V_i$ is a representation of $\Gamma$.

In the expressions given for $w_{V_i}$ in various cases before Theorem~\ref{T:intro-measure-general}, the only dependence on $s_H$ is through the subgroup  
$H^2( H \rtimes \Gamma, V_i)^{\mathcal L, s_H}$ of $H^2( H \rtimes \Gamma, V_i)$, and only in cases when the characteristic of $V_i$ divides $n$.
%
%
Thus it suffices to check $H^2( H \rtimes \Gamma, V_i)^{\mathcal L, s_H}$ depends on $s_H$ only through $\psi$ for each $V_i$ with $\ell:=\operatorname{char} V_i$ dividing $n$.

We first observe that $H^2(H\rtimes \Gamma, V_i) = H^2(H, V_i)^\Gamma = ( H^2(H, \mathbb F_\ell) \otimes V_i)^\Gamma$ since the action of $H$ on $V_i$ is trivial. Inside $H^2(H, \mathbb F_\ell)$ we have a subspace $H^2(H, \mathbb F_\ell)^{B}$ consisting of classes in the image of the Bockstein map $B_{m,1}: H^1( H , \mathbb Z/\ell^m) \to H^2 (H, \mathbb F_\ell)$ for at least one (equivalently, all sufficiently large) $m$. The subspace $H^2(H, \mathbb F_\ell)^{B}$ is clearly $\Gamma$-invariant so $(H^2(H, \mathbb F_\ell)^{B} \otimes V_i)^\Gamma$ is a well-defined subspace of $( H^2(H, \mathbb F_\ell) \otimes V_i)^\Gamma$. Defining $H^2( H \rtimes \Gamma, V_i)^{\mathcal L} $ to consist of classes whose corresponding $\Gamma$-extension
(via Lemma~\ref{L:diffext})
 is in $\mathcal L$, we next check that
\[ H^2( H \rtimes \Gamma, V_i)^{\mathcal L} \subseteq (H^2(H, \mathbb F_\ell)^{B} \otimes V_i)^\Gamma .\]  
To do this, observe that a class $\alpha\in H^2( H \rtimes \Gamma, V_i)^{\mathcal L}$ can be viewed as a ($\Gamma$-invariant) $\dim V_i$-tuple of classes $\alpha_j \in H^2(H, \mathbb F_\ell)$, with the corresponding extension of $H$ by $V_i$ being the fiber product of the $\dim V_i$ extensions of $H$ by $\mathbb F_\ell$, and with $\alpha \in (H^2(H, \mathbb F_\ell)^{B} \otimes V_i)^\Gamma$ if and only if each $\alpha_j \in H^2(H, \mathbb F_\ell)^B$. Since the extension corresponding to $\alpha$ lies in $\mathcal L$, it is certainly abelian, hence all the extensions in the fiber product are abelian. So it suffices to check that if a class $\alpha_i \in H^2(H, \mathbb F_\ell)$ corresponds to an abelian extension, then it lies in $H^2(H, \mathbb F_\ell)^{B}$.

From the definition of the Bockstein homomorphism $H^1( H , \mathbb Z/\ell^m) \to H^2 (H, \mathbb F_\ell)$,  the elements in its image are represented by symmetric cochains, and thus corresponding to abelian extensions.  
Abelian extensions of $H$ by $\mathbb F_\ell$ are classified by $\operatorname{Ext}^1(H, \mathbb F_\ell)$, which has size $|H[\ell]|$. 
We choose an $m$ large enough so that the image of $H^1( H , \mathbb Z/\ell^m) \to H^2 (H, \mathbb F_\ell)$ is $H^2(H, \mathbb F_\ell)^{B}$ and so that $\ell^mH$ has no $\ell$-torsion.
We consider the exact sequence
$$
H^1( H , \mathbb Z/\ell^{m+1}) \ra H^1( H , \mathbb Z/\ell^m) \to H^2 (H, \mathbb F_\ell).
$$
One can check that $\ell^m H$ having no $\ell$-torsion implies that the image of the first map is $\ell H^1( H , \mathbb Z/\ell^m) \sub H^1( H , \mathbb Z/\ell^m)$, and thus 
$H^2(H, \mathbb F_\ell)^{B}$ also has size $|H[\ell]|$.  Thus all abelian extension classes in $H^2(H, \mathbb F_\ell)$ are in $H^2(H, \mathbb F_\ell)^{B}$.


We now handle the case that $V_i$ is trivial. In this case $ (H^2(H, \mathbb F_\ell)^{B} \otimes V_i)^\Gamma=  (H^2(H, \mathbb F_\ell)^{B})^\Gamma$ is the image of $H^1(H, \mathbb Z/\ell^m)^\Gamma$ (since $\Gamma$ acts semisimply on abelian groups of order a power of $\ell$). However $H^1(H, \mathbb Z/\ell^m)^\Gamma=0$ since $H_\Gamma$ is trivial and any class in $H^1(H, \mathbb Z/\ell^m)^\Gamma$ would give a $\Gamma$-invariant quotient. Thus $H^2( H \rtimes \Gamma, V_i)^{\mathcal L, s_H}=0$ regardless of $s_H$ in this case.

In the case that $V_i$ is nontrivial, $H^2(H\rtimes \Gamma, V_i)^{\mathcal L, s_H}$ consists of classes $\alpha \in H^2(H\rtimes \Gamma, V_i)^{\mathcal L}$ such that $s_H( \alpha \cup \beta) =0 $ for all $\beta \in H^1( H\rtimes \Gamma, \Hom(V_i,\Z/n))$.  We claim that for such $\alpha$ and $\beta$, the value of
$s_H( \alpha \cup \beta)$ depends on $s_H$ only through $\psi$.   One can easily reduce the proof of this claim to the case that $n$ is a power of $\ell$, so we now assume $n=\ell^v$.   
We may compute $s_H( \alpha \cup \beta)$ by pulling back $\alpha$ to $H^2(H , V_i) \cong  H^2(H, \mathbb F_\ell)^{\dim V_i}$ and pulling $\beta$ back to $H^1(H, V_i^\vee) \cong H^1(H, \mathbb F_\ell)^{\dim V_i}$, where they become tuples of classes $\alpha_1,\dots, \alpha_{\dim V_i} \in H^2(H, \mathbb F_\ell)^B $ and $\beta_1,\dots, \beta_{\dim V_i} \in H^1(H, \mathbb F_\ell)$ respectively (after choosing a basis of $V_i$ and the dual basis of $V_i^\vee$). 
Let $t_{v-1}:  H^*( H, \F_\ell) \ra H^*( H, \Z/\ell^v)$ be the map induced by the inclusion $\F_\ell \sub \Z/\ell^v$ which multiplies elements by $\ell^{v-1}$.
Thus it suffices to show that $s_H(t_{v-1}(\alpha_i\cup\beta_j))$ depends on $s_H$ only through $\psi$.  

Let $t_0: H^*( H, \Z/\ell^v) \ra H^*( H, \F_\ell)$ be the map induced by the quotient $\Z/\ell^v\ra \F_\ell$ sending $1$ to $1$.
For $u\in H^2(H,\Z/\ell^v)$ and $v\in H^1(H,\F_\ell)$,  it is straightforward to check that $t_{v-1}(t_0(u)\cup v  )=u\cup t_{v-1}(v)$.  
By definition of $H^2(H, \mathbb F_\ell)^B$, we have $\alpha_i =B_{m,1}(a)$ for some $a \in H^1(H, \mathbb Z/\ell^m)$.
As in the proof of Lemma~\ref{h3-abelian-classification}, we have that $B_{m,1}=t_0\circ B_{m,v}$ and so $\alpha_i =t_0(B_{m,v}(a))$.
Thus
$$
s_H(t_{v-1}(\alpha_i\cup\beta_j))=s_H(t_{v-1}(t_0(B_{m,v}(a))\cup\beta_j))=s_H(B_{m,v}(a)\cup t_{v-1}(\beta_j) )
=
a (\psi (t_{v-1}(\beta_j))).
$$
Thus we conclude that $s_H( \alpha \cup \beta)$ and hence $H^2( H \rtimes \Gamma, V_i)^{\mathcal L, s_H}$ depends on $s_H$ only through $\psi$, and have proven the lemma.
\end{proof}

\subsection{The nonabelian $\omega$-invariant, after Liu}\label{SS:Liu}

Let $H$ be a finite group and $s_H\in H_3(H, \mathbb Z/n)$ an element. We define the \emph{lifting invariant associated to $s_H$} to be the image of $s_H$ in $H_2(H,\mathbb Z)$ under the connecting homomorphism  associated to the short exact sequence $0 \to \mathbb Z \to \mathbb Z \to \mathbb Z/n \to 0$. The lifting invariant defines a linear form on $H^2(H, \mathbb Q/\mathbb Z)$. Dually this linear form may be expressed as the composition \[ H^2(H, \mathbb Q/\mathbb Z) \to H^3(H, \mathbb Z/n) \to \mathbb Z/n \to \mathbb Q/\mathbb Z\] with the first arrow the connecting homomorphism $B'_n$ associated to the short exact sequence $0 \to \mathbb Z/n \to \mathbb Q/\mathbb Z \to \mathbb Q/\mathbb Z \to 0$, the second arrow $s_H$, and the third arrow division by $n$.

We will check that this is compatible (up to a constant factor) with previously-defined notions of the lifting invariant. 

In this section, we fix compatibly a generator $\xi_m$ of $\mu_m(\overline{\mathbb F}_q)$ for each $m$ coprime to $q$, thus giving an isomorphism $\mathbb Q/\mathbb Z_\ell \to \mu_{\ell^\infty}$.  
 We now see how in the function field setting, the lifting invariant comes from Poincar\'e duality on the geometric curve.

\begin{lemma}\label{lifting-is-trace} Let $K$ be the function field of a connected, smooth, complete curve $C$ over a finite field $\mathbb F_q$ of characteristic $p$. 
Let $L/K$ be a finite Galois extension that is unramified everywhere, with $H=\Gal(L/K)$ of order prime to $p$.
Let $n$ be a divisor of $q-1$ and let $s_H \in  H_3(H, \mathbb  Z/n)$ be the Artin-Verdier trace for $L/K$, using $\xi=\xi_n$. 
Then the lifting invariant linear form is the composition
\[ H^2(H, \mathbb Q/\mathbb Z) \to H^2( C, \mathbb Q/\mathbb Z) \to H^2(C_{\overline{\mathbb F}_q},\mathbb Q/\mathbb Z) \to  \mathbb Q/\mathbb Z\]
where the first map is $\phi_{L/K}$ from Lemma~\ref{L:BHpullback},  the second map is the usual pullback, 
and  the last map is the sum over $\ell \neq p$ of $- \frac{q-1}{n}$ times the composite 
$H^2(C_{\overline{\mathbb F}_q},\mathbb Q_\ell/\mathbb Z_\ell)\to H^2(C_{\overline{\mathbb F}_q}, \mu_{\ell^\infty}) \to  \mathbb Q_\ell/\mathbb Z_\ell$
of the map induced by the isomorphism $\mathbb Q_\ell/\mathbb Z_\ell \to \mu_{\ell^\infty}$ sending $\ell^{-m}\mapsto \xi_m$ and the 
trace map associated to Poincar\'e duality in \'{e}tale cohomology. 
\end{lemma} 

We could include a term for $\ell=p$ in the sum above but since $H^2(C_{\overline{\mathbb F}_q},\mathbb Q_p/\mathbb Z_p)=0$ the zero map is the only possible choice of map $H^2(C_{\overline{\mathbb F}_q},\mathbb Q_p/\mathbb Z_p) \to  \mathbb Q_p/\mathbb Z_p.$  

\begin{proof} Every class in $H^2(H, \mathbb Q/\mathbb Z)$ breaks into a sum of $\ell$-power torsion classes for primes $\ell\neq p$, and every $\ell$-power torsion class in $H^2(H, \mathbb Q/\mathbb Z)$ arises from $H^2(H, \mathbb Z/\ell^m)$ for some positive integer $m$, where the map $H^2(H, \mathbb Z/\ell^m)\to H^2(H, \mathbb Q/\mathbb Z)$ is by division by $\ell^m$. It suffices to prove the stated identity for each class arising from $H^2(H, \mathbb Z/\ell^m)$. Let $v$ be the $\ell$-adic valuation of $n$ and let $ \tilde{n}$ be the inverse of $n/ \ell^{v}$ modulo $\ell^{m+v}$. We have a commutative diagram of short exact sequences
\[ \begin{tikzcd}  \mathbb Z/\ell^{v} \arrow[r , "\times \ell^m" ] \arrow[d, "\times \frac{n  \tilde{n}}{\ell^{v}}"] &  \mathbb Z/\ell^{m+{v}} \arrow[d, "\times \frac{\tilde{n} }{  \ell^{m+v }}"] \arrow[r, "\times 1"] & \mathbb Z/\ell^m \arrow[d, "\times \frac{1}{\ell^m}"] \\ \mathbb Z/n \arrow[r, "\times \frac{1}{n}"] & \mathbb Q/\mathbb Z \arrow[r, "\times n "]  & \mathbb Q/\mathbb Z \end{tikzcd} .\]

The lifting invariant on classes that arise from $H^2(H, \mathbb Z/\ell^m)$ is the composition of the following maps 
\[ \begin{tikzcd} 
H^2(H, \mathbb Z/\ell^m) \arrow[r, "1\mapsto \frac{1}{\ell^m}"]  & 
H^2(H,\Q/\Z) \arrow[r, "B'_n"] &
 H^3(H,\Z/n) \arrow[r, "\phi_{L/K}"] & 
  H^3(C,\Z/n) \arrow[r, "1\mapsto \xi"]  &
  H^3(C,\G_m) \arrow[r, "A"] &
  \Q/\Z,  
    \end{tikzcd} \]
    where $A$ is the basic Artin-Verdier isomorphism of Equation~\eqref{E:basic-AV}.
    By the commutative diagram above, we can replace the first two maps to obtain another description of the lifting invariant as
    \[ \begin{tikzcd} 
H^2(H, \mathbb Z/\ell^m) \arrow[r, "B_{m,v}"]  & 
H^3(H,\mathbb Z/\ell^v) \arrow[r, "1\mapsto \frac{n\tilde{n}}{\ell^v}"] &
 H^3(H,\Z/n) \arrow[r, "\phi_{L/K}"] & 
  H^3(C,\Z/n) \arrow[r, "1\mapsto \xi"]  &
  H^3(C,\G_m) \arrow[r, "A"] &
  \Q/\Z,  
    \end{tikzcd}\]
    where $B_{m,v}$ is the connecting homomorphism from the top row of the commutative diagram above.
By Lemma~\ref{L:BHpullback}, we can replace the second, third, and fourth maps to get another description of the lifting invariant as
    \[ \begin{tikzcd} 
H^2(H, \mathbb Z/\ell^m) \arrow[r, "B_{m,v}"]  & 
H^3(H,\mathbb Z/\ell^v) \arrow[r, "\phi_{L/K}"] &
 H^3(C,\Z/\ell^v)  \arrow[r, "1\mapsto \xi^{\frac{n\tilde{n}}{\ell^v}}"]  &
  H^3(C,\G_m) \arrow[r, "A"] &
  \Q/\Z, 
    \end{tikzcd} \]
and then replace the first and second maps to get another description of the lifting invariant as
    \[ \begin{tikzcd} 
H^2(H, \mathbb Z/\ell^m) \arrow[r, "\phi_{L/K}"]  & 
H^2(C, \mathbb Z/\ell^m) \arrow[r, "B_{m,v}"] &
 H^3(C,\Z/\ell^v)  \arrow[r, "1\mapsto \xi^{\frac{n\tilde{n}}{\ell^v}}"]  &
  H^3(C,\G_m) \arrow[r, "A"] &
  \Q/\Z.
    \end{tikzcd} \]
    Let $\zeta_m$ be the image of $\xi^{-\frac{n}{\ell^v}}=\xi_{\ell^v}^{-1}$ in $\Hf^1 (C , \mu_{\ell^m})$ under the connecting homomorphism of the exact sequence $\mu_{\ell^m} \to \mu_{\ell^{m+v}} \to \mu_{\ell^v}.$ By Lemma \ref{exact-sequence-duality}, we get another description of the lifting invariant as
    \[ \begin{tikzcd} 
H^2(H, \mathbb Z/\ell^m) \arrow[r, "\phi_{L/K}"]  & 
H^2(C, \mathbb Z/\ell^m) \arrow[r, "\cup \zeta_m^{\tilde{n}}"] &
  H^3(C,\G_m) \arrow[r, "A"] &
  \Q/\Z.
    \end{tikzcd} \]



By \cite[Proof of Lemma 7.5]{Lipnowski2020},  using that 
$\xi_{n}=\xi$ and the compatibility of the $\xi_{k}$,
we can replace the last two maps and the lifting invariant is
    \[ \begin{tikzcd} 
H^2(H, \mathbb Z/\ell^m) \arrow[r, "\phi_{L/K}"]  & 
H^2(C, \mathbb Z/\ell^m) \arrow[r] &
H^2(C_{\overline{\mathbb F}_q}, \mathbb Z/\ell^m) \arrow[r, "1\mapsto \xi_{\ell^m}"] &
H^2(C_{\overline{\mathbb F}_q}, \mu_{\ell^\infty}) \arrow[r, "P.D. \times{\frac{\tilde{n}(1-q)}{\ell^v}}"] &
  \Q_\ell/\Z_\ell.
    \end{tikzcd} \]    
Since the map lands in $\ell^m$ torsion, and $n\tilde{n}$ is $1$ mod $\ell^{m+v}$, we can replace the last multiplication by ${\frac{\tilde{n}(1-q)}{\ell^v}}=\frac{\tilde{n}n}{\ell^v}{\frac{(1-q)}{n}}$ with multiplication by ${\frac{(1-q)}{n}}$.

The lemma now follows from the fact that the maps induced by the map of coefficients $\Z/\ell^m \ra \Q/\Z$ sending $1\mapsto \ell^{-m}$ commutes with
$\phi_{L|K}$ (Lemma~\ref{L:BHpullback}) and with the pull back induced by $C_{\overline{\mathbb F}_q}\ra C$.
\end{proof}

This lemma allows us to compare our lifting invariant with the notion of the lifting invariant defined by Liu~\cite{Liu2022}, which contains essentially the same information \cite[Proposition 4.3]{Liu2022} as the lifting invariant defined by the second author~\cite{Wood2021} based on work of Ellenberg, Venkatesh, and Westerland, as well as the earlier lifting invariant for quadratic extensions \cite[Theorem 3.13]{Wood2019}.

 We next review Liu's notion. 

Liu considers a Galois extension $K/\mathbb F_q(t)$ with Galois group $\Gamma$, corresponding to a curve $X/\mathbb P^1_{{\mathbb F}_q}$ of positive genus. She defines $K_{\emptyset}' $ to be the maximal extension of $K \overline{\mathbb F}_q$ that is everywhere unramified and a limit of finite Galois extensions of degree prime to $\abs{\Gamma} q$.  Liu first defines an isomorphism
\[ \omega_K \colon \hat{\mathbb Z}(1)_{ (q \abs{\Gamma})'} \to H_2( \Gal ( K_\emptyset'/ \overline{\mathbb F}_q(t)), \mathbb Z)_{ (q \abs{\Gamma})'} \]
which is characterized \cite[Corollary 2.12]{Liu2022} by the following property: the inverse map $  \omega_K^{-1} \colon H_2( \Gal ( K_\emptyset'/ \overline{\mathbb F}_q(t)), \mathbb Z)_{ (q \abs{\Gamma})'}\to \hat{\mathbb Z}(1)_{ (q \abs{\Gamma})'} $ corresponds to a class in $H^2( \Gal ( K_\emptyset'/ \overline{\mathbb F}_q(t)) , \hat{ \mathbb Z}(1))_{(q \abs{\Gamma})'}$ whose image in the \'etale cohomology $H^2(X_{\overline{\mathbb F}_q}, \hat{\Z}(1)_{(q |\Gamma|)'})$ has Poincar\'{e} duality trace $-\abs{\Gamma}\in \hat{\Z}_{(q \abs{\Gamma})'}$. 
(The subscripts denote a pro-prime-to-$q|\Gamma|$ completion.)

For $L$ a Galois extension of $K$, also Galois over $\F_q(t)$, that is everywhere unramified, split completely over $\infty$, and has degree prime to $q \abs{\Gamma}$,  Liu 
\cite[Definition 2.13, Remark 2.16(2)]{Liu2022}
defines a homomorphism
\[ \omega_{L/K} \colon  \hat{\mathbb Z}(1)_{ (q \abs{\Gamma})'} \to_{\omega_K}  H_2( \Gal ( K_\emptyset'/ \overline{\mathbb F}_q(t)), \mathbb Z)_{ (q \abs{\Gamma})'}  \to H_2(\Gal(L/ \mathbb F_q(t)), \mathbb Z) _{(q\abs{\Gamma})'} \]
where the second map is induced by the group homomorphism $\Gal(K_\emptyset'/ \overline{\mathbb F}_q(t) ) \to \Gal(L/\mathbb F_q(t))$ arising from the inclusions $\mathbb F_q(t) \subset \overline{\mathbb F}_q(t)$ and $L  \overline{\mathbb F}_q(t) \subseteq  K_{\emptyset}'$. 

\begin{lemma}\label{yuan-liu} 
Let $K/\mathbb F_q(t)$ be a  Galois extension with Galois group $\Gamma$, corresponding to a curve $X/\mathbb P^1_{{\mathbb F}_q}$ of positive genus. 
Let $L$ be a Galois extension of $K$, also Galois over $\F_q(t)$,  that is everywhere unramified, split completely over $\infty$, and has degree prime to $q \abs{\Gamma}$. Let $n$ be the maximal divisor of $q-1$ prime to $\abs{\Gamma}$. Our fixed system of roots of unity gives an isomorphism $\hat{\mathbb Z}_{ (q\abs{\Gamma})'} \to \hat{ \mathbb Z} (1)_{ (q\abs{\Gamma})'}$.  Let $\underline{\xi}$ be the image of $1$ under this isomorphism, with image $\xi_m$ in $\mu_m$ for all $m$.  Assume $\xi_n=\xi$. 

  Then $\frac{ q-1}{n} \abs{\Gamma} \omega_{L/K}(\underline{\xi}) $  is equal to the image of the lifting invariant in  $H_2(\Gal(L/K), \mathbb Z)$ of the Artin-Verdier trace in $H_3(\Gal(L/K), \mathbb Z/n)$ under the natural map $ H_2(\Gal(L/K),\mathbb Z)\to H_2(\Gal(L/ \mathbb F_q(t)), \mathbb Z)$.
   \end{lemma}

\begin{proof}
For an abelian group $A$, we write $A^{{(q\abs{\Gamma})'}}$ to denote the subgroup of $A$ of all elements of finite order that is relatively prime to $q|\Gamma|$.  
 To check that two elements of $H_2(\Gal(L/ \mathbb F_q(t)), \mathbb Z) _{(q\abs{\Gamma})'} $ agree, it suffices to check that their pairings with  an element $\alpha$ of $H^2( \Gal(L/\mathbb F_q(t)), \mathbb Q/\mathbb Z)^{(q\abs{\Gamma})'}$ agree.  The pairing of $\omega_{L/K}(\underline{\xi})$ and $\alpha$ may be calculated by pulling back $\alpha$ to $\hat{\alpha}\in H^2 ( \Gal( K_\emptyset' / \overline{\mathbb F}_q(t)) , \mathbb Q/\mathbb Z)^{(q \abs{\Gamma})'}$ and then pairing with $\omega_K(\underline{\xi})$.

Let $\epsilon \in H^2( \Gal ( K_\emptyset'/ \overline{\mathbb F}_q(t) ), \hat{ \mathbb Z}(1)_{(q \abs{\Gamma})'})$ be a class corresponding to $\omega_K^{-1}$ as in \cite[Remark 2.7]{Liu2022}.  In other words, this class has the property that for $a \in  \hat{\mathbb Z}(1)_{ (q \abs{\Gamma})'} $ we have $\langle \omega_K(a), \epsilon \rangle = a$ where $\langle, \rangle $ denotes the pairing $H_2(\mathbb Z) \times H^2( A ) \to A$ for any abelian group $A$. Hence if $f\colon\hat{ \mathbb Z}(1)_{(q \abs{\Gamma})'} \to  \mathbb Q/\mathbb Z$  is a homomorphism, we have \[ \langle \omega_K(a), f_*(\epsilon) \rangle = f ( \langle \omega_K(a),\epsilon \rangle)=f(a).\]

Since $H^2 ( \Gal( K_\emptyset' / \overline{\mathbb F}_q(t)) , \mathbb Q/\mathbb Z)^{(q \abs{\Gamma})'}$ is a direct limit of group cohomology of finite groups with $\Q/\Z$ coefficients, each of which is a direct limit of group cohomology with finite coefficients, we have that $\hat{\alpha}$ is in the image
of $H^2( \Gal( K_\emptyset' / \overline{\mathbb F}_q(t)), \mu_m)$ for some positive integer $m$, relatively prime to $q|\Gamma|$, and injective map $\mu_m\ra\Q/\Z$.  By \cite[Proof of Lemma 2.2]{Liu2022}, we have that the map induced by the natural quotient on coefficients $H^2( \Gal( K_\emptyset' / \overline{\mathbb F}_q(t)), \hat{ \mathbb Z}(1)_{(q \abs{\Gamma})'})\ra H^2( \Gal( K_\emptyset' / \overline{\mathbb F}_q(t)), \mu_m)\isom \Z/m\Z$ is surjective.  By \cite[Corollary 2.12]{Liu2022}, since $\epsilon$ has trace which is a unit in $\hat{ \mathbb Z}_{(q \abs{\Gamma})'}$, it is a topological generator of $H^2( \Gal ( K_\emptyset'/ \overline{\mathbb F}_q(t) ), \hat{ \mathbb Z}(1)_{(q \abs{\Gamma})'})$.  Thus, for some integer $k$, the map on $f$ coefficients that is the natural quotient  $\hat{ \mathbb Z}(1)_{(q \abs{\Gamma})'}\ra \mu_m$ followed by multiplication by $k$ and then our map $\mu_m\ra \Q/\Z$, induces a map on cohomology such that  
$f_*(\epsilon)=\hat{\alpha}$.

Thus the pairing of $\alpha$ with $\omega_{L/K}(\underline{\xi})$  is equal to $\langle \omega_K(\underline{\xi}), \hat{\alpha}\rangle=\langle \omega_K(\underline{\xi}), f_*(\epsilon)\rangle$ and thus to $f(\underline{\xi})$ for this value of $f$. To calculate $f(\underline{\xi})$, we use that $\epsilon$ has Poincar\'{e} duality trace $-\abs{\Gamma}$ so that the image of $f_*(\epsilon)\in H^2( \Gal ( K_\emptyset'/ \overline{\mathbb F}_q(t) ), \mathbb Q/\mathbb Z)$,
under the map on coefficients $\phi: \Q/\Z \ra \prod_{\ell\nmid q |\Gamma|} \mu_{\ell^\infty}$ sending $1/m$ to $\xi_m$,  has 
Poincar\'{e} duality trace $-|\Gamma|f(\underline{\xi})$ for the element $\underline{\xi}\in \hat{\Z}(1)_{(q|\Gamma|)'}$ with image $\xi_m$ in $\mu_m$ for all $m$.
In other words, $\phi_*(\hat{\alpha})$ has Poincar\'{e} duality trace $-|\Gamma|f(\underline{\xi})$.
Thus the pairing of $\alpha$ with $\frac{ q-1}{n} \abs{\Gamma} \omega_{L/K}(\underline{\xi})$  is  $\frac{ q-1}{n} \abs{\Gamma}f(\underline{\xi}), $ which is $-\frac{ q-1}{n}$ times the Poincar\'{e} duality trace of $\phi_*(\hat{\alpha})$.

By Lemma \ref{lifting-is-trace},  $-\frac{ q-1}{n}$ times the Poincar\'{e} duality trace of $\phi_*(\hat{\alpha})$ is precisely the evaluation of the lifting invariant linear form on the image of $\alpha$ in $H^2( \Gal(L/K), \mathbb Q/\mathbb Z)$.  This evaluation is also equal to the evaluation on $\alpha$ itself of
the linear form $H^2(\Gal(L/ \mathbb F_q(t)), \Q/\Z)\ra \Q/\Z$ corresponding to the image of the lifting invariant under 
the natural map $ H_2(\Gal(L/K),\mathbb Z)\to H_2(\Gal(L/ \mathbb F_q(t)), \mathbb Z)$.   
We conclude that $\frac{ q-1}{n} \abs{\Gamma} \omega_{L/K}(\underline{\xi})$ is the  image of the lifting invariant under 
the natural map $ H_2(\Gal(L/K),\mathbb Z)\to H_2(\Gal(L/ \mathbb F_q(t)), \mathbb Z)$, since they have the same pairing with any $\alpha$.
\end{proof}

Since $\frac{ q-1}{n} \abs{\Gamma}$ is divisible only by primes dividing $\Gamma$, multiplying by it is invertible on any $(q\abs{\Gamma})'$-module, so Lemma \ref{yuan-liu} shows that the image under $\omega_{L/K} $ of a generator, and hence $\omega_{L/K}$ itself, is uniquely determined by the Artin-Verdier trace.

In fact, even more precisely, we have the following.

\begin{lemma}\label{liu-domain-comparison} In the setup of Lemma \ref{yuan-liu}, the natural map \[H_2(\Gal(L/K), \mathbb Z)^{\Gamma} \to \Hom ( \hat{\mathbb Z}(1)_{ (q \abs{\Gamma})'}, H_2(\Gal(L/ \mathbb F_q(t)), \mathbb Z) _{(q\abs{\Gamma})'}  )\] that sends a class in $H_2(\Gal(L/K),\mathbb Z)$ to the unique homomorphism sending a fixed topological generator of  $\hat{\mathbb Z}(1)_{ (q \abs{\Gamma})'}$ to the image of that class under the inflation map $H_2(\Gal(L/K), \mathbb Z) \to H_2(\Gal(L/ \mathbb F_q(t)), \mathbb Z) _{(q\abs{\Gamma})'} $ is an isomorphism. \end{lemma}

In particular, there is an isomorphism between the set in which our lifting invariant lies and the set in which $\omega_{L/K}$ lies that sends one to the other.

\begin{proof} Since $\abs{\Gal(L/K)}$ is prime to $q\abs{\Gamma}$, the order of $H_2(\Gal(L/K), \mathbb Z)$ is prime to $q\abs{\Gamma}$ and so  $H_2(\Gal(L/K), \mathbb Z)_{(q\abs{\Gamma})'} \to H_2(\Gal(L/K), \mathbb Z)$ is an isomorphism.  The inflation map $H_2(\Gal(L/K), \mathbb Z)^{\Gamma}_{(q\abs{\Gamma})'}  \to H_2(\Gal(L/ \mathbb F_q(t)), \mathbb Z) _{(q\abs{\Gamma})'} $ is an isomorphism by taking the Lyndon-Hochschild-Serre spectral sequence computing $ H_2(\Gal(L/ \mathbb F_q(t)), \mathbb Z)$  from $H^a (\Gamma, H^b ( \Gal(L/K), \mathbb Z) $ and observing that all terms with $a>0$ are $\abs{\Gamma}$-torsion so after taking $q \abs{\Gamma}'$-parts the spectral sequence degenerates on the second page.

Finally, for $A$ a finite group of order prime to $q \abs{\Gamma}$, the map $ A \to \Hom (  \hat{\mathbb Z}(1)_{ (q \abs{\Gamma})'} , A)$ sending an element to the unique homomorphism sending a generator to that element is an isomorphism since it is constructed as the inverse of the map $ \Hom (  \hat{\mathbb Z}(1)_{ (q \abs{\Gamma})'} , A) \to A$ given by evaluating at a fixed generator which is an isomorphism since every element of $A$ generates a finite cyclic group of order prime to $q\abs{\Gamma}$ and this finite cyclic group is a quotient of  $\hat{\mathbb Z}(1)_{ (q \abs{\Gamma})'} $.\end{proof}

\subsection{2-class field towers of cyclic cubic fields, after Boston-Bush}\label{ss-bb}

In this subsection, we give an example of how to derive explicit probabilities from our conjectures. Our goal is to compute probabilities that the Galois group of the $2$-class field tower (i.e. the composition of all everywhere unramified Galois extensions of degree a power of $2$) of a cyclic cubic field is equal to a given $2$-group. In particular, we will do this for two of the simplest possible $2$-groups, the Klein four group and the eight-element quaternion group. Our choice of this example is motivated by work of Boston and Bush~\cite{Boston2021a}, who collected numerical data and theoretical results on this problem.  In particular, we observe that our conjectural probabilities, based on function field theorems and the moment method, agree with their conjectural probabilities, based on numerical evidence in the number field case. This provides some basic numerical evidence for our conjectures. In addition, we explain the relationship between our theorems and conjectures and some existential theorems and conjectures made by Boston and Bush.

We begin with the following result, which uses our calculation of the probability of obtaining a given oriented $\Gamma$-group to give a formula for the probability of obtaining a given $\Gamma$-group without orientation. In the following statement we add the parenthetical notation $(s_H)$ to $w_{V_i}$ and $w_{N_i}$ to make clear that they depend on $s_H$, though in the rest of the paper we work with a fixed $(H,s_H)$ and thus drop this notation for simplicity.
 
\begin{lemma}\label{prob-without-orient} Let $U$ be a multiset of elements of $\Gamma$. Assume all nonzero representations of $\Gamma$ of characteristic dividing $n$ contain some nontrivial vector fixed by at least one element of $U$. Let $\nu$ be the measure of Theorem \ref{T:measure}. Let $\cL$ be a level in the category of finite $\Gamma$-groups.  Let $H$ be a finite $\Gamma$-group in $\cL$ with $H_\Gamma= 1$. Then \[ \nu( \{ X \mid X^{\mathcal L} \cong H\}) = \frac{ |H^2(H \rtimes \Gamma, \mathbb Z/n) |}{ |H^3(H \rtimes \Gamma, \mathbb Z/n) | |\Aut_\Gamma(H)| H^{\cdot U} } \sum_{ s_H \in H_3 (H , \mathbb Z/n)^{\Gamma}}\prod_{ i=1}^{r} w_{V_i}(s_H)
 \prod_{i=1}^s  w_{N_i}(s_H).\] 
If $H_\Gamma\ne 1$, then $\nu( \{ X \mid X^{\mathcal L} \cong H\})=0$.
\end{lemma}

\begin{proof}  Let $\cW$ be the set of isomorphism classes of finite $n$-oriented $\Gamma$-groups whose underlying $\Gamma$-group is in $\cL$. Then $X^{ \cL } \cong H$ if and only if $X^{\cW} \cong (H, s_H)$ for some $s_H \in H_3 (H , \mathbb Z/n)^{\Gamma}$. Since the condition $X^{\cW} \cong (H, s_H)$ depends only on the isomorphism class of $(H,s_H)$, we have by Theorem \ref{T:measure}  and Corollary \ref{C:vfactor}
\[ \nu ( \{X \mid X^{\cL} \cong H \}) =\nu \Bigl( \bigcup_{ s_H\in H_3 (H , \mathbb Z/n)^{\Gamma} / \Aut_\Gamma(H)} \{ X \mid X^{\cW} \cong (H,s_H) \}\Bigr )\] \[ = \sum_{ s_H\in H_3 (H , \mathbb Z/n)^{\Gamma} / \Aut_\Gamma(H)} \nu( \{ X \mid X^{\cW} \cong (H,s_H) \}  ) \] \[=\sum_{ s_H\in H_3 (H , \mathbb Z/n)^{\Gamma} / \Aut_\Gamma(H)} v_{\cW, (H,s_H)}  = \sum_{ s_H\in H_3 (H , \mathbb Z/n)^{\Gamma} / \Aut_\Gamma(H)} \frac{M_\bH}{|\Aut(H,s_H)|}  
\prod_{ i=1}^{r} w_{V_i}(s_H )
 \prod_{i=1}^s  w_{N_i}(s_H ).\]
 
 When $H_\Gamma\ne 1$, then $M_{\bH} =0$, and otherwise
 \[ M_{\bH} =  \frac{ |H^2(H \rtimes \Gamma, \mathbb Z/n) |}{ |H^3(H \rtimes \Gamma, \mathbb Z/n) |H^{\cdot U} } \] depends only on $H$ and therefore can be brought outside the sum over $s_H$. Furthermore we have 
 \[ \sum_{ s_H\in H_3 (H , \mathbb Z/n)^{\Gamma} / \Aut_\Gamma(H)} \frac{1}{|\Aut(H,s_H)|}  
\prod_{ i=1}^{r} w_{V_i}(s_H )
 \prod_{i=1}^s  w_{N_i}(s_H ) \] \[= \frac{1}{|\Aut_\Gamma(H)|} \sum_{ s_H\in H_3 (H , \mathbb Z/n)^{\Gamma}}  
\prod_{ i=1}^{r} w_{V_i}(s_H )
 \prod_{i=1}^s  w_{N_i}(s_H ) \]
 since $\Aut(H,s_H)$ is the stabilizer of $s_H$ in $\Aut_\Gamma(H)$ so by the orbit-stabilizer theorem each $\Aut(H)$-orbit contributes \[\frac{| \Aut_\Gamma(H)| }{|\Aut(H,s_H)|}  
\prod_{ i=1}^{r} w_{V_i}(s_H )
 \prod_{i=1}^s  w_{N_i}(s_H )\] to the sum $ \sum_{ s_H\in H_3 (H , \mathbb Z/n)^{\Gamma}}  
\prod_{ i=1}^{r} w_{V_i}(s_H )
 \prod_{i=1}^s  w_{N_i}(s_H)$. Combining these, we obtain the statement. \end{proof}

 The following lemma lets us see what terms appear in the products.
 
 \begin{lemma}\label{p-group-simplifications} Let $p$ be a prime not dividing $|\Gamma|$. Let $\cL$ be a level in the category of finite $\Gamma$-groups consisting only of $p$-groups with an action of $\Gamma$.  Let $H$ be a finite $\Gamma$-group in $\cL$. Then:
 
 \begin{enumerate}
 
 \item  There are no admissible nonabelian simple $[H \rtimes \Gamma]$-groups $N_i$ and the only abelian simple $H \rtimes \Gamma$-modules $V_i$ are the groups of the form $\mathbb F_p^d$ with a trivial $H$-action and an irreducible action of $\Gamma$.
 
 \item A group $V_i$ of the form $\mathbb F_p^d$ with  a trivial $H$-action and an irreducible action of $\Gamma$ is an admissible finite simple $H \rtimes \Gamma$-group if and only if $V_i $ is in $\cL$. 
 
\item  An admissible finite simple $H \rtimes \Gamma$-group $V_i$ is never intermediate.\end{enumerate} \end{lemma}
 
 \begin{proof} For a nonabelian simple $[H \rtimes \Gamma]$-group $N_i$ to be admissible, an extension of the form $H \times_{ \Out(N_i) } \Aut(N_i)$ of $H$ by $N_i$ must be in $\cL$ and in particular must be a $p$-group, so $N_i$ must be a $p$-group. But the commutator subgroup of a non-abelian $p$-group is a nontrivial proper normal subgroup which is characteristic and hence $[H \rtimes \Gamma]$-invariant, so this contradicts simplicity. 
 
 For an abelian simple $H \rtimes \Gamma$-group $V_i$ to be admissible, $V_i \rtimes H$ must be a $p$-group, so $V_i$ must be a $p$-group. If $H$ acts nontrivially on $V_i$ then the maximal subgroup on which $H$ acts trivially is a nontrivial proper subgroup of $V_i$, fixed as a set under the action of $H\rtimes \Gamma$, contradicting simplicity, so $H$ acts trivially and the action of $H \rtimes \Gamma$ factors through $\Gamma$. This verifies (1).

If $V_i$ is an admissible finite simple $H\rtimes\Gamma$-group,  then $V_i \rtimes H$ in $\cL$.  If the action of $H$ on $V_i$ is trivial, we have $V_i \rtimes H=V_i \times H$. Since $ V_i$ is a quotient of $V_i \times H$, we have that 
$V_i \times H \in 
 \cL$ implies
$V_i \in \cL$.   On the other hand, if $V_i\in\cL$, then since $V_i\times H$ is a fiber product of $V_i$ and $H$, we have $V_i\times H \in \cL$ as long as $H \in \cL$, which is true by assumption. This verifies (2).

 The intermediate case occurs only if $p=2$, and if the natural homomorphism $H \rtimes \Gamma\to \Aut(V_i)$ factors through the symplectic group but not the affine symplectic group. Since the homomorphism $H \rtimes \Gamma\to \Aut(V_i)$ factors through $\Gamma$, this can happen only if $\Gamma \to \Aut(V_i)$ factors through the symplectic group but not the affine symplectic group. The affine symplectic group is an extension of the symplectic group by an abelian $2$-group, and thus every map from a group with odd order to the symplectic group lifts to the affine symplectic group. Since $p=2$, we have that $|\Gamma|$ is odd, and thus every map from $\Gamma$ to the symplectic group lifts to the affine symplectic group. This verifies (3).
 
 \end{proof}
 
 The next lemma lets us apply our formulas to calculate the probability that the maximum $p$-group quotient of a random group is a given group.  
 
 \begin{lemma}\label{maximal-p-group-quotient} Let $p$ be a prime not dividing $|\Gamma|$. Let $H$ be a finite $p$-group with an action of $\Gamma$. Let $\cL$ be a level in the category of $\Gamma$-groups consisting only of $p$-groups with an action of $\Gamma$, including $H$ and every extension of $H$ by $V_i$ for $V_i$ any group of the form $\mathbb F_p^d$ with a trivial action of $H$ and an irreducible action of $\Gamma$. Then for $X$ a profinite $\Gamma$-group, the maximal $p$-group quotient of $X$ is isomorphic to $H$ if and only if $X^{\cL} \cong H$. 
 
 Furthermore, in this case, for any group $V_i$ of the form $\mathbb F_p^d $ with an irreducible action of $\Gamma$, we have  $H^2(H \rtimes \Gamma, V_i)^{\cL} = H^2(H \rtimes \Gamma, V_i)$.  \end{lemma}
 
 \begin{proof} Certainly if the maximal $p$-group quotient of $X$ is isomorphic to $H$, $\cL$ consists only of $p$-groups, and $H \in \cL$, then the maximal quotient of $X$ in $\cL$ is isomorphic to $H$. Conversely, if $X^{\cL} \cong H$ but $H$ is not the maximal $p$-group quotient of $X$, there exists another finite $p$-group with an action of $\Gamma$ that admits a surjection from $G$ and a surjection to $H$ with nontrivial kernel. Let $G'$ be the quotient of $G$ by a maximal proper normal subgroup of the kernel. Then the kernel of $G' \to H$ is a simple $G'\rtimes \Gamma$-group, which since $G$ and the kernel are both $p$-groups implies that the kernel is a group of the form $\mathbb F_p^d$ with an irreducible action of $\Gamma$. Hence $G'\in \cL$ which contradicts $X^{\cL} \cong H$.
 
 A class is in $H^2(H \rtimes \Gamma, V_i)^{\cL}$ if and only if the associated extension of $H$ by $V_i$ is in $\cL$, which is true for all classes in $H^2(H \rtimes \Gamma, V_i)$ by assumption. \end{proof}

Boston and Bush~\cite{Boston2021a} studied the statistics of the $2$-class field towers of cyclic cubic fields. In other words, these are the composition of all everywhere unramified Galois extensions of degree a power of $2$.  They consider two variants, one in which $\mathbb C/\mathbb R$ is treated as an unramified extension (the narrow class field tower) and one in which $\mathbb C/\mathbb R$ is treated as ramified (the wide class field tower). In both variants, they proved theorems restricting which groups can appear as the Galois group of the $2$-class field tower, made conjectures about the set of groups that appear, and made conjectures, based on numerical data, about the frequency with which certain small finite groups appear. These conjectures are focused on the special case where the Galois group of the $2$-class field has exactly two generators (equivalently, the class group of the cyclic cubic field has $2$-rank two).

Motivated  by this work, let us specialize further to $\Gamma=\mathbb Z/3$, $n=2$, and $U= \{ 1\}$. These are the correct values to match $k=\mathbb Q$, since $\mathbb Q$ has $2$ roots of unity and the only possible choice of $\gamma_\infty$ in $\mathbb Z/3$ of order $1$ or $2$ is $1$.  For $k=\mathbb Q$ and $\Gamma=\mathbb Z/3$, the $\Gamma$-extensions $K/k$ we consider are cyclic cubic fields. 

The groups of the form $\mathbb F_2^d$ with an irreducible action of $\Gamma$ are $V_1 := \mathbb F_2$ with a trivial action of $\Gamma$ and $V_2 := \mathbb F_2^2$ with the unique-up-to-isomorphism nontrivial action of $\Gamma$. 
Let $H$ be a finite $\Gamma$-group $2$-group such that $H_\Gamma=1$.
Let $\cL$ be a level containing only of $2$-groups and  containing $H$ and every $\Gamma$-group extension of $H$ by $V_1$ or $V_2$. 

The formula of Lemma \ref{prob-without-orient} specializes to 

 \begin{equation}\label{cct-prob} \nu( \{ X \mid X^{\mathcal L} \cong H\}) = \frac{ |H^2(H \rtimes \Gamma, \mathbb F_2) | |H^\Gamma| }{ |H^3(H \rtimes \Gamma,  \mathbb F_2) | |\Aut_\Gamma(H)| |H|  } \sum_{ s_H \in H_3 (H ,  \mathbb F_2)^\Gamma} w_{V_1}(s_H) w_{V_2} (s_H) .\end{equation} 

We now calculate $w_{V_1}$ and $w_{V_2}$. 

We have that $V_1 = \mathbb F_2$ is trivial of characteristic dividing $n$. Thus we apply Lemma \ref{L:wntrivial} to calculate $w_{V_1} (s_H) $. We have $\kappa_1=\mathbb F_2$ and $q_1=2$. We have \[V_1^{\cdot U'} = \frac{ |\mathbb F_2|}{ |\mathbb F_2^{H\rtimes \Gamma}|}=\frac{2}{2}=1.\] Thus for $u_1=0$ we have $q_1^{u_1} = V_1^{\cdot U'} $.  Letting \[ z_1(s_H)  = \dim_{\kappa_1} H^2( H \rtimes \Gamma, V_1)^{\cL,s_H}= \dim_{\mathbb F_2} H^2(H \rtimes \Gamma, \mathbb F_2)^{\cL, s_H}\] we have 
\begin{equation}\label{wv1} w_{V_1}(s_H)  = \prod_{k=0}^{z_1(s_H) -1} ( 1- q_1^{k-u_1})= \prod_{k=0}^{z_1(s_H) -1} (1-2^k) = \begin{cases} 1 & \textrm{if } z_1(s_H) =0 \\ 0 & \textrm{otherwise} \end{cases} = \begin{cases} 1 & \textrm{if } H^2(H \rtimes \Gamma, \mathbb F_2)^{\cL, s_H}=0 \\ 0 & \textrm{otherwise} \end{cases}.\end{equation}

We now consider $V_2$. By Lemma \ref{p-group-simplifications}(3),  $V_2$ is not intermediate and thus not anomalous.  We have $\kappa_2=\mathbb F_4$ and $q_2=4$. The representation $V_2$ is unitary and $\mathbb F_2$-orthogonal, so $\epsilon_{V_2}=0$. We have \[ V_2^{ \cdot U'} = \frac{ |V_2|}{|V_2^{H \rtimes \Gamma}|}=\frac{4}{1}=4.\] Thus for $u_2=1$ we have $q_2^{u_2} = V_2^{\cdot U'}$. Letting \[ z_2(s_H)  = \dim_{\kappa_2}  H^2( H \rtimes \Gamma, V_2)^{\cL,s_H}= \dim_{\mathbb F_4} H^2(H \rtimes \Gamma, \mathbb F_2)^{\cL, s_H}\] we have
by Lemma~\ref{odd-selfdual-nontrivial-w}
\[ w_{V_2} (s_H) = \prod_{k=0}^{\infty} (1+ q_2^{ -k - \frac{\epsilon_{V_2} +1 }{2} - u_2})^{-1} \prod_{k=0}^{z_2(s_H) -1} (1-q_2^{k-u_2}) = \prod_{k=0}^{\infty} (1+ 4^{ -k - \frac{1}{2} - 1})^{-1}\prod_{k=0}^{z_2(s_H) -1} (1-4^{k-1})\] \begin{equation}\label{wv2} = \prod_{k=0}^{\infty} (1+ 4^{ -k - \frac{3}{2} })^{-1} \times \begin{cases} 1 & \textrm{if } z_2(s_H) =0 \\ \frac{3}{4} & \textrm{if } z_2(s_H) =1 \\ 0 & \textrm{otherwise} \end{cases} .\end{equation}

Combining \eqref{cct-prob}, \eqref{wv1}, and \eqref{wv2}, we see that the probability under our conjectured distribution that the $2$-class field tower of a cyclic cubic field has Galois group $H$ is positive if and only if  there exists $s_H\in H_3(H , \mathbb F_2)^\Gamma$ such that two conditions are satisfied: First that the composition \[ H^2( H \rtimes \Gamma, \mathbb F_2) \to H^3(H \rtimes \Gamma, \mathbb F_2) \to \mathbb F_2\] of the Bockstein homomorphism with the map induced $s_H$ is injective, and second that the kernel in $ H^2( H\rtimes \Gamma, V_2) $  of the pairing \[H^1(H \rtimes \Gamma, V_2) \times H^2( H\rtimes \Gamma, V_2) \to H^3( H \rtimes \Gamma, \mathbb F_2) \to \mathbb F_2 \] arising from cup product and $s_H$ has size at most $4$.

The first of these conditions, in different language, was highlighted in the work of Boston and Bush.

\begin{lemma}\label{schur-is-v1} Let $\Gamma = \mathbb Z/3$ and let $H$ be a $2$-group with an action of $\Gamma$ such that $H_\Gamma=1$.

There exists $s_H \in H_3( H, \mathbb F_2)^\Gamma$ such that the composition $H^2( H \rtimes \Gamma, \mathbb F_2) \to H^3(H \rtimes \Gamma, \mathbb F_2) \to \mathbb F_2$ of the Bockstein homomorphism with the pairing with $s_H$ is injective if and only if the Schur multiplier of $H \rtimes \Gamma$ has order $1$ or $2$. \end{lemma}

\begin{proof} Since $H_*( H\rtimes \Gamma, \mathbb F_2)$ and $H^*(H\rtimes \Gamma, \mathbb F_2)$ are dual vector spaces,  there exists $s_H$ such that the composition $H^2( H \rtimes \Gamma, \mathbb F_2) \to H^3(H \rtimes \Gamma, \mathbb F_2) \to \mathbb F_2$ is injective if and only if $H^2( H \rtimes \Gamma, \mathbb F_2) $ has order at most $2$ and the natural map $H^2( H \rtimes \Gamma, \mathbb F_2) \to H^3(H \rtimes \Gamma, \mathbb F_2)$ is injective. Since the Bockstein sequence is exact, this occurs if and only if $H^2( H \rtimes \Gamma, \mathbb F_2) $ has order at most $2$ and the natural map $H^2( H \rtimes \Gamma, \mathbb Z/4 ) \to H^2( H \rtimes \Gamma, \mathbb F_2) $ is zero. Since $H_\Gamma=1$ we have $H_1(H\rtimes \Gamma, \mathbb Z) = \mathbb Z/3$ and thus for every finite group $A$ of order prime to $3$, $H^2( H \rtimes \Gamma, A) = \Hom ( H_2(H \rtimes \Gamma, \mathbb Z), A)$. 

The map $ \Hom ( H_2(H \rtimes \Gamma, \mathbb Z), \mathbb Z/4)\to  \Hom ( H_2(H \rtimes \Gamma, \mathbb Z), \mathbb F_2)$ is zero if and only if the $2$-Sylow subgroup of $ H_2(H \rtimes \Gamma, \mathbb Z)$ is $2$-torsion, and additionally $\Hom ( H_2(H \rtimes \Gamma, \mathbb Z), \mathbb F_2)$ has order at most $2$ if and only if the $2$-Sylow subgroup of  $ H_2(H \rtimes \Gamma, \mathbb Z)$ has order at most $2$. But for $p\neq 2,3$ the $p$-Sylow subgroup of $ H_2(H \rtimes \Gamma, \mathbb Z)$ is trivial since $H\rtimes \Gamma$ has order prime to $p$, and the $3$-part of $H_2(H\rtimes \Gamma, \mathbb Z)$ is the $3$-part of $H_2(\Gamma, \mathbb Z)$ which vanishes, so this occurs if and only if $ H_2(H \rtimes \Gamma, \mathbb Z)$ has order at most $2$.\end{proof}

In \cite[Theorem 5.2]{Boston2021a}, Boston and Bush prove that every 2-generator $2$-group that arises as the wide $2$-class tower group of a cyclic cubic field is $2$-select. Here $2$-select is a somewhat complicated group-theoretic notion. One first defines special groups~\cite[Definition 4.3]{Boston2021a} as 2-generator $2$-groups $H$ with an action of $\mathbb Z/3$ such that the Schur multiplier of $H \rtimes \Gamma$ is trivial. By \cite[Lemma 4.5]{Boston2021a}, one can then define $2$-special groups as special groups with relation rank at most $4$. Finally one defines $2$-select groups~\cite[Definition 5.1]{Boston2021a} as all groups arising as quotients of a $2$-special group by the normal $\Gamma$-invariant subgroup generated by an element of order at most $2$.

In \cite[Theorem 5.3]{Boston2021a}, Boston and Bush prove that for a 2-select group $H$, the group $H\rtimes \Gamma$ has Schur multiplier $1$ or $2$.  Combining this with \cite[Theorem 5.2]{Boston2021a}, one obtains in particular that for $H$ a wide $2$-class tower group of a cyclic cubic field, $H \rtimes \Gamma$ has Schur multiplier $1$ or $2$. Combining Lemma \ref{schur-is-v1} with Lemma \ref{ne-d-finite} (with $V=V_1$), we obtain the same conclusion by different means. This gives the overlap between our work and Boston and Bush's.

Boston and Bush further conjecture~\cite[Conjecture 5.4]{Boston2021a} that every $2$-select group in fact appears as a $2$-class field tower.

This raises two questions, which we are not currently able to answer.  Do all 2-select groups have positive probability under our conjectured distribution? Are all groups which have positive probability 2-select? The above characterization of positive probability makes these questions purely group-theoretic. Positive answers to these questions would demonstrate a compatibility between \cite[Conjecture 5.4]{Boston2021a}  and our conjectures

%
%
%

Finally, we calculate our predicted probabilities for the first two $2$-select $2$-groups with the numerical evidence of Boston-Bush. Note that they compute frequencies only within the set of cyclic cubic fields whose 2-class field tower group has two generators, so we first calculate the probability that the $2$-class field tower group has two generators. To do this, we can take $\mathcal L = \mathcal L^{\textrm{elem}} $ to consist of elementary abelian $2$-groups with an action of $\Gamma$ and compute the probability in our distribution  that the maximum quotient in $\mathcal L^{\textrm{elem}}$ of the Galois group is $V_2$.    Since $(V_1^2)_\Gamma\ne 1$, we have $\nu( \{ X \mid X^{\mathcal L} \cong (V_1)^2
\} ) =0$.
Since $H^2( V_2 \rtimes \Gamma, \mathbb F_2)^{\mathcal L} = H^2(V_2\rtimes \Gamma, V_2)^{\mathcal L}=0$, we have $w_{V_1}(s_H) = 1$ and $w_{V_2}(s_H) = 
\prod_{k=0}^{\infty} (1+ 4^{ -k - \frac{3}{2} })^{-1}$ independently of $s_H$. Thus the sum over $s_H$  in \eqref{cct-prob} cancels with the $H_3$ term and we obtain
\[ \nu( \{ X \mid X^{\mathcal L} \cong V_2\} ) = \frac{ |H^2(V_2 \rtimes \Gamma, \mathbb F_2) | |V_2^\Gamma| }{ |\Aut_\Gamma(V_2)| |V_2|  } \prod_{k=0}^{\infty} (1+ 4^{ -k - \frac{3}{2} })^{-1}  = \frac{ 2 \cdot 1}{ 3 \cdot 4}  \prod_{k=0}^{\infty} (1+ 4^{ -k - \frac{3}{2} })^{-1} = \frac{1}{6}  \prod_{k=0}^{\infty} (1+ 4^{ -k - \frac{3}{2} })^{-1}  .\]

This agrees with the probability that a cyclic cubic field has class group of $2$-rank two conjectured by Malle~\cite[Equation (1) on p. 2827]{Malle2008}, and with the probability conjectured in our prior work~\cite[Conjecture 1.1, Theorem 1.2]{Sawin2023}. (This was guaranteed since the conjecture in this paper agrees with \cite{Sawin2023} for all abelian groups, and \cite{Sawin2023} was already checked to agree with Malle in many cases, including this one.)

We must divide the probability of obtaining a given $2$-generator $2$-group by $ \frac{1}{6}  \prod_{k=0}^{\infty} (1+ 4^{ -k - \frac{3}{2} })^{-1}$ to obtain its expected frequency among $2$-generator $2$-groups.
 
 We now calculate the conjectural probability that the $2$-class field tower is the Klein four group, before doing the same for the eight-element quaternion group. In each case, we begin with a lemma about the cohomology of that group.
 
 \begin{lemma}\label{klein-four-cohomology} Let $V_2$ be the Klein four group, viewed as a nontrivial representation of $\Gamma = \mathbb Z/3$. 
 
 \begin{itemize} 
 
 \item For $i\leq 3$ we have \[ H^i ( V_2 \rtimes \Gamma, \mathbb F_2) = \begin{cases} \mathbb F_2 & \textrm{if } i=0 \\ 0 & \textrm{if } i=1 \\ \mathbb F_2 & \textrm{if } i=2 \\ \mathbb F_2^2 & \textrm{if } i=3 \end{cases}. \] 
 
 \item For $i\leq 2$ we have \[ H^i (V_2 \rtimes \Gamma, V_2) = \begin{cases} 0 & \textrm{if } i=0 \\ \mathbb F_2^2 &  \textrm{if } i=1 \\ \mathbb F_2^2  & \textrm{if } i=2 \end{cases} \] 
 
 \item The Bockstein map $H^2 ( V_2 \rtimes \Gamma , \mathbb F_2) \to H^3( V_2 \rtimes \Gamma, \mathbb F_2)$ is nontrivial.
 
 \item There exist isomorphisms of  of abelian groups  from $H^3( V_2 \rtimes \Gamma, \mathbb F_2), H^2( V_2 \rtimes \Gamma, V_2),$ and $ H^1 ( V_2 \rtimes \Gamma, V_2^\vee) $ to $\mathbb F_4$ 
 such that the  cup product map $H^2( V_2 \rtimes \Gamma, V_2) \times H^1 ( V_2 \rtimes \Gamma, V_2^\vee) \to H^3( V_2 \rtimes \Gamma, \mathbb F_2) $ is given by a nontrivial sesquilinear form $\mathbb F_4 \times \mathbb F_4 \to \mathbb F_4 $. 
 
 \end{itemize}
 
 \end{lemma}
 
 \begin{proof} These follow from the following more basic facts. We can express $H^* ( V_2, \mathbb F_2)$ as the free ring over $\mathbb F_2$ on generators $x,y $ in degree $1$. (This follows from the K\"unneth formula and the analogous description of the cohomology of $\mathbb F_2$.)  A fixed generator $\gamma$ of $\Gamma$ acts by sending $x$ to $y$ and $y $ to $x+y$. 
So $H^1(V_2,\F_2)\isom V_2$ as $\Gamma$-representations. 
 The Bockstein map satisfies  the Leibnitz rule and the identities $B(x)= x^2, B(y)= y^2$. Finally, for any representation $V$ of $\Gamma$ we have $H^i ( V_2 \rtimes \Gamma, V) = H^i (V_2, V)^{\Gamma} = ( H^i (V_2, \mathbb F_2) \otimes V)^\Gamma\isom(\Sym^i V_2\tensor V)^\Gamma$.  From this we can easily compute the dimensions in the first two claims, though we will find elements more explicitly in order to prove the last two claims.
 
 In particular, $H^i ( V_2 \rtimes \Gamma, \mathbb F_2)$ is the space of $\Gamma$-invariant polynomials of degree $i$ over $\mathbb F_2$. 
 The Bockstein map applied to the generator $x^2 + xy+ y^2$  of the invariants in degree $2$ gives $2x^3 + x^2y+ xy^2 + 2y^3 = x y (x+y)$ which is a nontrivial invariant in degree $3$, proving the third claim.
 
 
We fix a generator $\omega$ of $\mathbb F_4$ satisfying $\omega^2=\omega+1$.  Let $\Gamma$ act on a one-dimensional $\mathbb F_4$-vector space $V=\F_4$, such that the generator $\gamma$ acts by multiplication by $\omega$.  As $\Gamma$-representations over $\F_2$, we have that $V$ is isomorphic to $V_2$.
The $\Gamma$-invariants in $H^*( V_2,  V)$ are a $\F_4$-vector space.  
We can check that $ \omega x+  y \in H^1( V_2 , V)$ is $\Gamma$-invariant and thus generates the $\Gamma$-invariants as an $\F_4$-vector space, and that $\omega x^2  + y^2$ generates the $\Gamma$-invariants in $H^2( V_2,  V)$ as an $\F_4$-vector space.

The cup product map  $ H^2( V_2 \rtimes \Gamma, V_2)\times H^1 ( V_2 \rtimes \Gamma, V_2^\vee)\to H^3( V_2 \rtimes \Gamma, \mathbb F_2) $, expressed in terms of  $\Gamma$-invariant polynomials over $\mathbb F_4$, may be obtained by multiplying one polynomial with the conjugate of the other and taking the trace from $\mathbb F_4$ to $\mathbb F_2$, since the $\Gamma$-invariant $\mathbb F_2$-bilinear form $\mathbb F_4 \times \mathbb F_4 \to \mathbb F_2$ is obtained by multiplying one element by the conjugate of the other and taking the trace.

Applied to $\alpha ( \omega x +y) $ and $ \beta ( \omega x^2 +y^2)$ we obtain \[ \alpha (\omega x+ y) \times \overline{\beta} (\overline{\omega} x^2 + y^2) = \alpha \overline{\beta} ( x^3 + \overline{\omega} x^2 y + \omega x y^2 + y^3 ) .\] If we choose an isomorphism between $\mathbb F_4$ and the ring of degree $3$ invariant polynomials over $\F_2$ that sends $1$ to \[ \operatorname{tr} ( x^3 + \overline{\omega} x^2 y + \omega x y^2 + y^3 ) =x^2y+xy^2= xy (x+y) \] and $\omega$ to  \[ \operatorname{tr} ( \omega ( x^3 + \overline{\omega} x^2 y + \omega x y^2 + y^3)  ) = x^3 + xy^2 +y^3 \] then this map sends $\alpha ( \omega x +y) $ and $ \beta ( \omega x^2 +y^2)$ to $\alpha \overline{\beta}$, verifying the final claim. \end{proof}

\begin{lemma}\label{klein-four-probability} For $X$ distributed according to the measure $\nu$, the probability that the maximal $2$-group quotient of $X$ is the Klein four group is $\frac{1}{12}  \prod_{k=0}^{\infty} (1+ 4^{ -k - \frac{3}{2} })^{-1}$ and the probability conditional on the abelianization of $X$ having $2$-rank $2$ is $\frac{1}{2}$. 

\end{lemma}

The probability $\frac{1}{2}$ agrees with \cite[Table 7, row 1]{Boston2021a}.

\begin{proof}  Let $\cL$ be a level containing every extension of $V_2$ by $V_1$ or $V_2$. By Lemma \ref{maximal-p-group-quotient}, the probability that the maximal $2$-group quotient of $X$ is the Klein four group is $ \nu( \{ X \mid X^{\mathcal L} \cong H\})$. We calculate this using the formula of \eqref{cct-prob}.

Since $H^3(H\rtimes \Gamma, \mathbb F_2) \cong \mathbb F_2^2$ by \ref{klein-four-cohomology}(1) there are four possible choices for $s_H$, corresponding to the four linear forms on $\mathbb F_2^2$. 

By Lemma \ref{maximal-p-group-quotient} we have $H^2( V_2 \rtimes \Gamma, V_1)^{\mathcal L}= H^2(V_2 \rtimes \Gamma, \mathbb F_2)$ so  $H^2( V_2 \rtimes \Gamma, V_1)^{\mathcal L, s_H} $ is the kernel of the Bockstein homomorphism $H^2(V_2 \rtimes \Gamma, \mathbb F_2)\to H^3(V_2 \rtimes \Gamma, \mathbb F_2)\to_{s_H} \mathbb F_2$. By Lemma \ref{klein-four-cohomology}(1,3), the source of this map is $\mathbb F_2$ so the kernel vanishes if and only if the composition is nontrivial, i.e. if and only if the linear form associated to $s_H$ vanishes on the image under the Bockstein homomorphism in $H^3(V_2 \rtimes \Gamma, \mathbb F_2)$ of the unique nontrivial element in $ H^2(V_2 \rtimes \Gamma, \mathbb F_2)$.  Thus the kernel vanishes for two possible values of $s_H$ and does not vanish for two other values of $s_H$, including $s_H=0$. By \eqref{wv1} we have $w_{V_1}(s_H)=1$ for the two values where the kernel vanishes and $w_{V_1}(s_H)=0$ for the two values where the kernel does not vanish.

By Lemma \ref{maximal-p-group-quotient} we have $H^2( V_2 \rtimes \Gamma, V_2)^{\mathcal L}= H^2(V_2 \rtimes \Gamma, V_2)$ so  $H^2( V_2 \rtimes \Gamma, V_2)^{\mathcal L, s_H} $ is the kernel of the cup product bilinear form $H^2(V_2 \rtimes \Gamma, V_2) \times H^1(V_2 \rtimes \Gamma,  V_2^\vee) \to H^3(V_2 \rtimes \Gamma, \mathbb F_2)\to_{s_H} \mathbb F_2$. As long as $w_{V_1}(s_H) \neq 0$, we have $s_H\neq 0$, so by Lemma \ref{klein-four-cohomology}(4)  no element of $H^2(V_2 \rtimes \Gamma, V_2)$ maps $H^1(V_2 \rtimes \Gamma,  V_2^\vee)$ into a $1$-dimensional subspace of 
$H^3(V_2 \rtimes \Gamma, \mathbb F_2)$.  Thus $s_H\ne 0$ implies the kernel is trivial. By \eqref{wv2} we have $w_{V_2}(s_H) = \prod_{k=0}^{\infty} (1+ 4^{ -k - \frac{3}{2} })^{-1}$ as long as $w_{V_1}(s_H)\neq 0$. 

Combining these observations, we see that

\[ \sum_{ s_H \in H_3 (V_2 \rtimes \Gamma,  \mathbb F_2)} w_{V_1}(s_H) w_{V_2} (s_H) =2 \cdot 1 \cdot \prod_{k=0}^{\infty} (1+ 4^{ -k - \frac{3}{2} })^{-1} .\]

Hence by \eqref{cct-prob} we have 

\[  \nu( \{ X \mid X^{\mathcal L} \cong V_2 ) = \frac{ |H^2(V_2  \rtimes \Gamma, \mathbb F_2) | |H^\Gamma| }{ |H^3(V_2  \rtimes \Gamma,  \mathbb F_2) | |\Aut_\Gamma(V_2)| |V_2|  }2 \cdot 1 \cdot \prod_{k=0}^{\infty} (1+ 4^{ -k - \frac{3}{2} })^{-1} \]
\[ = \frac{2 \cdot 1}{ 4 \cdot 3 \cdot 4} 2 \cdot 1 \cdot \prod_{k=0}^{\infty} (1+ 4^{ -k - \frac{3}{2} })^{-1} =\frac{1}{12}  \prod_{k=0}^{\infty} (1+ 4^{ -k - \frac{3}{2} })^{-1} .\]

Dividing by $\frac{1}{6}  \prod_{k=0}^{\infty} (1+ 4^{ -k - \frac{3}{2} })^{-1}  $ to compute the conditional probability, we obtain $\frac{1}{2}$. \end{proof}

  \begin{lemma}\label{quaternion-cohomology} Let $Q_8$ be the eight-element quaternion group, with the unique-up-to-isomorphism nontrivial action of $\Gamma = \mathbb Z/3$. 
 
 \begin{itemize} 
 
 \item For $i\leq 3$ we have \[ H^i ( Q_8 \rtimes \Gamma, \mathbb F_2) = \begin{cases} \mathbb F_2 & \textrm{if } i=0 \\ 0 & \textrm{if } i=1 \\ 0 & \textrm{if } i=2 \\ \mathbb F_2 & \textrm{if } i=3 \end{cases}. \] 
 
 \item For $i\leq 3$ we have \[ H^i (Q_8 \rtimes \Gamma, V_2) = \begin{cases} 0 & \textrm{if } i=0 \\   \mathbb F_2^2 & \textrm{if } i=1 \\ \mathbb F_2^2  & \textrm{if } i=2 \\ 0 & \textrm{if }i=3 \end{cases} \] 
  
 \item The  cup product map $H^2( Q_8 \rtimes \Gamma, V_2) \times  
H^1 ( Q_8 \rtimes \Gamma, V_2^\vee) 
 \to H^3( Q_8 \rtimes \Gamma,\mathbb F_2) $ is a non-degenerate $\mathbb F_2$-bilinear form. 
 \end{itemize}
 
 \end{lemma}
 
 \begin{proof} These follow from the fact that $ Q_8 \rtimes \Gamma$ is a finite subgroup of $SU(2)$, being the inverse image of the group of rotational symmetries of the tetrahedron, and hence acts freely on $S^3$.  Thus for $V$ a representation of $Q_8 \rtimes \Gamma$  we have a natural map $H^i ( Q_8 \rtimes \Gamma, V) \to H^i ( S^3/ Q_8 \rtimes \Gamma, V)$ where $S^3/ Q_8 \rtimes \Gamma$ is a manifold.
 
 Let us check that the natural map $H^i ( Q_8 \rtimes \Gamma, V) \to H^i ( S^3/ Q_8 \rtimes \Gamma, V)$ is an isomorphism for $i\leq 3$ and any irreducible representation $V$ of characteristic $2$.  Since $S^3$ has cohomology only in degrees $0$ and $3$, the spectral sequence computing  $H^{p+q} ( S^3/ Q_8 \rtimes \Gamma, V)$ from $H^p ( Q_8 \rtimes \Gamma, H^q(S^3, V))$  degenerates to a long exact sequence
 \[ H^{i-4} (Q_8 \rtimes \Gamma, V)   \to H^i ( Q_8 \rtimes \Gamma, V) \to H^i ( S^3/ Q_8 \rtimes \Gamma, V)  \to H^{i-3} (Q_8 \rtimes \Gamma, V)  .\]
 For $i\leq 3$ we have $ H^{i-4} (Q_8 \rtimes \Gamma, V)  =0$ and we have $H^{i-3} (Q_8 \rtimes \Gamma, V)  =0$ unless $i=3$, in which case  $H^{i-3} (S^3/Q_8 \rtimes \Gamma, V)  = V^{Q_8 \rtimes \Gamma}$ vanishes unless $V$ is the trivial repesentation. In the case that $V$ is the trivial representation $\mathbb F_2$, the map $ H^3 ( S^3/ Q_8 \rtimes \Gamma, \mathbb F_2)  \to H^{0} (Q_8 \rtimes \Gamma, \mathbb F_2)\to \mathbb F_2$ represents integration against the fundamental class of $S^3$, which is trivial as the fundamental class of $S^3$ pushes forward to $|Q_8 \rtimes \Gamma|=24 \equiv 0 \bmod 2$ times the fundamental class of $S^3/ Q_8 \rtimes \Gamma$ and thus has trivial pairing with the mod $2$ cohomology of $S^3/ Q_8 \rtimes \Gamma$.  Because the map $H^i ( S^3/ Q_8 \rtimes \Gamma, V)  \to H^{i-3} (Q_8 \rtimes \Gamma, V)  $ is zero regardless of whether $V \cong \mathbb F_2$, the natural map $H^i ( Q_8 \rtimes \Gamma, V) \to H^i ( S^3/ Q_8 \rtimes \Gamma, V)$ is an isomorphism, as desired.

Having checked this, for any representation $V$, by Poincar\'e duality we have $H^2( Q_8 \rtimes \Gamma, V) \cong H^1 ( Q_8 \rtimes \Gamma, V^\vee)^\vee$ with the isomorphism arising from the cup product bilinear form $H^2( Q_8 \rtimes \Gamma, V) \times  H^1 ( Q_8 \rtimes \Gamma, V^\vee)\to H^3( Q_8 \rtimes \Gamma, \mathbb F_2)$. This immediately proves the third claim and reduces the calculation of $H^2$ to the calculation of $H^1$, which is immediate from the fact that the abelianization of $Q_8$ is $\mathbb F_2^2$, with $\Gamma$ acting nontrivially, and reduces the calculation of $H^3$ to the even easier calculation of $H^0$.
 \end{proof}

 \begin{lemma}\label{quaternion-probability} For $X$ distributed according to the measure $\nu$, the probability that the maximal $2$-group quotient of $X$ is the quaternion group $Q_8$ is $\frac{7}{2^5 \cdot 3 }  \prod_{k=0}^{\infty} (1+ 4^{ -k - \frac{3}{2} })^{-1} $ and the probability conditional on the abelianization of $X$ having $2$-rank $2$ is $\frac{7}{16}$. 

\end{lemma}

The probability $\frac{7}{16}$ agrees with \cite[Table 7, row 2]{Boston2021a}.

\begin{proof}  Let $\cL$ be a level containing every extension of $Q_8$ by $V_1$ or $V_2$. By Lemma \ref{maximal-p-group-quotient}, the probability that the maximal $2$-group quotient of $X$ is $Q_8$ is $ \nu( \{ X \mid X^{\mathcal L} \cong H\})$. We calculate this using the formula of \eqref{cct-prob}.

Since $H^3(Q_8\rtimes \Gamma, \mathbb F_2) \cong \mathbb F_2$ by \ref{quaternion-cohomology}(1) there are two possible choices for $s_H$, one corresponding to a trivial linear form and one corresponding to a nontrivial form.

By Lemma \ref{quaternion-cohomology}(1) we have $H^2( Q_8 \rtimes \Gamma, \mathbb F_2) =0$ so certainly $H^2( Q_8 \rtimes \Gamma, \mathbb F_2)^{\mathcal L,s_H}=0$. Hence by \eqref{wv1} we have $w_{V_1}(s_H)=1$ for either choice of $s_H$. 

By Lemma \ref{maximal-p-group-quotient} we have $H^2( Q_8  \rtimes \Gamma, V_2)^{\mathcal L}= H^2(Q_8  \rtimes \Gamma, V_2)$ so  $H^2( Q_8 \rtimes \Gamma, V_2)^{\mathcal L, s_H} $ is the kernel of the cup product bilinear form $H^2(Q_8 \rtimes \Gamma, V_2) \times H^1(Q_8 \rtimes \Gamma, V_2^\vee) \to H^3(Q_8 \rtimes \Gamma, \mathbb F_2)\to_{s_H} \mathbb F_2$. By Lemma \ref{quaternion-cohomology}(3) this bilinear form is nondegenerate if $s_H$ is nonzero and zero if $s_H$ is zero. If $s_H$ is nonzero, this makes the kernel trivial, so by \eqref{wv2} we have $w_{V_2}(s_H) = \prod_{k=0}^{\infty} (1+ 4^{ -k - \frac{3}{2} })^{-1} $ . If $s_H$ is zero, this makes the kernel $2$-dimensional over $\mathbb F_2$, and hence $1$-dimensional over $\mathbb F_4$, so by \eqref{wv2} we have $w_{V_2}(s_H) =\frac{3}{4} \prod_{k=0}^{\infty} (1+ 4^{ -k - \frac{3}{2} })^{-1} $.

Combining these observations, we see that

\[ \sum_{ s_H \in H_3 (Q_8 \rtimes \Gamma,  \mathbb F_2)} w_{V_1}(s_H) w_{V_2} (s_H) =\prod_{k=0}^{\infty} (1+ 4^{ -k - \frac{3}{2} })^{-1}+ \frac{3}{4} \prod_{k=0}^{\infty} (1+ 4^{ -k - \frac{3}{2} })= \frac{7}{4} \prod_{k=0}^{\infty} (1+ 4^{ -k - \frac{3}{2} })^{-1}.\]

Hence by \eqref{cct-prob} we have 
\[  \nu( \{ X \mid X^{\mathcal L} \cong Q_8\} ) = \frac{ |H^2(Q_8 \rtimes \Gamma, \mathbb F_2) | |H^\Gamma| }{ |H^3(Q_8   \rtimes \Gamma,  \mathbb F_2) | |\Aut_\Gamma(Q_8)| |Q_8|  } \frac{7}{4} \prod_{k=0}^{\infty} (1+ 4^{ -k - \frac{3}{2} })^{-1}  \]
\[ = \frac{ 1 \cdot 2 }{ 2 \cdot 3 \cdot 8}  \frac{7}{4}  \cdot \prod_{k=0}^{\infty} (1+ 4^{ -k - \frac{3}{2} })^{-1}  =\frac{7}{2^5 \cdot 3 }  \prod_{k=0}^{\infty} (1+ 4^{ -k - \frac{3}{2} })^{-1}  .\]

Dividing by $\frac{1}{6}  \prod_{k=0}^{\infty} (1+ 4^{ -k - \frac{3}{2} })^{-1}  $ to compute the conditional probability, we obtain $\frac{7}{16}$. \end{proof}

\section{Evidence in the function field case}\label{S:evidence}

The goal of this section is to prove Theorem~\ref{ff-theorem}, a function field analogue of Conjecture~\ref{intro-moments-conjecture}. 
We fix an $n$-oriented $\Gamma$-group $\mathbf H$ throughout this section.
Our function field analogue estimates a certain sum over $\Gamma$-extensions $K$ of $\mathbb F_q(t)$, split at $\infty$, of the number of surjections from the Galois group of the maximal unramified, split at $\infty$, extension of $K$ to $\mathbf H$. Similar sums have been considered in prior work \cite{Liu2024,Liu2022,LandesmanLevy,LiuWillyard}. Such a sum was estimated in the large $q$ limit by Liu, Wood, and Zurieck-Brown~\cite[Theorem 1.4]{Liu2024}, assuming that $\abs{H}$ is prime to $q-1$, in which case the orientation is always trivial. Liu~\cite{Liu2022} removed this assumption, introduced the $\omega$-invariant, which we have already seen includes some, but not all, of the information in the orientation, and proved a large $q$ estimate keeping track of the $\omega$-invariant. A similar estimate in the more difficult case of fixed large $q$ was proven by Landesman and Levy~\cite[Theorem 1.3.2]{LandesmanLevy}. Finally, Liu and Willyard~\cite[Theorem 1.1]{LiuWillyard} relaxed the assumption that the $\Gamma$-extension be split at infinity. 

These arguments all rely on certain topological results. For $q\to\infty$ limits, it suffices to have results describing the set of connected components of Hurwitz space. To handle the fixed large $q$ case, Landesman and Levy~\cite{LandesmanLevy} proved breakthrough results controlling the low-degree homology groups of Hurwitz space, generalizing and strengthening earlier work of Ellenberg, Venkatesh, and Westerland~\cite{Ellenberg2016}.

Our strategy will therefore be similar to the strategy of these prior works, except accounting for the orientation. This requires us to take a certain covering of Hurwitz space, and prove results on the set of components of this covering and the low-degree homology groups of this covering. The description of the set of components of this covering is essentially a large monodromy result, and follows from classical results in low-dimensional topology. The computation of the low-degree homology groups of the cover of Hurwitz space proceeds by observing that this cover is itself covered by another Hurwitz space, whose cohomology groups can be controlled by applying the results of \cite{LandesmanLevy}.

We begin with the statement of our main counting theorem.
We follow with a brief general notation Subsection \ref{SS:AGnotation}, and then Subsection~\ref{SS:DefHur} defining the Hurwitz schemes we will use in our proof.  In Subsection~\ref{SS:PropHur} we establish the basic properties of Hurwitz schemes and prove Theorem~\ref{ff-theorem}  conditional on a result comparing the number of $\mathbb F_q$-points on a component of Hurwitz space and a cover of that component. The remaining sections prove that result, with Subsection~\ref{SS:toparg} proving a topological result about the braid fundamental class, Subsection~\ref{SS:homology} proving a topological result about the homology of covers of Hurwitz space, and then Subsection~\ref{SS:Hurcount} applying these results to Hurwitz spaces and their covers over finite fields.  

For $K$ a $\Gamma$-extension of $\mathbb F_q(t)$, we let $\rDisc K$ be the product of the norms of the places of $\mathbb F_q(t)$ that ramify in $K$ (the norm of a place $v$ of $\mathbb F_q(t)$ being $q^{\deg v}$). Let $K^{ \operatorname{un}, \abs{\Gamma}'} $ be the composition of all finite Galois extensions of $K$ with degree prime to $|\Gamma|$ that are everywhere unramified and split over $\infty$. 
If $\mathbb F_q(t)$ contains the $n$th roots of unity, i.e. $n \mid q-1$, then since $K^{ \operatorname{un}, \abs{\Gamma}'}/K$ is an extension with Galois group $\Gal( K^{ \operatorname{un}, \abs{\Gamma}'} /K) $,  it has an associated Artin-Verdier fundamental class $s\in H_3(\Gal( K^{ \operatorname{un}, \abs{\Gamma}'} /K)  ,\Z/n)^\Gamma$.
This makes $\bGal( K^{ \operatorname{un}, \abs{\Gamma}'} /K)$ naturally an $n$-oriented $\Gamma$-group.
(The $\Gamma$ action is by using a homomorphic section of $\Gal( K^{ \operatorname{un}, \abs{\Gamma}'} /\F_q(t))\ra\Gal(K/\F_q(t))=\Gamma$
and conjugating. By the Schur--Zassenhaus theorem \cite[Theorem 2.3.15]{Ribes2010}, this section is well-defined up to conjugation by elements of $\Gal( K^{ \operatorname{un}, \abs{\Gamma}'} /K)$, and thus the resulting  $n$-oriented $\Gamma$-group is well-defined up to isomorphism.)

 Let $E_\Gamma( q^m, \mathbb F_q(t) )$ be the set of $\Gamma$-extensions
  $K$ of $F$, split at $\infty$, such that  $\rDisc K = q^m$.
 
 The goal of this section will be to prove the following theorem.

\begin{theorem}\label{ff-theorem} 
Let $\Gamma$ be the finite group and $n$ the positive integer coprime to $\abs{\Gamma}$ fixed throughout the paper, and 
$\mathbf H$ the finite $n$-oriented $\Gamma$-group fixed throughout this section.  Suppose that $H_\Gamma=1$. Let $q$ be a prime power. As long as $q$ satisfies the congruence conditions that $ (q, |\Gamma| |H|)=1$ and $ q \equiv 1 \bmod n$ and $(q-1, |H|)= (n,|H|)$ and $q$ is sufficiently large depending on $\Gamma, H$ we have

\[ \lim_{b \to\infty}  \frac{\sum_{m\leq b} \sum_{K \in E_\Gamma( q^m , \mathbb F_q(t))}  \Sur ( \bGal ( K^{ \operatorname{un}, \abs{\Gamma}'} / K), \mathbf H)  }{ \sum_{m\leq b}\abs{ E_\Gamma( q^m , \mathbb F_q(t))}} =  \frac{ \abs{H^\Gamma} \abs{ H^2( H\rtimes \Gamma , \mathbb Z/n) } }{  \abs{H} \abs{ H^3( H\rtimes \Gamma , \mathbb Z/n) }  }.\]\end{theorem}

We begin with the following analogue of \cite[Lemma 9.3]{Liu2024} used to translate surjection counting to extension counting.

\begin{lemma}\label{count-to-count} 
Let $q \equiv 1\bmod n$ be a prime power. 
 Let $N(\mathbf H, \Gamma, q^m, \mathbb F_q(t))$ be the number of surjections $\Gal ( \overline{\mathbb F_q(t)}/\mathbb F_q(t)) \to H \rtimes \Gamma$ such that the corresponding $(H \rtimes \Gamma)$-extension $K$ has $\rDisc K = q^m$, the associated $H$-extension is unramified everywhere, split completely at all places lying above $\infty$,  the induced map $\Gal ( K^{ \operatorname{un}, \abs{\Gamma}'} / K) \ra H$ gives a homomorphism of $n$-oriented groups $\bGal ( K^{ \operatorname{un}, \abs{\Gamma}'} / K) \ra \bH$, and the associated $\Gamma$-extension $K^H/\mathbb F_q(t)$ is split completely above $\infty$. Then

\[ \sum_{K \in E_\Gamma( q^m , \mathbb F_q(t))}  \Sur ( \bGal ( K^{ \operatorname{un}, \abs{\Gamma}'} / K), \mathbf H)= \frac{\abs{H^\Gamma}}{ \abs{H}}N(\mathbf H, \Gamma, q^m, \mathbb F_q(t)).\]  \end{lemma}

\begin{proof} 
From
\cite[Lemma 9.3]{Liu2024}, we have
\[ \sum_{K \in E_\Gamma( q^m , \mathbb F_q(t))}  \Sur ( \Gal ( K^{ \operatorname{un}, \abs{\Gamma}'} / K),  H)= \frac{\abs{H^\Gamma}}{ \abs{H}}N( H, \Gamma, q^m, \mathbb F_q(t))\] where the surjections on the left-hand side need not respect the orientation and $N( H, \Gamma, q^m, \mathbb F_q(t))$ is defined the same way as $N(\mathbf H, \Gamma, q^m, \mathbb F_q(t))$ but without the orientation condition.  The proof of \cite[Lemma 9.3]{Liu2024} gives a $\frac{\abs{H}}{\abs{H^\Gamma}}$-to-one surjective map from the set counted by $N( H, \Gamma, q^m, \mathbb F_q(t))$ to the union over $K \in E_\Gamma( q^m , \mathbb F_q(t))$ of the set counted by $\Sur ( \Gal ( K^{ \operatorname{un}, \abs{\Gamma}'} / K),  H)$. To check our statement, it suffices to check that the constructed surjection sends orientation-compatible maps to surjections whose associated $H$-extension has compatible orientation. This is clear since the map sends a surjection $\Gal ( \overline{\mathbb F_q(t)}/\mathbb F_q(t)) \to H \rtimes \Gamma$ associated to an extension $K$ to a pair consisting of the $\Gamma$-extension $K^H$ and the surjection associated to the $\Gamma$-extension $K/K^H$, so the orientations considered on each side are the same.\end{proof}

\subsection{Algebraic geometry notation}\label{SS:AGnotation}
We now set some algebraic geometry notation to be used in Section~\ref{S:evidence}.  For a scheme $X$ over $\Spec \Z$ and a ring $R$, we write $X_R$ for the pullback to $X$ from $\Spec \Z$ to $\Spec R$.  
We write $\Z[\mu_n]$ for the ring $\Z[x]/\Phi_n(x)$, where $\Phi_n$ is the $n$th cyclotomic polynomial.  
All cohomology of schemes is \'etale cohomology, and for a scheme $X$ we write $\pi_1(X)$ for the algebraic fundamental group (with a choice of geometric basepoint suppressed).  For a scheme $X$ over $\Spec \bbC$, when we write $X(\bbC)$, we mean the complex points of $X$ as a topological space with their analytic topology.  For a scheme $X$ over $\F_q$, we write $\Frob_q$ for the automorphism of $X_{\bar{\F}_q}$ given by the ring map $a\mapsto a^q$ on affine pieces. 
 For $X$ a scheme defined over $\F_q$, and $x\in X(\F_q)$, we write
 $\Frob_{q,x}$ for the image of Frobenius in the map $\pi_1(\Spec \F_q)\ra \pi_1(X)$ associated to $x$, which is only well-defined up to conjugation since we do not choose base-points. 
 
For a scheme $X$ and a finite group $G$,  the definition of the algebraic fundamental group associates to
each homomorphism $\pi_1(X)\ra G$  a degree $|G|$ \'etale cover $Y\ra X$ with a given action of $G$ on $Y$ equivariant with the trivial action on $X$.  In other words, $Y/X$ is a $G$-torsor.

\subsection{Definition of Hurwitz schemes}\label{SS:DefHur}

We now define the Hurwitz schemes that are relevant to the proof, in particular because they will be moduli spaces for the kinds of extensions we want to count. 
For $G$ a finite group and $m$ a nonnegative integer, let $H_{G,*}^m$ be the scheme over $\mathbb Z [\abs{G}^{-1}]$ whose $S$-points are given by isomorphism classes of tame Galois covers of $\mathbb P^1$ with branch locus of degree $m$, unramified above $\infty$, with a marked point on the fiber over $\infty$,  and with an identification of their automorphism group with $G$,  defined as ``$\Hur_{G,*}^m$'' after Remark 11.3 in \cite{Liu2024}. 
We call a Galois cover of $\P^1$ with an identification of its automorphism group with $G$ a \emph{$G$-covering}.
  (See also  \cite[Section 11.1]{Liu2024} for definitions of these terms.)
The fact that $H_{G,*}^m$ exists (or, if the language is preferred, that the stack $H_{G,*}^m$ is a scheme) is \cite[Lemma 11.4]{Liu2024}. For $c$ a union of conjugacy classes of $G$, stable under the operation of raising a conjugacy class to the $d$th power for $d$ prime to $|G|$, let $H_{G,c}^m$ be the union of connected components of $H_{G,*}^m$ parameterizing covers all whose inertia groups are generated by elements of $c$, as in \cite[Section 11.4]{Liu2024}.   We have fixed $n$ throughout the paper.  We will work over $\Z[\mu_n]$, and we define $\Hur_{G,c}^m$ to be $(H_{G,c}^m)_{\Z[\mu_n]}.$
Even though this differs from the use of the notation in \cite{Liu2024}, there will be little conflict, since we are mostly focused on the base change to a finite field, and for any scheme $X$ and prime power $q\equiv 1 \pmod {n}$, we have $(X)_{\F_q}= (X_{\Z[\mu_n]})_{\F_q}$.
For the rest of this section, we will only ever consider prime powers $q\equiv 1 \pmod n$.

Let $c_\Gamma=  \Gamma \setminus \{1\}$ and let $c_{H  \rtimes \Gamma}$ be the set of all nontrivial elements of $H \rtimes \Gamma$ conjugate to an element of $\Gamma$.  

We will mostly consider $\Hur_{G,c}^m$ in the case when $(G,c)$ is $(H 
\rtimes \Gamma,c_{H 
\rtimes \Gamma})$.

The denominator in Theorem \ref{ff-theorem} was already computed in \cite{Liu2024} in terms of $\Hur_{\Gamma, c_\Gamma}^m$. We review this computation now. 
For $q$ relatively prime to $\abs{\Gamma}$ and $q\equiv 1 \pmod {n}$ we have
\begin{equation}\label{eq-count-to-geom}\abs{ E_\Gamma( q^m , \mathbb F_q(t))} = \abs{ \Hur_{\Gamma, c_\Gamma}^m(\mathbb F_q)} \end{equation} by \cite[Lemma 10.2]{Liu2024}. This enables the computation of $\abs{ E_\Gamma( q^m , \mathbb F_q(t))} $ via the \'{e}tale cohomology of $\Hur_{\Gamma, c_\Gamma}^m$.



To derive Theorem \ref{ff-theorem} we must obtain a similar result about $N(\mathbf H, \Gamma, q^m, \mathbb F_q(t))$. Without the orientation on $H$, \cite[Lemma 10.2]{Liu2024} shows that the analogous count is the number of $\mathbb F_q$-points of the Hurwitz space $\Hur_{H\rtimes \Gamma, c_{H\rtimes \Gamma}}$. We will show that $N(\mathbf H, \Gamma, q^m, \mathbb F_q(t))$ is proportional to the number of $\mathbb F_q$-points on a covering of the Hurwitz space. We will define this covering component by component.

  Let $F\colon Y \to X$ be a family of curves (i.e. a smooth proper morphism of schemes of relative dimension $1$ with connected geometric fibers) defined over $\mathbb Z[ \mu_n, 1/n]$. Then we have an ``integration along the fibers'' map $\int_F \colon  H^{i+2}(Y, \mu_n) \to H^i(X,\mathbb Z/n)$, defined as follows. We have $R^j F_*  \mu_n =0$ for $j>2$ so there is a map of complexes  $R F_* \mu_n \to R^2 F_* \mu_n [-2] $ inducing a map $H^{i+2}( X, R F_* \mu_n) \to H^i(X, R^2 F_* \mu_n)$.  The isomorphism $H^{i+2}(Y, \mu_n) \to H^{i+2}( X, R F_* \mu_n)$ arises from the derived category version of the Leray spectral sequence (i.e. the derived functor of a composition being the composition of derived functors). The trace morphism \cite[XVIII, Theorem 2.9]{sga4-3} defines a morphism $R^2 F_* \mu_n \to \mathbb Z/n$. (Much more generally, for $F$ flat proper of relative dimension $d$ and $\mathcal F$ a sheaf on $X$ torsion of order invertible on $x$, \cite[XVIII, Theorem 2.9]{sga4-3} defines a map $R^{2d} F_* R^*  \mathcal F (d) \to \mathcal F$. The case $d=1$ and $\mathcal F =\mathbb Z/n$ is what we use). By functoriality of cohomology, this induces a map $H^i ( X, R^2 F_* \mu_n ) \to H^i(X, \mathbb Z/n )$. Composing these
  \[ H^{i+2}(Y, \mu_n) \to  H^{i+2}( X, R F_* \mu_n) \to H^i(X, R^2 F_* \mu_n)\to H^i(X, \mathbb Z/n )\]
 we obtain the integration along the fibers map. 
  
 Over $\Z[\mu_n]=Z[x]/\Phi_n(x)$, there is a natural map of sheaves $\Z/n\ra\mu_n$, which can pull back to any scheme over $\Z[\mu_n]$, including $Y$ above,
to obtain  $H^{i+2}(Y, \Z/n) \ra H^{i+2}(Y, \mu_n)$. We  compose this with the above to obtain an integration along the fibers map
$$
H^{i+2}(Y, \Z/n) \ra H^i(X, \mathbb Z/n ).
$$
For each $q\equiv 1 \pmod{n}$, we choose, for use throughout this paper, a primitive $n$th root of unity in $\F_q$.  When we base change a scheme over $ \Z[\mu_n]$ to $\F_q$, it is always with the map $ \Z[\mu_n]$ to $\F_q$ sending $x$ to that chosen root of unity.  When we consider a global function field $K$ over $\F_q$ and compute an Artin-Verdier trace, we always use that same choice of root of unity in $\F_q$.
  When applying this map, we will always identify $\mu_n$ with $\mathbb Z/n$ by the same choice of $n$th root of unity we use when defining the Artin-Verdier trace.

We now state some lemmas giving fundamental properties of the integration along the fibers map which we will use later. First, that pullback on cohomology is compatible with integration along the fibers.
  \begin{lemma}\label{integration-vs-pullback}  For $g \colon X_1 \to X_2$ a morphism with $Y_2 \to X_2$ a family of curves
defined over $\mathbb Z[ \mu_n, 1/n]$
   and $Y_1 = Y_2 \times_{X_2} X_1$ its pullback, we have the commutative diagram where the vertical arrows are integration along the fibers and the horizontal arrows are pullback
   \[ \begin{tikzcd} H^{i+2}(Y_2, \mu_n) \arrow[r]\arrow[d] & H^{i+2} ( Y_1,\mu_n) \arrow[d] \\ H^i (X_2, \mathbb Z/n) \arrow[r] & H^i (X_1, \mathbb Z/n). \end{tikzcd}\]
    \end{lemma}

   We also have a ``projection formula'' showing that integration along the fibers is compatible with the cup product.

   \begin{lemma}\label{projection-formula} Let $F\colon Y \to X$ be a family of curves defined over $\mathbb Z[ \mu_n, 1/n]$   and $\alpha$ a class in $H^j(X, \mathbb Z/n)$, then for $j\leq 1$ we have a commutative diagram
   
      \[ \begin{tikzcd} H^{i+2}(Y, \mu_n) \arrow[r,"\cup F^* \alpha" ]\arrow[d] & H^{i+j+2} ( Y,\mu_n) \arrow[d] \\ H^i (X, \mathbb Z/n) \arrow[r,"\cup \alpha"] & H^{i+j}(X, \mathbb Z/n) \end{tikzcd}\]
      
      \end{lemma}

   Both Lemma \ref{integration-vs-pullback} and Lemma \ref{projection-formula} will be proved later in \S\ref{ss-long-proofs}.

  Let $X$ be a connected component of  $\Hur^m_{H\rtimes \Gamma ,c_{H \rtimes \Gamma} } $. Let $Z$, over $X$, be the universal family of $(H \rtimes \Gamma)$-coverings of $\mathbb P^1$, so that $Z \to X \times \mathbb P^1$ is an $H \rtimes \Gamma$-covering. Let $Y$ be the universal family of induced $\Gamma$-coverings of $\mathbb P^1$, i.e. the quotient of $Z$ by $H$. 
  Then $Z \to Y $ is a finite \'{e}tale $H$-covering 
  and thus,  by Lemma~\ref{L:BHpullback}, determines a map from group cohomology to \'{e}tale cohomology $H^3(H, \mathbb Z/n) \to H^3(Y, \mathbb Z/n)$. This, together with the integration map and the explicit evaluation of \'{e}tale $H^1$, defines a composed homomorphism
 \[ H^3(H, \mathbb Z/n) \stackrel{\phi_{Z/Y}}{\to} H^3(Y, \mathbb Z/n) \to H^1( X, \mathbb Z/n) \cong \Hom ( \pi_1(X) ,\mathbb Z/n) \]
which we denote by $\tilde{e}$.   Since $H_3(H,\mathbb Z/n)$ is dual to $H^3(H,\mathbb Z/n)$ (see Lemma~\ref{L:homdual}),  
we let $e$ be the induced map $\pi_1(X)\to \Hom ( H^3(H,\mathbb Z/n), \mathbb Z/n) \cong H_3(H,\mathbb Z/n)$.

Since $Z$ is an $H\rtimes \Gamma$-covering of $X\times \P^1$, the morphism $H^3(H, \mathbb Z/n) \to H^3(Y, \mathbb Z/n) $ is $\Gamma$-equivariant for the 
$\Gamma$-group action on $H$ and the action of $\Gamma$ on $Y$ by the base-change claim in Lemma~\ref{L:BHpullback}.
Since $Y \to X$ is $\Gamma$ invariant,
the integration map $H^3(Y, \mathbb Z/n) \to H^1( X, \mathbb Z/n)$ is $\Gamma$-invariant. Hence the composition $\tilde{e}$ is $\Gamma$-invariant, and thus $e$ is $\Gamma$-invariant as well. Thus the image of $e$ lands in $H_3(H, \mathbb Z/n)^\Gamma$.
 
Hence $e$ defines a finite \'{e}tale $H_3(H, \mathbb Z/n)^\Gamma$-torsor over $X$. 
We let $X^0$ be this \'etale covering and let $\Hur^{m,0}_{H\rtimes \Gamma ,c_{H \rtimes \Gamma} }$ be the union of $X^0$ over all connected components $X$ of $\Hur^m_{H\rtimes \Gamma ,c_{H \rtimes \Gamma} } $.  

For $q\equiv 1\pmod n$ and $s_H \in H_3(H, \mathbb Z/n)^\Gamma$, we now define a covering of a component $X_q$ of $(\Hur^m_{H\rtimes \Gamma ,c_{H \rtimes \Gamma} } )_{\F_q}$
 using the homomorphism $\pi_1(X_q) \to H_3(H, \mathbb Z/n)^\Gamma$ that sends an element $\sigma \in \pi_1(X)$ to $ e(\sigma) - s_H \deg \sigma$, 
 where $\deg \sigma \in \hat{\mathbb Z}$ is the image of $\sigma$ under the natural homomorphism $\pi_1(X_q) \to \Gal(\bar{\F}_q/\mathbb F_q) \cong \hat{\mathbb Z}$ that sends the Frobenius element to $1$. Let $X_q^{ s_H, \mathbb F_q}$ be the finite \'{e}tale $H_3(H, \mathbb Z/n)^\Gamma$-torsor over $X_q$
  associated to this homomorphism and let $\Hur^{m,s_H, \mathbb F_q }_{H\rtimes \Gamma ,c_{H \rtimes \Gamma} }$ be the union of $X_q^{s_H, \mathbb F_q}$ over all connected components $X_q$ of $(\Hur^m_{H\rtimes \Gamma ,c_{H \rtimes \Gamma} } )_{\mathbb F_q} $.  (We put $\mathbb F_q$ in the exponent because this construction is not compatible with base change: For $q'$ a power of $q$ we do not necessarily have $(\Hur^{m,s_H, \mathbb F_q }_{H\rtimes \Gamma ,c_{H \rtimes \Gamma} })_{\mathbb F_{q'}} \cong \Hur^{m,s_H, \mathbb F_{q'} }_{H\rtimes \Gamma ,c_{H \rtimes \Gamma} }$.)  Note that when we base change our scheme $\Hur^m_{H\rtimes \Gamma ,c_{H \rtimes \Gamma} } $ to $\F_q$, we have implicitly chosen a map $\Z[\mu_n]\ra \F_q$.

\subsection{Properties and applications of Hurwitz schemes}\label{SS:PropHur}


The next lemma, proved in \S\ref{ss-long-proofs}, relates the Artin-Verdier fundamental class to the map $e$.

 \begin{lemma}\label{homomorphism-to-class} 
   Let $m$ be a positive integer and $X$ be a connected component of  $\Hur^m_{H\rtimes \Gamma ,c_{H \rtimes \Gamma} } $. Let $Z/X$ be the universal family of $(H \rtimes \Gamma)$-coverings of $\mathbb P^1$ and let $Y/X$ be the universal family of induced $\Gamma$-coverings of $\mathbb P^1$. Let $q$ be a prime power that is $1$ mod $n$ and is relatively prime to $|H||\Gamma|$.
 For $x \in X(\F_q)$ and $Z_x$ and $Y_x$ the fibers of $Z$ and $Y$ at $x$, respectively,  we have that  $e(\Frob_{q,x})$ is the Artin-Verdier fundamental class for $Z_x/Y_x$. 
 \end{lemma}

 Recall 
$N(\mathbf H, \Gamma, q^m, \mathbb F_q(t))$ was defined in Lemma~\ref{count-to-count} to count certain $H$-extensions of $\F_q(t)$. Using Lemma \ref{homomorphism-to-class}, we relate $N(\mathbf H, \Gamma, q^m, \mathbb F_q(t))$ to the count of points on $ \Hur^{m,s_H, \mathbb F_q }_{H\rtimes \Gamma ,c_{H \rtimes \Gamma} }$.

\begin{lemma}\label{lem-count-to-geom} For $q\equiv 1\bmod n $ relatively prime to $|H|\abs{\Gamma}$ and $\mathbf H = (H, s_H)$ we have, for every positive integer $m$,
\[\abs{ N(\mathbf H, \Gamma, q^m, \mathbb F_q(t)) } = \frac{1}{ \abs{ H^3(H \rtimes \Gamma, \mathbb Z/n)} } \abs{ \Hur^{m,s_H, \mathbb F_q }_{H\rtimes \Gamma ,c_{H \rtimes \Gamma} }(\mathbb F_q)}.\]
\end{lemma}

\begin{proof} The argument in \cite[Proof of Lemma 10.2 on p. 54]{Liu2024} gives a bijection between $\mathbb F_q$-points of $\Hur^m_{H\rtimes \Gamma, c_{H \rtimes \Gamma}}$ and surjections $\Gal ( \overline{\mathbb F_q(t)}/\mathbb F_q(t)) \to H \rtimes \Gamma$ such that the corresponding $H \rtimes \Gamma$-extension $K$ has $\rDisc K = q^m$, the associated $H$-extension is unramified everywhere and split completely at all places lying above $\infty$, and the associated $\Gamma$-extension $K^H/\mathbb F_q(t)$ is split completely above $\infty$.  Furthermore $N(\mathbf H, \Gamma, q^m, \mathbb F_q(t))$ counts such surjections with the additional condition that the  natural orientation on $H$ arising from the Artin-Verdier fundamental class on the associated $H$-extension is equal to $s_H$. 

It suffices to prove that the bijection of \cite[Proof of Lemma 10.2 on p. 54]{Liu2024} sends a surjection satisfying that additional condition to an $\mathbb F_q$-point of  $\Hur^m_{H\rtimes \Gamma, c_{H \rtimes \Gamma}}$ that admits exactly $\abs{ H^3(H \rtimes \Gamma, \mathbb Z/n)} $ lifts to $\mathbb F_q$-points of $\Hur^{m,s_H, \mathbb F_q }_{H\rtimes \Gamma ,c_{H \rtimes \Gamma}}$, and sends a surjection not satisfying that additional condition to an  $\mathbb F_q$-point of  $\Hur^m_{H\rtimes \Gamma, c_{H \rtimes \Gamma}}$ that does not lift to a $\mathbb F_q$-point of $\Hur^{m,s_H, \mathbb F_q }_{H\rtimes \Gamma ,c_{H \rtimes \Gamma}}$.

We can check whether a point lifts component-by-component and so fix a component $X_q$ of $(\Hur^m_{H\rtimes \Gamma, c_{H \rtimes \Gamma}})_{\mathbb F_q}$. By definition of $\Hur^{m,s_H, \mathbb F_q }_{H\rtimes \Gamma ,c_{H \rtimes \Gamma}}$ as the disjoint union of the $X_q^{s_H, \mathbb F_q}$, the number of lifts of $x \in X_q \subseteq (\Hur^m_{H\rtimes \Gamma, c_{H \rtimes \Gamma}})_{\mathbb F_q}$ to $\Hur^{m,s_H, \mathbb F_q }_{H\rtimes \Gamma ,c_{H \rtimes \Gamma}}(\mathbb F_q)$ is equal to the number of lifts of $x$ to $X_q^{s_H, \mathbb F_q}$. Since $X_q^{s_H, \mathbb F_q}$ is a finite \'etale  cover of $X_q$ of degree $|H_3( H, \mathbb Z/n)^\Gamma|$,
 an $\mathbb F_q$-point $x$ of $X_q$ has exactly $\abs{H_3( H, \mathbb Z/n)^\Gamma}=\abs{H^3( H\rtimes \Gamma, \mathbb Z/n)}$
(Lemmas~\ref{L:GammaH} and \ref{L:homdual}) 
  lifts to $\mathbb F_q$-points of  $X_q^{s_H, \mathbb F_q}$ if the image of $\Frob_{q,x}$ in the associated homomorphism $\pi_1(X_q) \to H_3( H, \mathbb Z/n)^\Gamma$ is trivial and $0$ otherwise.

By definition, this homomorphism sends $\Frob_{q,x}$ to \[ e(\Frob_{q,x}) -  s_H \deg \Frob_{q,x} = e(\Frob_{q,x})-s_H \] since $\deg \Frob_{q,x}=1$. 
This image is trivial if and only if $ e(\Frob_{q,x})=s_H$ which by Lemma \ref{homomorphism-to-class} happens if and only if the Artin-Verdier fundamental class for $Z_x/Y_x$ is $s_H$. Under the bijection of \cite[Proof of Lemma 10.2 on p. 54]{Liu2024}, $K^H$ is the function field of $Y_x$ and $K$ is the function field of the associated $H$-covering, so this occurs if and only if the surjection satisfies the additional condition, as desired. \end{proof}

It remains to count points on  $\Hur^{m,s_H, \mathbb F_q }_{H\rtimes \Gamma ,c_{H \rtimes \Gamma} }$.

 Liu \cite{Liu2022} defines an $\omega$-invariant associated to an $H \rtimes \Gamma$ extension. Liu's $\omega$-invariant is a certain homomorphism, but we checked in Lemmas \ref{yuan-liu} and \ref{liu-domain-comparison} 
that the set of homomorphisms is isomorphic to $H_2(H, \mathbb Z)^\Gamma$ by an isomorphism sending Liu's $\omega$-invariant to our lifting invariant. We will apply this isomorphism, treating Liu's $\omega$-invariant as an element of $H_2(H, \mathbb Z)^\Gamma$ 
and thus as being equal to the lifting invariant. This will allow us to make use of prior work \cite{Liu2022,LandesmanLevy} which counts points on Hurwitz spaces with a given lifting invariant.

Liu \cite[\S4.2]{Liu2022} shows that the $\omega$-invariant is determined by the lifting invariant defined by Ellenberg-Venkatesh-Westerland and the second author. Since this lifting invariant is constant on components of  $\Hur^m_{ H \rtimes \Gamma, c_{H \rtimes \Gamma}}$ by \cite[Theorem 6.1]{Wood2021}, 
Liu's $\omega$ and our lifting invariant are determined by the component. This, together with Lemma \ref{homomorphism-to-class}, has the following consequence. Let $\omega_H$ be the image of $s_H$ in $H_2(H, \mathbb Z)^\Gamma$.  

\begin{lemma}\label{empty-component-observation} Let $m$ be a positive integer, let $q$ be a prime power that is $1$ mod $n$, and is relatively prime to $|H||\Gamma|$.  Let $X_q$ be a connected component of $(\Hur^m_{H\rtimes \Gamma ,c_{H \rtimes \Gamma} } )_{\mathbb F_q} $. If the cover $X_q^{s_H, \mathbb F_q}$ of $X_q$ has a $\mathbb F_q$-rational point, then the lifting invariant of $X_q$ is $\omega_H$. \end{lemma}

\begin{proof} An $\mathbb F_q$-point of $X_q^{s_H, \mathbb F_q}$ lies over an $\mathbb F_q$-point $x\in X_q$. By definition, Frobenius acts on the fiber over $x$ by adding $e( \Frob_{q,x})-s_H$. For the fiber to contain a $\mathbb F_q$ point, we must have $e( \Frob_{q,x})-s_H=0$. Thus by Lemma \ref{homomorphism-to-class}, $x$ parameterizes a $H\rtimes \Gamma $-covering of $\mathbb P^1$ with Artin-Verdier fundamental class $s_H$. Thus its lifting invariant is $\omega_H$, and so the lifting invariant of the component $X_q$ is $\omega_H$. \end{proof}

Motivated by Lemma \ref{empty-component-observation}, let $\Hur^{m,\omega}_{ H \rtimes \Gamma, c_{H \rtimes \Gamma}}$ be the union of the components of $\Hur^m_{ H \rtimes \Gamma, c_{H \rtimes \Gamma}}$ with lifting invariant $\omega$.

Using prior work, we may now control the number of points on ``bad'' components of  $\Hur^{m,\omega}_{ H \rtimes \Gamma, c_{H \rtimes \Gamma}}$.

\begin{lemma}\label{component-bounding-lemma} Suppose that $H_\Gamma=1$.
Let $q$ be a prime power that is $1$ mod $n$, and is relatively prime to $|H||\Gamma|$. 
 Fix $\omega \in H_2(H, \mathbb Z)^\Gamma [q-1]$.  Let $a$ be a positive integer. Then for $b \geq 0$ an integer,
if $q$ and $b$ are sufficiently large given $H$ and $\Gamma$,  
  the fraction of $\mathbb F_q$-points on $\bigcup_{  m\leq b}  \Hur^{m,\omega}_{ H \rtimes \Gamma, c_{H \rtimes \Gamma}}$ that  lie on a component parameterizing curves where at least one conjugacy class in $c$ has $< a$ branch points of that conjugacy class is $O(a/b)$ where the implicit constant may depend on $q,H,$ and $\Gamma$. 
\end{lemma}

\begin{proof}
In this proof ``sufficiently large'' always meaning only depending on $H,\Gamma$, and the implicit constants in $O$ notation depend on $q,H,\Gamma$.

Let $d_{ \Gamma}(q)$ be the number of orbits of $q$th powering on nontrivial conjugacy classes in $\Gamma$. Since each conjugacy class in $c_{H\rtimes \Gamma}$ contains a unique nontrivial conjugacy class of $\Gamma$, this is also the number of orbits of $q$th powering on nontrivial conjugacy classes in $c_{H \rtimes \Gamma}$.

We will show that the number of points of  $\bigcup_{  m\leq b}  \Hur^{m,\omega}_{ H \rtimes \Gamma, c_{H \rtimes \Gamma}}$ is, for $q$ and $b$ sufficiently large, greater than a positive constant depending on $q,H,\Gamma$ times  $  b^{ d_\Gamma(q)-1}q^b $. We will then show that the number of points that  lie on a component parameterizing curves where at least one conjugacy class in $c$ has $\leq a$ branch points of that conjugacy class is $O( a b^{ d_\Gamma(q)-2}q^b)$. Combining these, we obtain the desired statement.

By \cite[Lemma 8.4.4]{LandesmanLevy}, for $q$ sufficiently large, 
each geometrically connected component of $\Hur^{m}_{ H \rtimes \Gamma, c_{H \rtimes \Gamma}}$ contains at least $q^{m}/2$
 points. The same is true for $\Hur^{m,\omega}_{ H \rtimes \Gamma, c_{H \rtimes \Gamma}}$ as it is a union of components of $\Hur^{m}_{ H \rtimes \Gamma, c_{H \rtimes \Gamma}}$. Following \cite[\S4.3]{Liu2022}, let $\pi^\omega_{H \rtimes \Gamma, c_{H\rtimes \Gamma}} (q,m) $ be the number of geometrically connected components of $\Hur^m_{ H \rtimes \Gamma, c_{H \rtimes \Gamma}}$ with lifting invariant $\omega$.

By $H_\Gamma=1$ and \cite[Lemma 4.4(1)]{Liu2022} (where we have put lines over Liu's variables $a$ and $b$ to avoid ambiguity with our $a$ and $b$) there is a positive integer $M_{\overline{a}}$, non-empty set $E_{\overline{a}}$ of residues mod $M_{\overline{a}}$, and positive number $r_{\overline{a},\overline{b}}'$ for $\overline{b} \in E_{\overline{a}}$, all depending only on $H$ and $ \Gamma$ and the residue class $\overline{a}$ of $q$ modulo $|H \rtimes \Gamma|^2$, such that if $m \equiv \overline{b} \in E_{\overline{a}}$ modulo $M_{\overline{a}}$ we have
\[\pi^\omega_{H \rtimes \Gamma, c_{H\rtimes \Gamma}} (q,m) = r'_{\overline{a},\overline{b}} m^{ d_{\Gamma}(q)-1} + O (m^{ d_{\Gamma}(q)-2}).\]

Letting $m$  be the greatest integer $\leq b$ whose residue mod $M_{\overline{a}}$ lies in $E_{\overline{a}}$, which is always $\geq b -M_{\overline{a}}$ and thus each component contributes at least $q^{ b-M_{\overline{a}}}/2$ points,  ensures that the number of components is, for $b$ sufficiently large, greater than a positive constant times $m^{ d_\Gamma(q)-1} $ and thus greater than a positive constant times $b^{ d_\Gamma(q)-1} $. Hence the total number of points of  $\bigcup_{  m\leq b}  \Hur^{m,\omega}_{ H \rtimes \Gamma, c_{H \rtimes \Gamma}}$, for $q$ and $b$ sufficiently large, is greater than a positive constant times $  b^{ d_\Gamma(q)-1}q^b $, as desired.

Again  by \cite[Lemma 8.4.4]{LandesmanLevy}, for $q$ sufficiently large, each geometrically connected component of $\Hur^{m}_{ H \rtimes \Gamma, c_{H \rtimes \Gamma}}$ contains at most $2q^m $ points. Following \cite{Liu2024}, let $\mathbb Z^{ c_{H \rtimes \Gamma}/ H\rtimes \Gamma}_{\equiv q}$ be the set of functions from conjugacy classes in $c$ to nonnegative integers which take the same value on each component and its $q$th power. Associated to each component is an element of  $\mathbb Z^{ c_{H \rtimes \Gamma}/ H\rtimes \Gamma}_{\equiv q}$: The function whose value on each conjugacy class is the number of branch points with ramification given by that conjugacy class. Furthermore, this element of $\mathbb Z^{ c_{H \rtimes \Gamma}/ H\rtimes \Gamma}_{\equiv q}$ is always a function that sums to $m$. We are interested in components corresponding to elements of $\mathbb Z^{ c_{H \rtimes \Gamma}/ H\rtimes \Gamma}_{\equiv q}$ which take some value $\leq a$. The number of such elements of $\mathbb Z^{ c_{H \rtimes \Gamma}/ H\rtimes \Gamma}_{\equiv q}$ is $O ( a m^{ d_{\Gamma}(q)-2})$ as they are determined by $d_{\Gamma}(q)$ distinct values, one of which is bounded by $a$ and the remainder of which are bounded by $m$, but with a linear relation such that any one is determined by the remaining ones. By \cite[Theorem 12.4]{Liu2024} and \cite[Lemma 3.3]{Ellenberg2005}, the number of components corresponding to each element of $\mathbb Z^{ c_{H \rtimes \Gamma}/ H\rtimes \Gamma}_{\equiv q}$ is $O(1)$. (This is proved over $\mathbb C$ but by \cite[Lemma 10.3]{Liu2024} with $i=0$ the same statement holds over $\bar{\F}_q$.) Thus the total number of components is $O ( a m^{ d_{\Gamma}(q)-2})$ and hence the total number of $\mathbb F_q$-points is $O ( a m^{ d_{\Gamma}(q)-2} q^m)$  which summed over $m\leq b$ is $O ( a b^{ d_{\Gamma}(q)-2} q^b)$, as desired. \end{proof}

The key new result we will need is the following:

\begin{theorem}\label{component-comparison} 
Let $\Gamma$ be the finite group and $n$ the positive integer coprime to $\abs{\Gamma}$ fixed throughout the paper, and 
$\mathbf H$ the finite $n$-oriented $\Gamma$-group fixed throughout this section.  Suppose that $H_\Gamma=1$. 
 Let $q$ be prime power $q\equiv 1\bmod n$ with $\gcd(q, \abs{\Gamma}\abs{H})=1$ and $\gcd(q-1, \abs{H})=\gcd(n,H)$,

Let $m$ be a positive integer and let $X_q$ be a component of $\Hur^{m,\omega_H}_{ H \rtimes \Gamma, c_{H \rtimes \Gamma}, \mathbb F_q}$. Let $a_{X_q}$ be the minimum over conjugacy classes in $c_{H\rtimes\Gamma}$ of the number of branch points of that conjugacy class in covers parameterized by $X_q$. 

There is $\delta>0$, only depending on $\Gamma,H,q$, such that as long as $q$ is sufficiently large in terms of $H, \Gamma$ we have
\[ | X_q^{ s_H, \mathbb F_q}( \mathbb F_q) |=  ( | H^2(H \rtimes \Gamma, \mathbb Z/n)|   +  O( e^{ - \delta a_{X_q} }) ) |X_q(\mathbb F_q)| ,\]
with the implicit constant in the $O$ depending on $\Gamma,H,q$.
 \end{theorem}

  We now show how Theorem \ref{component-comparison} implies Theorem \ref{ff-theorem}.

\begin{proof}[Proof of Theorem \ref{ff-theorem}]
Through this proof, the implicit constant in  $O$ notation depends on $\Gamma,H,q$.

By Lemmas \ref{count-to-count} and \ref{lem-count-to-geom} we have 
\begin{equation}\label{num-to-geom} \sum_{m\leq b} \sum_{K \in E_\Gamma( q^m , \mathbb F_q(t))}  \Sur ( \bGal ( K^{ \operatorname{un}, \abs{\Gamma}'} / K), \mathbf H)= \frac{\abs{H^\Gamma} }{ \abs{H} \abs{ H^3 (H \rtimes \Gamma, \mathbb Z/n)}} \sum_{m\leq b} \abs{ \Hur_{H\rtimes \Gamma, c_{H \rtimes \Gamma}}^{m,s_H, \mathbb F_q}(\mathbb F_q)}\end{equation}
and by \eqref{eq-count-to-geom} we have
\begin{equation}\label{denom-to-geom}\sum_{m\leq b}\abs{ E_\Gamma( q^m , \mathbb F_q(t))} = \sum_{m\leq b} \abs{ \Hur_{\Gamma, c_\Gamma}^m(\mathbb F_q)}.\end{equation}

We have \cite[(9.3)]{LandesmanLevy}
\begin{equation}\label{denom-to-simplified-num} \frac{\sum_{m\leq b} \abs{ \Hur_{H\rtimes \Gamma, c_{H \rtimes \Gamma}}^{m,\omega_H}(\mathbb F_q)} }{\sum_{m\leq b} \abs{ \Hur_{\Gamma, c_\Gamma}^m(\mathbb F_q)}} = 1 + O(1/b). \end{equation}

Combining \eqref{num-to-geom}, \eqref{denom-to-geom} and \eqref{denom-to-simplified-num}, we see that to establish Theorem \ref{ff-theorem} it suffices to show
\begin{equation}\label{ff-theorem-1} \frac{\sum_{m\leq b} \abs{ \Hur_{H\rtimes \Gamma, c_{H \rtimes \Gamma}}^{m,s_H, \mathbb F_q}(\mathbb F_q)}} {\sum_{m\leq b} \abs{ \Hur_{H\rtimes \Gamma, c_{H \rtimes \Gamma}}^{m,\omega_H} (\mathbb F_q)}} = \abs{H^2(H \rtimes \Gamma, \mathbb Z/n) } + O(1/b),  \end{equation}
or equivalently
\begin{equation}\label{ff-theorem-2} 
\sum_{m\leq b} \abs{ \Hur_{H\rtimes \Gamma, c_{H \rtimes \Gamma}}^{m,s_H, \mathbb F_q}(\mathbb F_q)} - \abs{H^2(H \rtimes \Gamma, \mathbb Z/n) } \sum_{m\leq b} \abs{ \Hur_{H\rtimes \Gamma, c_{H \rtimes \Gamma}}^{m,\omega_H} (\mathbb F_q)} = O ( \frac{1}{b} \sum_{m\leq b} \abs{ \Hur_{H\rtimes \Gamma, c_{H \rtimes \Gamma}}^{m,\omega_H} (\mathbb F_q)})  . \end{equation}

We have (by Lemma \ref{empty-component-observation})
\begin{equation}\label{by-components-1} \sum_{m\leq b} \abs{ \Hur_{H\rtimes \Gamma, c_{H \rtimes \Gamma}}^{m,s_H, \mathbb F_q}(\mathbb F_q)}  = \sum_{m\leq b}  \sum_{\substack{ X_q \subseteq  \Hur_{H\rtimes \Gamma, c_{H \rtimes \Gamma}}^m \\ \textrm{component} \\ \textrm{lifting invariant } \omega_H} }  \abs{ X_q^{\mathbb F_q, s_H}( \mathbb F_q)} \end{equation}
and (by definition)
\begin{equation}\label{by-components-2} \sum_{m\leq b} \abs{ \Hur_{H\rtimes \Gamma, c_{H \rtimes \Gamma}}^{m,\omega_H}(\mathbb F_q)}=  \sum_{m\leq b}  \sum_{\substack{ X_q \subseteq  \Hur_{H\rtimes \Gamma, c_{H \rtimes \Gamma}}^m \\ \textrm{component} \\ \textrm{lifting invariant } \omega_H} }  \abs{ X_q( \mathbb F_q)}\end{equation}
 so the left hand side of \eqref{ff-theorem-2} is
 
\begin{equation}\label{lhs}  \sum_{m\leq b}  \sum_{\substack{ X_q \subseteq  \Hur_{H\rtimes \Gamma, c_{H \rtimes \Gamma}}^m \\ \textrm{component} \\ \textrm{lifting invariant } \omega_H} }  \Bigl(   \abs{ X_q^{\mathbb F_q, s_H}( \mathbb F_q)} -  \abs{H^2(H \rtimes \gamma, \mathbb Z/n) }  \abs{ X_q( \mathbb F_q)} \Bigr).\end{equation} 
 
For each $X_q$ we apply Theorem \ref{component-comparison} and get 
 \begin{equation}\label{bounded-lhs}  \abs{ X_q^{\mathbb F_q, s_H}( \mathbb F_q)} -  \abs{H^2(H \rtimes \Gamma, \mathbb Z/n) }  \abs{ X_q( \mathbb F_q)}   = O ( e^{ -\delta a_{X_q}}   \abs{ X_q( \mathbb F_q)})\end{equation} so \eqref{lhs} is  \[  O\Bigl(\sum_{m\leq b}  \sum_{\substack{ X_q \subseteq  \Hur_{H\rtimes \Gamma, c_{H \rtimes \Gamma}}^m \\ \textrm{component} \\ \textrm{lifting invariant } \omega_H} }   e^{ -\delta a_{X_q}} \abs{ X_q( \mathbb F_q)}\Bigr).\]
By Lemma \ref{component-bounding-lemma}, the terms with $a_{X_q} =a$ on the right hand side of \eqref{bounded-lhs} contribute 
\[O \Bigl(  e^{- \delta a}  \cdot \frac{a}{b}  \sum_{m\leq b} \abs{ \Hur_{H\rtimes \Gamma, c_{H \rtimes \Gamma}}^{m,\omega_H}(\mathbb F_q)} \Bigr)\]
so summing over $a$, \eqref{lhs} is
\[ O \Bigl(  \sum_{a=1}^\infty \frac{a}{b}  e^{- \delta a} \sum_{m\leq b} \abs{ \Hur_{H\rtimes \Gamma, c_{H \rtimes \Gamma}}^{m,\omega_H}(\mathbb F_q)} \Bigr) = O \Bigl( \frac{1}{b}  \sum_{m\leq b} \abs{ \Hur_{H\rtimes \Gamma, c_{H \rtimes \Gamma}}^{m,\omega_H}(\mathbb F_q)}\Bigr)\]
since $\sum_{a=0}^\infty a e^{-\delta a}=O(1)$. Hence \eqref{lhs} is bounded by the right hand side of \eqref{ff-theorem-2}, so we have \eqref{ff-theorem-2}, which we already saw suffices.\end{proof}

 
  
  Theorem \ref{component-comparison} is trivial when $X_q$ is not geometrically irreducible as in that case neither $X_q$ nor $X_q^{s_H, \mathbb F_q}$ has any $\mathbb F_q$-points, and also trivial if $a(X_q)$ is small. In the remaining case, we will prove Theorem \ref{component-comparison} by showing that $X_q^{s_H, \mathbb F_q}$ has $\abs{H^2( H \rtimes \Gamma, \mathbb Z/n)}$ components and each component has $ ( 1+O( e^{ - \delta m }) ) |X_q(\mathbb F_q)| $ $\mathbb F_q$-points. These will follow from analogous statements over the complex numbers, one about the image of the map $e$  from the fundamental group to $H_3(H, \mathbb Z/n)$, and one about the homology of connected components of $\Hur^{m, 0}_{H \rtimes \Gamma, c_{H \rtimes \Gamma}}$. We will prove these by topological methods.

\subsection{The braid fundamental class}\label{SS:toparg}

We now explain the topological analogue of our geometric setup.

Consider a configuration of $m$ points in the (topological) plane $P$.   For psychological convenience we may take the points to be in a line.
Let $P^\circ$ be the plane minus those points.
The braid group $\Br_m$ on $m$ strands naturally maps to the mapping class group of  $P^\circ$ \cite[Section 9.1.4]{Farb2011}, and thus also acts on the 
fundamental group $\pi_1(P^\circ)$ (taking a base point near infinity).
Given a finite group $G$, the braid group $\Br_m$ then acts on the finite set of
homomorphisms $\pi_1(P^\circ)\ra G.$
This action preserves the set of homomorphisms where the conjugacy class of the local monodromy around each point lies in a subset $c$ of $G$ that is closed under taking $G$-conjugates.  

Given a braid $\sigma\in \Br_m,$ we have a geometric realization of the braid as $m$ sections $s_1,\dots,s_m$ of the map $P\times [0,1] \ra [0,1]$, 
such that the sections do not intersect,  and the set $s_i(0)$ are our $m$ chosen points in $P$,  as is the set of $s_i(1)$.  We can then complete  $P\times [0,1]$ to
$S^2\times [0,1]$, and identify $S^2\times 0$ and $S^2\times 1$.  Then in the compact $3$-manifold $S^2\times S^1$, our braids, along with $\infty\times S^1$, describe a link $L_\sigma$.

For our given $\Gamma$ and $H$, fix a surjective homomorphism $f:\pi_1(P^\circ) \ra H\rtimes \Gamma$ such that 
the local monodromy around each missing point is in $c_{H\rtimes \Gamma}$, and the local monodromy around $\infty$ is either in $c_{H\rtimes \Gamma}$ or trivial. 
Let $\sigma\in \Stab_f\sub \Br_m.$
Then we can check that there is a unique homomorphism $\tilde{f} : \pi_1(S^2\times S^1\setminus L_\sigma) \ra H\rtimes \Gamma$
that restricts to $f$ on $P^\circ\times 0$ and is trivial on the loop $\infty \times S^1$.
We let $\tilde{f}_\Gamma$ be the composition of $\tilde{f}$ with $H\rtimes \Gamma \ra \Gamma$.
We let $M_\sigma$ be the compact oriented $3$-manifold such that $\tilde{f}_\Gamma$ corresponds to a degree $|\Gamma|$ branched covering $u_\Gamma: M_\sigma \ra S^2\times S^1$, with an action of $\Gamma$.
We note by our monodromy condition that the restriction of $\tilde {f} : \pi_1(M_\sigma \setminus u_\Gamma^{-1} L_\sigma)\ra H$
factors through a unique map $ \pi_1(M_\sigma)\ra H$.
We can define a map from the stabilizer $\Stab_f\sub \Br_m$ of $f$ to $H_3(H,\Z)^\Gamma$, by sending $\sigma$
to the image of the fundamental class of $M_\sigma$ in $H_3(H,\Z)^\Gamma$.  In more detail,
the fundamental class of $M_\sigma$ is an element of $H_3(M_\sigma,\Z)$, which
maps to  $H_3(\pi_1(M_\sigma),\Z)$ by the usual map, and then we compose with $H_3(\pi_1(M_\sigma),\Z)\ra H_3(H,\Z)$ using the map $ \pi_1(M_\sigma)\ra H$ above, and finally we note 
that the image of the fundamental class of $M_\sigma$ in $H_3(H,\Z)$ is $\Gamma$-invariant.

We call this map $\Stab_f\ra H_3(H,\Z)^\Gamma$  the \emph{braid fundamental class map}. 

We now explain the relevance of the braid fundamental class map to our problem.

Let $X$ be a connected component of $(\Hur^m_{H\rtimes \Gamma ,c_{H \rtimes \Gamma} } )_{\mathbb C}$. We can express the topological fundamental group $X(\mathbb C)$ in terms of the braid group on $m$ strands.  
Let $x$ be a point of a component $V$ of  $\Hur^m_{H\rtimes \Gamma ,c_{H \rtimes \Gamma} } (\bbC)$, such that $x$ corresponds to a cover unramified over $P^\circ$.   
 Let $f: \pi_1(P^\circ) \ra H\rtimes \Gamma$ be the map corresponding to that cover.
Then the isomorphism $\Br_m \ra \pi_1(\Conf^m \A^1(\bbC))$ takes $\Stab_f$ to $\pi_1(V)$. (See \cite[Section 1.4]{Fried} and \cite[Proof of Theorem 12.4]{Liu2024}.) We have the following lemma, proven later in \S\ref{ss-long-proofs}.

\begin{lemma}\label{b-vs-bfc} Let $m$ be a nonnegative integer and $f$ and $X$ as above. Under the above identification of $\Stab_f$ and $\pi_1(X(\mathbb C))$, for any braid $\sigma$ in the stabilizer of $f$, the reduction mod $n$ of the braid fundamental class of $\sigma$ is the product of $e(\sigma)$ by an element of $(\Z/n)^\times$. \end{lemma}

It should be possible to check that the element of $(\Z/n)^\times$ mentioned in Lemma \ref{b-vs-bfc} is $1$ for any braid as long as the fixed isomorphism $\mu_n \to \mathbb Z/n$ sends $e^{2\pi i/n}$ to $1$, but this would be unnecessary and add additional technical details so we do not attempt it.

The goal of this subsection is to prove the following result, which, because of Lemma \ref{b-vs-bfc}, will allow us to control the image of $e$:
\begin{proposition}\label{topological-braid-surjective} 
Let $m$ be a nonnegative integer, and consider the (topological)  plane $P^\circ$ punctured at $m$ points.  
Let $f$ be a surjective homomorphism $f: \pi_1(\P^\circ) \ra H \rtimes \Gamma$.   Assume that the local monodromy of $f$ around each puncture is in $c_{H\rtimes \Gamma}$, and the local monodromy around $\infty$ is either in $c_{H\rtimes \Gamma}$ or trivial.

Assume that for each conjugacy class of $H\rtimes \Gamma$, the number of points in the configuration space with local monodromy in that conjugacy class (i.e a small loop around that point goes to that conjugacy class)
  is sufficiently large depending $H$ and $\Gamma$.

Then the braid fundamental class map from the stabilizer $\Stab_f$ in the braid group 
 to $H_3( H,\Z)^\Gamma$ is surjective. \end{proposition}

To prove Proposition \ref{topological-braid-surjective} we define a more general fundamental class. Let $M$ be an oriented compact 3-manifold,
 $L$ an oriented link in $M$, and $f \colon \pi_1(M \setminus L) \to H \rtimes \Gamma $ that, for each knot of $L$, sends a small loop around that knot to an element of a conjugate of $\Gamma$. Let $\tilde{M}$ be the induced 
oriented compact 3-manifold such that $\tilde{M}\ra L$ is a branched 
 $\Gamma$-covering. 
 We have a morphism $\pi_1(\tilde{M})\to H$ and hence a map $\tilde{M} \to  BH$ where $BH$ is the classifying space of $H$. The fundamental class of $\tilde{M}$ gives an element of $H_3(H,\Z)$ which is manifestly $\Gamma$-invariant. Let $ [M, L, f] \in H_3(H,\Z)^\Gamma$ be the fundamental class of $\tilde{M}$.

\begin{lemma} \label{L:braid-fundamental-homomorphism}
Let $m$ be a nonnegative integer.
Let $P^\circ$ be the topological plane punctured at $m$ points, as above.
Let $f$ be a surjective homomorphism $f: \pi_1(P^\circ) \ra H \rtimes \Gamma$ such that 
the local monodromy around each puncture is in $c_{H\rtimes \Gamma}$ and the local monodromy at $\infty$ is either in $c_{H \rtimes \Gamma}$ or trivial. 
The braid fundamental class is a homomorphism from the stabilizer of $f$ in $\Br_m$ to $H_3( H,\Z)^\Gamma$. \end{lemma}

\begin{proof} The proof is essentially the same as \cite[Lemma 3.2]{Sawin2024}.
We use \cite[Lemma 2.13]{Sawin2024} to interpret the fundamental class of $M_\sigma$ in $H_3(H,\Z)$ as the  class of $M_\sigma$ in the bordism group of $BH$.

We check that the braid fundamental class of a product of two braids is the sum of the braid fundamental classes.

  The braid fundamental class is constructed by starting with a braid and producing a link in $S^2 \times S^1$. Given two braids whose product is a third braid, we can construct a surface in the product of $S^2$ with a pair of pants, whose restriction to each boundary component is one of the three links in $S^2 \times S^1$.

To construct the surface, we view each braid as an element of the fundamental group of the configuration space of points in $S^1$.
We then use the fact that when the product of two fundamental group elements of some space is a third, we can find a map from a pair of pants to the space whose boundary  consists of three loops representing the three elements of the fundamental group.

The complement of this surface in the product of $S^2$ with a pair of pants is a fibration over the pair of pants with fiber a punctured sphere. Hence its fundamental group is an extension of the fundamental group of the pair of pants by the fundamental group of the punctured sphere. The fundamental group of the pair of pants is a free group on two generators, each corresponding to one of the boundary loops, and thus each corresponding to a braid that lies in the stabilizer of $f$ in $\Br_m$. We can lift these two generators explicitly to the fundamental group of the surface complement by taking the product of the corresponding loop in the pair of pants with the point $\infty$ in $\mathbb P^1$.  Since each generator of the free group stabilizes the  homomorphism $f $ from the fundamental group of the punctured sphere to $H \rtimes \Gamma$, we can extend the homomorphism $f$ uniquely to a homomorphism from the fundamental group of the surface complement to $H \rtimes \Gamma $ which is trivial on the two generators of the braid group. 

This gives a $\Gamma$-covering of the product of $S^2$ with a pair of pants branched at this surface, and a map from the branched covering to $BH$.  This gives a four-manifold with boundary witnessing the desired relation in the homology of $BH$ between the three braid fundamental classes.\end{proof}

\begin{lemma}\label{represented-by-3} Every class in $H_3( H,\Z)^\Gamma$ is of the form $[M,L,f]$
 for some $3$-manifold $M$, link $L$ in $M$, and $ f \colon \pi_1(M \setminus L) \to H \rtimes \Gamma$. In fact, we may take $L$ empty. \end{lemma}

\begin{proof} In the case $L$ is empty, $[M]$ represents a class in $H_3( H \rtimes \Gamma,\Z)$, from which $[M, L,f]$ is obtained via the trace map $H_3( H \rtimes \Gamma,\Z) \to H_3(H,\Z)^\Gamma$. Because $H$ and $\Gamma$ have coprime orders, the trace map $H_3( H \rtimes \Gamma,\Z) \to H_3(H,\Z)^\Gamma$ is surjective, so it suffices to show every class in $H_3( H \rtimes \Gamma,\Z) $ is represented by a 3-manifold, which follows from \cite[Lemma 2.13]{Sawin2024}. \end{proof}

\begin{lemma}\label{represented-by-link} The group $H_3(H,\Z)^\Gamma$ is generated by classes of the form $[S^3, L, f]$ where $L$ is a link in $S^3$.  \end{lemma}

\begin{proof} Fix an element of $H_3(H)^\Gamma$, which by Lemma \ref{represented-by-3} arises from a triple $[M, \emptyset, f]$. By the Lickorish-Wallace theorem, $M$ arises from surgery on a link in $S^3$. Thus, there exists a link $L$ in $S^3$ such that when we remove a tubular neighborhood of $L$ and glue in a new solid torus at each boundary component of the neighborhood, we obtain $M$.

Furthermore, this is surgery with slope $\pm 1$, so the meridian of the new solid torus is the sum of the meridian and standard longitude of the removed solid torus. However, it will be more convenient for us to choose coordinates in $B^2 \times S^1$ for the removed solid torus where the new meridian is simply the longitude, which is possible as adding any integer multiple of the meridian to a longitude gives a new longitude.

This isomorphism gives a homomorphism from the fundamental group of $S^3$ minus the tubular neighborhood of $L$ to $H \rtimes \Gamma$. We can extend this to the complement of $L$, but there is no reason for the local monodromy elements to lie in $\Gamma$. However, since $H_\Gamma=1$, the conjugates of $\Gamma$ generate $H \rtimes \Gamma$, and thus we can write the local monodromy around each component $K_i$ of $L$ as a product of $n_i$ conjugates of elements of $\Gamma$.  We define a link $L'$ which for each component $K_i$ of $L$ consists of $n_i$ loops in the tubular neighborhood of $K_i$, parallel to each other and to the meridian. (Inside $B^2 \times S^1$, this is the product of $n_i$ points with $S^1$).

 We can then extend $f$ to $f' \colon \pi_1( S^3 \setminus L') \to H \rtimes \Gamma$, in such a way that the local monodromy around each new component is the corresponding conjugate of an element  of $\Gamma$. Indeed, the fundamental group of $B^2 \times S^1$ minus the product of $n_i$ points with $S^1$ is the product of the free group on $n_i$ letters with $\mathbb Z$, where the loop around each point goes to the corresponding letter, the meridian goes to the product of the letters, and the longitude goes to a generator of $\mathbb Z$.  Restricted to the fundamental group of the torus boundary, $f$ sends the product of the letters to the local monodromy of $f$ along $K_i$ and the generator of $\mathbb Z$ (the meridian curve in $M$) to the identity, so by writing the local monodromy as a product, we can extend the homomorphism.
 
 Let us now calculate the difference $[M, \emptyset, f] - [S^3, L' ,f']$. The difference between the fundamental class of $\tilde{M}$ and the fundamental class of $\tilde{S}^3$ is the fundamental class of the $\Gamma$-covering of the tubular neighborhood of $L$ we removed minus the $\Gamma$-covering of the solid tori we added. In other words, it is the fundamental class of the $\Gamma$-covering of the manifold obtained by gluing the tubuluar neighborhood of $L$ to the solid tori along their common boundary. Since the union of these two tori is $S^3$, this manifold is itself a disjoint union of $\Gamma$-coverings of copies of $S^3$ containing links, so we can write $[M,\emptyset,f]$ as a sum and difference of classes of links, as desired. \end{proof}
 
 \begin{lemma}\label{link-vs-bfc} 
Let $p$ and $m$ be nonnegative integers.
 Fix a configuration of $p$  points in the plane $P$ and a homomorphism $f$ from the fundamental group of their complement in $P$
 to $H \rtimes \Gamma$, with local monodromy around each point in $c_{H\rtimes \Gamma}$ and local monodromy around $\infty$ either in $c_{H\rtimes \Gamma}$ or trivial.
 
 Fix a disc $D$ in the plane, containing $m$ of the points. Consider a braid $\sigma$ in the braid group on $m$ strands, shifting only the points in $D$. Let $\sigma'$ be the braid in the braid group on $p$ strands obtained by extending $\sigma$ so that the points outside $D$ are unmoved. Assume that $\sigma'$ fixes $f$. Let $L(\sigma)$ be the link in $D\times S^1\subset S^3$  obtained by closing up the braid $\sigma$. 
 
 There exists a unique homomorphism $\tilde{f}$ from the fundamental group of complement of $L(\sigma)$ in $S^3$ to $H \rtimes \Gamma$ whose restriction to $D$ minus $m$ points is equal to the restriction of $f$.
 
 Then $[S^3, L(\sigma), \tilde{f}] $  is the braid fundamental class of $\sigma'$ (with respect to $f$). 
 \end{lemma}
 
 \begin{proof} To show that $\tilde{f}$ exists and is unique, we use the classical fact that the fundamental group of the complement of $L(\sigma)$ is the quotient of the fundamental group of the $m$-times-punctured disc, i.e. the free group on generators $x_1,\dots,x_m$  by the relations $x_i =\sigma(x_i)$ for $i$ from $1$ to $m$. This fact itself may be checked using the van Kampen theorem applied to the decomposition of $S^3$ into $D \times S^1$ and a complementary solid torus, whose intersection is homotopic to a torus. The fundamental group of the complement of $L(\sigma)$ in $D \times S^1$ is the semidirect product $\langle x_1,\dots, x_m \rangle \rtimes \langle \sigma \rangle$, and adding the solid torus produces the relation $\sigma=1$ which gives the stated fundamental group. Since $f$ is $\sigma'$-invariant, the restriction of $f$ to the fundamental group of the $m$-times punctured disc is $\sigma$-invariant, and hence factors uniquely through the fundamental group of the complement of $L(\sigma)$ in $S^3$.
 
  We must show two manifolds, one a branched $\Gamma$-covering of $ S^3$ and the other a branched $\Gamma$-covering of $S^2\times S^1$, with provided $H$-coverings,  have identical fundamental classes in $H_3(H,\Z)$. To do this, we consider maps from each manifold to the classifying space $BH$ of $H$. The fundamental classes in $H_3(H,\Z)$ of these manifolds may be represented by the images of these maps, viewed as chains. Since the $H \rtimes \Gamma$-coverings of $S^3$ and $S^2 \times S^1$ become isomorphic when restricted to $D \times S^1$, we can choose the maps from the $\Gamma$-coverings to $BH$ to agree on the $\Gamma$-coverings of $D \times S^1$.
  
  Having done this, the difference of the two fundamental classes is represented by a map from a 3-manifold given by gluing a $\Gamma$-covering of the complement of $D \times S^1$ in $S^3$ to a $\Gamma$-covering of the complement of $D \times S^1$ in $S^2 \times S^1$. The relevant covering of $S^3$ is unbranched in the complement of $D \times S^1$, which is a solid torus. The relevant covering of $S^2 \times S^1$ is branched at the product of $p-m$ points outside $D$ with $S^1$, and possibly branched also at the product of the point at $\infty$ of $S^1$. Thus the difference of the two fundamental classes is represented by a map from a $3$-manifold given by gluing an unbranched $\Gamma$-covering of the complementary solid torus to $D \times S^1$ in $S^3$ to a $\Gamma$-covering of $(S^2 \setminus D) \times S^1$ branched at the product of either $p-m$ or $p+1-m$ points with $S^1$.  
  
  We claim that this $\Gamma$-covering of $S^3$ may be obtained by taking a number of copies of $S^3$, removing balls, and gluing them together along their sphere boundaries. It suffices to break $S^3$ into pieces, with boundary $S^2$, such that the $\Gamma$-covering of each piece is a disjoint union of copies of $S^3$ with balls removed, as we can then glue the pieces together by identifying some of their boundary components, all of which are spheres.
  
  To do this, we observe that the product of either $p-m$ or $p+1-m$ points in $S^2 \setminus D$  with $S^1$ inside $(S^2 \setminus D) \times S^1 \subset S^3$ is a union of parallel circles inside the unknotted torus $(S^2 \setminus D) \times S^1$ and thus a union of unknotted, unlinked circles. We may place each of these circles into a disjoint ball.
  
  Let $C$ be the complement of all these balls. Certainly $C$ may be obtained from $S^3$ by removing balls. Hence $C$ is simply connected. The induced $\Gamma$-cover of $C$ is unbranched. Hence the induced $\Gamma$-cover of $C$ is a union of copies of $C$, and thus a union of complements of balls in $S^3$.
  
  Restricted to a ball containing a single circle, the $\Gamma$-covering is branched only at that circle. The fundamental group of the complement of a single unknotted circle in a ball is $\mathbb Z$, so the $\Gamma$-covering must be a disjoint union of cyclic coverings, where the $n$-fold cyclic covering is homeomorphic to the complement of $n$ balls in $S^3$.
  
  Thus the $\Gamma$-covering of $S^3$ may be obtained by taking a number of copies of $S^3$, removing balls, and gluing them together along their sphere boundaries. Any such manifold has fundamental group a free group by the groupoid van Kampen theorem. Hence the natural map from the homology of this $\Gamma$-covering to $H_3(H,\Z)$ factors through the $H_3$ of the free group and thus vanishes, so the fundamental class of this $\Gamma$-covering of $S^3$ vanishes, as desired.  \end{proof}

 

 \begin{lemma}\label{braid-in-braid} Let $p$ and $m$ be nonnegative integers. Let $f$ be a surjective homomorphism from the fundamental group of an $p$-times punctured plane to $H\rtimes \Gamma$, with local monodromies at each puncture in $c_{H\rtimes \Gamma}$ and with sufficiently many (depending on $m$) punctures with each possible conjugacy class of local monodromy. 
 
 Let $f'$ be a homomorphism from the fundamental group of an $m$-times punctured disc $D$ to $H \rtimes \Gamma$, with local monodromies in $c_{H\rtimes \Gamma}$. Then there is a disc in the plane containing $m$ punctures  which is homeomorphic to $D$ by a homeomorphism sending the restriction of $f$ to $f'$. \end{lemma}
 
 \begin{proof} It suffices to find a homomorphism from the fundamental group of the $(p-m)$-punctured disc $D^*$ such that the homomorphism from the fundamental group of the $p$-punctured plane obtained by gluing together $D$ and $D^*$ lies in the same braid group orbit as $f$, as then we can choose a braid in the braid group sending this homeomorphism to $f$, take an associated homeomorphism from the punctured plane to itself, and choose the image of $D$ under this homeomorphism to be our disc in the plane.
 
 This follows from the statement \cite[Theorem 3.1]{Wood2021} that braid group orbits of covers of the $p$-punctured plane with sufficiently many punctures of each local monodromy class are classified by their lifting invariant, and that lifting invariants are additive when gluing discs: We must subtract the lifting invariant of $f'$ from the lifting invariant of $f$ and confirm a braid with this lifting invariant exists, which happens if the difference has coordinates sufficiently large, which occurs if the lifting invariant of $f$ has coordinates sufficiently large, where these coordinates are the number of punctures with local monodromy realizing each conjugacy class of $c$. \end{proof}
  
 \begin{proof}[Proof of Proposition \ref{topological-braid-surjective}]  Fix finitely many links $L_i$ in $S^3$ and homomorphisms $f_i \colon \pi_1( S^3 \setminus L_i) \to H \rtimes \Gamma$, whose local monodromy around each component of $L$ is conjugate to an element of $\Gamma$, such that $[S^3, L_i , f_i]$ generate $H_3(H)^\Gamma$ (using Lemma \ref{represented-by-link}). For each of these, by  Alexander's theorem, $L_i$ can be made by closing up a braid $B_i$ in the disc, in which case $f_i$ restricts to a homomorphism from the fundamental group of the punctured disc.
 
If we now take a configuration + surjective homomorphism $f$ with sufficiently many local monodromy elements in each conjugacy class, we can find by Lemma \ref{braid-in-braid} a disc $D_i$ in the plane, avoiding the $k$ marked strands, where $f$ restricts to the homomorphism $f_i$ from the fundamental group of a punctured disc to $H \rtimes \Gamma $. We can then extend $B_i$ to a braid $B_i'$ involving all the points in the plane by leaving all the points outside $D_i$ is unmoved. By Lemma \ref{link-vs-bfc}, the braid fundamental class of $B_i'$ is $[S^3, L_i , f_i]$.

Composing the $B_i'$s, we thus can generate an arbitrary element in $H_3(H)^\Gamma$.
  \end{proof}
  
  \subsection{Homology of covers of Hurwitz space}\label{SS:homology}
  
  Fix a connected component $X$ of $\operatorname{Hur}_{H \rtimes \Gamma, c_{H \rtimes \Gamma}, \mathbb C}^m$. We will compute the homology of the connected components of the cover $X^0(\mathbb C)$ of $X(\mathbb C)$ arising from the homomorphism $\pi_1(X(\mathbb C)) \to \pi_1(X) \to H_3(H, \mathbb Z/n)^\Gamma$, where the first map is the map from the topological fundamental group to the algebraic fundamental group, and the second map is the $e$ defined above.  We also call this composition $e:\pi_1(X(\mathbb C)) \to H_3(H, \mathbb Z/n)^\Gamma$. 
   To do this, we will show that this cover is itself covered by a connected component of another Hurwitz space. We first introduce the notation that allows us to define this Hurwitz space.
  
  For a group $G$ and a union $c$ of conjugacy classes of $G$, let $\operatorname{Hur}^{m,\infty}_{G,c}(\mathbb C)$ be the topological Hurwitz space allowing arbitrarily ramification at $\infty$: That is, $\operatorname{Hur}^{m,\infty}_{G,c}(\mathbb C)$ is the covering of the configuration space of $m$ points in the plane corresponding to the action of the braid group on the set of $m$-tuples of elements of conjugacy classes in $c$ that generate $G$. (We could define this space as a scheme following \cite{LandesmanLevy}, but don't need to as we only need the manifold.)
  
  
   Let $A=H_3(H \rtimes \Gamma, \mathbb Z/n)$. We have a short exact sequence of $H\rtimes \Gamma$-modules
   \[ 0 \to A \to \Ind^{H\rtimes \Gamma}_1 A \to \Ind^{H\rtimes \Gamma}_1 A/A\to 0,\] 
   where $A$ is the trivial module embedding diagonally in $\Ind^{H\rtimes \Gamma}_1$.
Since $H^3( H \rtimes \Gamma,  \Ind^{H\rtimes \Gamma}_1 A)=0$, we thus have a surjection
  \[ H^2(H \rtimes \Gamma,  \Ind^{H\rtimes \Gamma}_1 A/A ) \to H^3(H \rtimes \Gamma , A ) \] and the universal coefficient theorem induces a surjection \[  H^3(H \rtimes \Gamma , A )  \to \operatorname{Hom} ( H_3(H \rtimes \Gamma, \mathbb Z), A).\]
  
 Fix $\alpha \in H^2(H \rtimes \Gamma,  \Ind^{H\rtimes \Gamma}_1 A/A ) $ whose image in $\operatorname{Hom} ( H_3(H \rtimes \Gamma, \mathbb Z), A)$ under the composition of these two surjections is the natural map $H_3(H \rtimes \Gamma, \mathbb Z)\to H_3(H \rtimes \Gamma, \mathbb Z/n)=A$
 (induced by $\Z\ra\Z/n$ on coefficients).
 
\begin{lemma}\label{covering-lemma}
Let $m$ be a nonnegative integer. Fix a connected component $X$ of $\operatorname{Hur}_{H \rtimes \Gamma, c_{H \rtimes \Gamma}, \mathbb C}^m$. Let $G$ be the extension of $H \rtimes \Gamma$ by $\Ind^{H\rtimes \Gamma}_1 A/A$ defined by the class $\alpha$ as defined above.  For each subgroup $G'$ of this $G$, let $c_{G'}$ be the set of nontrivial elements of $G'$ whose order divides the order of $\Gamma$. 

There is a subgroup $G'$ of $G$, whose image inside $H \rtimes \Gamma$ is all of $H \rtimes \Gamma$, such that the natural map $\operatorname{Hur}^{m,\infty}_{G',c_G'}(\mathbb C) \to   \operatorname{Hur}_{H \rtimes \Gamma, c_{H \rtimes \Gamma}}^{m,\infty}(\mathbb C) $ has image containing the component $X(\mathbb C)$. 

Furthermore, we can choose a component $Z$ of $\operatorname{Hur}^{m,\infty}_{G',c_G'}(\mathbb C) $ whose image in $\operatorname{Hur}_{H \rtimes \Gamma, c_{H \rtimes \Gamma}}^m(\mathbb C)$  is $X(\mathbb C)$, such that for each component $Y(\mathbb C)$ of $X^0(\mathbb C)$, 
there exists a covering space map $Z\to Y(\mathbb C)$ such that the composition $Z\to Y(\mathbb C)\to X (\mathbb C)$ is
the given map $Z(\mathbb C) \to X (\mathbb C)$.
\end{lemma}

\begin{proof} The topological space $ \operatorname{Hur}_{H \rtimes \Gamma, c_{H \rtimes \Gamma}}^m(\mathbb C)$ is the covering of the configuration space of $m$ points in the plane corresponding to the action of the braid group on the set of $m$-tuples of elements in conjugacy classes in $c_{H \rtimes \Gamma}$ that multiply to $1$ and generate $H \rtimes \Gamma$ (\cite[Section 11.3]{Liu2024}). Each connected component corresponds to a braid group orbit on the set of $m$-tuples. Let $g_1,\dots, g_m$ be a tuple in the orbit corresponding to $X(\mathbb C)$. To prove the first claim, we must check that the orbit of $g_1,\dots,g_m$ is the image of some orbit of $c_{G'}$-tuples for some $G'$.

Each $g_i$ is conjugate to a nontrivial element of $\Gamma$ and so has order dividing $|\Gamma|$. Lift each element $g_i$ from $H \rtimes \Gamma$ to an element $g_i'$ of $G$. Since the kernel of $G \to H \rtimes \Gamma$ has order prime to $|\Gamma|$, we can choose the lift to have the same order. Let $G'$ be the subgroup generated by $g_1',\dots, g_m'$. Then $g_1',\dots,g_m'$ lie in $c_{G'}$ and therefore the orbit of the tuple $g_1',\dots, g_m'$ defines a component of $\operatorname{Hur}^{m,\infty}_{G',c_G'}(\mathbb C)$ whose image is $X(\mathbb C)$. Take $Z$ to be this component.

We now prove the second part of the claim. The topological space $X^0(\mathbb C)$ is the $H_3(H, \mathbb Z/n)^\Gamma$-torsor over $X^0(\mathbb C)$ corresponding to the map $e \colon \pi_1(X(\mathbb C)) \to H_3(H, \mathbb Z/n)^\Gamma$.
 Hence all the connected components of  $X^0(\mathbb C)$ are isomorphic to each other
  and the natural map $Z(\mathbb C) \to X(\mathbb C)$ lifts to one such component if the composition \[ \pi_1(Z) \to \pi_1(X(\mathbb C)) \to  H_3(H, \mathbb Z/n)^\Gamma\] vanishes.
  Thus it suffices to check this vanishing. Let $\sigma$ be an element of $\pi_1(Z)$, which we may view as an element of the braid group stabilizing a tuple $g_1',\dots, g_m'$ of elements of conjugacy classes in $c_{G'}$. We must check that $e(\sigma)=0$.

By Lemma \ref{b-vs-bfc}, it suffices to check that the image of the braid fundamental class of $\sigma$ under  $H_3(H, \mathbb Z) \to H_3(H, \mathbb Z/n)$  vanishes. By Lemma \ref{link-vs-bfc} with $p =m$, this braid fundamental class is given by $[S^3, L, f]$ where $L$ is the closure of the braid $\sigma$ and $f$ is a homomorphism $\pi_1 ( S^3 \setminus L) \to H \rtimes \Gamma$ whose restriction to an $m$-punctured disc in $S^3$ is the homomorphism corresponding to a point of $X({\mathbb C})$. The fundamental group of the complement of the braid closure in $S^3$ is the free group on generators $x_1,\dots,x_m$ modulo the relation $x_i =\sigma(x_i)$ for all $i$, and $f$ is the homomorphism sending $x_i$ to $g_i$.

 Thus the fact that $\sigma$ fixes the tuple $g_1',\dots,g_m'$ implies that $f$ lifts to a homomorphism $f'\colon \pi^1(S^3 \setminus L) \to G'$ sending each $x_i$ to $g_i'$.  Let $H'$ be the kernel of the natural map $G' \to H\rtimes \Gamma \to \Gamma$.  Let $\tilde{S}^3$ be the $\Gamma$-covering of $S^3$ corresponding to $f$. Then $f'$ corresponds to a $H'$-covering of $\tilde{S}^3$ lifting the $H$-covering corresponding to $f$. Since the image of a loop around each component of $L$ under $f'$ is an element of $g_i'$ and thus has order dividing $\Gamma$ and hence has order prime to $H'$, the $H'$-covering must be unbranched. Hence the class $[S^3, L, f]$, which by definition is the image of the fundamental class of $\tilde{S}^3$ in $H_3(H, \mathbb Z)$, lifts to $H_3(H', \mathbb Z)$.  Since this $H'$-covering admits an action of $\Gamma$,
 this lift is $\Gamma$-invariant. Thus the braid fundamental class
  lies in the image of $H_3(H',\mathbb Z)^\Gamma \to H_3(H, \mathbb Z)^\Gamma$.
 
The class $\alpha$ vanishes when pulled back to from $H\rtimes \Gamma$ to $G$, and hence to $H'$, since the extension $\alpha$ splits over $G$. Since pullback is compatible with the connecting homomorphism, the pullback from $H^3( H \rtimes \Gamma, A) $ to  $H^3( H', A) $ of the class corresponding to $\alpha$ vanishes. Since pullback is compatible with the universal coefficient theorem, the composition of the natural homomorphism $H_3(H \rtimes \Gamma, \mathbb Z)\to  H_3(H \rtimes \Gamma, \mathbb Z/n) = A$ with $H_3(H ' , \mathbb Z) \to H_3(H \rtimes \Gamma, \mathbb Z)$ is zero. Thus the map along the bottom row of the below commutative diagram vanishes.

\[ \begin{tikzcd}  H_3(H', \mathbb Z)^\Gamma \arrow[r]\arrow[d]& H_3(H, \mathbb Z)^\Gamma \arrow[d]\arrow[r] & H_3(H, \mathbb Z/n)^\Gamma \arrow[d] \\ H_3(H ' , \mathbb Z)\arrow[r] &H_3(H \rtimes \Gamma, \mathbb Z) \arrow[r] & H_3 (H \rtimes \Gamma, \mathbb Z/n) \end{tikzcd}\]
Since the rightmost arrow $H_3(H, \mathbb Z/n)^\Gamma  \to H_3 (H \rtimes \Gamma, \mathbb Z/n) $ is an isomorphism, the composition of the top row $H_3(H', \mathbb Z)^\Gamma\to H_3(H, \mathbb Z)^\Gamma\to H_3(H, \mathbb Z/n)^\Gamma$ vanishes.
We saw that the image of the braid fundamental class in $H_3(H ,\mathbb Z/n)^\Gamma$ lies in the image of this composition, and hence vanishes, as desired. \end{proof}

\begin{proposition}\label{cohomology-bound} There exist positive integers $I,J,K$ depending only on $H$ and $\Gamma$ such that for each positive integer $m$, each connected component $X(\mathbb C)$ of $\operatorname{Hur}_{H \rtimes \Gamma, c_{H \rtimes \Gamma}}^m(\mathbb C)$, each component $Y(\mathbb C)$ of $X^0(\mathbb C)$, and each nonnegative integer $i$ we have:

\begin{enumerate}

\item $ \dim H_i (Y(\mathbb C), \mathbb Q) \leq K^{i+1}$.
\item If, for each conjugacy class in $c$, $X(\mathbb C)$ parameterizes covers with at least $I i +J$ branch points of that conjugacy class, then the natural map $H_i(Y(\mathbb C), \mathbb Q) \to H_i (X(\mathbb C), \mathbb Q)$ is an isomorphism. \end{enumerate}

\end{proposition} 

\begin{proof} Choose a subgroup $G'$ of $G$ and a component  $Z$ of  $\operatorname{Hur}^{m,\infty}_{G',c_G'}(\mathbb C) $ satisfying the conclusion of  Lemma \ref{covering-lemma}.

Let us check that there are finitely many possibilities for $G'$ for a given $H, \Gamma$, even independently of $n$:  $A = H_3(H \rtimes \Gamma, \mathbb Z/n)$ has size bounded by $|H_3(H\rtimes \Gamma, \mathbb Z)| |H_2(H \rtimes \Gamma, \mathbb Z)|$ and thus has finitely many possible values, then there are finitely many extensions of $H \rtimes \Gamma$ by $A$, and finitely many subgroups $G'$ of each extension.

 Let $g_1,\dots,g_v$ be representatives of the conjugacy classes of $c_{H \rtimes \Gamma}$ and let $m_1,\dots,m_v$ be the number of branch points of each conjugacy class in covers parameterized by $X(\mathbb C)$. Then we have a sequence of maps
\[  Z\to Y(\mathbb C) \to X(\mathbb C)\to \operatorname{Conf}_{m_1,\dots,m_v}(\mathbb C) \] where $\operatorname{Conf}_{m_1,\dots,m_v} (\mathbb C)$ parameterizes configurations of points in the plane with $m_1$ points labeled by $g_1$, $m_2$ points labeled by $g_2$, and so on, and the map $ X (\mathbb C)\to \operatorname{Conf}_{m_1,\dots,m_v} (\mathbb C)$ labels each branch point by a representative of its conjugacy class.

We have the induced maps on homology
\[  H_i(Z,\mathbb Q)\to H_i(Y(\mathbb C),\mathbb Q) \to  H_i(X(\mathbb C),\mathbb Q) \to H_i(\operatorname{Conf}_{m_1,\dots,m_v}(\mathbb C),\mathbb Q). \]   
All our maps between spaces are covering maps since the individual spaces are all connected covering spaces of $\operatorname{Conf}_m$. Thus the induced maps on homology with rational coefficients are injective.

Part (1) then follows from  \cite[Lemma 8.4.2]{LandesmanLevy} by
\[ \dim H_i (Y(\mathbb C), \mathbb Q) \leq \dim H_i (Z, \mathbb Q) \leq K^{i+1}\] where we choose $K$ large enough to be valid for any of the finitely many possible values of $G'$.

For part (2), we first check that the natural map from the set $c_G'/G'$ of conjugacy classes in $c_{G'}$ to the set $c_{H \rtimes \Gamma}/ H \rtimes \Gamma$ of conjugacy classes in $c_{H \rtimes \Gamma}$ is a bijection. This follows because $G'$ is the extension of $\Gamma$ by a group $H'$ of order prime to $|\Gamma|$, hence is $H' \rtimes \Gamma$ by Schur-Zassenhaus, and thus every conjugacy class of $G'$ of order dividing $|\Gamma|$ is in fact the conjugacy class of an element of $\Gamma$. Hence the natural map $c_G'/G' \to c_\Gamma/\Gamma$ is a bijection, and the same is true for $c_{H \rtimes \Gamma}/ H\rtimes \Gamma \to c_\Gamma / \Gamma$, so the natural map $c_{G'}/G'\to c_{H \rtimes \Gamma}/H \rtimes \Gamma$ is a bijection as well.

It follows from \cite[Theorem 1.4.2]{LandesmanLevy} that as long as all the $m_j$ are at least $Ii+J$ for some $I,J$ depending only on $G'$ that the natural map $H_i(Z,\mathbb Q) \to H_i(\operatorname{Conf}_{m_1,\dots,m_v}(\mathbb C),\mathbb Q)$ is an isomorphism. Thus, because all the maps are injective, all the maps are isomorphisms, as desired. We can choose a single $I,J$ that works for all the finitely many possible values of $G'$.\end{proof}

\subsection{From characteristic $0$ to characteristic $p$}\label{SS:Hurcount}

We will use Proposition~\ref{topological-braid-surjective} and Proposition~\ref{cohomology-bound} to prove the lemmas which will be the final ingredients in the proof of Theorem~\ref{component-comparison}. To do this, we relate the components of Hurwitz space in characteristic $0$ and characteristic $p$. While this was done in prior work, we do it in a slightly different way that directly keeps track of the fundamental group of the components. Our key lemma, Lemma \ref{general-pi1-comparison}, proves several useful facts under the assumption that a map between the spaces induces a surjection on their fundamental groups. Because of this, we will need Lemma  \ref{conf-pi1-surjective}, which shows that natural maps between configurations spaces induce a surjection on fundamental groups, which we prove as a consequence of two general results on fundamental groups. It would likely be possible to give a shorter proof avoiding fundamental groups, but we do not do this as we think the argument via fundamental groups is natural, and the results we prove about fundamental groups could be useful later.


Using Proposition~\ref{topological-braid-surjective} and Lemma~\ref{b-vs-bfc}, we obtain the following.

\begin{lemma}\label{complex-pi1-image} 
Let $X$ be a connected component of $(\Hur^m_{H\rtimes \Gamma ,c_{H \rtimes \Gamma} } )_{\mathbb C}$. 
The image of $e\colon\pi_1(X)\ra H_3(H,\Z/n)^\Gamma$ is
equal to the image of $H_3(H,\Z)^\Gamma$ under the natural map  $H_3(H,\Z)^\Gamma\to H_3(H, \mathbb Z/n)^\Gamma$, 
as long as the component parametrizes covers with sufficiently many (given $H,\Gamma$) branch points of each conjugacy class in $c_{H \rtimes \Gamma}$. \end{lemma}

\begin{proof} The fundamental group $\pi_1(X)$ is the completion of the topological fundamental group of $X(\bbC)$. Composition gives a homomorphism $e: \pi_1(X(\bbC)) \ra
H_3(H, \mathbb Z/n)$ and it suffices to show that this composition has image the image of $H_3(H,\Z)^\Gamma$.

Using Lemma \ref{b-vs-bfc}, since the braid fundamental class map has image in $H_3(H,\Z)^\Gamma$, we can immediately conclude that the image of the geometric fundamental group is contained in the image of $H_3(H,\Z)^\Gamma$.  By Proposition \ref{topological-braid-surjective}, the image of the braid fundamental class map is equal to the image of $H_3(H,\Z)^\Gamma$ as long as the component parameterizes covers with sufficiently many branch points in each conjugacy class. It then follows that the image of the geometric fundamental group is equal to the image of $H_3(H,\Z)^\Gamma$: Indeed, for $\alpha$ in the image of $H_3(H,\Z)^\Gamma$, $\alpha$ is the braid fundamental class of some $\sigma$, hence by Lemma \ref{b-vs-bfc} an integer multiple of $e(\sigma)$, hence equal to $e(\sigma^d)$ for some $d$.\end{proof}

\begin{lemma}\label{general-pi1-comparison} Let $A \to B$ be a map of schemes such that $\pi_1(A) \to \pi_1(B)$ is surjective, and let $C$ be a finite \'{e}tale cover of $B$. Let $g$ be the map from the set of connected components of $A \times_B C$ to the set of connected components of $C$ that sends a component $X_A$ of $A \times_B C$ to the unique connected component of $C$ containing the image of $X_A$ in $C$.

Then \begin{enumerate}
\item $g$ is a bijection.

\item  For $X_A$ a connected component of $A\times_B C$ and $X_B =g(X_A)$ the corresponding component of $C$, the map $\pi_1(X_A) \to \pi_1(X_B)$ induced by the map $X_A \to X_B$ is surjective. in particular, any homomorphism from $\pi_1(X_B)$ to a fixed group $G$ has the same image as its composition with $\pi_1(X_A) \to \pi_1(X_B)$. 

\item For $X_A$ a connected component of $A \times_B C$ and $X_B = g(X_A)$ the corresponding component of $C$, we have $X_A = A \times_B X_B$.\end{enumerate}\end{lemma}

\begin{proof} The only fact about the \'{e}tale fundamental group we need for this is that the map of fundamental groups $\pi_1(A) \to \pi_1(B)$ is surjective if and only if the pullback of each connected finite \'{e}tale cover of $B$ to $A$ is connected \cite[V, Proposition 6.9]{Grothendieck2003}.

Indeed, $C$ is a finite union of connected finite \'{e}tale covers $C_1,\dots, C_k$ of $B$, those being the connected components of $C$, and $A \times_B C$ is the union of the pullbacks $A \times_B C_i$ of these components from $B$ to $A$. Since the pullback $A \times_B C_i$ are connected, these are the connected components of $A\times_B C$, and this gives a bijection between the connected components of $A \times_B C$ and the connected components of $C$, which by definition has the property (3). The image of a component of $A\times_B C_i$ in $C$ is certainly contained in the component $C_i$, so this bijection is the same as $g$. For $C_i$ a component of $C$, a connected finite \'{e}tale cover of $C_i$ is itself a connected finite \'{e}tale cover of $B$, and its pullback from $C_i$ to $A\times_B C_i$ agrees with its pullback from $B$ to $A$ and hence is connected, showing that $\pi_1(A \times_B C_i) \to \pi_1(C_i)$ is surjective. Finally, precomposition of $\pi_1(C_i)\to G$ with the surjective morphism $\pi_1(A \times_B C_i) \to \pi_1(C_i)$ certainly preserves the image in $G$. \end{proof} 

\begin{lemma}\label{surjectivity-nc} Let $R$ be a strict Henselian local ring of generic characteristic $0$. Let $s$ be a geometric special point of $R$ and $\eta$ a geometric generic point. Let $X$ be a scheme over $R$ which is the complement in a smooth proper scheme $\overline{X}$ over $R$ of a normal crossings divisors $D$. Then the natural maps $\pi_1(X_{\eta} ) \to \pi_1(X)$ and $\pi_1(X_s) \to \pi_1(X)$ are surjective.  \end{lemma}

\begin{proof}Indeed, note first that since $R$ has generic characteristic $0$, the generic point of each component of $D$ has characteristic $0$, so the local rings at these points have residue characteristic $0$, and thus all finite \'{e}tale extensions over the field of fractions of these local rings are tamely ramified. Hence by   \cite[XIII, Definition 2.1.1]{Grothendieck2003}, all finite \'{e}tale covers of $X$ are tamely ramified and thus by \cite[XIII, Definition 2.1.3]{Grothendieck2003} the natural map from the usual fundamental group of $X$ to the tame fundamental group $\pi^t_1(X)$ is an isomorphism. It is shown in \cite[XIII, p. 48, l. 3-4]{Grothendieck2003} that $ \pi_1^t ( X_s) \to \pi_1^t(X)$ is an isomorphism which combined with the surjectivity of $ \pi_1 (X_s) \to  \pi_1^t ( X_s) $ \cite[XIII, 2.1.3]{Grothendieck2003} gives that $\pi_1 (X_s ) \to \pi^t_1(X)= \pi_1(X)$ is surjective. 

For the morphism $\pi_1 (X_\eta) \to \pi_1(X)$, it suffices to show that for each finite group $G$, the natural map from isomorphism classes of $G$-torsors over $X$ to isomorphism classes of $G$-torsors over $X_\eta $ is injective. Since all $G$-torsors over $X$  are tame, it suffices to show the injectivity of the analogous map on tame torsors. This is the content of \cite[XIII, Corollary 2.8]{Grothendieck2003} for $\mathcal F$ the constant sheaf $G$, taking $s_1=\eta$ and $s_2$ the geometric special point associated to the algebraic closure of the residue field, once we make several observations: the condition that $X \to \operatorname{Spec} R$ is locally acyclic is satisfied by \cite[XIII, Remark 1.17]{Grothendieck2003}. The fact that $(R^1_t f_* \mathcal F)_{\eta}$ agrees with the set of tame $G$-torsors over $ X_{\eta}$ follows from the stated compatibility of that pushforward with arbitrary change of base applied to the change of base $\eta \to \operatorname{Spec} R$. Since $R$ is strict Henselian so $\operatorname{Spec} R$ is the strict localization of $\operatorname{Spec} R$ at $s_2$, by \cite[XIII, (2.1.2.1)]{Grothendieck2003}  $(R^1_t f_* \mathcal F)_{s_2}$ agrees with the set of tame $G$-torsors over $\operatorname{Spec} R$.\end{proof}

\begin{lemma}\label{surjectivity-descent} Let $R$ be a ring and $z$ a point of $\operatorname{Spec} R$. Let $X$ be a scheme over $R$, let $G$ be a finite group acting freely on $X$, and let $Y=X/G$ be the quotient. If the natural map $\pi_1(X_z) \to \pi_1(X)$ is surjective, then the natural map $\pi_1(Y_z)\to\pi_1(Y)$ is surjective. \end{lemma}

\begin{proof} Note that $X\to Y$ is a finite \'{e}tale cover with Galois group $G$.  Thus for $C$ a connected finite \'{e}tale cover of $Y$,  $C\times_Y X$ is a finite \'{e}tale cover of $Y$ with an action of $G$, where $G$ acts transitively on the connected components. Since the inverse image of each connected component of $C\times_{ Y}X$ in the pullback 
\[C\times_{Y} X\times_{ X} X_z  =C \times_Y X_z =  C\times_{ Y } Y_z \times_{ Y_z }  X_z \]
is connected by \cite[V, Proposition 6.9]{Grothendieck2003}, $G$ also acts transitively on the connected components of $C\times_{ Y } Y_z \times_{ Y_z }  X_z $, and thus $C\times_{ Y } Y_z $ is connected, showing that $\pi_1 (Y_z) \to \pi_1(Y)$ is surjective by \cite[V, Proposition 6.9]{Grothendieck2003}, as desired.\end{proof}

\begin{lemma}\label{conf-pi1-surjective} Let $R$ be the ring of Witt vectors of $\bar{\F}_q$, and fix an embedding $R\to \mathbb C$. The maps $\pi_1 ( \Conf^m_{\bar{\F}_q}) \to \pi_1(\Conf^m_{R})$ and  $\pi_1 ( \Conf^m_{\bbC}) \to \pi_1(\Conf^m_{R})$ are each surjective. \end{lemma}

\begin{proof} $R$ is a strict Henselian local ring. The fixed embedding $R  \to \mathbb C$ defines a geometric generic point of $R$, and $R \to \bar{\F}_q$ gives a special point. Since $\operatorname{PConf}^m_R$ is the complement of a normal crossings divisor in the smooth projective scheme $\overline{\mathcal M}_{0,m}(\mathbb P^1,1)_R$, we may apply Lemma \ref{surjectivity-nc} to get that $\pi_1 ( \operatorname{PConf}^m_{\bar{\F}_q}) \to \pi_1(\operatorname{PConf}^m_{R})$ and  $\pi_1 ( \operatorname{PConf}^m_{\bbC}) \to \pi_1(\operatorname{PConf}^m_{R})$ are each surjective. Because $\Conf^m$ is the quotient of $\operatorname{PConf}^m$ by a free action of $S_m$, we may apply Lemma \ref{surjectivity-descent} to both points to obtain the desired conclusion.\end{proof}

\begin{lemma}\label{homology-comparison} Let $U$ be a finite \'{e}tale cover of $\operatorname{Conf}^m_R$. Then \[ \dim H^{2m-i}_c(U_{\bar{\F}_q}, \mathbb Q_\ell)= \dim H_i ( U_{\mathbb C}(\mathbb C), \mathbb Q) .\]\end{lemma}

\begin{proof} This follows from \cite[proof of Lemma 10.3 on p. 55]{Liu2024}: This proof is stated for a specific finite \'{e}tale covering of $\Conf_m$, but works in general. \end{proof}

\begin{lemma}\label{geometric-pi1-image} 
Let $q\equiv 1 \pmod{n}$ be a prime power relatively prime to $|H||\Gamma|$.
Let $X_{\bar{\F}_q}$ be a connected component of $(\Hur^m_{H\rtimes \Gamma ,c_{H \rtimes \Gamma} } )_{\bar{\F}_q}$. 
The image of $e\colon \pi_1(X_{\bar{\F}_q})\ra H_3(H,\Z/n)^\Gamma$ is
equal to the image of $H_3(H,\Z)^\Gamma$ under the natural map  $H_3(H,\Z)^\Gamma\to H_3(H, \mathbb Z/n)^\Gamma$, 
as long as $X_{\bar{\F}_q}$ parametrizes covers with sufficiently many (given $H,\Gamma$) branch points of each conjugacy class in $c_{H \rtimes \Gamma}$. \end{lemma}

\begin{proof} In view of Lemma \ref{complex-pi1-image} it suffices to show for each connected component $X_{\bar{\F}_q}$ of  $(\Hur^m_{H\rtimes \Gamma ,c_{H \rtimes \Gamma} } )_{\bar{\F}_q}$ that there exists a connected component $X_\mathbb C$ of  $(\Hur^m_{H\rtimes \Gamma ,c_{H \rtimes \Gamma} } )_{\bbC}$, with the same number of branch points of each conjugacy class, such that the image of $\pi_1(X_\mathbb C)$ in $H_3(H,\Z)^\Gamma$ is equal to the image of $\pi_1(X_{\bar{\F}_q})$, as we may then apply Lemma \ref{complex-pi1-image} to that component.

The space $(\Hur^m_{H\rtimes \Gamma ,c_{H \rtimes \Gamma} } )_{\bar{\F}_q}$ is the fiber product of $(\Hur^m_{H\rtimes \Gamma ,c_{H \rtimes \Gamma} } )_{R}$ over  $\Conf^m_{R}$ with   $\Conf^m_{\bar{\F}_q}$. Similarly, $(\Hur^m_{H\rtimes \Gamma ,c_{H \rtimes \Gamma} } )_{\bbC}$ is the fiber product of $(\Hur^m_{H\rtimes \Gamma ,c_{H \rtimes \Gamma} } )_{R}$ over  $\Conf^m_{R}$ with   $\Conf^m_{\bbC}$ 

By Lemma \ref{conf-pi1-surjective} the maps $\pi_1 ( \Conf^m_{\bar{\F}_q}) \to \pi_1(\Conf^m_{R})$ and  $\pi_1 ( \Conf^m_{\bbC}) \to \pi_1(\Conf^m_{R})$ are each surjective, as then applying Lemma \ref{general-pi1-comparison} twice gives us two bijections, first between components of $(\Hur^m_{H\rtimes \Gamma ,c_{H \rtimes \Gamma} } )_{\bar{\F}_q}$ and components of $(\Hur^m_{H\rtimes \Gamma ,c_{H \rtimes \Gamma} } )_{R}$, and then between components of $(\Hur^m_{H\rtimes \Gamma ,c_{H \rtimes \Gamma} } )_{R}$ and components of $(\Hur^m_{H\rtimes \Gamma ,c_{H \rtimes \Gamma} } )_{\mathbb C}$, which each preserve the image of the homomorphism $e$ to $H_3(H,\Z)^\Gamma$. Since the number of branch points of each conjugacy class is constant on components of $(\Hur^m_{H\rtimes \Gamma ,c_{H \rtimes \Gamma} } )_{R}$, both bijections preserve this quantity.\end{proof}

\begin{lemma}\label{fq-cohomology-bound}
There exist positive integers $I,J,K$ depending only on $H$ and $\Gamma$ such that for each $m$, each geometrically connected component $X_q$ of $\operatorname{Hur}_{H \rtimes \Gamma, c_{H \rtimes \Gamma}, \F_q}^m$, each connected component $Y_{\bar{\F}_q}$ of $(X_{q}^{s_H, \mathbb F_q})_{\bar{\mathbb F}_q}$, and each nonnegative integer $i$ we have:

\begin{enumerate}

\item $ \dim H^{2m-i}_c (Y_{\bar{\F}_q}, \mathbb Q_\ell) \leq K^{i+1}$.
\item If, for each conjugacy class in $c$, $X_{q}$ parameterizes covers with at least $I i +J$ branch points of that conjugacy class, then the natural map $ H^{2m-i}_c((X_q)_{\bar{\F}_q},\mathbb Q_\ell) \to  H^{2m-i}_c(Y_{\bar{\F}_q}, \mathbb Q_\ell) $ is an isomorphism. \end{enumerate}

\end{lemma}

\begin{proof}The proof is similar to the proof of Lemma \ref{geometric-pi1-image}.

Given a component $X_q$,  the homomorphisms defining the coverings $ X^{s_H, \mathbb F_q} _{q}$, for varying $s_H$, all 
 agree on restriction to the geometric fundamental group of $X_q$, and thus the coverings agree after base change from $\F_q$ to ${\bar{\mathbb F}_q}$. In other words, we have \[ (X^{s_H, \mathbb F_q} _q)_{\bar{\F}_q}\cong (X^{0, \mathbb F_q}_q)_{\bar{\F}_q} = (X^0_q)_{\bar{\F}_q}.\]  So it suffices to prove the same result for components $Y$ of $ (X^0_q)_{\bar{\F}_q}$, which are all  connected components of $\operatorname{Hur}_{H \rtimes \Gamma, c_{H \rtimes \Gamma}, \bar{\mathbb F}_q}^{m,0}$.

Lemmas \ref{general-pi1-comparison} and \ref{conf-pi1-surjective} give a bijection between connected components $Y_{\bar{\F}_q}$ of $\operatorname{Hur}_{H \rtimes \Gamma, c_{H \rtimes \Gamma}, \bar{\mathbb F}_q}^{m,0}$ and  connected components $Y_\mathbb C$ of $\operatorname{Hur}_{H \rtimes \Gamma, c_{H \rtimes \Gamma}, \mathbb C}^{m,0}$. By Lemma \ref{general-pi1-comparison}(3), components $Y_{\bar{\F}_q}$  and $Y_{\mathbb C}$ related by this bijection are respectively the base changes to $\bar{\F}_q$ and $\mathbb C$ of a single component $Y_{R}$ of $\operatorname{Hur}_{H \rtimes \Gamma, c_{H \rtimes \Gamma},R}^{m,0}$, which is a finite \'{e}tale cover of $\operatorname{Conf}^m_R$.

By Lemma \ref{homology-comparison} it follows that $\dim H^{2m-i}_c( Y_{\bar{\F}_q},\mathbb Q_\ell) = \dim H^i ( Y_{\mathbb C} (\mathbb C),\mathbb Q)$.

We apply Proposition \ref{cohomology-bound} to control the homology groups of $Y_{\mathbb C} (\mathbb C)$. Part (1), bounding the dimension, immediately translates to a bound on $\dim H^{2m-i}_c( Y_{\bar{\F}_q},\mathbb Q_\ell) $. For part (2), we have to check that the number of branch points of each conjugacy class is the same, but this follows because  the number of branch points of each conjugacy class is constant on components of $(\Hur^m_{H\rtimes \Gamma ,c_{H \rtimes \Gamma} } )_{R}$, and then we only get from Proposition \ref{cohomology-bound} and Lemma \ref{homology-comparison} that 
\[ \dim  H^{2m-i}_c(Y_{\bar{\F}_q}, \mathbb Q_\ell) =\dim H^{2m-i}_c((X_q)_{\bar{\F}_q},\mathbb Q_\ell).\]

However, because $\pi \colon Y \to X_q$ is a finite \'{e}tale cover between geometrically connected varieties, the induced map on compactly supported cohomology $ \pi^\circ \colon H^{2m-i}_c(X_q,\mathbb Q_\ell) \to  H^{2m-i}_c(Y_{\bar{\F}_q},, \mathbb Q_\ell) $ is an injection. (Since $\pi$ is finite, we have $\pi_! =\pi_*$ and  $R^j \pi_! \mathbb Q_\ell =0$ for $j>0$ so the Leray spectral sequence with compact supports gives $H^{2m-i}_c(Y_{\bar{\F}_q}, \mathbb Q_\ell) = H^{2m-i}_c((X_q)_{\bar{\F}_q}, \pi_* \mathbb Q_\ell ) $. In this perspective, the natural map $\pi^\circ$ arises from the adjunction map $\mathbb Q_\ell \to \pi_* \mathbb Q_\ell$ which has a left inverse given by the trace map $\pi_* \mathbb Q_\ell \to \mathbb Q_\ell$ (after division by the degree of $\pi$) by \cite[XVIII, Theorem 2.9 (Var 4)(I)]{sga4-3}, inducing a left inverse to $\pi^\circ$, showing that $\pi^\circ$ is injective.) Thus because the source and target of $\pi^\circ$ have the same dimension, $\pi^\circ$ is an isomorphism, completing the proof of part (2).\end{proof}

We can now prove Theorem \ref{component-comparison}.

\begin{proof}[Proof of Theorem \ref{component-comparison}]

Fix a component $X_q$ of $\Hur_{H \rtimes \Gamma, c_{H \rtimes \Gamma}, \mathbb F_q}$. If $X_q$ has no $\mathbb F_q$-points then the same is true for $X_q^{ s_H, \mathbb F_q}$, and the theorem is automatically true. Thus we assume that $X_q$ has a $\mathbb F_q$-point. Since $X_q$ is smooth, it further follows that $X_q$ is geometrically irreducible.

If $a(X_q)$ is lower than any fixed function of $\Gamma,H$, we can make the statement true by taking the implicit constant sufficiently large, and using the trivial bound $ | X_q^{s_H, \mathbb F_q} (\mathbb F_q)| = O(|X_q(\mathbb F_q)|)$, so we may assume $a(X_q)$ is larger than any desired function of $\Gamma,H$.

By definition $X^{ s_H, \mathbb F_q}_q$ is the finite \'{e}tale Galois cover 
 of $X_q$ associated to the homomorphism that sends an element $\sigma \in \pi_1(X_q)$ to $ e(\sigma) - s_H \deg \sigma$. 
 Restricted to the geometric fundamental group, this homomorphism is $e$ and hence by Lemma \ref{geometric-pi1-image}, since 
we can assume 
 $a(X_q)$ is sufficient large in terms of $H$ and $\Gamma$, has image the image 
 of $H_3(H,\mathbb Z)^\Gamma$ inside $H_3(H, \mathbb Z/n)^\Gamma$. The components of $(X_q^{ s_H, \mathbb F_q})_{\overline{\mathbb F_q}}$ 
 are naturally a principal homogenous space for
the quotient of $H_3(H,\mathbb Z/n)^\Gamma$ by this image, which by the long exact sequence is $H_2(H, \mathbb Z)^\Gamma [n]$. The long exact sequence associated to $0 \to \mathbb Z \to \mathbb Z \to \mathbb Z/n$ expresses $H_2( H, \mathbb Z/n)$ as the extension of
$H_1(H, \mathbb Z) [n]$ by
 $H_2(H,\mathbb Z)/n H_2(H,\mathbb Z)$, and $(|H|,|\Gamma|)=1$ implies that the same is true if we take $\Gamma$-invariants of everything.  Further,
 $H_\Gamma=1$ implies that $H_1(H, \mathbb Z)^\Gamma=0$.  Thus 
$$\abs{ H_2(H, \mathbb Z)^\Gamma [n]}= \abs{ H_2(H, \mathbb Z)^\Gamma /n H_2(H, \mathbb Z)^\Gamma }=
\abs{H_2( H, \mathbb Z/n)^\Gamma}=\abs{H^2( H\rtimes \Gamma, \mathbb Z/n)},
$$
 using Lemmas~\ref{L:GammaH} and \ref{L:homdual} in the last step. Thus the number of components over $\overline{\mathbb F_q}$ of $X^{ s_H, \mathbb F_q}_q$ is $\abs{H^2( H\rtimes \Gamma, \mathbb Z/n)}$.

Next we check that the action of $\Frob_q$ on the set of components is trivial. The action of $\Frob_q$ on the set of components is the same as the action of the Frobenius element $\sigma$  at an arbitrary $\mathbb F_q$-point of $X_q$.

It follows from Lemmas \ref{yuan-liu},  \ref{liu-domain-comparison}, and \ref{homomorphism-to-class} that the projection of $e(\sigma)$ to $H_2(H, \mathbb Z)^\Gamma [n]$ is the $\omega$-invariant of the component, which we have assumed is $\omega_H$. The degree of a Frobenius element is 1 so the projection of $s_H \deg \sigma$ to $H_2(H, \mathbb Z)^\Gamma [n]$ is the projection of $s_H$ to $H_2(H, \mathbb Z)^\Gamma [n]$, which by definition is $\omega_H$. 

 These two terms cancel and thus the action of $\Frob_q$ on the set of components is trivial, so each of these components is defined over $\mathbb F_q$.
 
 Fix a component $Y$ of $X_q^{s_H, \mathbb F_q}$, necessarily geometrically connected. The Lefschetz fixed point formula gives
 
 \begin{equation}\label{x-counting-formula}| X_q(\mathbb F_q) | = \sum_{i=0}^{2m}(-1)^i  \operatorname{tr}(\Frob_q, H^{2m-i}_c( (X_q)_{\bar{\F}_q}, \mathbb Q_\ell)) \end{equation}
 and
 \begin{equation}\label{z-counting-formula} | Y(\mathbb F_q) | = \sum_{i=0}^{2m}(-1)^i  \operatorname{tr}(\Frob_q, H^{2m-i}_c( Y_{\bar{\F}_q}, \mathbb Q_\ell) )\end{equation}
 
 By Lemma \ref{fq-cohomology-bound}(2), for $i \leq  \frac{ a(X_q)-J}{I}$ the Frobenius-equivariant natural map $H^{2m-i}_c( (X_q)_{\bar{\F}_q}, \mathbb Q_\ell)\to H^{2m-i}_c( Y_{\bar{\F}_q}, \mathbb Q_\ell)$ is an isomorphism and thus the terms with $i \leq  \frac{ a(X_q)-J}{I}$ in the above two formulas agree. Hence
 
 \[ || Y(\mathbb F_q) | - | X_q(\mathbb F_q) | |\leq \sum_{i > \frac{ a(X_q)-J}{I}} ( | \operatorname{tr}(\Frob_q, H^{2m-i}_c( (X_q)_{\bar{\F}_q}, \mathbb Q_\ell)) | + |\operatorname{tr}(\Frob_q, H^{2m-i}_c( Z_{\bar{\F}_q}, \mathbb Q_\ell))|)\]
 \[ \leq \sum_{i > \frac{ a(X_q)-J}{I}} q^{m-\frac{i}{2}} ( | \dim H^{2m-i}_c( (X_q)_{\bar{\F}_q}, \mathbb Q_\ell) | + |\dim ( H^{2m-i}_c( Z_{\bar{\F}_q}, \mathbb Q_\ell))|)\]
 \[ \leq \sum_{i > \frac{ a(X_q)-J}{I}} q^{m-\frac{i}{2}} 2 K^{i+1}  \leq 2 K^{\frac{ a(X_q)-J}{I}} q^{m - \frac{ a(X_q)-J}{2I}} \sum_{i=0}^\infty  K^{i+1} q^{ - \frac{i}{2} } \] \[= 2K^{\frac{ a(X_q)-J}{I}} q^{m - \frac{ a(X_q)-J}{2I}}  \frac{K}{ 1- Kq^{-1/2} } = O (e^{ - \delta a(X_q)  } ) q^m  .  \]
 by Deligne's Riemann hypothesis and Lemma \ref{fq-cohomology-bound}(1) (whose $n=1$ case handles $X_q$), taking $\delta \leq \frac{1}{I}\left( \frac{\log q}{2} - \log K  \right)$, which we can always do as long as $q$ is sufficiently large in terms of $H,\Gamma$. 
 
 Since this is true for every component and there are  $\abs{H^2( H\rtimes \Gamma, \mathbb Z/n)}$ components, we get
\[ |X_q^{s_H, \mathbb F_q}(\mathbb F_q)| -  \abs{H^2( H\rtimes \Gamma, \mathbb Z/n)} |X_q(\mathbb F_q)| = O ( \abs{H^2( H\rtimes \Gamma, \mathbb Z/n)} e^{ - \delta a(X_q)}    q^m) = O(  e^{ - \delta a(X_q)} |X_q(\mathbb F_q)| ),\]
since $|X_q(\mathbb F_q)|  \geq q^{m}/2$ for 
$q$ sufficiently large in terms of $\Gamma,H$.
 by \cite[Lemma 8.4.4]{LandesmanLevy} and we can absorb the factors of $2$ and $\abs{H^2( H\rtimes \Gamma, \mathbb Z/n)} $ into the big $O$. This is the desired statement. \end{proof}

  \subsection{Proofs of lemmas about integration along the fibers}\label{ss-long-proofs}
 
 We first prove the compatibility of integration along fibers with pullbback, and then the projection formula.
  
     \begin{proof}[Proof of Lemma \ref{integration-vs-pullback}] In this setting \cite[XVIII, Theorem 2.9 (Var 2)]{sga4-3} gives the commutative diagram, where $F' \colon Y_1= Y_2 \times_{X_2} X_1 \to X_1$ is the second projection and $g' \colon Y_1 = Y_2 \times_{X_2} X_1 \to Y_2$ is the first projection,
   
  \[ \begin{tikzcd} R^2 F'_* \mu_n \arrow[d, "\operatorname{Tr}_{F'}"]  &  R^2 F'_*  g^{'*} \mu_n \arrow[l]   & g^* R^2 F_* \mu_n \arrow[d, "g^* \operatorname{Tr}_F"] \arrow[l] \\   \mathbb Z/n & &g^* \mathbb Z/n \arrow[ll] \end{tikzcd}\]
  where the top-right arrow is a base change map, the bottom arrow is the isomorphism $g^* \mathbb Z/n \cong \mathbb Z/n$, and the top-right arrow arises from the isomorphism $ g^{'*} \mu_n  \cong \mu_n$. Applying the functor $H^i (X_1,\cdot) $ to this diagram, we obtain the commutative diagram
    \[ \begin{tikzcd}  H^i ( X_1, R^2 F'_* \mu_n \arrow[d, "H^i(\operatorname{Tr}_{F'})"] )  & H^i(X_1, R^2 F'_*  g^{'*} \mu_n ) \arrow[l] &  H^i(X_1, g^* R^2 F_* \mu_n \arrow[d, "H^i(g^* \operatorname{Tr}_F)"] ) \arrow[l] \\   H^i(X_2, \mathbb Z/n) && H^i(X_1, g^* \mathbb Z/n) \arrow[ll] \end{tikzcd}\]
    
    We can expand this diagram as follows, where ``pullback'' always denotes the natural map from the cohomology of a sheaf to the cohomology of the pullback of the sheaf.
        \[ \begin{tikzcd}   H^{i+2} (Y_1, \mu_n) \arrow[d, "\textrm{Leray}"]  & H^{i+2} (Y_1, g^{'*}  \mu_n) \arrow[d, "\textrm{Leray}"]  \arrow[l]  &&  H^{i+2}( Y_2, \mu_n) \arrow[d,"\textrm{Leray}"] \arrow[ll, "\textrm{pullback}"] \\
                H^{i+2}( X_1, R F'_* \mu_n \arrow[d] ) &   H^{i+2}(X_1, R F'_*  g^{'*} \mu_n ) \arrow[l]  \arrow[d] & H^{i+2}(X_1, g^* R F_* \mu_n \arrow[d] ) \arrow[l]  & H^{i+2} ( X_2, R F_* \mu_n) \arrow[d,] \arrow[l,"\textrm{pullback}" ]\\ 
        H^i ( X_1, R^2 F'_* \mu_n \arrow[d, "H^i(\operatorname{Tr}_{F'})"] ) &   H^i(X_1, R^2 F'_*  g^{'*} \mu_n ) \arrow[l] & H^i(X_1, g^* R^2 F_* \mu_n \arrow[d, "H^i(g^* \operatorname{Tr}_F)"] ) \arrow[l]  & H^i ( X_2, R^2 F_* \mu_n) \arrow[d,"H^i(\operatorname{Tr}_F)"] \arrow[l,"\textrm{pullback}" ]\\ 
          H^i(X_1, \mathbb Z/n) && H^i(X_1, g^* \mathbb Z/n) \arrow[ll] & H^i(X_2, \mathbb Z/n) \arrow[l,"\textrm{pullback}"] \end{tikzcd}\]
          It suffices to check that this diagram is commutative, since the left and right sides are by definition the integration along the fibers map and the top and bottom sides are the maps giving functoriality of cohomology in the space. Since we have checked the bottom-left rectangle is commutative, it remains to check the other squares and rectangles are commutative.
          
          In the bottom-right square, the horizontal arrows are the map on derived functors induced by a natural transformation on the original functors, that being the natural map from the global sections of a sheaf to the global sections of the pullback. The vertical arrows are the maps on derived functors introduced by a map of objects. Hence, this diagram expresses that a natural transformation on functors induces a natural transformation on the derived functors, which is \cite[Lemma 1.3.9]{Hovey99}.
          
          In the top-left square and center-left square, the horizontal arrows are maps induced by an isomorphism of sheaves, with the vertical arrows induced by applying the same construction to each of the two isomorphic sheaves. The cohomology, the derived pushforward, and the maps relating them may all be computed in terms of an injective resolution of the sheaf. Since the two sheaves are isomorphic, we may simply choose the same complex as the injective resolution of both, making the vertical maps identical and the horizontal maps identities so that the diagram commutes.
          
        In the top-right rectangle, all the arrows are natural transformations of derived functors, or compositions of derived functors, induced from natural transformations of the original functors, or compositions of functors. To check that the derived diagram commutes, it suffices to check that the original diagram commutes (since a commutative diagram of functors between abelian categories induces a commutative diagram of functors on the category of chain complexes and then we can apply \cite[Theorem 1.4.3 and Definition 1.4.2.6]{Hovey99}). The original diagram states that if we take a global section of a sheaf $\mathcal F$ on $Y_2$, pull back to a global section of $g^{'*}\mathcal F$, then view as a global section of $F^{'*} g^{'*} \mathcal F$, we obtain the same global section as if we view the original section as a section of $F_* \mathcal F$, pull back to a section of $g^*F_* \mathcal F$, and apply a base change map $g^*F_* \mathcal F \to F^{'*} g^{'*} \mathcal F$. This can be checked directly using the definition of the base change map.
          
      In the center square, all the arrows are maps on cohomology of $X_1$ with coefficients in a complex of sheaves induced by maps between complexes of sheaves, so it suffices to show the following diagram of maps between complexes of sheaves commutes. 
      \[    \begin{tikzcd} R F'_*  g^{'*} \mu_n \arrow[d] & g^* R F_* \mu_n \arrow[d] \arrow[l]  \\ 
        R^2 F'_*  g^{'*} \mu_n [-2] &X_1, g^* R^2 F_* \mu_n [-2] \arrow[l]  \end{tikzcd}\] 
Since the complexes have cohomology in degree $\leq 2$, we can equivalently express the vertical maps as the canonical map from a complex to its canonical truncation in degree $\geq 2$, so it suffices to check commutativity of the diagram
      \[    \begin{tikzcd} R F'_*  g^{'*} \mu_n \arrow[d] & g^* R F_* \mu_n \arrow[d] \arrow[l]  \\ 
        \tau^{\geq 2}  RF'_*  g^{'*} \mu_n  &  \tau^{\geq 2}  g^* R^2 F_* \mu_n  \arrow[l]  \end{tikzcd}\] 
Since the bottom horizontal arrow is the truncation of the top horizontal arrow, this can be checked directly from the definition of the canonical truncation of a complex. 

In the center-right square, the vertical arrows both arise from the map of complexes $g^* R F_* \mu_n  \to g^* R^2 F_* \mu_n [-2]$ while the horizontal arrows arise from the map from the cohomology of a complex on $X_2$ to the cohomology of its pullback to $X_1$. Thus the square expresses that the map from a cohomology of a complex on $X_2$ to the cohomology of its pullback to $X_1$ is a natural transformation. This again follows from \cite[Lemma 1.3.9]{Hovey99}. \end{proof}

      \begin{proof}[Proof of Lemma \ref{projection-formula}]
              The assumption that $j\leq 1$ is only used in one step of the proof, and can probably be removed. For this reason we have written the rest of the proof for arbitrary $j$, though we only need the $j=1$ case.

The relevant commutative diagram, where we have broken the cup product map into two steps for clarity, is as follows:
      
      \[ \begin{tikzcd}  H^{i+2} (Y,\mu_n) \otimes H^j(X, \mathbb Z/n)  \arrow[r, "\cup"] \arrow[dd,"\textrm{Leray}"] & H^{i+j+2} ( Y, \mu_n \otimes F^* \mathbb Z/n ) \arrow[r] \arrow[d,"\textrm{Leray}"] & H^{i+j+2} (Y, \mu_n)  \arrow[d,"\textrm{Leray}"] \\ 
& H^{i+j+2} (X, RF_* (\mu_n \otimes F^* \mathbb Z/n))\arrow[r] &  H^{i+j+2} (X, RF_* \mu_n) \arrow[dd] \\
       H^{i+2} (X, RF_* \mu_n)  \otimes H^j(X,\mathbb Z/n) \arrow[d]  \arrow[r,"\cup"]  & H^{i+j+2} (X, RF_* \mu_n \otimes \mathbb Z/n)  \arrow[ur] \arrow[d] 
  &  \\
       H^{i}(X, R^2 F_* \mu_n) \otimes H^j(X,\mathbb Z/n)  \arrow[d, "H^i(\operatorname{Tr}_F)"] \arrow[r,"\cup"] &  H^{i+j}(X, R^2 F_* \mu_n \otimes \mathbb Z/n) \arrow[r] \arrow[d, "H^{i+j}(\operatorname{Tr}_F\otimes id)"] & H^{i+j}(X,  R^2 F_* \mu_n) \arrow[d,"H^{i+j} (\operatorname{Tr}_F)"] \\
         H^i(X, \mathbb Z/n)  \otimes H^j(X,\mathbb Z/n) \arrow[r,"\cup"]  &  H^{i+j}(X, \mathbb Z/n \otimes \mathbb Z/n) \arrow[r] & H^{i+j}(X, \mathbb Z/n) 
       \end{tikzcd} \]
      All tensor products are over $\mathbb Z/n$. Here all the horizontal and diagonal arrows on the right side of the diagram are induced by isomorphisms of the corresponding complexes arising from the fact that the tensor product (over $\mathbb Z/n$) of a complex of $n$-torsion groups with $\mathbb Z/n$ recovers the original complex. 
       
            Let us see why the individual polygons commute. We begin with the bottom-right square and the quadrilateral above it. Both of these diagrams arise from commutative diagrams of complexes of sheaves on $X$, where the vertical arrows on the right arise from certain maps of complexes and the vertical arrows on the left arise from the same maps tensored with $\mathbb Z/n$. So these diagrams express the statement that isomorphism between the tensor product of a complex with $\mathbb Z/n$ and the complex itself is a natural transformation, which follows by \cite[Lemma 1.3.9]{Hovey99} from deriving the natural transformation of functors from the category of $n$-torsion sheaves to the category of sheaves given by $ A \otimes \mathbb Z/n \to A $.
            
            All the arrows in the top-right square are natural transformations between compositions of derived functors induced by natural transformations of the underlying functors (each functor is applied to the sheaf $\mu_n$). By \cite[Theorem 1.4.3 and Definition 1.4.2.6]{Hovey99}, to check the commutativity of the diagrams it suffices to check the commutativity of the corresponding diagrams for sheaves. This can be done directly from the definitions. 
            
       For the left-hand side, it is convenient to think of tensor product as a derived functor in the model category sense of the tensor product functor from the product of two copies of the category of complexes, with the product of the usual model structure, to the category of complexes. The cup product is then the natural transformations on derived functors induced by the natural transformation $\Gamma(\mathcal F) \otimes \Gamma(\mathcal G) \to \Gamma(F \otimes \mathcal G)$. Since replacing both complexes of sheaves with injective resolutions is a right Quillen replacement functor, this agrees with the definition of \cite{Swan1999}.
       
        That being done, all the arrows in the nonconvex hexagon are natural transformations of derived functors of maps on categories of complexes induced by natural transformations of functors on the original abelian categories. However, since the tensor product is a left derived functor and not a right derived functor, we cannot directly apply \cite[Theorem 1.4.3]{Hovey99}. It may be possible to do this with the theory of \cite{Shulman2011}, but we prove commutativity by hand. The key thing that makes this calculation easier is that given a complex of sheaves, for example any resolution of a sheaf, not necessarily injective, a cocycle in the global sections of the complex represents a class in the hypercohomology of the complex, or, in the example, in the cohomology of the sheaf, which may be expressed in terms of an injective resolution by mapping this resolution to the injective one. We start with a cohomology class $\alpha$ in $H^{i+2}(Y,\mu_n)$ and a cohomology class $\beta$ in $H^j(X,\mu_n)$. To represent them, we take an injective resolution $A^\cdot$ of $\mu_n$  on $Y$ and a resolution $B^\cdot$ of $\mathbb Z/n$ on $X$. Rather than choosing $B^\cdot$ to be an injective resolution, we choose it to be a finite resolution of flat $\mathbb Z/n$-modules such that $H^0(X, B^j)$ admits a cocycle representing $\beta$. The fact that it is a finite complex of flat $\mathbb Z/n$-modules means that tensor product with this complex computes the derived tensor product. We may do this for $j\leq 1$ since if $j=0$ we can take the complex to simply be $\mathbb Z/n$ and if $j=1$ then $\beta$ represents an extension $E$ of $\mathbb Z/n$ by $\mathbb Z/n$ as sheaves on $X$ and then $E \to \mathbb Z/n$ gives the desired complex. 
        
        Then $F^* B^\cdot$ gives a resolution of $\mathbb Z/n$ on $Y$ where the same cycle represents the pullback of $\beta$. The top horizontal arrow is obtained by taking the cocycle $\alpha \otimes \beta$ in the tensor product of resolutions $A^\cdot \otimes  F^* B^\cdot$, which is a resolution of $\mu_n\otimes F^* \mathbb Z/n$. We can map $\alpha \otimes \beta$ to a class in an injective resolution $C^\cdot$ of $\mu_n\otimes F^* \mathbb Z/n$.
         The top-right downwards Leray arrow then corresponds to taking the induced section in $H^0(X, F_* C^{i+j+2})$. Since $F_* C^{\cdot}$ is a complex representing $RF_* (\mu_n\otimes F^* \mathbb Z/n)$, this gives a cohomology class in $H^{i+j+2}(X, RF_* (\mu_n \otimes F^* \mathbb Z/n))$. Equivalently, this is induced by taking the section $\alpha \otimes \beta$ of $F_* ( A \otimes F^* B)$ and mapping along the map $F_* ( A \otimes F^* B)\to F_* C$ induced by $A \otimes F^* B \to C$. The rightward arrow on the right-hand side can be obtained by using the isomorphism $\mu_n \otimes F^* \mathbb Z/n \to \mu_n$ to view $C$ as a resolution of $\mu_n$, and interpreting the same cocycle as a class in the cohomology of $F_* \mu_n$. 
         
         For the left vertical arrow, we map $\alpha$ to an element $F_* \alpha$ of $H^0(X, F_* A^{i+2})$, where $F_* A^{i+2}$ is a complex representing $H^{i+2}(X, RF_* \mu_n)$. For the bottom horizontal arrow, we observe that $F_* \alpha \otimes \beta$ gives an element in $H^0(X,  (F_* A \otimes B)^{i+j+2})$ which induces a class in the cohomology of $R F_* \mu_n \otimes \mathbb Z/n$. For the bottom-right diagonal arrow, we observe that $F_* A \to (F_* A) \otimes B$ is a quasi-isomorphism so $F_* A\otimes B$ is also a resolution of $F_* A$, and hence $\alpha \otimes \beta$ induces a class in the hypercohomology of $F_*A $, i.e. in $H^* ( X, RF_* \mu_n)$.
         
         In summary, the upward transit involves taking the classes $\alpha$ and $\beta$, obtaining a cocycle $F_*(\alpha \otimes F^* \beta)$ in $F_*(A  \otimes F^* B)$ and obtaining a cocycle in $F_* C$, which since the composition of the two vertical arrows on the right side of the below diagram is a quasi-isomorphism (being the functor $F_*$ applied to a quasi-isomorphism between complexes of injectives) gives a class in the cohomology of $F_* A$. The downward transit involves obtaining a cocycle $F_* \alpha \otimes \beta$ in $F_* A \otimes B $,  which since the vertical arrow on the left side of the below diagram is a quasi-isomorphism gives a class in the cohomology of $F_A$.

         \[ \begin{tikzcd} 
         F_* A \arrow[r,"id"] \arrow [d] &  F_* A \arrow[d]  \\
         F_* A \otimes B \arrow[r]  & F_*(A  \otimes F^* B)  \arrow [d] \\  & F_* C \end{tikzcd}\]
         The horizontal arrow   $(F_* A \otimes B)  \to F_* (A \otimes F^* B)$ may be constructed for arbitrary complexes of sheaves $A, B$ as part of the projection formula. It is defined to send the tensor product of a section $s_1$ of $F_* A$ and a section $s_2$ of $B$ on the open set $U$, with $s_1$ arising from a section $s_1'$ of $A$ on $F^{-1}(U)$, to the section $s_1' \otimes s_2$ of $A \otimes F^* B$ on $F^{-1}(U)$. To check the two cocycles constructed from the two transits induce the same class in cohomology, we must check that they produce the same cocycle in $F_* C$, which reduces us to checking that the map on sheaves $F_* A \otimes B \to F_* ( A \otimes F^* B)$ sends $F_*\alpha \otimes  \beta$ to $F_*(\alpha \otimes F^* \beta)$. This is true by definition.
       
       For the center-left square, the vertical arrows arise from the same map of complexes $R F_* \mu_n \to R^2 F_* \mu_n[-2]$ (and the identity map $\mathbb Z/n\to\mathbb Z/n$) while the horizontal arrows are the derived functors of the natural transformation $\Gamma(\mathcal F) \otimes \Gamma(\mathcal G) \to \Gamma(F \otimes \mathcal G)$. The commutativity of this quadilateral follows from the fact that the derived functor of a natural transformation is a natural transformation \cite[Lemma 1.3.9]{Hovey99}.
       
       A similar argument work for the bottom-left square, as both vertical arrows arise from the map $( \operatorname{Tr}_F , \textrm{id}) \colon ( R^2f_* \mu_n, \mathbb Z/n) \to (\mathbb Z/n, \mathbb Z/n)$ in the product of two derived categories, and the horizontal arrow is the derived functor of the natural transformation  $\Gamma(\mathcal F) \otimes \Gamma(\mathcal G) \to \Gamma(F \otimes \mathcal G)$. Naturality of the derived functor gives commutativity of the diagram.\end{proof}
       
     Before proving Lemma \ref{homomorphism-to-class} relating the Artin-Verdier trace and the map $e$, we need the following calculational tool:

 \begin{lemma}\label{calculating-near-to-global} Let $\mathcal F$ be a sheaf in the \'{e}tale topology on a smooth projective curve $Y_x$ over a field.  Let $y$ be a point of $Y_x$ and let $j \colon Y_x\setminus y\to Y_x$ be the open immersion. Let $i \geq 0$ be an integer. The composition 
 \begin{equation}\label{near-point-to-global} H^i( k_y, \mathcal F) \to H^{i+1}_c ( Y_x\setminus y, j^* \mathcal F) \to H^{i+1}_c(Y_x, \mathcal F)\end{equation}
 of maps defined in \cite[II, Section 2]{Milne2006} may be computed on Cech cocycles in the following way:
 
 Let $\alpha$ be a cocycle in $C^{i+1} ( Y_x, j_! j^* \mathcal F)$. Let $\beta$ be a cochain in $C^i ( Y_x \setminus y, j^* \mathcal F)$ with $d \beta = j^* \alpha$. Let $\gamma$ be a cochain in $C^i( \mathcal O_{k_y}, j_! j^* \mathcal F)$ with $ d \gamma$ equal to the pullback of $\alpha$ to $C^{i+1}( \mathcal O_{k_y}, j_! j^* \mathcal F)$. Then we can pull both $\beta$ and $\gamma$ back to $C^i( k_y, j^* \mathcal F)$ where $ d \beta = d\gamma$ so $\beta-\gamma$ is a cocycle. 
 
 The composition \eqref{near-point-to-global} sends the class of $\beta-\gamma$ to the image of the class of $\alpha$ under the counit map $C^{i+1} ( Y_x, j_! j^* \mathcal F)\to C^{i+1} ( Y_x, \mathcal F)$. \end{lemma}
 
 \begin{proof} This can be proved by examining how the maps making up \eqref{near-point-to-global} are defined in \cite[II, Section 2]{Milne2006}. The map from compactly supported cohomology to ordinary cohomology is obtained by composing the isomorphism \cite[Proposition 2.3(d)]{Milne2006} between the compactly supported cohomology  group $H^{i+1}_c ( Y_x\setminus y,j^* \mathcal F)$ and the cohomology of the extension by zero $H^{i+1}(Y_x, j_! j^* \mathcal F)$ (n the reference both groups are compactly supported, but since $Y_x$ is proper, we may drop the $c$ there) and the map $H^{i+1}(Y_x, j_! j^* \mathcal F) \to H^{i+1}(Y_x, \mathcal F)$ appearing in the cohomology long exact sequence, which is the map induced  from the counit $j_! j^* \mathcal F \to \mathcal F$.  So it suffices to check that the composition 
 \[H^i( k_y, \mathcal F) \to H^{i+1}_c ( Y_x\setminus y, j^* \mathcal F) \to H^{i+1}(Y_x,  j_! j^* \mathcal F)\] sends the class of $\beta-\gamma$ to the class of $\alpha$.
 
 The compactly supported cohomology group $H^{i+1}_c ( Y_x\setminus y,j^* \mathcal F)$ is defined as the $i+1$st cohomology of the mapping cone of the map from $C^{*}(Y_x, j^*\mathcal F) $ to $C^{*}(k_y, j^*\mathcal F)$. To map this to the cohomology of the extension by zero, Milne uses \cite[Lemma 2.4(d)]{Milne2006}, which in our case gives a long exact sequence
 \[ \dots \to H^i(k_y, j^* \mathcal F) \to   H^{i+1}(Y_x,  j_! j^* \mathcal F) \to H^{i+1}(Y_x \setminus y,  j^* \mathcal F)   \to H^{i+1}(k_y, j^* \mathcal F) \to \dots \] 
 Since $H^{i+1}_c(Y_x,  j^* \mathcal F)$ is defined to be the $i+1$st cohomology of the mapping cone of  $
C^{*}(Y_x,  j_! j^* \mathcal F) \to C^{*}(Y_x \setminus y,  j^* \mathcal F) $, which is quasi-isomorphic to $C^{*} (Y_x,  j_! j^* \mathcal F)$ by the above exact sequence, this gives an isomorphism between $H^{i+1}_c ( Y_x\setminus y,j^* \mathcal F)$ and $H^{i+1}(Y_x,  j_! j^* \mathcal F)$. The isomorphism constructed this way sends the natural map from $H^i(k_y, j^* \mathcal F)$ to $H^{i+1}(Y_x,  j_! j^* \mathcal F)$ to the map $H^i(k_y, j^* \mathcal F) \to   H^{i+1}(Y_x,  j_! j^* \mathcal F) $ from the exact sequence above, since the natural map to the cohomology of the mapping cone of two complexes from the cohomology of one complex is one of the maps in the long exact sequence associated to the mapping cone. So it suffices to check that the map \[ H^i(k_y, j^* \mathcal F) \to   H^{i+1}(Y_x,  j_! j^* \mathcal F) \] of \cite[Lemma 2.4(d)]{Milne2006} sends the class of $\beta-\gamma$ to the class of $\alpha$.
 

The cohomology of the pair $H^*_y(Y_x,\cdot)$ is defined as the derived functor of the functor that takes a sheaf to its global sections on $Y_x$ that vanish on the complement of $y$. One can check that this derived functor is represented by the functor that takes a sheaf $\mathcal G$ to the cohomology of the mapping cone of $C^{*}( Y_x, \mathcal G) \to C^*(Y_x \setminus y, j^* \mathcal G)$ since this mapping cone has the correct value in degree $0$, its higher cohomology vanishes for injective sheaves, and it is compatible with short exact sequences of sheaves since $C^{*}( Y_x, \mathcal G) $ and $ C^*(Y_x \setminus y, j^* \mathcal G)$ both are.

The long exact sequence of the pair  \cite[III, Proposition 1.25]{Milne1980} \[   H^{i} (Y_x\setminus y, j^* \mathcal G) \to H^{i+1}_y (Y_x, \mathcal G)\to H^{i+1} (Y_x, \mathcal G) \to H^{i+1} (Y_x \setminus y, \mathcal G)\to \] includes maps $H^{i} (Y_x\setminus y, j^* \mathcal G) \to H^{i+1}_y (Y_x, \mathcal G)$ and $H^{i+1}_y (Y_x, \mathcal G)\to H^{i+1} (Y_x, \mathcal G) \to H^{i+1} (Y_x \setminus y, \mathcal G)$, defined as maps of $\operatorname{Ext}$ groups induced by a short exact sequence of sheaves. We can check that these agree with the natural maps between the cohomology of the mapping cone of two complexes and the cohomology of the original complexes by checking that the natural maps of the mapping cone are the correct maps on $H^0$ and compatible with connecting homomorphisms of short exact sequences of sheaves, both of which are straightforward.

In particular, the map $  H^{i+1}_y( Y_x, j_! j^* \mathcal F) \to H^{i+1}( Y_x, j_! j^* \mathcal F)$ from the long exact sequence of the pair, expressed in terms of the mapping cone, is equivalent to forgetting $ C^*(Y_x \setminus y, j^* \mathcal F)$, so the pair $(\alpha,\beta)$ represents a class in $H^{i+1}_{y} ( Y_x, j_! j^* \mathcal F) $, whose image in $H^{i+1}( Y_x, j_! j^* \mathcal F) $ is $\alpha$. So it suffices to check that the class represented by $(\alpha,\beta) $ in $H^{i+1}_{y} ( Y_x, j_! j^* \mathcal F) $ is sent by the isomorphism $H^{i+1}_{y} ( Y_x, j_! j^* \mathcal F)  \cong H^i( k_y,  j^* \mathcal F)$ to $\beta-\gamma$. 

The isomorphism $H^{i+1}_{y} ( Y_x, j_! j^* \mathcal F)  \cong H^i( k_y,  j^* \mathcal F)$ is obtained by composing an isomorphism $H^{i+1}_{y} ( Y_x, j_! j^* \mathcal F) \to H^{i+1}_y ( \mathcal O_{k_y}, j_! j^* \mathcal F)$ arising from excision \cite[III.1.28]{Milne1980} with an isomorphism $H^{i+1}_y ( \mathcal O_{k_y}, j_! j^* \mathcal F)\to H^{i}_y ( j_y ,  j^* \mathcal F) $ of \cite[II, Proposition 1.1(a)]{Milne2006}. The map proved to be an isomorphism in \cite[III, Corollary 1.28]{Milne1980} is simply the pullback map, as one can tell from \cite[III, proof of Corollary 1.28 and statement of Proposition 1.27]{Milne1980}. We can check that pullback commutes with the identification of relative cohomology with a mapping cone since pullback is a natural transformation of derived functors induced by the natural transformation of original functors, natural transformations between derived functors which are compatible with the connecting homomorphism of a short exact sequence are uniquely determined by their value on $R^0$, and the pullback map on the cohomology of the mapping cone is compatible with the connecting homomorphism of a short exact sequence and takes the correct value on $R^0$. Thus the pullback map $H^{i+1}_{y} ( Y_x, j_! j^* \mathcal F) \to H^{i+1}_y ( \mathcal O_{k_y}, j_! j^* \mathcal F)$ sends the class represented by $\alpha$ and $\beta$ to the class represented by the pullback of $\alpha$ and the pullback of $\beta$. The isomorphism $H^{i+1}_y ( \mathcal O_{k_y}, j_! j^* \mathcal F)\to H^{i}_y ( k_y ,  j^* \mathcal F) $ is defined to be the inverse of  the connecting morphism $H^{i}_y ( k_y ,  j^* \mathcal F) \to H^{i+1}_y ( \mathcal O_{k_y}, j_! j^* \mathcal F)$ of the long exact sequence of the pair
\[  H^{i}( \mathcal O_{k_y}, j_! j^* \mathcal F)\to H^{i}_y ( k_y ,  j^* \mathcal F) \to H^{i+1}_y ( \mathcal O_{k_y}, j_! j^* \mathcal F) \to H^{i+1}( \mathcal O_{k_y}, j_! j^* \mathcal F) \to \dots \]
which is an isomorphism since $H^{i}( \mathcal O_{k_y}, j_! j^* \mathcal F)=0$ for all $i$ by \cite[II, Proposition 1.1(b)]{Milne2006}. The connecting homomorphism sends the class represented by $\beta-\gamma$ to the class represented by $(0,\beta-\gamma)$ (as we checked above) which is equivalent to the class represented by $(\alpha, \beta)$ since $(\alpha,\beta)- (0,\beta-\gamma)= (\alpha, \gamma) = (d\gamma, \gamma) = d(\gamma,0)$ is a coboundary in the mapping cone. \end{proof}

\begin{proof}[Proof of Lemma \ref{homomorphism-to-class}] It suffices to check that for each $\alpha \in H_3( H, \mathbb Z/n)$, 
its image  $\tilde{e} (\alpha)\in  \Hom ( \pi_1(X) ,\mathbb Z/n) $ takes value on $\Frob_{q,x}$ equal to $\AV_{Z_x/Y_x}(\alpha)$, that is, the Artin-Verdier trace for $Z_x/Y_x$ applied to $\alpha$.  
Let $f$ be the map $Y_x\ra\Spec \F_q$.  
Because $\phi_{Z/Y}$, the integration map and the explicit calculation of \'{e}tale $H_1$  are compatible with pullback,
 we can compute $\tilde{e}(\alpha)(\Frob_{q,x})$ upon pullback to $Y_x$, i.e.  taking the image of $\alpha$ under the composite map
 \[ H^3 ( H, \mathbb Z/n) \stackrel{\phi_{Z_x/Y_x}}{\to} H^3 ( Y_x , \mathbb Z/n) \stackrel{\int_f}{\to}  H^1( \mathbb F_q , \mathbb Z/n) \cong \Hom ( \Gal( \mathbb F_q)  , \mathbb Z/n)
\stackrel{\operatorname{eval}_{\Frob_q}}{\longrightarrow} \Z/n .
  \] 
 
 We have two maps $H^3 ( H, \mathbb Z/n)\ra \Z/n$
 (the above composite and $ \AV_{Z_x/Y_x}(\alpha)$),
  both factoring through $H^3 ( H, \mathbb Z/n) \stackrel{\phi_{Z_x/Y_x}}{\to} H^3 ( Y_x , \mathbb Z/n)$,
so we can check they are the same just by checking that the corresponding maps  $H^3 ( Y_x , \mathbb Z/n)\ra \Z/n$ are the same.
Since $H^3(Y_x,\mathbb Z/n) \cong \mathbb Z/n$ (by \cite[Proposition 2.6]{Milne2006} and the fact that on $Y_x$ we have $\Z/n=\mu_n$),
 it suffices to check this on a generator. 
 
The Kummer map $H^1(Y_x, \mathbb G_m) \to H^2(Y_x, \mu_n)$ applied to the class of a degree $1$ line bundle $L$  in $H^2(Y_x,\mu_n)$, gives, using our identification $\mu_n\isom \Z/n\Z$, a class 
$[L]\in H^2(Y_x,\Z/n)$. 
  We consider the cup product of  $[L]$ with the pullback from $\mathbb F_q$ of the class $\beta$ of a $\mathbb Z/n$-torsor on $\mathbb F_q$ on which $\Frob_q$ acts by $+1$.

  We can check $\int_f [L] = \deg L = 1$ as follows. By compatibility of integration with pullback, because $H^0(\mathbb F_q, \mathbb Z/n) \to H^0 ( \bar{\mathbb F}_q, \mathbb Z/n)$ is an isomorphism, we can reduce to the case where the base is $\operatorname{Spec}  \bar{\mathbb F}_q$.  When the base is the spectrum of a algebraically closed field, the integration-along the fibers map reduces to the trace map from the second cohomology of the fiber to $\mathbb Z/n$. When the base is the spectrum of an algebraically closed field and the fiber is a connected curve, this map is given explicitly in  \cite[XVIII, 1.1.3.4]{sga4-3}. This explicit map is given by the inverse of the Kummer map $H^1(Y_x, \mathbb G_m) \to H^1(Y_x,\mu_n)$ composed with the identification between the Picard group of a curve modulo $n$ and $\mathbb Z/n$ arising in \cite[IX, Corollary 4.7]{sga4-3} from tensoring with $\mathbb Z/n$ the exact sequence \cite[IX, (4.8)]{sga4-3} $ 0 \to \operatorname{Pic}^0(Y_x) \to \operatorname{Pic}(Y_x) \to \mathbb Z \to 0$. Since the map $  \operatorname{Pic}(Y_x) \to \mathbb Z$ in that exact sequence is the degree map (this is not quite specified in \cite[IX, (4.8)]{sga4-3}, which only says to recall this exact sequence, but seems to be what is meant), it follows that the trace map sends $[L]$ to the degree of $L$, modulo $n$. The fact that this trace agrees with the more general notion of trace is  \cite[XVIII, Proposition 2.10]{sga4-3}, which states that in the case the that fiber is a curve, the trace map of \cite[XVIII, Theorem 2.9]{sga4-3} agrees with the trace map of \cite[XVIII, Proposition 1.1.6]{sga4-3}, which is itself defined to agree when the base is the spectrum of an algebraically closed field  with the trace map of \cite[XVIII 1.1.3.4]{sga4-3}.
    
  By the projection formula (Lemma \ref{projection-formula})
   \[ \int_f ( [L] \cup f^*\beta) =( \int_f [L ]) \cup \beta = 1 \cup \beta= \beta.\]
   In particular,  $\operatorname{eval}_{\Frob_q}(\int_f([L] \cup f^*\beta))=1$.
   Also, from this it follows that $[L] \cup f^*\beta$ is a generator of $H^3(Y_x,\mathbb Z/n)$.

Thus, we will have proven the lemma if we can show that the image of $[L]\cup f^*\beta$ under the following map used in defining the Artin-Verdier trace
\begin{equation}\label{E:maptocheck}
H^3 ( Y_x, \mathbb Z/n)\to H^3(Y_x ,\mu_n) \stackrel{D}{ \to} \mathbb Z/n
\end{equation}
 is 1.  
It suffices to check, for each point $y$, that the image of $[\mathcal O_{Y_x}(y) ] \cup f^* \beta $ under \eqref{E:maptocheck} is $\deg y$, as we can take a point $y_1$ of degree $d$ and $y_2$ of degree $d+1$ for $d$ sufficiently large, and then choose $L = \mathcal O_{Y_x}(y_2-y_1)$.

The definition of the map $D$ comes from the isomorphism $
H^3_c(Y_x,\G_m)\isom \Q/\Z
$ constructed in \cite[Chapter II, Section 3]{Milne2006}.  In particular it follows from this construction that, for a point $y \in Y_x$ with local field $k_y$, the composition of $D$ with $H^2( k_y, \mathbb \mu_n)\to H^3_c( Y_x \setminus y, \mathbb \mu_n) \to H^3 (Y_x, \mathbb \mu_n)$ agrees with the invariant map on the Brauer group of $k_y$. We have a diagram
\[\begin{tikzcd} H^2( k_y, \mu_n)\arrow[r] &  H^3_c( Y_x \setminus y, \mu_n) \arrow[r] &  H^3 (Y_x, \mu_n) \\ H^1( k_y , \mu_n)\arrow[r]  \arrow[u] &  H^2_c(Y_x \setminus y, \mu_n) \arrow[r] \arrow[u] &  H^2 (Y_x, \mu_n)\arrow[u] \\ H^0(k_y, \mathbb G_m) \arrow[r] \arrow[u] & H^1(Y_x \setminus y, \mathbb G_m) \arrow[u]\arrow[r] & H^1(Y_x, \mathbb G_m)\arrow[u] \end{tikzcd}\] where the vertical arrows in the top row are cup products with $f^*\beta$ and the vertical arrows in the bottom row are the Kummer map.

We first check that the horizontal arrows in the bottom row send the uniforrmizer in $H^0(k_y, \mathbb G_m) = k_y^\times$ to $\mathcal O_{Y_x}(y)\in H^1(Y_x, \mathbb G_m)$. However, the diagram does not quite commute. Instead, we check that the class obtained by sending the uniformizer along the bottom row and then the right column differs from the class obtained by sending the uniformizer along the left column and top row by a sign of $-1$. From this it follows that 
 the image of $[\mathcal O_{Y_x}(y)] \cup f^* \beta$ 
under the map in \eqref{E:maptocheck}
is the Brauer invariant of the cup product of the Kummer class of the \emph{inverse} of the uniformizer with $f^*\beta$, which is $\deg y$ by \cite[Lemma 6.1]{Lipnowski2020} since the Frobenius at a degree $\deg y$ point is the $\deg y$th power of $\operatorname{Frob}_q$. This was the desired statement, so when we check this, the proof is complete.

To calculate the horizontal arrows, we use Lemma \ref{calculating-near-to-global}. Let $\pi \in k(Y_x)$ vanish to order $1$ at $y$ and therefore function as a uniformizer. Let $U_1 = Y_x \setminus y$ and let $U_2$ be obtained from $Y_x$ by removing all zeroes and poles of $\pi$ other than $y$, so that $U_1 \cup U_2$ is a Zariski open cover of $Y_x$.

We start with the Cech cocycle in $a_1 \in C^1( Y_x,  j_!  \mathbb G_m)$ with underlying covering $U_1 \cup U_2$ defined by the section that sends $U_1 \times U_1$ and $U_2 \times U_2$ to $1$, $U_1\times U_2$ to $\pi$, and $U_2 \to U_1$ to $\pi^{-1}$. Here we express sections of $j_! j^* \mathbb G_m$ as sections of $\mathbb G_m$, i.e. invertible functions, that happen to be trivial, i.e. equal to $1$, on $U_1 \times U_1$, the only open set under consideration containing $y$.  The image of $a_1$ under $C^1( Y_x,  j_!  \mathbb G_m) \to C^1( Y_x,   \mathbb G_m) $ is expressed using the same cocycle. One immediately checks this is the cocycle arising from the line bundle $\mathcal O(y_x)$, with the trivialization on $U_1$ given by the natural map $\mathcal O \to \mathcal O(y)$ and the trivialization on $U_2$ given by the natural map $\mathcal O \to \mathcal O(y)$ divided by $\pi$.

We can find $b_0 \in C^0( Y_x \setminus y, \mathbb G_m)$ with  $d b_0 = j^* a_1$ with underlying covering the pullback of $U_1 \cup  U_2$, that sends the pullback of $U_1$ to $1$ and the pullback of $U_2$ to $\pi$. 

We can find $c_0$ in $C^0( \mathcal O_{k_y}, j_!\mathbb G_m)$ with $d c_0 $ the pullback of $a_1$ with underlying covering the pullback of $U_1 \cup U_2$, that sends the pullback of $U_1$ to $\pi^{-1}$ and the pullback of $U_2$ to $1$.

Then $b_0-c_0$ is the cocycle in $C^0({k_y}, \mathbb G_m)$ with underlying covering the pullback of $U_1 \cup U_2$, that sends both $U_1$ and $U_2$ to $\pi$. This represents the class $\pi$, so by Lemma \ref{calculating-near-to-global} the bottom row sends $\pi$ to $\mathcal O_{Y_x}(y)$.

We now apply the Kummer map to all this data. To do this, we obtain by \cite[III, Lemma 2.19]{Milne1980} \'{e}tale coverings $V_1$  and $V_2$ of $U_1$ and $U_2$ such that $\pi$ admits an $n$'th root $\pi^{1/n}$ on $V_1 \times V_2$.  We have a cochain in $ a_1 /n \in  C^1( Y_x,  j_!  \mathbb G_m)$ with underlying covering $V_1 \cup V_2$, defined by the function that sends $V_1 \times V_1$ and $V_2 \times V_2$ to $1$, $V_1 \times V_2$ to $\pi^{1/n}$ and $V_2 \times V_1$ to $\pi^{-1/n}$. The cochain $a_1/n$ is no longer a cocycle, but its  coboundary $d(a_1/n)$ takes values in $\mu_n$, and therefore defines a cocycle $a_2 \in C^2 ( Y_x, j_!\mu_n)$. The projection of $a_2$ to $C_2(Y_x, \mu_n)$ is obtained from $a_1$ by exactly the snake lemma construction that is used to define the connecting homomorphism of the Kummer exact sequence and therefore is the Kummer class $[\mathcal O_{Y_x}(y)]$ corresponding to $\mathcal O_{Y_x}(y)$.

We take $b_0/n \in C^0(Y_x\setminus y, \mathbb G_m)$, with underlying covering the refinement of the pullback of $V_1 \cup V_2$ given by the pullback of $V_1 \cup V_1 \times V_2$, defined by the function that sends $V_1$ to $1$ and $V_1 \times V_2$ to $\pi^{1/n}$. Then $(a_1/n) - d (b_0/n)  $ takes values in $\mu_n$, and therefore defines a cocycle $b_1 \in C^1 (Y_x\setminus y, \mu_n)$.  We have
\[d b_1= d ( (a_1/n) - d (b_0/n)) = d (a_1/n) - d^2(b_0/n) = a_2 - 0 = a_2 .\]

We take $c_0/n$ in $C^0( \mathcal O_{k_y}, j_!\mathbb G_m)$ with underlying covering the refinement of the pullback of $V_1 \cup V_2$ given by the pullback of $V_1 \times V_2 \cup V_2$, defined by the function that sends $V_1 \times V_2$ to $\pi^{1/n}$ and $V_2$ to $1$. Then $(a_1/n) - d ( c_0/n)$ takes values in $\mu_n$, and therefore defines a cocycle $c_1 \in C^1( \mathcal O_{k_y}, j_!\mathbb G_m)$. Identical reasoning gives $d c_1 =a_2$. 

Then $b_1- c_1 = d (c_0/n) - d(b_0/n)$ is obtained from $c_0-b_0$ by exactly the snake lemma construction that is used to define the connecting homomorphism of the Kummer exact sequence and is therefore the Kummer class corresponding to $\pi^{-1}$.

We now take the cup product of everything with $f^* \beta$. To do this, we can choose a Cech cocycle representing $f^* \beta$, arising from a finite \'{e}tale covering of $\operatorname{Spec} \F_q$, and use the formula \cite[07MB]{stacks} for the cup product of Cech cochains to cup $a_2$ with $f^* \beta$, obtaining $a_3$, and similarly cup $b_1$ and $c_1$ with $f^* \beta$, obtaining $b_2$ and $c_2$. We have $d b_2 = a_3 $ and $d c_2 = a_3$ since the differential is a derivation with respect to cup product \cite[01FP]{stacks}. Then $a_3$ represents the cup product of $a_2$ with $f^* \beta$ (by \cite[Corollary 3.10]{Swan1999}) and $b_2 - c_2$ represents the cup product of $b_1-c_1$ with $f^*\beta$, and thus Lemma \ref{calculating-near-to-global} shows that the cup product of the Kummer class of $\pi^{-1}$ with $f^* \beta$ is sent to the cup product of the Kummer class of $\mathcal O_{Y_x}(y)$ with $f^* \beta$, completing the proof. \end{proof}

\begin{proof}[Proof of Lemma \ref{b-vs-bfc}] Let $Z$, over $X$, be the universal family of $(H \rtimes \Gamma)$-coverings of $\mathbb P^1_\mathbb C$, so that $Z \to X \times \mathbb P^1$ is an $H \rtimes \Gamma$-covering. Let $Y$ be the universal family of induced $\Gamma$-coverings of $\mathbb P^1$, i.e. the quotient of $Z$ by $H$. 
  Then $Z \to Y $ is a finite \'{e}tale $H$-covering 
  and thus,  by Lemma~\ref{L:BHpullback}, determines a map from group cohomology to \'{e}tale cohomology $H^3(H, \mathbb Z/n) \to H^3(Y, \mathbb Z/n)$, and from \'{e}tale cohomology to singular cohomology $H^3(Y(\mathbb C), \mathbb Z/n)$. Since the integration along the fibers map is compatible with pullbacks, the map $e: \pi_1(X) \to H_3(H, \mathbb Z/n)$ is induced by the composition   \[ H^3(H, \mathbb Z/n) \to H^3 ( Y, \mathbb Z/n) \to H^1( X, \mathbb Z/n) \to \Hom (\pi_1(X(\mathbb C)), \mathbb Z/n) .\]
  
  Let us discuss how to define an ``integration along the fibers'' map in the singular cohomology setting. This can be defined as a composition
 \[ H^{i+2} (Y (\mathbb C), \mu_n) \to H^{i+2} ( X(\mathbb C), R F'_* \mu_n) \to H^i(X(\mathbb C), R^2 F'_* \mu_n) \to H^i (X(\mathbb C) , \mathbb Z/n) \]  where $F'$ is the map $Y(\mathbb C) \to X(\mathbb C)$ induced by the map $Y \to X$ once we have a trace map $R^2 F'_* \mu_n \to \mathbb Z/n$ of sheaves on the topological space $X(\mathbb C)$. We choose this trace map to be compatible with our original trace map, in the sense of the commutative diagram 
   \begin{equation}\label{pullback-to-analytic-diagram} \begin{tikzcd} R^2 F'_* \mu_n \arrow[d, "\operatorname{Tr}_{F'}"]  &  R^2 F'_*  g^{'*} \mu_n \arrow[l]   & g^* R^2 F_* \mu_n \arrow[d, "g^* \operatorname{Tr}_F"] \arrow[l] \\   \mathbb Z/n & &g^* \mathbb Z/n \arrow[ll] \end{tikzcd}\end{equation}
where $g$ is the map from the site of open sets on the topological space $X(\mathbb C)$ to the \'{e}tale site of $X$ and $g'$ is the analogous map for $Y$.

We can always choose a trace map to make this diagram commute, and furthermore, this trace map is an isomorphism. To check this, it suffices to check that every other arrow in the diagram is an isomorphism. For the top-left arrow and bottom arrow, this is immediate from the definition. For the top-right arrow this is a case of the comparison theorem \cite[XVI, Theorem 4.1]{sga4-3} between \'{e}tale cohomology and singular cohomology. For the right vertical arrow, we must check that the trace map in \'{e}tale cohomology is an isomorphism. It suffices to check it is an isomorphism on stalks. To do this, we use \cite[XVIII, Proposition 2.10]{sga4-3}, which states that in the case the that fiber is a curve, the trace map of \cite[XVIII, Theorem 2.9]{sga4-3} agrees with the trace map of \cite[XVIII, Proposition 1.1.6]{sga4-3}, which is itself defined to agree when the base is the spectrum of an algebraically closed field (in our case, $\mathbb C$) with the trace map of \cite[1.1.3.4]{sga4-3}.

So it suffices to consider the trace map of \cite[XVIII, 1.1.3.4]{sga4-3}. This is defined, in the case of a projective curve $Y_x$, by taking the Kummer isomorphism between $H^2( Y_x, \mu_n)$ and $\operatorname{Pic}( Z_x)/n$, observing that $\operatorname{Pic}( Y_x)/n$ is isomorphic to $(\mathbb Z/n)^c$, $c$ the number of components, and taking a sum map to $\mathbb Z/n$, weighted by the multiplicity of the components. When the curve is smooth projective irreducible, as it is in our case, there is only one component with multiplicity one, so the weighted sum is simply the identity map $\mathbb Z/n\to \mathbb Z/n$, and thus the trace map is an isomorphism.

  Using the ``integration along the fibers'' map, we can define an analogue of the map $e$ in singular cohomology the same way we defined the original map, as the map  $\pi_1(X(\mathbb C))\to H_3(H,\mathbb Z/n)$ induced by the composition
  \[ H^3(H, \mathbb Z/n) \to H^3 ( Y(\mathbb C), \mathbb Z/n) \to H^1( X(\mathbb C), \mathbb Z/n) \to \Hom (\pi_1(X(\mathbb C)), \mathbb Z/n) .\]
inducing a homomorphism $\pi_1(X(\mathbb C) )\to H_3(H, \mathbb Z/n)$.

Let us check that the analogue agrees with $e$ after composition with the map $\pi_1(X(\mathbb C))\to\pi_1(X)$ from the topological to the \'{e}tale fundamental group. To prove this, we need to check the integration-along-the-fibers map is compatible with the comparison between \'{e}tale cohomology in characteristic zero and singular cohomology. This follows by an argument identical to Lemma \ref{integration-vs-pullback}, except that we use pullback from the \'{e}tale site to the analytic site instead of pullback between different schemes, and we use \eqref{pullback-to-analytic-diagram} as our starting point.


For $\sigma$ in the fundamental group of a component of $\Hur^m_{H\rtimes \Gamma ,c_{H \rtimes \Gamma} } (\bbC)$ corresponding to a braid, a loop in $\Hur^m_{H\rtimes \Gamma ,c_{H \rtimes \Gamma} } (\bbC)$ 
representing $\sigma$ is a lift from $\Conf^m(\mathbb C)$  to $\Hur^m_{H\rtimes \Gamma ,c_{H \rtimes \Gamma} } (\bbC)$ of a  geometric realization of that braid. The pullback of the universal family $Y$ to that loop is exactly the covering $\tilde{M}$ of $S^2 \times S^1$ constructed from that braid.

Let $S^1$ be the loop, $\tilde{M}$ be the pullback of the universal family to that loop, and $\tilde{F} \colon \tilde{M} \to S^1$ the map between them. Let $\tilde{g} \colon  S^1\to X(\mathbb C)$ and $\tilde{g} '\colon \tilde{M} \to Y(\mathbb C)$ be the maps induced by the embedding of the loop $S^1$ into $X$. Again, we can define a trace map $R^2 \tilde{F}_* \mu_n \to \mathbb Z/n$ to make the diagram
   \begin{equation}\label{pullback-to-circle-diagram} \begin{tikzcd} R^2 \tilde{F}_* \mu_n \arrow[d, "\operatorname{Tr}_{\tilde{F}}"]  &  R^2\tilde{F}_*  \tilde{g}^{'*} \mu_n \arrow[l]   & \tilde{g}^* R^2 F'_* \mu_n \arrow[d, "\tilde{g}^* \operatorname{Tr}_F'"] \arrow[l] \\   \mathbb Z/n & &\tilde{g}^* \mathbb Z/n \arrow[ll] \end{tikzcd}\end{equation}
   commute, and again this map is an isomorphism: The right vertical arrow is an isomorphism by what we checked before, that the top-left horizontal arrow and bottom horizontal arrow are isomorphisms follows from the definitions, and the top-right vertical arrow is an isomorphism because of the proper base change theorem in the setting of singular cohomology. Again this definition gives an integration-along-the-fibers map and a compatibility of integration-along-the-fibers with pullback.  Hence  $e(\sigma)$ is given by the composition
\[ H^3(H, \mathbb Z/n) \to H^3( \tilde{M}, \mathbb Z/n) \to  H^1(S^1, \mathbb Z/n) \cong \mathbb Z/n \]
where we now consider the fibration $ \tilde{M} \to S^1$ and its associated integration-along-the-fibers map.

The braid fundamental class of $\sigma$, on the other hand, is given by composing $H^3(H, \mathbb Z/n) \to H^3( \tilde{M}, \mathbb Z/n) $ with the fundamental class isomorphism $H^3( \tilde{M}, \mathbb Z/n) \to \mathbb Z/n$. Hence to check our desired statement that the braid fundamental class is equal to $e(\sigma)$ up to multiplication by an element of $(\mathbb Z/n)^\times$, it suffices to check that the integration-along-the-fibers map $ H^3( \tilde{M}, \mathbb Z/n) \to  H^1(S^1, \mathbb Z/n) $ is an isomorphism. 

This map is the composition
\[ H^3 (\tilde{M}, \mathbb Z/n) \to H^3(S^1, R \tilde{F}_* \mathbb Z/n) \to H^1  (S^1, R^2 \tilde{F}_* \mathbb Z/n) \to H^1(S^1, \mathbb Z/n).\]
The first map is an isomorphism by the derived category version of the Leray spectral sequence. The second map is an isomorphism since, in the spectral sequence computing $H^{p+q}(S^1, R \tilde{F}_* \mathbb Z/n) $ from $ H^p  (S^1, R^q \tilde{F}_* \mathbb Z/n) $, every differential to or from $H^1  (S^1, R^2 \tilde{F}_* \mathbb Z/n) $ after the second page must vanish since   $ H^p  (S^1, R^q \tilde{F}_* \mathbb Z/n) =0$ for $p>0$ or $q>0$. The third map is an isomorphism since we already saw that the trace map is an isomorphism, and we are done.\end{proof}

\section{Group theory preliminaries}\label{S:GroupTheory}

In this section, we give many results in group theory that we will use repeatedly to determine our distribution from its moments.  These results may be considered well-known, but we give them here for completeness and to fix notation.

\subsection{Alternating tensor powers}



\begin{lemma}\label{L:Homwed}
For a finite abelian group $A$ and positive integer $n$, there is a homomorphism
\begin{align*}
\Delta: &\Hom(A\tensor A,\Z/n) &\ra &\Hom(A\tensor A, \Z/n)\\
&f &\mapsto &(a\tensor b \mapsto f(b,a)-f(a,b)),
\end{align*}
whose image $\Delta\Hom(A\tensor A,\Z/n)$ is the set of $g\in \Hom(A\tensor A, \Z/n)$ such that $g(a\tensor a)=0$ for all $a\in A$ (\emph{alternating} maps).
Also, $\ker\Delta$ is the set of $g\in \Hom(A\tensor A, \Z/n)$ such that $g(a\tensor b)=g(b\tensor a)$ for all $a,b\in A$ (\emph{symmetric} maps).
\end{lemma}

\begin{proof}
Clearly all maps in $\Delta \Hom(A\tensor A, \Z/n)$ are alternating.
 We can check, using bases, that every alternating $g\in \Hom(A\tensor A, \Z/n)$  is in  $\Delta\Hom(A\tensor A,\Z/n).$ 
\end{proof}

We write $\wedge\Hom(A\tensor A,\Z/n)$ for the quotient $\Hom(A\tensor A,\Z/n)/\ker\Delta$ of $\Hom(A\tensor A,\Z/n)$ by the symmetric maps $\ker\Delta$.   Note that Lemma~\ref{L:Homwed} implies we have an isomorphism
$$
\Delta:  \wedge\Hom(A\tensor A,\Z/n) \isom \Delta\Hom(A\tensor A,\Z/n).
$$

\subsection{Representations over finite fields}
We explain in detail the relationships between the different types of self-dual representations $V$ and the properties of an invariant element of $V^\vee\tensor V^\vee$.


\begin{lemma}\label{L:typeV}
Let $V$ be an irreducible self-dual representation of a finite group $\Pi$ over a field $\F_p$ for some prime $p$.
Let $\kappa=\End_\Pi(V)$.
Let $\omega \in (V^\vee\tensor V^\vee)^\Pi\setminus 0.$
Then there is a unique  automorphism $\sigma$ of $\kappa$,
 such that for all $k\in\kappa$, we have $(k\tensor 1) \omega =(1\tensor \sigma(k) )\omega$.  Further, $\sigma^2=1$.  
 Let $\omega^t$ be the image of $\omega$ under switching factors.
There is a $\lambda\in \kappa$ such that $\lambda \sigma(\lambda)=1$ and $\omega^t=(\lambda \tensor 1) \omega$.

Then 
\begin{itemize}
\item 
$V$ is unitary if and only if $\sigma$ is non-trivial,
\item 
for odd $p$: 
$V$ is symmetric if and only if $\sigma$ is trivial and $\lambda=1$,
\item 
for $p=2$: 
$V$ is symmetric if and only if $\sigma$ is trivial and $V$ is trivial, 
\item 
for odd $p$: 
$V$ is symplectic if and only if $\sigma$ is trivial and $\lambda=-1$, and
\item 
for $p=2$: 
$V$ is symplectic if and only if $\sigma$ is trivial and $V$ is non-trivial. 
\end{itemize}

Also,
\begin{itemize}
\item if $V$ is unitary, then $\dim_{\F_p} (\wedge_2 V^\vee)^\Pi=\frac{1}{2}\dim_{\F_p} \kappa$, 
\item if $V$ is symmetric, then $\wedge_2 V^\vee=0$, and
\item if $V$ is symplectic, then $\dim_{\F_p} (\wedge_2 V^\vee)^\Pi=\dim_{\F_p} \kappa$.
\end{itemize}
\end{lemma}
\begin{proof}

Since $V$ is irreducible and self-dual, we have  $0\neq \omega\in (V^\vee\tensor V^\vee)^\Pi$ and 
$(\kappa \tensor 1)\omega =(1 \tensor \kappa)\omega = (V^\vee\tensor V^\vee)^\Pi$.   Thus or a $k\in \kappa$, we have $(k\tensor 1) \omega =(1\tensor \sigma(k) )\omega$, for
some permutation $\sigma$ of $\kappa$, and we have $\omega^t=(\lambda \tensor 1) \omega$ for some $\lambda\in \kappa$. It is immediate from the bilinearity of $\omega$ that $\sigma$ is an automorphism.
By comparing the transposes of $(k\tensor 1)\omega$ and $(1\tensor \sigma(k))\omega$,
\details{$((1\tensor k)\omega)^t=(\sigma^{-1}(k)\tensor 1)\omega^t =(\lambda \sigma^{-1}(k)\tensor 1)\omega$ 
and $((1\tensor \sigma(k))\omega)^t=(\sigma(k)\tensor 1)\omega^t=(\lambda \sigma(k)\tensor 1)\omega$}
we find that $\sigma^2=1$, and by considering the fact that $(\omega^t)^t=\omega$, we have that $\lambda \sigma(\lambda)=1$.

We have isomorphisms of $\Pi$-representations
\begin{align*}
&V^\vee \tensor V^\vee&\isom & \Hom(V,V^\vee) &\isom & \Hom(V,V^\vee)\\
&u\tensor v\mapsto & &(a \mapsto u(a)v ) & &(a \mapsto v(a)u )\\
& & &f\mapsto  &&(a\mapsto (b\mapsto f(b)(a) ))
\end{align*}
where for $M\in V^\vee \tensor V^\vee$, we write $M_L$ for the image in the first $\Hom(V,V^\vee)$ and $M_R$ for the image in the second $\Hom(V,V^\vee)$,
and for $M\in \Hom(V,V^\vee)$, we write $M^t$ for the image in the other $\Hom(V,V^\vee)$.
For $\lambda\in \{\pm 1\}$, we have that $M\in V^\vee \tensor V^\vee$ has $M^t=\lambda M$  if and only if $M_L=\lambda M_R$, or equivalently $M_L^t=\lambda M_L$.
There is an entirely analogous situation when $V^\vee$ is replaced by $\Hom_\kappa(V,\kappa)$,  and the tensor products and $\Hom$'s are over $\kappa$.
The trace isomorphism of $H$-representations
$$
\Tr :  \Hom_\kappa(V,\Hom_\kappa(V,\kappa)) \ra \Hom_\kappa(V,V^\vee)
$$
has the property that
for $\lambda\in \{\pm 1\}$ and $M\in \Hom_\kappa(V,\Hom_\kappa(V,\kappa))$, we have
 $M^t=\lambda M $ if and only if $(\Tr M)^t=\lambda \Tr M$.

When $\sigma$ is not the identity, then $\omega_L$ gives an isomorphism from  $V$ to $V^\vee$ as representations of $\Pi$ over $\F_p$.
Moreover, all such isomorphisms are $k \omega_L$ for some $k\in \kappa$, and hence none are $\kappa$-equivariant, and $V$ and $V^\vee$ are not isomorphic as $\Pi$-representations over $\kappa$.  Thus $V$ and $\Hom_\kappa(V,\kappa)$ are not isomorphic representations over $\kappa$, and hence $(V \tensor_\kappa V)^H=0$, and $V$ is unitary.

Next we consider the case when $\sigma$ is the identity. 
Then $\omega_L$ gives a $\kappa$-equivariant and $\Pi$-equivariant isomorphism $V\ra V^\vee$, and since 
$V^\vee\isom \Hom_\kappa(V,\kappa)$, this implies $V$ is self-dual over $\kappa.$  
By the remarks above,  there is a non-zero $\Pi$-invariant element $\Tr^{-1} \omega_L \in \Hom_\kappa(V,\Hom_\kappa(V,\kappa))$ such that
$(\Tr^{-1} \omega_L)^t=\lambda \Tr^{-1} \omega_L$, and hence,  a non-zero element in $(\Hom_\kappa(V,\kappa)\tensor_\kappa \Hom_\kappa(V,\kappa))^\Pi$ with the same property.
Using that $V$ is self-dual over $\kappa$,  we have a non-zero element 
$M\in (V\tensor_\kappa V)^\Pi$ such that $M^t=\lambda M$.

When $p$ is odd, we have an isomorphism of $\Pi$-representations $V \tensor_\kappa V\isom \Sym^2_\kappa V \times \wedge_\kappa^2 V$,
where $\Sym^2_\kappa V$ are the fixed points of the factor switching and  $\wedge_\kappa^2 V$ is the $-1$ eigenspace of the factor switching.
Since $V$ is irreducible, if $V \tensor_\kappa V$ has any $\Pi$-invariants, it has a one dimensional $\kappa$ vector space of $H$-invariants, which is either in 
$\Sym^2_\kappa V$ or $\wedge_\kappa^2 V$. Thus, for odd $p$, when $\lambda=1$ we have that $V$ is symmetric and when $\lambda=-1$ we have that $V$ is symplectic.

When $p=2$,  we have an exact sequence of representations
$$
0 \ra V \ra \Sym^2_\kappa V \ra \wedge^2_\kappa V \ra 0
$$
sending $v\in V$ to $v\tensor v \in \Sym^2_\kappa V$, and $[v\tensor w]\mapsto [v\tensor w]$
for the second map.  Thus if $V$ is a non-trivial representation, self-dual over $\kappa$, then it is symplectic, and if $V$ is trivial then it is symmetric.

For the claims about $\dim_{\F_p} (\wedge_2 V^\vee)$,  we have
$\wedge_2 V^\vee \isom \wedge^2 V$ as $H$-representations. 
We have $\dim_{\F_p} (\wedge^2 V)^\Pi=\dim_{\kappa} (\wedge^2_\kappa (V\tensor_{\F_p} \kappa))^\Pi$. 
We have $V\tensor_{\F_p} \kappa$ is the sum of $\dim_{\F_p} \kappa$ absolutely irreducible representations over $\kappa$,
which are dual in pairs if $V$ is unitary, and self-dual and symplectic if $V$ is symplectic, and the claims about $\dim_{\F_p} (\wedge_2 V^\vee)$ follow.\end{proof}

\begin{lemma}\label{L:ortho}
Let $\Pi$ be a finite group.
Let $V$ be an irreducible $\F_2$-orthogonal representation of $\Pi$ over $\F_2$.
Then the map induced by $\Delta$
$$(\Sym^2 V)^\Pi \ra (\wedge_2 V)^\Pi$$
is surjective.
\end{lemma}
\begin{proof}
We have an exact sequence
$$
0 \ra V \ra \Sym^2 V \stackrel{\Delta}{\ra} \wedge_2 V\ra  0,
$$
where the first map sends $v\mapsto [v\tensor v]$.  If $V$ is a trivial representation, then the lemma follows.
If $V$ is non-trivial, then $V^\Pi=0$ implies that $(\Sym^2 V)^\Pi \ra (\wedge_2 V)^\Pi$
is an injection.  So let $Q$ be a non-zero element of $(\Sym^2 V)^\Pi$, and then $\omega=\Delta Q$ is a non-zero element of
$(\wedge_2 V)^H.$  For $k\in\kappa$, we have that $(k\tensor k)Q\in (\Sym^2 V)^\Pi$ and maps under $\Delta$ to $(k\tensor k)\omega$.
We will show the elements $(k\tensor k)\omega$ are all the elements of $(\wedge_2 V)^\Pi$.  In the language of Lemma~\ref{L:typeV},
if $\sigma$ is non-trivial,  then we are considering the elements $(k\sigma(k)\tensor 1)\omega,$ which gives $|\kappa|^{1/2}$ elements of
 $(\wedge_2 V)^\Pi$, which by Lemma~\ref{L:typeV} must be all of them.
 \details{Elements of the form $k\sigma(k)$ are of the form $k^{p^{d/2}+1}$ and we can count these, they are the $p^{d/2}-1$ roots of unity and $0$.
If $(n\tensor 1 )\omega =(n'\tensor 1 )\omega $, then $((n-n')\tensor 1)\omega=0$, but $n-n'$ is invertible, so these are all different.  }
If $\sigma=1$, then we are considering the elements $(k^2\tensor 1)\omega,$ which gives $|\kappa|$ elements of
 $(\wedge_2 V)^\Pi$, which by Lemma~\ref{L:typeV} must be all of them.
\end{proof}

\subsection{Group homology and cohomology}

When we write group homology or cohomology with $\Z/n$ coefficients, we always mean for the trivial action of the group on $\Z/n$.

\begin{lemma}\label{L:GammaH}
For a finite group $\Gamma$ and a finite $\Gamma$-group $G$,  and finite abelian group $A$ of order relatively prime to $|\Gamma|$,
the pull-back map $H^k (G\rtimes \Gamma,A)\ra H^k (G,A)^\Gamma$ is an isomorphism for all $k\geq 0$.
\end{lemma}
\begin{proof}
In Lyndon-Hochschild-Serre spectral sequence computing $H^{p+q} (G\rtimes \Gamma,A)$ from $H^p(\Gamma, H^q(G,A))$
the only non-zero terms are when $p=0$, and so the edge maps to those terms are isomorphisms.  
\end{proof}

\begin{lemma}[{\cite[03VI]{stacks}}]\label{L:homdual}
Let $n$ be a positive integer.  There is an natural isomorphism of functors on the category of profinite groups from $H_n(-,\Z/n)$ to $\Hom(H^n(-,\Z/n),\Z/n)$.
\end{lemma}

Via Lemma~\ref{L:homdual}, we often view elements of $H_n(G,\Z/n)$ as homomorphisms  $H^n(-,\Z/n)\ra \Z/n.$

\begin{lemma}\label{L:Delta}
If $A$ is an abelian group and $n$ a positive integer, then there is an $\Aut(A)$-equivariant homomorphism
\begin{align*}
\Delta: &H^2(A,\Z/n) &\ra &\Hom(A\tensor A, \Z/n)\\
&[f] &\mapsto &(a\tensor b \mapsto f(b,a)-f(a,b)),
\end{align*}
where $f$ is any normalized $2$-cocycle.  
The image of the map $\Delta$ above is the same as $\Delta\Hom(A\tensor A,\Z/n)$ from Lemma~\ref{L:Homwed}.
\end{lemma}

\begin{proof}
We claim for a normalized cocycle $f$, the map $a\tensor b \mapsto f(b,a)-f(a,b)$ is a homomorphism
from $A\tensor A$ to $\Z/n$.
This can be checked by noting it is the map that sends $a\tensor b$ to the commutator of lifts of $a,b$ in the extension corresponding to $f$.  This will also be confirmed in the proof of Lemma~\ref{L:d202gen} below.

The map in the lemma is well-defined on $H^2$ because any coboundary is a symmetric function of the two inputs, and thus is mapped to $0$.  Given that it is well-defined, it is clearly a homomorphism.  

Since an element of $\Hom(A\tensor A,\Z/n) $ is a cocycle, the image from $H^2(A,\Z/n)$ is at least as large as 
$\Delta \Hom(A\tensor A,\Z/n)$.  The image from $H^2(A,\Z/n)$ is a subset of the set alternating $g$, and we conclude the final statement of the lemma, using Lemma~\ref{L:Homwed}.
\end{proof}

\subsection{Extensions by center free groups}

 We see there are unique extensions by center free $[H]$-groups.

\begin{lemma}\label{L:nonabrig}
Let $H$ be a   group and let $N$ be a   $[H]$-group with trivial center. 
The conjugation map $N\ra \Aut(N)$ gives an exact sequence
$1\ra N \ra \Aut(N) \times_{\Out(N)} H \ra H \ra 1$,  and the induced $[H]$-structure on $N$ agrees with the given one.

Let $1 \ra N \ra G \ra H \ra 1$ be an exact sequence of groups,  such that the induced $[H]$-structure on $N$ agrees with the given one.
Then there is a unique isomorphism $G\ra \Aut(N) \times_{\Out(N)} H$ compatible with the identity maps on $N$ and $H$.
\end{lemma}

\begin{proof}
Since $N$ has trivial center $N \ra \Aut(N) \times_{\Out(N)} H$ is an injection, and the quotient by $N$ is $H$, which proves the first claim.

For the second claim, we send $g\in G$ to $\alpha(g):=(c_g, \bar{g})\in \Aut(N) \times_{\Out(N)} H$, where $c_g$ is the automorphism of $N$ given by conjugation by $g$ and 
$\bar{g}$ is the image of $g$ in $H$.  They must have the same image in $\Out(N)$ by the condition on the exact sequence.
This map is compatible with the identity maps on $N$ and $H$ and it follows (via the short five lemma in the category of groups) that $\alpha:G\ra \Aut(N) \times_{\Out(N)} H$ is an isomorphism.

Any other compatible automorphism $G \ra \Aut(N)\times_{\Out(N)} H$ must have the form $\alpha \circ \beta$ for  $\beta: G\ra G$ an isomorphism compatible with the identity maps on $N$ and $H$.
Then for $g\in G$ and $n\in N$, we have $gng^{-1}=\beta(gng^{-1})=\beta(g) n \beta(g)^{-1}$.  Also $\bar{g}=\bar{\beta(g)}$.  Hence $\alpha(g)=\alpha(\beta(g))$ so $\alpha = \alpha \circ \beta$.
\end{proof} 
 
\subsection{Extensions of $\Gamma$-groups}
Let $\Gamma$ be a group, $H$ a $\Gamma$-group, and $F$ a $[H\rtimes\Gamma]$-group.  
 We consider an exact sequence $$ 1 \to  F \to G \stackrel{\pi}{\to} H \to 1$$ of $\Gamma$-groups,  in which the $[H\rtimes \Gamma]$ structure on $F$ from the exact sequence agrees with the given $[H\rtimes \Gamma]$ structure.   We call such a thing an \emph{extension} of $\Gamma$-groups of $H$ by $F$. 
  Two such extensions 
$ 1 \to  F \to G \stackrel{\pi}{\to} H \to 1$ and 
$ 1 \to  F \to G' \stackrel{\pi'}{\to} H \to 1$ 
  are isomorphic if there is an isomorphism of $\Gamma$-groups from $G$ to $G'$ that restricts to the identity of $F$ and $H$.  We write $\Ext_{\Gamma}(H,F)$ for the set of
  isomorphism classes of such extensions.  We often write $(G,\pi)$ for such an extension (even though the map $F\ra G$ is also part of the data, we leave it implicit).
We write $\Aut_{F,H}(G)$ for the isomorphisms of such an extension, i.e. $\Gamma$-automorphisms of $G$ that restrict to the identity  on $F$ and $H$.

\subsubsection{Abelian extensions}

If $\Gamma$ were trivial, then isomorphism classes of extensions of a group $H$ by an abelian $H$-group $A$
exactly correspond to classes of $H^2(G,A)$.  We seek a similar characterization of $\Gamma$-extensions by cohomology, which turns out to be simplified when $|\Gamma|$ is relatively prime to $|H|$ and $|F|$.

\begin{lemma}\label{L:diffext}
Let $\Gamma$ be a group, $H$ a $\Gamma$-group, and $A$ an abelian $(H\rtimes\Gamma)$-group.  
Suppose $|\Gamma|$ is relatively prime to $|H|$ and $|A|$.  There is a bijection between $\Ext_\Gamma(H,A)$ and $H^2(H\rtimes \Gamma,A)$
taking $A\ra G\ra H$ to $A\ra G\rtimes \Gamma\ra H\rtimes \Gamma$, with the natural maps.
\end{lemma}
\begin{proof}
We consider an inverse to the given map.  Let $A\ra E \stackrel{\pi}{\ra} H\rtimes \Gamma$ be an extension of groups giving the correct $(H\rtimes \Gamma)$-action on $A$.
Let $G=\ker(E\ra \Gamma)$.
  Then there is a splitting $s:\Gamma\ra E$ by the Schur-Zassenhaus theorem.
  Moreover, by that theorem the image of the splitting in $H\rtimes \Gamma$ must be a conjugate of the trivial splitting,
  and so  by conjugating our original splitting, we can assume $s$ lifts the trivial splitting of $H\rtimes \Gamma\ra \Gamma$.

We then have group homomorphisms $A\ra G \ra H$.  The $\Gamma$ action on $G$ given by conjugation by $s(\gamma)$ in $E$ makes the map $G\ra H$ $\Gamma$-equivariant. 
Moreover, $A$ has the intended $H\rtimes\Gamma$ action.
 Thus we have given a map from $H^2(H\rtimes \Gamma,A)$ to $\Ext_\Gamma(H,A)$, though we should check it is well-defined.  
Writing $E=G\rtimes \Gamma$ using our original section,  
if an alternate section also agrees with the trivial section to $H\rtimes \Gamma$,  it must land in $A\rtimes\Gamma$. Hence by 
Schur-Zassenhaus
we have an $a\in A$ such that our alternate section sends $\gamma \mapsto (a\gamma(a)^{-1},\gamma)$.
Then, we can check that conjugation by $a$ gives a map $G\ra G$,  that is the identity on $A$ and $H$, taking the original $\Gamma$-action 
to the $\Gamma$-action that is conjugation by this alternate section.  So we have constructed the same isomorphism class in $\Ext_\Gamma(H,A)$.

It is straightforward to check that these two maps are inverses, and hence we've given a bijection from $\Ext_\Gamma(H,A)$ to $H^2(H\rtimes \Gamma,A)$.
 \end{proof}
 
In the setting of Lemma~\ref{L:diffext}, we then often consider  $H^2(H\rtimes \Gamma,A)$ as parametrizing the elements of $\Ext_\Gamma(H,A)$.

For a set $\cL$ of isomorphism classes of $\Gamma$-groups,  we write $H^2(H\rtimes \Gamma,A)^\cL$ for the elements of 
$H^2(H\rtimes \Gamma,A)$ that correspond to $\Gamma$-group extensions of $H$ by $A$ whose underlying $\Gamma$-group is in $\cL$.
The following lemma is straightforward to check.

\begin{lemma}\label{L:basicH2}
Let $\Gamma$ be a group and $H$ a $\Gamma$-group.
Let $\cL$ be a set of isomorphism classes of $\Gamma$-groups,  closed under taking fiber products and quotients, 

Let $A$ be a simple abelian $(H\rtimes \Gamma)$-group.  Then $\kappa:=\End_{H\rtimes \Gamma}(A)$ (i.e.  $(H\rtimes \Gamma)$-group morphisms from $A$ to $A$) is a field.  Moreover, $H^2(H\rtimes \Gamma,A)$ is naturally a $\kappa$-vector space through the $\kappa$ action on $A$ and $H^2(H\rtimes \Gamma,A)^\cL$ is a $\kappa$-subspace of $H^2(H\rtimes \Gamma,A)$ as long as it is nonempty.

Let $B$ and $C$ be abelian $(H\rtimes \Gamma)$-groups.  The isomorphism $H^2(H\rtimes \Gamma, B\times C)\ra H^2(H\rtimes \Gamma, B) \times H^2(H\rtimes \Gamma, C)$ takes
$H^2(H\rtimes \Gamma, B\times C)^\cL$ to $H^2(H\rtimes \Gamma, B)^\cL \times H^2(H\rtimes \Gamma, C)^\cL$.
\end{lemma}

\begin{proof} $\kappa$ is a field since any non-zero non-invertible element of $\kappa$ would have image a nontrivial proper subspace, which is impossible, so $\kappa$ is a division algebra, and every finite division algebra is a field.

The bijection between $H^2(H\rtimes \Gamma,A)$ and $\Gamma$-group extensions of $H $ by $A$ sends addition on $H^2(H\rtimes \Gamma,A)$ to the Baer sum, defined by taking fiber products of $H$ and quotienting by the antidiagonally embedded $A$. Since $\cL$ is closed under taking fiber products and quotients, thet Baer sum is in $\cL$ if the individual extensions are. Thus $H^2(H\rtimes \Gamma,A)^\cL$ is a subspace of a finite group closed under addition, hence it is a subgroup unless it is empty. The action of invertible elements of $\kappa$ on $H^2(H\rtimes \Gamma,A)$ sends an extension to an extension which is isomorphic as a $\Gamma$-group but with the map from $A$, part of the data of the extension, composed with an invertible endomorphism of $A$. Since the extension is isomorphic as a $\Gamma$-group to the original, it is in $\cL$ if the original extension is. Thus $H^2(H\rtimes \Gamma,A)^\cL$ is $\kappa$-stable, and hence a $\kappa$-subspace, as long as $H^2(H\rtimes \Gamma,A)^\cL$ is nonempty.

The isomorphism $H^2(H\rtimes \Gamma, B\times C)\ra H^2(H\rtimes \Gamma, B) \times H^2(H\rtimes \Gamma, C)$ sends an extension of $H$ by $B \times C$ to a pair consisting of its quotient by $C$ and its quotient by $H$. The inverse map sends a pair of extension by $B$ and $C$ respectively to their fiber products. Since $\cL$ is stable under fiber products and quotients, the extension by $B \times C$ is in $\cL$ if and only if both the extension by $B$ and the extension by $C$ are. \end{proof}

 \begin{lemma}\label{L:GGamma}
 If $\Gamma$ is a group, and $H$ is a $\Gamma$-group with $H_\Gamma=1$, and $A$ is an abelian $(H\rtimes\Gamma)$-group, 
and $A\ra G\ra H$ is an extension of $\Gamma$-groups such that the induced $H\rtimes \Gamma$ action on $A$ is as given,
then $G_\Gamma$ is a quotient of
$A_{H\rtimes\Gamma}$.
 \end{lemma}
 \begin{proof}
 Since $H_\Gamma=1$, it follows that the normal subgroup $A$ of $G$ surjects onto $G_\Gamma$, and in particular that $G_\Gamma$ is abelian.  
 Hence any $a\in A$, and lift $\tilde{h}$ to $G$ of an element $h\in H$, have images that commute in $G_\Gamma$.
 So the image of $A$ in $G_\Gamma$ factors through $A_{H\rtimes \Gamma}$, and is all of $G_\Gamma$.    
 \end{proof}

\subsubsection{Non-abelian extensions of $\Gamma$-groups}

\begin{lemma}\label{L:nonabGam}
Let $\Gamma$ be a group, 
$H$ be a $\Gamma$-group,  and $N$ a $[H\rtimes \Gamma]$-group with trivial center.
Let $1\ra N \ra G \ra H \ra 1$ be an exact sequence of groups with induced $[H]$-group structure on $N$ agreeing with the given $[H]$-group structure.

There is a bijection from (1) $\Gamma$ actions on $G$ such that the given map $G\ra H$ is $\Gamma$-equivariant and the induced $\Gamma$-group structure on $N$ reduces to the given $[\Gamma]$-structure on $N$, and (2) lifts of the given map $\Gamma\ra \Out(N)$ to $\Gamma\ra \Aut(N)$ (respecting the map $\Aut(N)\ra \Out(N)$).  This bijection is given by taking the induced $\Gamma$-group structure on $N$.
\end{lemma}
\begin{proof}
By Lemma~\ref{L:nonabrig}, we may assume $G=\Aut(N) \times_{\Out(N)} H$.  We give an inverse to the given map.
Given a map $\Gamma\ra \Aut(N)$, lifting $\Gamma\ra\Out(N)$, we get a homomorphism $\Gamma \ra \Aut(N) \times_{\Out(N)} (H\rtimes\Gamma)$, using the identity map of $\Gamma$ for the second coordinate.  The group $G$ is a normal subgroup of $\tilde{G}:=\Aut(N) \times_{\Out(N)} (H\rtimes\Gamma)$, and thus
we get an action of $\Gamma$ on $G$ by conjugation in $\tilde{G}$.  We can check that this action makes $G\ra H$ a $\Gamma$-equivariant map, and the induced $\Gamma$-group structure on $N$ is compatible with the given $[\Gamma]$-group structure.

We now check these constructions are inverses.
Let $\alpha: \Gamma \ra \Aut(G)$ be an action as in (1).  If $\alpha(\gamma)(c_n,1)=(c_{\beta(\gamma)(n)},1)$ for $\gamma\in \Gamma$ and $n\in N$, then $\beta: \Gamma \ra \Aut(N)$ is the 
given map in (2).  Then we use $\beta$ to obtain $\alpha': \Gamma \ra \Aut(G)$ such that for $\gamma\in \Gamma$ and $\sigma\in\Aut(N)$ and $h\in H$, we have 
 $$\alpha'(\gamma)(\sigma,h )=(\beta(\gamma),\gamma )(\sigma,h )(\beta(\gamma)^{-1},\gamma^{-1} )=(\beta(\gamma)\sigma\beta(\gamma)^{-1},\gamma(h)).$$
 If $\sigma=c_n$ for some $n\in N$, note that $\beta(\gamma)c_n\beta(\gamma)^{-1}=c_{\beta(\gamma)(n)}$
We claim $\alpha(\gamma)$ and $\alpha'(\gamma)$ are the same element in $\Aut(G)$.  This is because they give the same action on $H$, and the same action on $N$,  and so  by Lemma~\ref{L:nonabrig} they must be the same.
For the other direction, we start with $\beta: \Gamma \ra \Aut(N)$, and obtain $\alpha'$ as above.  We've already seen above that the action coming from $\alpha'$ on $N$ agrees with $\beta$.  This shows that the given construction is a bijection.
\end{proof}

\section{Analysis of $E^3_{0,2}$ in the Lyndon-Hochschild-Serre spectral sequence}\label{S:LHS}

In this section, we give a detailed analysis of the $E^3_{0,2}$ term in the Lyndon-Hochschild-Serre spectral sequence, especially for an extension of $H$
by a semisimple abelian $H$-group.
This will be crucial for our determination of a distribution from its moments, but strictly speaking, the analysis we do is not restricted to the context of this problem.  

\subsection{$H^2(A,\Z/n)$}

First, we must carefully understand $H^2(A,\Z/n)$ for an abelian group $A$ of exponent dividing $n$.
Recall the definition of $\Delta$ and $\wedge\Hom(A\tensor A,\Z/n)$ from Lemma~\ref{L:Homwed}.

\begin{lemma}\label{L:H2struc}
Let $n$ be a positive integer and let $A$ be a finite abelian group of exponent dividing $n$,
and such that $4$ does not divide the exponent of $A$.  Then there is an 
$\Aut(A)$-equivariant injective Bockstein homomorphism
$$
B: \Hom(A,\Z/n) \ra H^2(A,\Z/n),
$$
and an $\Aut(A)$-equivariant  morphism using bilinear forms as cochains
$
C: \Hom(A\tensor A,\Z/n) \ra H^2(A,\Z/n).
$

When $n$ is odd or $4\mid n$, these maps induce an isomorphism
$$
B\times \bar{C}: \Hom(A,\Z/n) \times \wedge\Hom(A\tensor A,\Z/n)  \isom H^2(A,\Z/n).
$$

If $n=2$, the map $C$ induces an isomorphism
$$
\Sym^2 \Hom(A,\Z/2) \isom  H^2(A,\Z/2).
$$

For any $n$,  the map $C$ induces an isomorphism
$$
\bar{\bar{C}}: \wedge\Hom(A\tensor A,\Z/n)  \ra H^2(A,\Z/n)/\im B
$$
and the map $\Delta$ induces an isomorphism
$$
\Delta: H^2(A,\Z/n)/\im B \ra \Delta\Hom(A\tensor A,\Z/n).
$$ 
\end{lemma}

\begin{proof}

If for each element of $s\in\Z/n$ we choose a lift $\tilde{s}\in\Z/n^2$, then for a element $\phi\in \Hom(A,\Z/n) $, we can send $v$ to a 
$2$-cocycle $f(a,b):=\frac{\widetilde{\phi(a)}+\widetilde{\phi(b)}-\widetilde{\phi(a+b)}}{n}$.  This induces the map $H^1(A,\Z/n)\ra H^2(A,\Z/n)$ from the Bockstein homomorphism for the short exact sequence $1\ra \Z/n\ra \Z/n^2\ra \Z/n\ra 1$, and thus $B$ is injective since the exponent of $A$ divides $n$.   

Any element of $\Hom(A\tensor A,\Z/n)$ 
naturally gives a $2$-cochain $A^2\ra\Z/n$, which is automatically a $2$-cocycle, giving the map $\Hom(A\tensor A,\Z/n) \ra H^2(A,\Z/n)$.
The $2$-cochains $A^2\ra\Z/n$ that are coboundaries are symmetric in the two $A$ coordinates.

Thus any bilinear form in the kernel of $C$ is symmetric since it is a coboundary of a 1-cocycle.
We claim that if $n$ is odd or $4|n$ then
any symmetric element of $\Hom(A\tensor A,\Z/n)$ is in the kernel of $\Hom(A\tensor A,\Z/n) \ra H^2(A,\Z/n)$, and that if $n\equiv 2 \pmod{4}$ then any symmetric element of $\Hom(A\tensor A,\Z/n)$ is in the kernel of $\Hom(A\tensor A,\Z/n) \ra H^2(A,\Z/n)/\im B$.
By factoring $\Z/n\Z$ we may assume $n$ is either odd or a power of $2$.  
Given $M_0\in \Hom(A\tensor A, \Z/n)$, then if we let $M(a\tensor b):=M_0(a\tensor b)+M_0(b\tensor a)$ for $a,b\in A$,  then 
$$
M(a\tensor b)= -M_0(a\tensor a)-M_0(b\tensor b) +M_0((a+b)\tensor (a+b)) 
$$
and so $M$ is a coboundary.  If $n$ is odd, for any 
symmetric element
$M\in \Hom(A\tensor A, \Z/n)$, we can take $M_0=(1/2)M$, and so $M$ is a coboundary.  
If $n$ is even,  let $A_2$ be the maximal abelian $2$-group quotient of $A$, which by assumption is elementary abelian.
We choose a basis of $A_2$, and if the matrix corresponding to $M$ is $0$ along the diagonal then we can choose $M_0$ to be the upper triangular part of $M$ and see that $M$ is a coboundary.
So the claim will follow if we can show that the $M_1\in\Hom(A_2\tensor A_2, \Z/n)$ whose matrix has a single $n/2$ in the upper left corner (and all other entries 
$0$) maps to 0 in $H^2(A,\Z/n)$ or $H^2(A,\Z/n)/\im B$ (as required).  If $4\mid n$, we let $g: A \ra \Z/n\Z$ be such that
$g(a)=n/4$ for $a\in A$ whose image in $A_2$, written in our chosen basis,  contains a copy of the first basis vector $e_1$, and $g(a)=0$ otherwise.  We can check, for all $a,b\in A$, that
$$
g(a)+g(b)-g(a+b)=M_1(a,b)
$$
by considering separately the cases that the images of $a,b$ in $A_2$ both contain $e_1$, neither contain $e_1$, or exactly one of them contains $e_1$.  Thus if $4\mid n$, we have that $M_1$ maps to $0$ in $H^2(A,\Z/n)$.
If $n=2$, then we can let $\phi\in\Hom(A,\Z/n)$ be the homomorphism such that $\phi(a)$ is $1$ if $a$ contains $e_1$ and $0$ otherwise.  We pick a lift $\Z/n \ra \Z/n^2\Z$ such that $\widetilde{0}=0$ and $\widetilde{1}=1$.
We can check, for all $a,b\in A$, that
$$
\frac{\widetilde{\phi(a)}+\widetilde{\phi(b)}-\widetilde{\phi(a+b)}}{n}=M_1(a,b)
$$
by considering separately the cases that the images of $a,b$ in $A_2$ both contain $e_1$, neither contain $e_1$, or exactly one of them contains $e_1$.   Thus if $2\mid n$, but $4\nmid n$,we have that $M_1$ maps to $\im B$ in $H^2(A,\Z/n)$.

If $n$ is odd or $4\mid n$, we have seen above that the kernel of $C$ is precisely the symmetric forms and hence 
$\bar{C}$ is injective.  Also $\im \bar{C}$ does not contain the class of any symmetric nontrivial $2$-cocycle and thus $\im\bar{C}\cap \im B=0$.  When $2\mid n,$ but $4\nmid n$,
we consider
$${C'}:\Hom(A\tensor A,\Z/n)  \ra H^2(A,\Z/n)/\im B.$$
Any $2$-cocycle representing an element trivial on the right-hand side must be symmetric, since it is a coboundary plus an element of $\im B$.  We have also seen that all the symmetric forms in $\Hom(A\tensor A,\Z/n)$ are in the kernel of $C'$, so the kernel must be precisely the symmetric forms in $\Hom(A\tensor A,\Z/n)$.  This proves that $\bar{\bar{C}}$ is injective.

For any $n$,
we can use the universal coefficient theorem to compute
$$
|H^2(A,\Z/n)|=|\Hom(A,\Z/n)|| \wedge\Hom(A\tensor A,\Z/n)|.
$$
When $n$ is odd or $4\mid n$,  we then see by counting that $ B \times  \bar{C}$ is an isomorphism.
For any $n$, we have by counting that $\bar{\bar{C}}$ is an isomorphism.  

We have from Lemma~\ref{L:Homwed} that the composite
$$
 \wedge\Hom(A\tensor A,\Z/n)  \stackrel{\bar{\bar{C}}}{\ra} H^2(A,\Z/n)/\im B
\ra H^2(A,\Z/n)/\im B \stackrel{\Delta}{\ra} \Delta\Hom(A\tensor A,\Z/n).
$$
is an isomorphism, which implies that $\Delta:H^2(A,\Z/n)/\im B \ra \Delta\Hom(A\tensor A,\Z/n)$ is an isomorphism.

When $n=2$,  we have that $\Sym^2 \Hom (A,\Z/2\Z)$ is the quotient of $\Hom(A\tensor A,\Z/2)=\Hom (A,\Z/2\Z) \tensor \Hom (A,\Z/2\Z)$ by $\wedge_2  \Hom (A,\Z/2\Z),$ the elements
corresponding to matrices over $\F_2$ that are symmetric and $0$ along the diagonal.  
We have seen above that  elements corresponding to matrices over $\F_2$ that are symmetric and $0$ along the diagonal are in $\ker C$.
Moreover, all elements in $\ker C$ correspond  to symmetric matrices, and it follows from that above that any non-zero $\rho\in \Hom (A,\Z/2\Z) \tensor \Hom (A,\Z/2\Z)$ corresponding to a non-zero diagonal matrix has $C(\rho)$ equal to 
a non-zero element of $\im B$.  Thus $\ker C=\wedge_2 \Hom (A,\Z/2\Z) $, which implies that $C$ induces an injection $\Sym^2 \Hom (A,\Z/2\Z) \tensor \Hom (A,\Z/2\Z) \ra H^2(A,\Z/2)$, and the injection must be an isomorphism since the two groups are the same size.
\end{proof} 

\subsection{Functoriality}

We will need to use the functoriality of the Lyndon-Hochschild-Serre spectral sequence.

\begin{lemma} \label{pullback-differential} 
Let $n,m$ be integers with a given map $\pi: \Z/m\ra \Z/n.$
Let \[ \begin{tikzcd} 1 \arrow[r] \arrow[d] & F \arrow[r] \arrow[d,"\rho"]& G \arrow[r]\arrow[d] & H \arrow[r] \arrow[d,"\psi"] & 1\arrow[d] \\ 
1 \arrow[r] & F'\arrow[r] & G' \arrow[r] & H' \arrow[r] & 1  \end{tikzcd}\] be a commutative diagram of groups with both rows exact.  Then the differentials in the Lyndon-Hochschild-Serre spectral sequences $(E,d),(E',d')$ computing $H^{p+q} (G, \Z/n)$ and $H^{p+q} (G', \Z/m)$, respectively, from $H^p ( H, H^q ( F, \Z/n))$ and $H^p ( H, H^q ( F', \Z/m))$, respectively, are compatible with the pullback map $\tau$ induced by $\rho$,$\psi,$ and $\pi$,
\[ \tau: H^p ( H' , H^q ( F', \Z/m)) \to H^p ( H, H^q ( F, \Z/n)),\]
i.e. for all $r\geq 2$,
$$
\tau (d')_r^{p,q} =d_r^{p,q} \tau,
$$  
where pullback maps $\tau$ on pages past the second page are well-defined by the commutativity of these diagrams on previous pages.
\end{lemma}

\begin{proof} The lemma follows in a straightforward way from the definition of the spectral sequence, e.g. see \cite[Lemma 7.11]{Sawin2024}.
\end{proof}

\subsection{Analysis of the differentials}
We will work to relate the differentials to concrete constructions.  The first two lemmas will describe $d_{2}^{0,2}$.

Recall a $2$-cochain $f$ is \emph{normalized} if $f(a,a')=0$ whenever $a=0$ or $a'=0$.

\begin{lemma}\label{L:ASpcon}
Let $H$ be a finite group,  $n$ a positive integer, and $A$ a finite $H$-module, of exponent dividing $n$.
Let $f:A^2\ra \Z/n\Z$ be a normalized 2-cocycle such that its class  $[f]\in H^2(A,\Z/n\Z)$ is $H$-invariant.

For each $h\in H$, let
 $c_h : A \ra \Z/n$ be a normalized cochain such that $f(h^{-1}(a),h^{-1}(b))-f(a,b)=c_h(a)+c_h(b)-c_h(a+b)$ for all $a,b\in A$, and
 $c_1(a)=0$ for all $a\in A$, and for $h,i\in H$ with the same action on $A$, we have $c_h=c_{i}$.
 
 Let $E$ be the extension of $A$ by $\Z/n$ corresponding to the class $[f]$.  Let $\mathcal{A}$ be the group of automorphisms
  of $E$ that fix $\Z/n$ pointwise and act on $A$ through the image $\bar{H}$ of $H$  in $\Aut(A)$.  Then $\mathcal{A}$ is an extension
  of $\bar{H}$ by $\Hom(A,\Z/n)$ with extension class in $H^2(\bar{H},\Hom(A,\Z/n))$
represented by the cocycle  $(h,i)\mapsto c_h+c_i\circ h^{-1}-c_{hi},$
with the action of $\bar{H}$ on  $\Hom(A,\Z/n)$ induced by the action of $H$ on $A$.
\end{lemma}

We will prove Lemma~\ref{L:ASpcon} in the course of the proof of Lemma~\ref{L:d202gen}.
We let 
$$
d_{2,2}^{0,2}: H^2(A,\Z/n\Z)^H \ra H^2(H,\Hom(A,\Z/n))
$$
denote the homomorphism sending $[f]$ to the cocycle $(h,i)\mapsto c_h+c_i\circ h^{-1}-c_{hi}$ from Lemma~\ref{L:ASpcon}.
Note that $d_{2,2}^{0,2}$ vanishes on classes represented by $H$-invariant cocycles.  Also, note that $d_{2,2}^{0,2}$ only depends on $H$, $A$, $n$, and the $H$-action on $A$.

\begin{lemma}\label{L:d202gen}
Let $n$ be a positive integer.
Let $1\ra A \ra G \ra H \ra 1$ be a short exact sequence of groups with $A$ abelian of exponent dividing $n$, 
and let $G$ have extension class in $H^2(H,A)$ 
represented by a cocycle $\alpha: H^2\ra A.$ 
Let $f:A^2\ra \Z/n\Z$ be a normalized 2-cocycle such that its class  $[f]\in H^2(A,\Z/n\Z)$ is $H$-invariant.
 We define a cochain $d_{2,1}^{0,2}(f): H^2\ra\Hom(A,\Z/n)$:
\begin{align*}
&d_{2,1}^{0,2}(f)(h,i)(a):=f(a,\alpha(h,i))-f(\alpha(h,i),a)=\Delta f(\alpha(h,i),a),
\end{align*}
which induces a homomorphism $d_{2,1}^{0,2}:H^2(A,\Z/n\Z)^H \ra H^2(H,\Hom(A,\Z/n))$.

In the Lyndon-Hochschild-Serre spectral sequence computing $H^\vee(G,\Z/n)$ from $H^p(H,H^q(A,\Z/n))$,  the differential
$$
d_2^{0,2} : H^2(A,\Z/n)^H \ra H^2(H, \Hom(A,\Z/n))
$$
satisfies
$$
d_2^{0,2}=d_{2,1}^{0,2}+d_{2,2}^{0,2}.
$$

If $\alpha=0$, then $d_3^{0,2}=0$.

Let
$\mathcal{B}: H^2(H,\Z/n)\ra H^3(H,\Z/n)$,
be the Bockstein map,
 which is defined as the
boundary map coming from the exact sequence $0\ra \Z/n \ra \Z/n^2\ra\Z/n\ra 0$ of coefficients. If we use the injection
$$
B:\Hom(A,\Z/n)^H  \ra H^2(A,\Z/n)^H
$$
from Lemma~\ref{L:H2struc}, 
 for $\phi\in\Hom(A,\Z/n)^H$,  we have  $d_2^{0,2}(B(\phi))=d_{2,1}^{0,2}(B(\phi))=d_{2,2}^{0,2}(B(\phi))=0$
 and
 $d_3^{0,2}([B(\phi)])$
is the image of $[\alpha]\in H^2(H,A)$ under $\phi_*:H^2(H,A) \ra H^2(H,\Z/n)$ and then $\mathcal{B}$.
\end{lemma}

\begin{remark}\label{R:e1bil}
We will identify $\Hom(A\tensor A,\Z/n\Z)$ with $\Hom(A,\Hom(A,\Z/n))$ such that
$f\in \Hom(A\tensor A,\Z/n\Z)$ sends $a$ to the map that sends $b$ to $f(a\tensor b)$.
Recall from Lemma~\ref{L:Delta}, for $f\in H^2(A,\Z/n)$, we have $\Delta f\in \Hom(A\tensor A,\Z/n)$,
which we now also view as an element of $\Hom(A,\Hom(A,\Z/n))$.

Then $[d_{2,1}^{0,2}(f)]\in H^2(H,\Hom(A,\Z/n))$ is the image of $\alpha$ under the induced map
$$
(\Delta f)_* :H^2(H,A) \ra H^2(H,\Hom(A,\Z/n)).
$$
\end{remark}

\begin{proof}
Huebschmann \cite{Huebschmann1981} describes some of the differentials in the spectral sequence in a more general situation than we have here.
Let
\begin{equation}
0\ra \Z/n \ra E \ra A \ra 1
\end{equation}
be the central extension corresponding to  $[f]\in H^2(A,\Z/n)^H$.  
Since $f$ is normalized, and we can write $E$, as a set, as $\Z/n\times A$, where for $j,k\in \Z/n$ and $a,b\in A$, we have
$$
(j,a)\cdot _E(k,b)=(j+k+f(a,b),a+b).
$$
Let $\bar{H}$
 be the subgroup of $\Aut(A)$ that is the image of the action of $H$ on $A$ from the given exact sequence. 
Using  Huebschmann's definition \cite[Section 2.1]{Huebschmann1981} of $\Aut(A,\Z/n)$, when $A$ acts trivially on $\Z/n$, we can see that $\Aut(A,\Z/n)=\Aut(A)\times \Aut(\Z/n).$
Also, Huebschmann's $\Aut_G(A,\Z/n)$ is $\bar{H}\sub \Aut(A)\times 1 \sub \Aut(A,\Z/n)$.  
Further, Huebschmann's $\Aut_{H}^{\Z/n}(E)$ is the automorphisms of $E$ that are the identity on $\Z/n$ and act on $A$ via $\bar{H}$.

The kernel of the map $\Aut_{H}^{\Z/n}(E)\ra \bar{H}$ 
(pushing forward the automorphism to $A$)
is $\Hom(A,\Z/n)$ (as in Lemma~\ref{automorphism-of-extension}),
\details{\\
First if $\phi\in \Hom(A,\Z/n)$, then we claim $(j,a)\mapsto(j+\phi(a),a)$ is in $\Aut_{H}^{\Z/n}(E)$.  (Then it is clearly in the kernel to 
$\bar{H}$.)  We see
\begin{align*}
(j,a)&\mapsto(j+\phi(a),a)\\
(k,b)&\mapsto(k+\phi(b),b)\\
(j+k+f(a,b),b)&\mapsto(j+k+f(a,b)+\phi(a+b),a+b)
\end{align*}
while
$$
(j+\phi(a),a)(k+\phi(b),b)= (j+k+\phi(a+b)+f(a,b),a+b).
$$
So   $(j,a)\mapsto(j+\phi(a),a)$ is in the kernel of the map $\Aut_{H}^{\Z/n}(E)\ra \bar{H}$. 
Moreover,  let $\sigma$ be an element of this kernel.
Define a set map $\phi: A \ra \Z/n\Z$ such that $(\phi(a),a):=\sigma(0,a)$.  Then
\begin{align*}
&(\phi(a+b),a+b)=\sigma(0,a+b) =\sigma(-f(a,b),0)\sigma(0,a)\sigma(0,b)\\
=&(-f(a,b),0)(\phi(a),a)(\phi(b),b)=(\phi(a)+\phi(b),a+b)
\end{align*}
proving that $\phi$ is an automorphism.  Note $\sigma(j,a)=(j,0)\sigma(0,a)$, so this implies that $\sigma$ is the automorphism we obtained from $\phi$ above.
\\}
so we have an exact sequence
$$
1\ra  \Hom(A,\Z/n)\ra \Aut_{H}^{\Z/n}(E)\ra \bar{H} \ra 1.
$$
Since $c_h$ only depends on the action of $h$ on $A$, we can extend the notation so that
for $\sigma\in\bar{H}$ we write $c_\sigma:=c_h$ for any $h\in H$ with the same action on $A$ as $\sigma$. 
Further, for $g\in G$, we let $c_g:=c_{\bar{g}}$, where $\bar{g}$ is the image of $g$ in $H$ (or $\bar{H}$).
We have a set-theoretic section $\bar{H}\ra  \Aut_{H}^{\Z/n}(E)$, taking $\sigma\in\bar{H}$ to 
the map taking $(j,a)\mapsto (j+  c_{\sigma}(\sigma(a))  ,\sigma(a))$
for all $(j,a)\in E$.  We can check this map is an automorphism of $E$ using the definition of $c_h$.
\details{\\
\begin{align*}
(j,a)&\mapsto (j+  c_{\sigma}(\sigma(a))  ,\sigma(a))\\
(k,b)&\mapsto (k+  c_{\sigma}(\sigma(b))  ,\sigma(b))\\
(j+k+f(a,b),a+b)&\mapsto (j+k+f(a,b)+  c_{\sigma}(\sigma(a+b))  ,\sigma(a+b))\\
\end{align*}
and
$$
(j+  c_{\sigma}(\sigma(a))  ,\sigma(a))(k+  c_{\sigma}(\sigma(b))  ,\sigma(b))=
(j+  c_{\sigma}(\sigma(a))+k+  c_{\sigma}(\sigma(b)) +f(\sigma(a),\sigma(b))  ,\sigma(a+b)).
$$
Into $f(h^{-1}(a),h^{-1}(b))-f(a,b)=c_h(a)+c_h(b)-c_h(a+b)$ we plug in $h(a)$ for $a$ and $h(b)$ for $b$ to obtain
$f(a,b)-f(h(a),h(b))=c_h(h(a))+c_h(h(b))-c_h(h(a+b))$, which exactly shows that the above map is an automorphism.  
\\
}
This allow us to write  $\Aut_{H}^{\Z/n}(E)$ as a set theoretic product $\Hom(A,\Z/n)\times \bar{H}$, and for $
\phi\in \Hom(A,\Z/n)$ and $\sigma\in \Aut_{H}$, we write $[\phi,\sigma]$ for the element of $\Aut_{H}^{\Z/n}(E)$ such that
$$
[\phi,\sigma] (j,a)=(j+ \phi(\sigma(a))+ c_{\sigma}(\sigma(a))  ,\sigma(a) )
$$
for all $(j,a)\in E$.  
For $\phi,\phi'\in \Hom(A,\Z/n)$ and $\sigma,\sigma'\in \Aut_{H}$, we can compute 
$[\phi,\sigma][\phi',\sigma']$ by seeing where it sends $(0,a)\in E$.  We have
\begin{align*}
[\phi,\sigma][\phi',\sigma'](0,a)&=[\phi,\sigma]( \phi'(\sigma'(a))+ c_{{\sigma'}}(\sigma'(a))  ,\sigma'(a) )\\
&=( \phi'(\sigma'(a))+ c_{{\sigma'}}(\sigma'(a)) +\phi(\sigma(\sigma'(a))) +c_{\sigma}(\sigma(\sigma'(a))) ,\sigma(\sigma'(a)) )\\
&=[\phi'\sigma^{-1}+\phi+c_{{\sigma'}}\sigma^{-1}+ c_{\sigma} -c_{{\sigma\sigma'}} ,\sigma\sigma' ](0,a),
\end{align*}
\details{\\
We check 
\begin{align*}
&(\phi'\sigma^{-1}+\phi+c_{{\sigma'}}\sigma^{-1}+ c_{\sigma} -c_{{\sigma\sigma'}} )(\sigma(\sigma'(a)))
+ c_{{\sigma\sigma' }}(\sigma\sigma' (a))\\
=&\phi'\sigma(a)+\phi\sigma(\sigma'(a))+c_{{\sigma'}}\sigma(a)+ c_{\sigma}\sigma(\sigma'(a)) -c_{{\sigma\sigma'}}(\sigma(\sigma'(a))) + c_{{\sigma\sigma' }}(\sigma\sigma' (a))
\end{align*}
\\}
and thus
$$
[\phi,\sigma][\phi',\sigma']=[\phi'\sigma^{-1}+\phi+c_{{\sigma'}}\sigma^{-1}+ c_{\sigma} -c_{{\sigma\sigma'}} ,\sigma\sigma' ].
$$ (This proves Lemma~\ref{L:ASpcon}.)

We next check how elements of $E$ act on $E$ by conjugation.
For $a\in A$ and $(0,a)\in E$, we have $(0,a)^{-1}=(-f(a,-a),-a)$ and so for $a,b\in A$, 
$$
(0,a)\cdot_E(0,b)\cdot_E(0,a)^{-1}=(f(a,b),a+b)\cdot_E (-f(a,-a),-a)=(f(a,b)-f(a,-a)+f(a+b,-a)  ,  b).
$$
By the cocycle condition and $f$ being normalized, 
$$
f(a,-a)-f(b+a,-a)-f(b,a)=0.
$$
Let $F_a(b):=f(a,b)-f(b,a)$ for all $a,b\in A$.
Then conjugation by $(0,a)$, written as an element of  $\Aut_{H}^{\Z/n}(E)$, is $[F_a,1]$.
(This proves the claim in Lemma~\ref{L:Delta} that $a\tensor b \mapsto f(b,a)-f(a,b)$ is a homomorphism.)
Since $\Z/n$ is central in $E$, conjugation by $(j,a)$ is the same as conjugation by $(0,a)$.

Huebschmann's $\Aut_G(e):= \Aut_{H}^{\Z/n}(E) \times_{\bar{H}} G$, 
is a set theoretic product $ \Hom(A,\Z/n)\times G$ 
and has elements $\{\phi,g \}$
for $\phi\in  \Hom(A,\Z/n)$ and $g\in G$, with multiplication
$$
\{\phi,g\}\{\phi',g'\}=\{\phi'g^{-1}+\phi+c_{{g'}}g^{-1}+ c_{{g}} -c_{{gg'}} ,gg' \}.
$$ This uses our earlier expression of $\Aut_{H}^{\Z/n}(E)$ as a set-theoretic product of $\Hom(A,\Z/n)$ and $H$.  The action of an element of $\Aut_G(e)$ on $E$ is given by the formula
\[ \{\phi, g\} \circ ( j,a) =  (j + \phi (g (a)) + c_g(g(a)), g(a) ) .\]

Using that $h_1=1$ and $c_1=0$, we can check that
$$
\{0,g\}^{-1}=\{-c_{{g}} g- c_{{g^{-1}}}   ,g^{-1}\}.
$$
\details{
\begin{align*}
\{-c_{{g}} g- c_{{g^{-1}}}   ,g^{-1}\}\{0,g\}&= \{ 
-c_{{g}} g- c_{{g^{-1}}}  + 
c_{{g}}g+ c_{{g^{-1}}} -c_{{1}}
 ,   1\}\\
&= \{ 
0 ,   1\}.
\end{align*}
}

Huebschmann's $\beta: E\ra  \Aut_G(e)$ sends $(j,a)\in E$ to  the element of
$\Aut_{H}^{\Z/n}(E) \times_{\bar{H}} G$ whose first coordinate is conjugation  by $(j,a)$ and whose second coordinate is $a$.
In our notation above, this means
$$
\beta(j,a)=\{F_a ,a\}.
$$
Huebschmann's $\Out_G(e):=\Aut_G(e)/\beta(E)$ is then the extension
$$
1\ra \Hom(A,\Z/n) \ra \Out_G(e) \ra H\ra 1
$$
that gives the image of $d_2^{0,2}$ according to \cite[Theorem 1]{Huebschmann1981}.

Let $s:H \ra G$ be a set theoretic section of the quotient $G\ra H$ such that for $h,i\in G$, we have $s(h)s(i)=\alpha(s,i)s(hi)$.
\details{\\
Suppose I change to $s'$ where $s'(h)=a(h)s(h)$ for some $a:H\ra A$. Then
$$
s'(h)s'(i)=a(h)s(h)a(i)s(i)=a(h) h(a(i))  s(h)s(i)=a(h) h(a(i)) \alpha(s,i) s(hi) =a(h) h(a(i)) \alpha(s,i)  a(hi)^{-1} s'(hi) 
$$
and so $\alpha'(h,i)=a(h)+h(a(i))-a(hi)+\alpha(s,i)$.
Thus I choose some section,  and if it gives $\alpha'$, since $\alpha'$ differs from $\alpha$ by $da$ for some $a:H\ra A$,
I can just change by section by that a to get $\alpha$.
\\}
Then for $h\in H$, we can lift it to $\{0,s(h) \}\beta(E)$ in $\Out_G(e)$.
Then we can compute a cocycle representing $d_2^{0,2}([f])$ by computing, for $h,i\in H$,
\begin{align*}
&\beta(E)\{0,s(h) \}\beta(E)\{0,s(i) \}\beta(E)\{0,s(hi) \}^{-1} \\
=&
\beta(E)\{c_{{i}}h^{-1}+ c_{{h}} -c_{{hi}} ,s(h)s(i) \}\{0,s(hi) \}^{-1}\\
=&\{c_{{i}}h^{-1}+ c_{{h}} -c_{{hi}} ,s(h)s(i) \}
\{-c_{{hi}} hi- c_{{(hi)^{-1}}}   ,s(hi)^{-1}\}\\
=&\beta(E)\{
-c_{{hi}} - c_{{(hi)^{-1}}}(hi)^{-1} +
c_{{i}}h^{-1}+ c_{{h}} -c_{{hi}}
+c_{{(hi)^{-1}}}(hi)^{-1}+ c_{{hi}} -c_{{\alpha(h,i)}}
,\alpha(h,i )\}
\\
=&\beta(E)\{
-c_{{hi}}  +
c_{{i}}h^{-1}+ c_{{h}} 
,\alpha(h,i )\}
\\
=&\beta(E)
\{F_{\alpha(h,i )}, \alpha(h,i )\}\{-F_{\alpha(h,i )}
-c_{{hi}}  +
c_{{i}}h^{-1}+ c_{{h}} 
,1\}\\
=&\beta(E)\{-F_{\alpha(h,i )}
-c_{{hi}}  +
c_{{i}}h^{-1}+ c_{{h}} 
,1\}.
\end{align*}
We used above that $\alpha(h,i)$ has image $1$ in $H$,  and so it acts trivially on $A$ and $c_{{\alpha(h,i)}}=0$.
We thus conclude $d_2^{0,2}=d_{2,1}^{0,2}+d_{2,2}^{0,2}$. 
 The fact that $d_2^{0,2}$ and $d_{2,2}^{0,2}$ are well-defined on classes in $H^2(A,\Z/n)^H$ (and not just cocycles)
implies the same for $d_{2,1}^{0,2}$.  The fact that $d_2^{0,2}$ and $d_{2,2}^{0,2}$ have image in 
$H^2(H,\Hom(A,\Z/n))$ implies the same for $d_{2,1}^{0,2}$.

Now suppose that $d_2^{0,2}([f])=0.$  Then we can choose a normalized cochain $H \ra \Hom(A,\Z/n)$ 
taking $h\mapsto v_h$ such that
for all $h,i\in H$, we have
$$
-v_ih^{-1}-v_h+v_{hi}=
-F_{\alpha(h,i )}
-c_{{hi}}  +
c_{{i}}h^{-1}+ c_{{h}}. 
$$
Then we have a section $\Psi:H\ra \Aut_G(e)$ taking $h\mapsto \{v_h,s(h)\}$.
We will check that the induced map $\bar{\Psi}: H\ra \Out_G(e):=\Aut_G(e)/\beta(E)$ is a homomorphism.
For $h,i\in H$, we have
\begin{align*}
\Psi(h)\Psi(i)=\{v_h,s(h)\}\{v_i,s(i)\}=\{v_h+v_ih^{-1}   +c_{i}h^{-1} +   c_{h} -c_{hi}
  ,s(h)s(i)\}
\end{align*}
whereas
\begin{align*}
\beta(0,\alpha(h,i))\Psi(hi)=\beta(0,\alpha(h,i))\{v_{hi},s(hi)\}&=\{F_{\alpha(h,i)},\alpha(h,i) \}\{v_{hi},s(hi)\}
\\
&=\{F_{\alpha(h,i)}+v_{hi}   ,\alpha(h,i)s(hi)\}.
\end{align*}
Thus $\Psi(h)\Psi(i)=\beta(0,\alpha(h,i))\Psi(hi)$ by our choice of $v$,  and $\bar{\Psi}$ is a homomorphism.

Using $\bar{\Psi}$, we can form a fiber product $\Aut_G(e)\times_{\Out_G(e)} H$, which fits in an exact sequence
$$
0\ra \Z/n\ra E \stackrel{\beta}{\ra} \Aut_G(e)\times_{\Out_G(e)} H\ra H \ra 1.
$$
The result \cite[Theorem 3]{Huebschmann1981}
shows that the above is a crossed $2$-fold extension, where $\Aut_G(e)\times_{\Out_G(e)} H$
acts on $E$ via the $\Aut_G(e)$ coordinate and its projection onto $\Aut^{\Z/n}_H(E)$. 
Further, by the correspondence between 
crossed $2$-fold extensions and $H^3(H,\Z/n\Z)$, 
\cite[Theorem 3]{Huebschmann1981} says the above
corresponds to $d_3^{0,2}([f])$.

There is a  cocycle description of the  class in $H^3(H,\Z/n)$ 
associated to any crossed $2$-fold extension
given in \cite[{Chapter 4, Section 5}]{Brown1982}.  
We use the section of $ \Aut_G(e)\times_{\Out_G(e)} H\ra H$ sending $h\mapsto (\Psi(h),h)$ as Brown's $s$.
Note that, for $h,i\in H$, in $ \Aut_G(e)\times_{\Out_G(e)} H$, we have
$$
(\Psi(h),h)(\Psi(i),i)=(\Psi(h)\Psi(i),hi)=(\beta(0,\alpha(h,i))\Psi(hi),hi)=\beta(0,\alpha(h,i))(\Psi(hi),hi)
$$
and so Brown's $F(h,i)$ can be taken to be $(0,\alpha(h,i))$.
 Then 
the cocycle $H^3\ra\Z/n$ associated to our extension above, and hence $d_3^{0,2}([f])$,  for $g,h,i\in H$ sends
\begin{align}
(g,h,i)\mapsto &(\{v_g,s(g)\} \circ (0,{\alpha}(h,i)) ) (0,{\alpha}(g,hi))  (0,{\alpha}(gh,i))^{-1} (0,{\alpha}(g,h))^{-1} \label{E:d3} \\
=&( v_g( g({\alpha}(h,i) )) +c_g(g({\alpha}(h,i))),g({\alpha}(h,i)))(0,{\alpha}(g,hi)) (0,{\alpha}(gh,i))^{-1} (0,{\alpha}(g,h))^{-1}.\notag
\end{align}
The individual terms of the above product are in $E$, but $\alpha$ being a $2$-cocycle from $H^2\ra A$ implies that the product has trivial image in $A$ and hence is in $\Z/n\Z$.  In particular, we see that if $\alpha=0$, then $d_3^{0,2}([f])=0$.

{\bf For $f$ in the image of $\Hom(A,\Z/n)^H$:}
If $f$ is the image of $\phi\in \Hom(A,\Z/n)^H$ under the map $
\Hom(A,\Z/n)^H  \ra H^2(A,\Z/n)^H
$
from Lemma~\ref{L:H2struc},  then note that we can take $c_h=0$ for all $h\in H$.
Thus $d_{2,2}^{0,2}(f)=0$.
  Further, we have $F_a=0$ for all $a\in A$.
Thus from the description of 
$d_2^{0,2}$ above we have
$d_2^{0,2}([f])=0$, and we can take $v_h=0$ for all $h\in H$.
Note an inverse formula for $E$ is given by $(j,a)^{-1}= (-j-f(a,-a),-a)$.
Then the cocycle \eqref{E:d3} sends $(g,h,i)$ to
\begin{align*}
&(0, g({\alpha}(h,i))) (0,\alpha(g,hi)) (0,{\alpha}(gh,i))^{-1}  (0,{\alpha}(g,h))^{-1} \\
=&(0, g({\alpha}(h,i))) (0,\alpha(g,hi))\\&\cdot  (-f({\alpha}(gh,i),-{\alpha}(gh,i)),-{\alpha}(gh,i))(-f({\alpha}(g,h),-{\alpha}(g,h)) ,-{\alpha}(g,h)).
\end{align*}
Let $a=g({\alpha}(h,i))$ and $b=\alpha(g,hi)$ and $c=-\alpha(gh,i)$ and $d=- \alpha(g,h)$.
Since the image of the above cocycle \eqref{E:d3} is in $\Z/n$, we have $a+b+c+d=0$.  
The image of the above cocycle \eqref{E:d3}
 has $\Z/n\Z$ coordinate  $f(a,b)+f(a+b,c)+f(a+b+c,d)-f(-c,c)-f(-d,d)$.
Since $f(x,y)=\frac{\widehat{\phi(x)} +\widehat{\phi(y)}-\widehat{\phi(x+y)}}{n}$ for some lift $k\mapsto \widehat{k}$ from $\Z/n\Z$ to $\Z/n^2\Z$ with $\widehat{0}=0$, we have
$$
f(a,b)+f(a+b,c)+f(a+b+c,d)-f(-c,c)-f(-d,d)=\frac{\widehat{\phi(a)}+\widehat{\phi(b)}-\widehat{\phi(-c)}-\widehat{\phi(-d)} }{n}.
$$
Since $\phi$ is $H$-invariant, we have $\phi(a)=\phi(\alpha(h,i))$.  Thus the cocycle associated to the crossed $2$-fold extension
is
\begin{equation}\label{E:3cocycle}
(g,h,k)\mapsto \frac{\widehat{\phi(\alpha(h,i))}+\widehat{\phi(\alpha(g,hi))}-\widehat{\phi(\alpha(gh,i))}-\widehat{\phi(- \alpha(g,h))} }{n}
\end{equation}

On the other hand, if we start with $\alpha:H^2 \ra A$, and then apply $\phi \in \Hom (A,\Z/n\Z)$ to obtain $\phi_*(\alpha)\in H^2(H,\Z/n\Z)$, and then apply the Bockstein connecting homomorphism to $\phi_*(\alpha)$, we obtain the cocycle
$$
(g,h,k)\mapsto  \frac{\widehat{\phi(\alpha(h,i))}-\widehat{\phi(\alpha(gh,i))}+\widehat{\phi(\alpha(g,hi))}-\widehat{\phi(\alpha(g,h))}}{n},
$$
agreeing with Equation~\eqref{E:3cocycle} and proving the claim about $d_3^{0,2}([f])$ for $f$ in the image of $\Hom(A,\Z/n)^H$. \end{proof}

We now relate the vanishing of $d_{2,2}^{0,2}$ to a group theoretic condition.
Recall the map from Lemma~\ref{L:Delta}, $\Delta: H^2(V,\Z/n\Z) \ra \Delta \Hom(V\tesnor V, \Z/n\Z)$.  If $V$ is an $\F_2$ vector space, then any homomorphism $V\tensor V \ra \Z/n\Z$ lands in $(n/2)\Z/n\Z\isom \F_2$, and thus we have an isomorphism  $\Delta \Hom(V\tensor V, \Z/n\Z)\isom \wedge_2 V^\vee.$  We use $\Delta$ to also denote the composite $H^2(V,\Z/n\Z) \ra \wedge_2 V^\vee$.

\begin{lemma}\label{L:PhiASp}
Let $H$ be a finite group, $V$ a finite-dimensional  representation of $H$ over $\F_2$, and let  $\omega \in (\wedge_2 V^\vee)^H$ be non-zero.
Hence $H$'s action on $V$ is through $\Sp(V)$.
Let $n$ be a power of $2$.
If $n=2$, we also assume $V$ is $\F_2$-orthogonal.
Let $\Delta^{-1}(\omega)$ be a preimage of $\omega$ in $H^2(V,\Z/n)^H$.
Then for the map 
$$
d_{2,2}^{0,2}: H^2(V,\Z/n\Z)^H \ra H^2(H,\Hom(V,\Z/n))
$$
defined above,
$d_{2,2}^{0,2}(\Delta^{-1}(\omega))=0$ if and only if the map $H\ra \Sp(V)$ factors through the affine symplectic group $\ASp(V)\ra \Sp(V)$.
\end{lemma}

\begin{remark}
If $V$ is an irreducible representation of $H$ over $\F_2$,  then for any two non-zero $\omega,\omega' \in (\wedge_2 V^\vee)^H$,
 there is an $H$-automorphism of $V$ taking $\omega$ to $\omega'$
(see the proof of  Proposition~\ref{P:dforbilinear4}).  Thus the question of whether 
$H\ra \Sp(V)$ factors through $\ASp(V)\ra \Sp(V)$ does not depend on the choice of $\omega$ for an irreducible representation $V$.
\end{remark}

\begin{proof}
Since $\im B$ is symmetric,  $d_{2,2}^{0,2}$ is trivial on $\ker \Delta=\im B$.  Thus
$d_{2,2}^{0,2}(\Delta^{-1}(\omega))$ does not depend on which preimage we take.  
If $4\mid n,$ by Lemma~\ref{L:H2struc},
we can take a choice
$[g]\in H^2(V,\Z/n)^H$ represented by a cocycle $g$ that is a bilinear form and $\Delta[g]=\omega$.
If $n=2$, then by Lemmas~\ref{L:H2struc} and \ref{L:ortho} 
we can take a choice
$[g]\in H^2(V,\Z/2)^H$ represented by a cocycle $g$ that is a bilinear form and $\Delta[g]=\omega$.
For any $n$, we can use the bilinear form $g$ to obtain a class $[g]\in H^2(V,\Z/4),$ which is necessarily $H$-invariant by Lemma~\ref{L:H2struc}, since
$\Delta[g]=\omega$.

We consider the Lyndon-Hochschild-Serre
spectral sequence for $V\ra V\rtimes H \ra H$.  By Lemma~\ref{L:d202gen}, in this case $d_2^{0,2}=d_{2,2}^{0,2}.$
By the functoriality of the spectral sequence (Lemma~\ref{pullback-differential}) in the $\Z/n$ term, we have that $d_{2,2}^{0,2}([g])\in H^2(H,V^\vee)$ does
not depend on which $n$ we take, and so we may assume $n=4$.
\details{There is a nice diagram in the comments here
}

By Lemma~\ref{L:d202gen}, 
$d_3^{0,2}=0$.
Thus  $[g]\in E_2^{0,2}= H^2(V,\Z/4)^H$ survives to  $E_\infty^{0,2}$ if and only if $d_{2,2}^{0,2}([g])=0$.  
By the edge maps of the spectral sequence,  this happens if and only if $[g]$ pulls back from $H^2(V\rtimes H,\Z/4),$
which is equivalent to the  extension $\Z/4\ra E \ra V$ associated to $[g]$ pulling back from an extension 
$\Z/4\ra \mathcal{E} \ra V\rtimes H.$  

If the extension $E$ pulls back from $\Z/4\ra \mathcal{E} \ra V\rtimes H,$ then $E$ is the preimage of $V$ in $\mathcal{E}$.
Then $H$ acts on $E$ by lifting to $\mathcal{E}$ and conjugating, which preserves $\Z/4$ pointwise and acts on $V$ through $\Sp(V)$.   
 This means the map $H\ra\Sp(V)$ factors through $\ASp(V)$.  

Conversely,  if 
$H\ra\Sp(V)$ factors through $\ASp(V)$,  then 
 the $H$ action on $V$ lifts to an action  on $E$ that fixes $\Z/4$ pointwise.  Hence  $\Z/n\ra E\rtimes H \ra V\rtimes H$ is an extension from which 
$\Z/4\ra E \ra V$ pulls back. 

In conclusion, $\Z/4\ra E \ra V$ pulls back from an extension 
$\Z/n\ra \mathcal{E} \ra V\rtimes H$ if and only if the map $H\ra \Sp(V)$ factors through $\ASp(V)\ra \Sp(V)$.
Thus $d_{2,2}^{0,2}(\Delta^{-1}(\omega))=0$ if and only if  the map $H\ra \Sp(V)$ factors through $\ASp(V)\ra \Sp(V)$.
\end{proof}

Next we show the vanishing of $d_3^{0,2}$ on classes from bilinear forms.

\begin{lemma}\label{L:d3is0}
 Let $V$ be a finite dimensional vector space over $\F_2$ and 
 and let $\omega$ be a full rank element of $\wedge_2 V^\vee$. Let $n$ be a positive integer with $4 \mid n$. 
Let $f\in\Hom(V\tensor V,\Z/n)$ be a bilinear form such that $\Delta f=\omega$ and  $[f] \in H^2(V, \mathbb Z/n)^{\Sp(V)}$.
We view $\ASp(V)$ as an extension $V\ra \ASp(V) \ra \Sp(V)$ using the isomorphism $V\ra V^\vee$ induced by $\omega$.
Then in the Lyndon-Hochschild-Serre spectral sequence computing $H^*(\ASp(V),\Z/n)$ from $H^p(\Sp(V),H^q(V,\Z/n))$ we have $d_3^{0,2}([f])=0.$
 \end{lemma}

\begin{proof}
We first reduce to the case $n=4$. The inclusion $\mathbb Z/4 \to \mathbb Z/n$ by multiplication by $\frac{n}{4}$ induces a map $H^2( A,\mathbb Z/4) \to H^2(A, \mathbb Z/n)$, which is compatible with $\Delta$  and by Lemma~\ref{L:H2struc} is an isomorphism on $\Sp(V)$ invariant classes in $\im C$.  
In particular any choice of $f'\in\Hom(V\tensor V,\Z/n)$ such that $\Delta f'=\omega$ and  $[f'] \in H^2(V, \mathbb Z/n)^{\Sp(V)}$ is $\frac{n}{4}f$, where $f\in\Hom(V\tensor V,\Z/4)$ such that $\Delta f=\omega$ and  $[f] \in H^2(V, \mathbb Z/4)^{\Sp(V)}$.
So by Lemma~\ref{pullback-differential},   $d_3^{0,2}([f'])$ in the general $n$ case is the image of $d_3^{0,2}([f])$ in the $n=4$ case under the homorphism $H^3(\Sp(V),\Z/4)\ra H^3(\Sp(V),\Z/n)$. So it suffices, for proving the lemma, to show that $d_3^{0,2}([f])=0$ in the $n=4$ case.

We will show that the class $[f]$ pulls back from a class in $H^2(\ASp(V),\Z/4)$, which will imply $d_3^{0,2}([f])=0$
by consideration of the edge maps of the spectral sequence.

In \cite[Theorem 1.3]{Gurevich2012}, Gurevich and Hadani construct a group $\operatorname{AMp}(V)$ as an extension of $\operatorname{ASp}(V)$ by $\mu_4$, giving an extension class in $H^2 ( \operatorname{ASp}(V), \mu_4)$. We will  show that this class is the desired class, when composed with either of the two isomorphisms $\mu_4 \cong \mathbb Z/4$. To begin the proof, we must explain the defining property of the group $\operatorname{AMp}(V)$ of $V$.

Let $E$ be the extension of $V$ by $\Z/4\Z$ given by $[f]$.
Let $\psi \colon \mathbb Z/4 \to \mathbb C^\times$ be a faithful character. The result \cite[Theorem 1.1]{Gurevich2012} states that there exists a unique faithful irreducible representation of $E$ with central character $\psi$, and denotes the underlying vector space of this representation by $\mathcal H$ and the homomorphism $E \to GL(\mathcal H)$ by $\pi$. Next \cite[\S1.3.2]{Gurevich2012} defines a homomorphism $\tilde{\rho} \colon \operatorname{ASp}(V) \to PGL(\mathcal H)$ by defining $\tilde{\rho} (g)$ as the unique-up-to-scalars element of $GL(\mathcal H)$ solving $\tilde{\rho}(g) \pi(h) \tilde{\rho}(g)^{-1} = \pi(g(h)) $ for all $h\in E$.

The statement of \cite[Theorem 1.3]{Gurevich2012} is that there exists a homomorphism $\rho \colon \operatorname{AMp}(V) \to GL(\mathcal H)$ lying over $\tilde{\rho}$. 

Also before giving the proof, we note that the pullback of the extension class $\alpha$ of $\operatorname{AMp}(V)$ to 
an extension class in $H^2 (V, \mu_4)$
is invariant under the conjugation action of $\operatorname{ASp}(V)/V^\vee$, which is the usual action of $\operatorname{Sp}(V)$.
Also $[f]\in H^2 (V, \Z/4)^{\Sp(V)}.$
 However, the group of $\operatorname{Sp}(V)$-invariant classes in $H^2(V, \mathbb Z/4)$ is isomorphic to $\mathbb F_2$ since,
by Lemma~\ref{L:H2struc}, 
  \[ H^2(A, \mathbb Z/4)^{\operatorname{Sp}(A)} \cong (\wedge^2 V^\vee \times V^\vee)^{\operatorname{Sp}(V)} = (\wedge^2 V^\vee)^{ \operatorname{Sp}(V)} \times (V^\vee)^{ \operatorname{Sp}(V)} \cong \mathbb F_2 \times 1 = \mathbb F_2.\]
Thus it suffices to show that the pullback of the extension class $\alpha$ to $H^2(V, \mu_4) \cong H^2(V, \mathbb Z/4)$ is not trivial. 

To show an extension class of abelian groups is not trivial, it is enough to show that two elements have a non-identity commutator, because a split extension would be abelian and so all commutators would be the identity. Since commutators are preserved by restrictions of extensions, it suffices to find two elements in the image of $V$ inside $\operatorname{ASp}(V)$ whose lifts to $\operatorname{AMp}(V)$ have nontrivial commutator.

Let $\overline{\pi} \colon V \to PGL(\mathcal H)$ be obtained from composing $\pi$ with the projection $GL(\mathcal H)\to PGL(\mathcal H)$ and using that elements in $\mathbb Z/4$ are sent to scalars by $\pi$. 

The restriction of $\tilde{\rho}$ to $V$ is simply $\overline{\pi}$, because for $g\in V$ and $g'$ a lift of $g$ to $E$,  $\pi(g')$ solves the equation $\pi(g') \pi(h) \pi(g')^{-1}=\pi(g(h))$, since the action of $g$ on $E$ is by conjugation, and thus modulo scalars we have $\tilde{\rho}(g) = \pi(g') =\overline{\pi}(g)$.

The commutator of two elements in $GL(\mathcal H)$ is invariant under multiplying each element by scalars and thus may be computed from their image in $PGL(\mathcal H)$.  Let $g$ and $h$ be two elements of $V$ on which the bilinear form $\omega$ is nontrivial, let $g'$ and $h'$ be their lifts to $\operatorname{AMp}(V)$, and let $g^\vee$ and $h^\vee$ be their lifts to $E$. Then
\[ \rho( [g',h']) = [\rho(g'), \rho(h') ] = [\tilde{\rho}(g), \tilde{\rho}(h)] =[\overline{\pi}(g),\overline{\pi}(h) ] =[\pi(g^\vee),\pi(h^\vee) ] = \pi ( [g^\vee,h^\vee] ).\]

Since $\omega([g,h])\neq 0$, we have $[g^\vee,h^\vee]$ a nonzero element of $\mathbb Z/4 \subset E$. Since $\pi$ restricted to $\mathbb Z/4$ is the faithful character $\psi$, this implies $  \pi ( [g^\vee,h^\vee] )\neq 1$ so $\rho( [g',h'])\neq 1$ and hence $[g',h']\neq 1$, as desired.\end{proof}

\begin{lemma}\label{L:d3is0-2}
 Let $V$ be a finite dimensional vector space over $\F_2$ 
 and let $q$ be an element of $\Sym^2 V^\vee$ whose associated (symplectic) bilinear form defined by $\omega(a,b)=q(a+b)-q(a)-q(b)$ is nondegenerate. Let $\operatorname{O}(V)$ be the group of automorphisms of $V$ preserving $q$, which maps to $\Sp(V)$ since automorphisms preserving the quadratic form $q$ preserve the symplectic form $\omega$. Let $\Ps(V)$ be the fiber product of $\ASp(V)$ with $\operatorname{O}(V)$ over $\Sp(V)$.
 
 Let $f\in\Hom(V\tensor V,\Z/2)$ be a bilinear form whose projection to $\Sym^2 V^\vee$ (as in Lemma~\ref{L:H2struc}) is $q$. Let $[f] \in H^2(V, \mathbb Z/2)^{O(V)}$ be the associated class.
 
We view $\Ps(V)$ as an extension $V\ra \Ps(V) \ra \operatorname{O}(V)$ using the isomorphism $V\ra V^\vee$ induced by $\omega$.
Then in the Lyndon-Hochschild-Serre spectral sequence computing $H^\vee(\Ps(V),\Z/2)$ from $H^p(\operatorname{O}(V),H^q(V,\Z/n))$ we have $d_3^{0,2}([f])=0.$
 \end{lemma}
 
\begin{proof}
We will show that the class $[f]$ pulls back from a class in $H^\vee(\Ps(V),\Z/2)$, which will imply $d_3^{0,2}([f])=0$
by consideration of the edge maps of the spectral sequence.

The group $\Ps(V)$ is given an explicit description in \cite[\S31]{Weil1964} as the group of pairs $\sigma, f$ with $\sigma \in \operatorname{O}(V)$ and $Q$ a quadratic form on $V$ satisfying $Q(a+b)-Q(a)-Q(b) = f (\sigma(a),\sigma(b))-f(a,b)$. Here the multiplication $(\sigma_1, Q_1) (\sigma_2,Q_2)=(\sigma_1\sigma_2, Q_2+ Q_1\circ \sigma_2)$. (Actually, this is the reverse order of Weil's multiplication convention -- we have switched this since Weil uses right multiplication for group actions but the switch is irrelevant since a group is isomorphic to its opposite group.) We check that this definition agrees with the definition as a fiber product. To do this, realize the extension $E$ of $V$ by $\mathbb Z/4$ associated to the bilinear form $f$ as the group of pairs $(v, t)$ with $v\in V$ and $t \in \mathbb Z_4$ under composition $(v_1,t_2)(v_2,t_2)= (v_1+v_2 , t_1+ t_2 + 2 f(v_1,v_2 ))$. Then let $(\sigma,Q)$ act on $E$ by $(\sigma, Q)(v,t)= (\sigma(v), t+ 2 Q(v))$. This is indeed an automorphism since
\[ (\sigma,Q)(v_1,t_2)  \cdot (\sigma,Q) (v_2,t_2) = (\sigma(v_1), t_1+ 2Q(v_1)) \cdot  (\sigma(v_2), t_2+ 2Q(v_2))\]\[=  (\sigma(v_1)+  \sigma(v_2) , t_1+ 2Q(v_1) + t_2 + 2Q (v_2) +2 f( \sigma(v_1),\sigma(v_2))) = ( \sigma(v_1+v_2), t_1+t_2 + 2Q(v_1+1_2) +2 f(v_1,v_2) ) \] \[= (\sigma,Q) (v_1+v_2, t_1+t_2 + 2f(v_1,v_2)) = (\sigma, Q)( (v_1,t_1)\cdot(v_2,t_2)) \]
and the map $\Ps(V) \to \Aut(E)$ is a homomorphism since 
\[ (\sigma_1,Q_1) ( ( \sigma_2,Q_2)(v,t))= (\sigma_1,Q_1) ( \sigma_2(v), t+ 2 Q_2(v)) =(\sigma_1(\sigma_2(v)), t+ 2 Q_2(v) + 2 Q_1(\sigma_2(v)))\] \[ = (\sigma_1\sigma_2, Q_2 + Q_1 \circ \sigma_2) (t,v)= ( (Q_1,\sigma_1) (Q_2,\sigma_2))(t,v) .\]
This homomorphism is an isomorphism since it is a homomorphism between extensions of $\operatorname{O}(V)$ by $V^\vee$ that induces isomorphisms on both $\operatorname{O}(V)$ and $V^\vee$.

In \cite[Proposition b.i. on p. 18]{Blasco1993}, Bl\'asco demonstrates the existence of a certain group extension $\widehat{\Ps} $ of $\Ps(V)$ by $\pm 1$, giving an extension class in $H^2(\Ps(V), \mathbb Z/2)$ (using the unique isomorphism $\mathbb Z/2\cong \pm 1$). We will show that this class is the desired class. We begin by describing the defining property of the group $\widehat{\Ps} $.

Bl\'asco~\cite[\S1.2]{Blasco1993} first defines the Heisenberg group $H(f)$ associated to the bilinear form $f$, an extension of $V$ by $\mathbb F_2$. Take $\psi$ the unique nontrivial character of $\mathbb F_2$. Then Bl\'asco defines a vector space $\mathcal V$ and homomorphism $\rho \colon  H(f) \to GL(\mathcal V)$ as the unique irreducible representation with central character $\psi$. (Bl\'asco also considers the case of a more general field, where there can be multiple choices of $\psi$. Since $\mathbb F_2$ suffices, we will drop the dependence on $\psi$ in all subsequent notation.) Next Blasco~\cite[\S2.1]{Blasco1993} defines a homomorphism $\overline{\omega} \colon \Ps(V) \to PGL(\mathcal V)$ by setting $\overline{\omega}(s)$ to be the unique-up-to-scalars isomorphism between $\rho$ and the composition of $\rho$ with the automorphism $s$ of $H(f)$, i.e. as the unique solution to the equation $\overline{\omega}(s) \rho(h) \overline{\omega}(s)^{-1} = \rho(s(h)) $. He next defines the extension $\widetilde{\Ps} $ as the fiber product of $\Ps(V)$ with $GL(\mathcal V)$ over $PGL(\mathcal V)$, thus an extension of $\Ps(V)$ by $\mathbb G_m(\mathbb C)$ that admits a homomorphism $\omega$ to $GL(\mathcal V)$ lifting $\overline{\omega}$.

Finally, \cite[Proposition b.i. on p. 18]{Blasco1993} states the existence of a subgroup  $\widehat{\Ps} $ of $\widetilde{\Ps}$ whose projection to $\Ps$ is surjective with kernel $\pm 1$.

Furthermore, before giving the proof, we observe that the pullback of the extension class of $\widehat{\Ps} $ to $V$ is invariant under the conjugation action of $\Ps(V)$, and thus invariant under the usual action of $O(V)$.  By Lemma~\ref{L:H2struc} we have  \[H^2( V, \mathbb Z/2)^{O(V)} = (\Sym^2 V^\vee)^{O(V)} \cong \mathbb F_2,\] and $[f]$ is a nontrivial class in $H^2( V, \mathbb Z/2)^{O(V)} $, hence to show that the pullback of the extension class is $[f]$ it suffices to show it is nontrivial. It suffices to show that two elements of $\widehat{\Ps}$ have nontrivial commutator. Since $\widehat{\Ps}$ is a subgroup of $\widetilde{\Ps}$ with surjective projection to $\Ps$, every element of $\widetilde{\Ps}$ is the product of an element of $\widehat{\Ps}$ with a (central) element of $\mathbb G_m$ and therefore it suffices to find two elements of $\widetilde{\Ps}$ with nontrivial commutator. Furthermore, it suffices to find two elements in the image of $V$ inside $\Ps(V)$ whose inverse images in $\widetilde{\Ps}(V)$ have nontrivial commutator.

Let $\overline{\rho} \colon V \to PGL(\mathcal V)$ be obtained by composing $\rho$ with the projection $GL(\mathcal V)\to PGL(\mathcal V)$ and using that elements of $\mathbb Z/2$ are sent to scalars by $\rho$.

The restriction of $\overline{\omega}$ to $V$ is simply $\overline{\rho}$, because for $s\in V$ and $s'$ a lift of $s$ to $H(f)$,  $\rho(s')$ solves the equation $\overline{\omega}(s) \rho(h) \overline{\omega}(s)^{-1} = \rho(s(h))$ since the action of $s$ on $H(f)$ is by conjugation, and thus modulo scalars we have $\overline{\omega}(s) =\rho(s') =\overline{\rho}(s)$.

The commutator of two elements in $GL(\mathcal H)$ is invariant under multiplying each element by scalars and thus may be computed from their image in $PGL(\mathcal H)$.  Let $g$ and $h$ be two elements of $V$ on which the bilinear form $\omega$ is nontrivial, let $g'$ and $h'$ be their lifts $\widetilde{\Ps}$, and let $g^\vee$ and $h^\vee$ be their lifts to $H(f)$. Then
\[ \omega( [g',h']) = [\omega(g'), \omega(h') ] = [\overline{\omega}(g), \overline{\omega}(h)] =[\overline{\rho}(g),\overline{\rho}(h) ] =[\rho(g^\vee),\rho(h^\vee) ] = \rho ( [g^\vee,h^\vee] ).\]

Since $\omega(g,h)\neq 0$, we have $[g^\vee,h^\vee]$ a nonzero element of $\mathbb Z/2 \subset H(f)$. Since $\rho$ restricted to $\mathbb Z/2$ is the faithful character $\psi$, this implies $\rho ( [g^\vee,h^\vee] )\neq 1$ so $\omega( [g',h'])\neq 1$ and hence $[g',h']\neq 1$, as desired

\end{proof}

Finally, we will prove the key information we need about $E_3^{0,2}$. 

\begin{definition}\label{D:im}
Let $V$ be an irreducible representation of a finite group $H$ over $\F_p$ for some prime $p$.  
Let $\kappa=\End_H(V)$.  
For $\alpha=(\alpha_1,\dots,\alpha_e)\in 
H^2(H,V^e)=H^2(H,V)\tensor \kappa^e, $ we let $\im \alpha\sub \kappa^e$ be the image of the corresponding map $\Hom_\kappa(H^2(H,V),\kappa) \ra\kappa^e$, or equivalently 
the $(a_1,\dots,a_e)\in\kappa^e$ such that for every $(k_1,\dots,k_e)\in \kappa^e$ with $\sum_i k_i\alpha_i=0$, we have
$\sum_i k_ia_i=0,$ or equivalently, 
the minimal $\kappa$-subspace $W$ of $\kappa^e$ such that the image of $\alpha$ is $0$ in $H^2(H,V^e) \tensor (\kappa^e/W)$.  
\details{
Let $e_i$ be a standard basis of $\kappa^e$ so that $\alpha=\sum_i \alpha_i \tensor e_i$.
If $\pi\in \Hom_\kappa(H^2(H,V),\kappa)$, then $\alpha$ sends it to $
\sum_i \pi_i(\alpha_i)  e_i=
(\pi(\alpha_1),\dots,\pi(\alpha_e))\in \kappa^e$.  
Suppose $\sum_i k_i\alpha_i=0$.  Then $\sum_i k_i\pi(\alpha_i)=0$.  So the image of $\alpha$ is a subset of these $(a_1,\dots,a_n).$
On the other hand, suppose that we have such an $(a_1,\dots,a_n)$.
Consider any linear form that is $0$ on the image, given by $j_1,\dots, j_e$ such that
$\sum_i j_i\pi(\alpha_i)=0$ for all $\pi$.
Then $\pi( \sum_i j_i\alpha_i)=0$ for all $\pi$ and hence  $\sum_i j_i\alpha_i=0$.  Thus $\sum_i j_i a_i=0$ (by assumption on the $a_i$, and any linear form $0$ on the image is $0$ on $(a_1,\dots,a_n)$., which implies these $(a_1,\dots,a_n)$ are in the image.
} We let $\sspan(\alpha)$ be the $\kappa$-linear span of the $\alpha_i$ in $H^2(H,V)$, or equivalently the 
image of the  map $\kappa^e \ra H^2(H,V)$ corresponding to $\alpha$.
\end{definition}

\begin{remark}\label{R:gete10} 
Let $\tau: V\ra V'$ be an isomorphism of representations of $G$ over $\F_p$.  Then $V'$ inherits a natural $\kappa$-action such that $\tau$ is an isomorphism over $\kappa$.  
Let $M\in \Hom(\kappa^e,\kappa^{e'})$.  Then $M$ with coordinates $m_{ij}$ gives a map of representations $M_\tau: V^e\ra (V')^{e'}$
sending the $j$th copy of $V$ to the $i$th copy of $V'$ via $m_{ij}\tau$ for all $i,j$.
This induces a map
$(M_\tau)_*:H^2(H, V^e)\ra H^2(H, (V')^{e'})$.  We can check that for $\alpha\in H^2(H,V^e)$, we have that $(M_\tau)_*(\alpha)=0$ if and only if $M(\im \alpha)=0$.  
\details{

Let $e_i$ be a standard basis of $\kappa^e$.  Then a general element of $H^2(H,V^e)$ is of the form $\alpha=\sum_j \alpha_j\tensor e_j$ for $\alpha_j\in H^2(H,V)$ and 
$$(M_\tau)_* (\sum_j \alpha_j\tensor e_j)=\sum_{i} \sum_j m_{ij}\tau_*(\alpha_j)  e_i =
\sum_{i} \tau_*( \sum_j m_{ij} \alpha_j)  e_i. 
$$ 
If $(M_\tau)_*(\alpha)=0$, then $\sum_j m_{ij} \alpha_j=0$ for all $i$, which means each row of $M$ is among the $(a_1,\dots,a_e)$ from the
detail above,  which means everything in $\im \alpha$ pairs to $0$ with them.   We expand on the above and note that the only $(j_1,\dots,j_e) $ that pair to $0$ with $\im \alpha$ are those that pair to $0$ with $\alpha$ (This is because $(W^\perp)^\perp=W$). If $M(\im \alpha)=0$, then all the rows of $M$ are among those $(k_1,\dots,k_e)$ that pair to $0$ with $\alpha$, and hence $(M_\tau)_* (\alpha)=0.$

}
\end{remark}

\begin{proposition}\label{P:dforbilinear4}
Let $p$ be a prime and $n$ be a positive integer divisible by $p$.  Let $e$ be a positive integer.
Let $V$ be an irreducible self-dual representation of a finite group $H$ over $\F_p$. 
Let $\kappa=\End_H(V)$.   Let $W=V^e$.
Let $1\ra W \ra G \ra H \ra 1$ be a short exact sequence of groups compatible with the given $H$ action on $W$.
 Let $G$ have extension class
$\alpha\in H^2(H,W)$  and let $r=\dim_{\kappa} \im \alpha$.
Consider the Lyndon-Hochschild-Serre spectral sequence computing $H^\vee(G,\Z/n)$ from $H^p(H,H^q(W,\Z/n))$.
If $V$ is non-anomalous,
then
$$
|\ker(d_2^{0,2})|=|V^H|^e|\kappa|^{\frac{(e-r)(e-r-\epsilon_V)}{2}}.
$$

If $V$ is intermediate (and if $4\nmid n$, also $\F_2$-orthogonal),  let $f\in H^2(V,\Z/n)^H$ with $\omega=\Delta f\in \wedge_2 V^\vee\sub \Hom(V,V^\vee)$, 
where $\omega\ne 0$,
and let $\Phi=d_{2,2}^{0,2}(f)\in H^2(H,V^\vee)$, which gives  $\omega_*^{-1}(\Phi)\in H^2(H,V)$.
Then 
we have
$$ 
|\ker(d_2^{0,2})|=
\begin{cases}
2|\kappa|^{\frac{(e-r)(e-r+1)}{2}} & \textrm{if }  \omega_*^{-1}(\Phi)\in\operatorname{span}(\alpha)\\
|\kappa|^{\frac{(e-r)(e-r-1)}{2}} & \textrm{if }  \omega_*^{-1}(\Phi)\not\in\operatorname{span}(\alpha).
\end{cases}
$$

Using the map $B$ from Lemma~\ref{L:H2struc},
$$
d_3^{0,2}( E_3^{0,2} )=d_3^{0,2}(B(\Hom(W,\Z/n )^H )).
$$
\end{proposition}

\begin{proof}
Since the spectral sequence entirely factors according to a factorization of $\Z/n\Z$ into cyclic groups of prime power order,  we can reduce to the case that $n$ is a power of $p$.

We change basis, relabelling $V^e$ so that $\im \alpha$ is the $\kappa^r\sub \kappa^e$ corresponding to the last $r$ coordinates.
Let $d=\dim_{\F_2} \kappa$.
We will identify $\Hom(W,\Z/n\Z)$ with $W^\vee$, and $\Hom(W\tensor W,\Z/n\Z)$ with $(W\tensor W)^\vee,$
which we identify with $W^\vee\tensor W^\vee$.
Then we identify the corresponding subgroups $\Delta \Hom(W\tensor W,\Z/n\Z)$ and $\wedge_2 W^\vee$.
We also identify $W^\vee\tensor W^\vee $ with $\Hom(W,W^\vee)$, where $u\tensor v$ is identified with the map that takes $w$ to $u(w)v$.

Lemma~\ref{L:d202gen} tells us that $d_2^{0,2}$, and indeed $d_{2,1}^{0,2}$ and $d_{2,2}^{0,2}$,
factor through $\Delta: H^2(W,\Z/n)^H \ra (\wedge_2 W^\vee)^H.$
We will next identify $(\wedge_2 W^\vee)^H$ with a certain set of matrices over $\kappa$, and then
compute $d_2^{0,2}$ in terms of these matrices.

{\bf I. Identification of $(\wedge_2 W^\vee)^H$ with matrices over $\kappa$:}

We take a $\omega \in (V^\vee\tensor V^\vee)^H$ as in Lemma~\ref{L:typeV}, and take $\omega$ in the subspace
$\wedge_2 V^\vee$ if possible (i.e. unless $V$ is symmetric).
The elements of $(W^\vee \tensor W^\vee)^H$ correspond to $e\times e $ matrices over $\kappa$ as follows.
A matrix $M$ with entries $m_{ij}\in\kappa$ corresponds to $\underline{m}:=\sum_{ij} (m_{ij}\tensor 1)\omega_{ij}\in (W^\vee \tensor W^\vee)^H$, where $\omega_{ij}$ is the copy of $\omega$ from the $j$th copy of $V^\vee$ on the left tensored with the $i$th copy of $V^\vee$ on the right.  As a map $W\ra W^\vee$, the element $\omega_{ij}$ maps the $j$th copy of $V$ to the $i$th copy of $V^\vee$ via $\omega:V\ra V^\vee$.
The switching factors action on $W^\vee \tensor W^\vee$ sends
$M\in M_{e\times e}(\kappa)$ to $\lambda \sigma(M^t).$

%

We can choose an $\F_p$ basis of  $V^\vee$ and a corresponding basis $w_i$ of $W^\vee$, and then 
$\wedge_2 W^\vee$ is the set of elements of the form $\sum_{i<j} a_{ij}(w_i\tensor w_j-w_j\tensor w_i)$, 
i.e. those that are alternating as matrices.  These matrices with $a_{ij}$ entries are block matrices, where each
$\dim_{\F_p} V\times \dim_{\F_p} V$ block corresponds to one entry of our matrix $M$ above.

If $V$ is symplectic, then by Lemma~\ref{L:typeV}, we have that $\kappa\omega=(V^\vee\tensor V^\vee)^H$, and since
$\wedge_2 V^\vee$ has the same size, we have $\kappa\omega=\wedge_2 V^\vee.$  Thus all of $\kappa\omega$
give alternating matrices in the $a_{ij}$ entries.
Also, $\lambda=-1$.
 Thus an
 an element $\underline{m}\in (W^\vee \tensor W^\vee)^H$ is in $(\wedge_2 W^\vee)^H$ if and only if $m_{ij}=m_{ji}$ for all $i,j$, or equivalently $M$ is symmetric.

 If $V$ is symmetric,  
by Lemma~\ref{L:typeV}, we have that  $\wedge_2 V^\vee=0$ and $\lambda=1$.
 An element $\underline{m}\in (W^\vee \tensor W^\vee)^H$ is in $(\wedge_2 W^\vee)^H$ if and only if $m_{ij}=-m_{ji}$ for all $i,j$ and $m_{ii}=0$ for all $i$, or equivalently corresponds to an alternating matrix $M$.
 
 If $V$ is unitary,  
we claim that for $k\in \kappa$, the element $(k\tensor k)\omega\in (\wedge_2 V^\vee)^H$.  
The element $(k\tensor k)\omega$ is $H$ invariant since $\kappa$ respects the $H$ action.  It is
in $\wedge_2 V^\vee$ since if $\omega=\Delta \psi$ for some $\psi\in V^\vee\tensor V^\vee$,
 then $(k\tensor k)\omega=\Delta(k\tensor k)\psi.$   Note that $(k\tensor k)\omega=(k\sigma(k)\tensor 1)\omega$.
 The elements in $\kappa$ of the form $k\sigma(k)$ for some $k\in\kappa$ are precisely the elements of the fixed field $\kappa^\sigma$.
 \details{we can check that every element of the fixed field $\kappa^\sigma$ is of the form
 $k\sigma(k)$ for some $k\in\kappa$.
 $\sigma(x)=x$ means that $x$ is a $p^{d/2}-1$th root of unity.  In a cyclic group of order $p^d-1$, this means that $x$ is 
 a $p^{d/2}+1$ power, i.e. $k\sigma(k)$.} 
 Thus,  by
 Lemma~\ref{L:typeV}, the elements of $(\wedge_2 V^\vee)^H$ are precisely those of the form
 $(m\tensor 1)\omega$ for $m\in \kappa^\sigma$.
 If $p\ne 2$, the fact that $\omega\in \wedge_2 V$ implies it is in the $-1$ eigenspace of the transpose action, and hence $\lambda=-1$.
 Thus for all $p$, we have $\lambda=-1$ in the unitary case.
 So an element $\underline{m}\in (W^\vee \tensor W^\vee)^H$ is in $(\wedge_2 W^\vee)^H$ if and only if $m_{ij}=\sigma(m_{ji})$ for all $i,j$ and $m_{ii}\in \kappa^\sigma$ for all $i$. 

When $n=2$ and $V$ is non-$\F_2$-orthogonal, then $H^2(W,\Z/n)^H=(\Sym^2 W^\vee)^H$ by Lemma~\ref{L:H2struc}, and
$(\Sym^2 W^\vee)^H$ is spanned by the image of $(V_i^\vee\tensor V_j^\vee)^H$ for $i<j$, where $V_k^\vee$ is the $k$th copy of $V^\vee$ in $W^\vee.$
From this it follows that  the image of $\Delta: H^2(W,\Z/n)^H\ra (\wedge^2 W^\vee)^H$ correspond precisely to the matrices $M$ with $m_{ii}=0$
and $m_{ij}=\sigma(m_{ji})$ for all $i,j$.

{\bf II. Computation of $d_2^{0,2}$ in terms of the matrix $M$:}

Now we will compute $d_2^{0,2}$
of  an element with a given image  $\underline{m}\in \wedge^2 W^\vee$ under $\Delta$,
 in terms of the matrix $M$ with entries $m_{ij}$.  Suppressing various maps, we call this $d_2^{0,2}(M)$,
 and the corresponding components $d_{2,1}^{0,2}(M)$ and $d_{2,2}^{0,2}(M)$.
Let $M_*:H^2(H,V)\ra H^2(H,V^\vee)$ be the map induced by $M:V \ra V^\vee$. 

For any $M$ with $m_{ii}=0$ for all $i$, we claim $d_{2,2}^{0,2}(M)=0$.
We can take $g=\sum_{i<j} (m_{ij}\tensor 1)\omega_{ij}\in (W^\vee \tensor W^\vee)^H$
as an $H$-invariant bilinear form such that $\Delta g=\underline{m}$
\details{
$$g-g^t=
\sum_{i<j} (m_{ij}\tensor 1)\omega_{ij}-(1\tensor m_{ij})\omega^t_{ji}
=\sum_{i<j} (m_{ij}\tensor 1)\omega_{ij}-( \lambda\sigma(m_{ij})\tensor 1)\omega_{ji}
=\sum_{i<j} (m_{ij}\tensor 1)\omega_{ij}+( \sigma(m_{ji})\tensor 1)\omega_{ji}$$
}
So for $[g]\in H^2(W,\Z/n)^H$ 
when computing $d_{2,2}^{0,2}([g])$,  we can take $c$ to be $0$ and thus $d_{2,2}^{0,2}([g])=0.$
Then $d_{2,2}^{0,2} (M)$ is a function only of the diagonal entries of $M$.

{\bf II.A $p$ odd: }If $p$ is odd, then we claim $d_{2,2}^{0,2}$ is $0$.  For any $f\in (\wedge_2 W^\vee)^H$, we can take
$f/2$ as an $H$-invariant bilinear form such that $\Delta (f/2)=f$.  So for $[f/2]\in H^2(W,\Z/n)^H$ 
when computing $d_{2,2}^{0,2}([f/2])$,  we can take $c$ to be $0$ and thus $d_{2,2}^{0,2}([f/2])=0.$
We then have $d_2^{0,2} (M)= M_*(\alpha)$ from Lemma~\ref{L:d202gen}.

{\bf II.B $p=2$: }
If $p=2$,  when either $4\mid n$ and $V$ is non-symmetric, or $n=2$ and $V$ is $\F_2$-orthogonal,  
using Lemmas~\ref{L:H2struc} and \ref{L:ortho},
we let $\Delta^{-1}(\omega)$ be some preimage of $\omega$ in $H^2(V,\Z/n)^H$, and let
$\Phi:=d_{2,2}^{0,2}(\Delta^{-1}(\omega))\in H^2(H,V^\vee)$ (in the setting $e=1$).
Note $d_{2,2}^{0,2}$ is the same as the differential $d_{2}^{0,2}$ in the $\alpha=0$ case.
From the functoriality of the differentials in the $\alpha=0$ case,
  for $k\in\kappa$, we have $d_{2,2}^{0,2}(\Delta^{-1}((k\sigma(k)\tensor  1)\omega))=d_{2,2}^{0,2}(\Delta^{-1}((k\tensor k)\omega))=k\Phi$ (when $V$ is  $\F_2$-orthogonal if $n=2$). 
If $\sigma\ne 1$,  there are $k\ne 1$ such that $k\sigma(k)=1$, i.e. any $2^{d/2}+1$ root of unity, for which we have $k\Phi=\Phi$.
 So when $V$ is unitary, we have $\Phi=0$. If $\sigma=1$ then applying functoriality in differentials to the projection to the $i$'th copy of $V$ followed by multiplication by $\sqrt{m_ii}$ gives
 \begin{equation}\label{sqrt-formula} d_{2,2}^{0,2}(M) = \sum_i \sqrt{m_{ii}} p_i( \Phi) .\end{equation}
If either $4\mid n$ and $V$ is symmetric or $n=2$ and $V$ is non-$\F_2$-orthogonal, let $\Phi=0$.
  Let $p_i:H^2(H,V^\vee)\ra H^2(H,W^\vee)$ be the map induced by the inclusion of the $i$th factor of $V^\vee$ into $W^\vee$. We see \eqref{sqrt-formula} holds in this case and the $\sigma\neq 1$ case as well since both sides are zero. Then Lemma~\ref{L:d202gen} tells us that
$$
d_2^{0,2} (M)= M_*(\alpha)+\sum_i \sqrt{m_{ii}} p_i( \Phi).
$$

{\bf III. Counting matrices in the kernel of $d_2^{0,2}$:}

Now, given these explicit expressions for $d_2^{0,2}$, we will count the number of matrices $M$ that correspond to elements of the image of $\Delta: H^2(W,\Z/n)^H\ra
(\wedge^\vee W)^H$
whose preimages in $H^2(W,\Z/n)^H$ are in the kernel of $d_2^{0,2}$.  If $n$ is odd, $4\mid n$, or $V$ is $\F_2$-orthogonal then
$\Delta: H^2(W,\Z/n)^H\ra (\wedge^\vee W)^H$ is surjective by Lemmas~\ref{L:H2struc} and \ref{L:ortho}.  (Otherwise we have already characterized above which $M$ correspond to elements in the image of $\Delta: H^2(W,\Z/n)^H\ra
(\wedge^\vee W)^H$.)

{\bf III.A $p$ is odd or $p=2$ and $\Phi=0$:}
We first consider the case when either $p$ is odd or $p=2$ and $\Phi=0$.
Then $M$ maps the $j$th copy of $V$ to the $i$th copy of $V^\vee$ via the map $m_{ij} \omega$ for all $i,j$. 
By  Remark~\ref{R:gete10},  
 $d_2^{0,2}(M)=0$ if and only if $\im \alpha\sub \kappa^e$ is in the kernel of $M$ (for matrix multiplication).  
  We have
$ d_2^{0,2}(M)=0$ exactly if $M$ has $0$'s in the last $r$ columns.
So, the number of $M$ representing elements of the image of  $\Delta: H^2(W,\Z/n)^H\ra
(\wedge^\vee W)^H$ such that $\ker d_2^{0,2}(M)=0$ is 
as given in the below table.
%

First in the case $p\neq 2$ we have

\begin{center}
\begin{tabular}{ |c|c|c| } 
 \hline
 $V$ & $M$ from & $\#\{ M|\ker d_2^{0,2}(M)=0\}$  \\\hline \hline
symplectic & symmetric matrices & $|\kappa|^{\frac{(e-r)(e-r+1)}{2}}$ \\\hline
 symmetric & alternating matrices & $|\kappa|^{\frac{(e-r)(e-r-1)}{2}}$\\\hline
  unitary & $\sigma(M^t)=M$ & $|\kappa|^{\frac{(e-r)(e-r)}{2}}$\\\hline
\end{tabular}
\end{center}

Next in the case $n=2$ we have

\begin{center}
\begin{tabular}{ |c|c|c| } 
 \hline
  $V$ & $M$ from & $\#\{ M|\ker d_2^{0,2}(M)=0\}$  \\\hline \hline
 $\F_2$-orthogonal and symplectic & symmetric matrices & $|\kappa|^{\frac{(e-r)(e-r+1)}{2}}$  \\\hline
 trivial  & alternating matrices & $|\kappa|^{\frac{(e-r)(e-r-1)}{2}}$\\\hline
 non-$\F_2$-orthogonal  & $m_{ii}=0$ and $m_{ij}=\sigma(m_{ji})$ & $|\kappa|^{\frac{(e-r)(e-r-1)}{2}}$\\\hline
 $\F_2$-orthogonal and unitary & $\sigma(M^t)=M$ & $|\kappa|^{\frac{(e-r)(e-r)}{2}}$\\\hline
\end{tabular}
\end{center}

Finally in the case $4\mid n$ we have

\begin{center}
\begin{tabular}{ |c|c|c| } 
 \hline
  $V$ & $M$ from & $\#\{ M|\ker d_2^{0,2}(M)=0\}$  \\\hline \hline
symplectic & symmetric matrices & $|\kappa|^{\frac{(e-r)(e-r+1)}{2}}$  \\\hline
 trivial  & alternating matrices & $|\kappa|^{\frac{(e-r)(e-r-1)}{2}}$\\\hline
unitary & $\sigma(M)=M$ & $|\kappa|^{\frac{(e-r)(e-r)}{2}}$\\\hline
\end{tabular}
\end{center}

We use the analysis above and
 $d_2^{0,2}((\im B)^H)=0$ to conclude that $|\ker d_2^{0,2}|$
is $|(\im B)^H|$ times the number of matrices $m$ counted above.  
We conclude the result in the cases that $p$ is odd or $p=2$ and $\Phi=0$, using Lemma~\ref{L:PhiASp} to see that when $p=2$ and $V$ is symplectic,
$\Phi=0$ if and only if $V$ is A-symplectic and the definition of $\epsilon_V$ in each case from \ref{s-notation}. 



{\bf III.B $p=2$ and $\Phi\ne0$:}
Now we consider the case when $p=2$ and $\Phi\ne 0$. 
So we have $\sigma=1$ and $V$ is non-trivial.
If we write $H^2(H,W)$ as column vectors $H^2(H,V)^e$, and similarly for $W^\vee$, then we have 
$$\alpha=\begin{bmatrix}\alpha_1\\\vdots\\\alpha_e\end{bmatrix},$$ where $\alpha_i\in H^2(H,V)$ and $\alpha_1=\dots=\alpha_{e-r}=0$, and $\alpha_{e-r+1},\dots\alpha_e$ are $\kappa$-linearly independent.
We write $\omega_*$ for the isomorphism $H^2(H,V) \ra H^2(H,V^\vee)$ induced by the isomorphism $\omega: V\ra V^\vee$.
Let $\beta_i=\omega_*(\alpha_i)$.
Note for $m \in \kappa$, the isomorphism ($m\tensor 1)\omega: V\ra V^\vee$ is the map that takes $v$ to $\omega(mv)=\sigma(m)\omega(v)$.
So for the matrix $M$ corresponding to $\underline{m}\in(\wedge_2 W^\vee)^H$,
$$
d_2^{0,2}(M)=\sigma(M)\begin{bmatrix}\beta_1\\  \vdots\\ \beta_e  \end{bmatrix}+
\begin{bmatrix}\sqrt{m_{11}}\Phi\\ \vdots\\ \sqrt{m_{ee}}\Phi  \end{bmatrix}=
\begin{bmatrix}
\sum_j \sigma(m_{1j}) \beta_j + \sqrt{m_{11}}\Phi\\
\vdots
\\
\sum_j \sigma(m_{ej}) \beta_j + \sqrt{m_{ee}}\Phi
\end{bmatrix}
.
$$
We can relabel $W$ as $V^e$ in a different way so that, if $\Phi$ is in the span of the $\beta_i$, then
$\beta_e=\Phi$.  

If  $\Phi\ne0$ is in the $\kappa$-span of the $\beta_i$, then we may change basis so that $\Phi = \beta_e$. With this convention,
$d_2^{0,2}(M)=0$ if and only if 
$m_{ie}=\sqrt{m_{ii}}$ for $1\leq i \leq e$ and $m_{ij}=0$ for $1\leq i \leq e $ and $e-r+1\leq j \leq e-1.$
Such an $M$ (that also corresponds to an element of $\wedge_2 W^\vee$) is precisely parameterized by a choice of $m_{ii}\in\kappa$ for $i\leq e-r$, and a choice $m_{ee}\in \{0,1\}$, and 
any choice of $m_{ij}$ for $i<j\leq e-r$.

Thus, the number of such $M$ and the size of $|\ker d_2^{0,2}|$ is $2|\kappa|^{\frac{(e-r)(e-r+1)}{2}}$.

If  $\Phi$ is not in the $\kappa$-span of the $\beta_i$, then
$d_2^{0,2}(M)=0$ if and only if $m_{ii}=0$ for all $i$ and $m_{ij}=0$ for $1\leq i \leq e $ and $r-e+1\leq j.$
Such $M$ (that also corresponds to an element of $\wedge_2 W^\vee$) are given by any choice of $m_{ij}$ for $i<j\leq e-r$.
The number of such $M$  and the size of $|\ker d_2^{0,2}|$ is $|\kappa|^{\frac{(e-r)(e-r-1)}{2}}$.

{\bf IV. Analysis of $d_3^{0,2}$:}

Now we consider the map 
$$d_3^{0,2}: \ker d_2^{0,2} \ra H^3(H,\Z/n).$$
Since it is clear that
$d_3^{0,2}(B(\Hom(W,\Z/n )^H ))\sub d_3^{0,2}( E_3^{0,2} )$, we only need to show that
$ d_3^{0,2}( E_3^{0,2} )\sub d_3^{0,2}(B(\Hom(W,\Z/n )^H )).$

Our relabelling of the factors of $W$, means that if we write $\alpha=(\alpha_1,\dots,\alpha_e)\in H^2(H,V^e)$,
we have $\alpha_1=\dots=\alpha_{r-e}=0$.  We take cocycles representing the $\alpha_i$ and write $\alpha_i$ for those cocycles as well, taking $\alpha_1=\dots=\alpha_{r-e}=0$.   
We can write the group $G$ as a set theoretic product $V^e\times H$ where for $w,w'\in W$ and $h,h'\in H$
$$
(w,h) \cdot_G (w',h')=(w+h(w')+\alpha(h,h'), hh').
$$
In this notation, there is a subgroup $G'$ of $G$ whose $W$ coordinate is trivial in the first $r-e$ $V$ coordinates.  
This subgroup $G'$ maps isomorphically onto the quotient of $G$ by the first $r-e$ factors of $V$, so we write $G'$ for this quotient as well.
 Let $(E')_r^{p,q}, (d')_r^{p,q}$ be the 
Lyndon-Hochschild-Serre spectral sequence computing $H^\vee(G',\Z/n)$ from $H^p(H,H^q(V^r,\Z/n))$.
We consider the maps $d_2^{1,1}, (d')_2^{1,1}$ which both map into $(E')_2^{3,0}=E_2^{3,0}=H^3(H,\Z/n)$.
The functoriality of the spectral sequence (Lemma~\ref{pullback-differential}) applied to the quotient $G\ra G'$
gives $\im (d')_2^{1,1} \sub \im d_2^{1,1}$, and applied to the map $G'\ra G$
gives $\im d_2^{1,1} \sub \im (d')_2^{1,1}.$  Thus $\im (d')_2^{1,1} = \im d_2^{1,1}$ and $(E')_3^{3,0}=E_3^{3,0}$.
Similarly, $\im (d')_3^{0,2} = \im d_3^{0,2}.$

{\bf IV.A $V$ is non-intermediate, or $V$ is intermediate and $\Phi$ is not in the span of the $\beta_i$, or $V$ is intermediate and non-$\mathbb F_2$-orthogonal and $n=2$:}
Now we consider the case when either $V$ is non-intermediate, or $V$ is intermediate and $\Phi$ is not in the span of the $\beta_i$.
We saw above that, in these cases, 
$(E')_3^{0,2}=\ker((d')_2^{0,2})=B(\Hom(V^r,\Z/n)^H)$.  
So $d_3^{0,2}( E_3^{0,2} )=(d')_3^{0,2}( E_3^{0,2} )=(d')_3^{0,2}(B(\Hom(V^r,\Z/n)^H))$, 
and by functoriality, this is equal to $d_3^{0,2}\circ B$ applied to the homomorphisms in $\Hom(V^e,\Z/n)^H$
that factor through the  last $r$ factors of $V$.  
Thus $d_3^{0,2}( E_3^{0,2} )\sub d_3^{0,2}(B(\Hom(V^e,\Z/n)^H))$.

{\bf IV.B $V$ is intermediate
and $\Phi$ is in the span of the $\beta_i$:}
Finally, when $V$ is intermediate
and $\Phi$ is in the span of the $\beta_i$,
let us first consider the case that $r=e=1$.  
Let  $\Sp(V) $ be the group of automorphisms of $V$ as a group (or $\F_2$-vector space) that preserve $\omega$. If $n=2$, so that $V$ is $\mathbb F_2$-orthogonal, let $\SO(V)$ be the group of automorphisms of $V$ as a group that preserve the unique $H$-invariant quadratic form lifting $\omega$. Let $\Aut_2(V)=\Sp(V)$ in the case $4\mid n$ and $\Aut(V)=\SO(V)$ in the case $n=2$.
Note that $d_{2,2}^{0,2}$ does not depend on the extension class of $G$, but only on $W$,  $H$, and the $H$ action on $W$.

We pick a preimage $\Delta^{-1} \omega\in H^2(V,\Z/n)^{\Aut_2(V)}$ represented by a bilinear form
 (which exists by  \eqref{L:H2struc}), and
use this same element as a choice of a preimage $\Delta^{-1} \omega\in H^2(V,\Z/n)^H$ for any $H$ acting on $V$ through $\Aut_2(V)$.
We consider $d_{2,2}^{0,2} (\Delta^{-1} \omega)$ in the case $H=\Aut_2(V)$, and see that the choice of $c$ we take here 
pulls back to a choice of $c$ for a general $H$ that acts on $V$ through $\Aut_2(V)$.  Thus the element
$\Phi=d_{2,2}^{0,2}(\Delta^{-1}\omega)\in H^2(\Aut_2(V),V^\vee)$ (for the spectral sequence when $H=\Aut_2(V))$ pulls back to 
$\Phi=d_{2,2}^{0,2}(\Delta^{-1}\omega)\in H^2(\Aut_2(V),V^\vee)$ (for the spectral sequence for a general $H$ acting through $\Aut_2(V)$).

When $H=\Aut_2(V)$, 
Lemma~\ref{L:ASpcon} implies that 
the class $\Phi=d_{2,2}^{0,2}(\Delta^{-1}\omega)\in H^2(\Aut_2(V),\Hom(V,\Z/n))$
is the extension class of $\ASp(V)\ra\Sp(V)$ (pulled back from $\Sp(V)$ to $\SO(V)$ in the case $n=2$).  The class $\alpha=\omega_*^{-1}(\Phi)$ (by our labelling)
is also the extension class of $\ASp(V)\ra\Sp(V)$, just with kernel labelled $V$ instead of $\Hom(V,\Z/n)$.
This implies we have a commutative diagram of extensions
$$
 \begin{tikzcd} 
 1 \arrow[r] \arrow[d] & V \arrow[r] \arrow[d,]& G \arrow[r]\arrow[d] & H \arrow[r] \arrow[d] & 1\arrow[d] \\ 
1 \arrow[r] & V\arrow[r] & \ASp(V) \arrow[r] & \Sp(V) \arrow[r] & 1.
 \end{tikzcd}
$$
By Lemmas~\ref{L:d3is0} when $H= \Aut_2(V)=\Sp(v)$ and~\ref{L:d3is0-2} when $H= \Aut_2(V)=\SO(V)$, we have that
when $H=\Aut_2(V)$, we have $d_3^{0,2}([\Delta^{-1}\omega])=0$, and thus by the functoriality of the spectral sequence we have
$d_3^{0,2}([\Delta^{-1}\omega])=0$ in the general $H$ case.

Now, we handle the general $r$ case similarly to how we handled it for non-intermediate $V$.
We saw above that 
$(E')_3^{0,2}=\ker((d')_2^{0,2})$ is generated by $B(\Hom(V^r,\Z/n)^H)$
and any element $\psi$ of $H^2(W,\Z/n)$ that under $\Delta$ maps to a element of $\wedge_2 W^\vee$ corresponding to the matrix with $m_{ee}=1$ and all other $m_{ij}=0$.
The pullback of $\Delta^{-1} \omega$ under the quotient $V^r\ra V$ onto the last factor of $V$
is such an element $\psi$.  
We have an analogous quotient map $G\ra G/V^{e-1}$ take takes the quotient of $G$ by the first $e-1$ factors of $V$.
Thus from the functoriality of the spectral sequence,  the fact that $d_3^{0,2}(\Delta^{-1}(\omega))=0$ when $r=1$ implies that $d_3^{0,2}(\psi)=0$ in our general $r$ case.  
Thus $d_3^{0,2}( E_3^{0,2} )=(d')_3^{0,2}( E_3^{0,2} )=d_3^{0,2}(B(\Hom(V^r,\Z/n)^H))$
and $d_3^{0,2}( E_3^{0,2} )\sub d_3^{0,2}(B(\Hom(V^e,\Z/n)^H))$.
\end{proof}

\begin{proposition}\label{P:dfordualpairsgen}
Let $n,e_1,e_2$ be positive integers.
Let $V_1,V_2$ be non-isomorphic dual representations of a finite group $H$ over $\F_p$ for some  prime $p\mid n$. 
Let $W= V_1^{e_1}\times V_2^{e_2} $.
Let $\kappa=\End_H(V_1)$.  
Let $1\ra W \ra G \ra H \ra 1$ be a short exact sequence of groups.
 Let $G$ have extension class
$(\alpha_1,\alpha_2)\in H^2(H,V_1^{e_1})\times H^2(H,V_2^{e_2})=H^2(H,W)$ (with the given $H$ action on $W$), and let $r_i=\dim_{\kappa} \im \alpha_i$.
In the Lyndon-Hochschild-Serre spectral sequence computing $H^\vee(G,\Z/n)$ from $H^p(H,H^q(W,\Z/n))$,  
$$
|\ker(d_2^{0,2})|=|\kappa|^{(e_1-r_1)(e_2-r_2)}
$$
and
$$d_3^{0,2}=0.$$
\end{proposition}
\begin{proof}
We change basis, relabelling $V_i^{e_i}$ so that $\im \alpha_i$ is the $\kappa^{r_i}\sub \kappa^{e_i}$ corresponding to the last $r_i$ coordinates for $i=1,2$.
Let $W^\vee=\Hom(W,\F_p)=V_1^{e_2}\times V_2^{e_1}$ be the dual representation.
Since, $V_1,V_2$ are non-self dual and hence non-trivial, $\Hom(W,\Z/n)^H=0.$

Thus by Lemma \ref{L:H2struc}, we have $\Delta: H^2(W,\Z/n)^H \isom ( \wedge_2 W^\vee)^H$ since $\Delta: H^2(W,\Z/n) \to ( \wedge_2 W^\vee)$ is a split surjection whose kernel has no $H$-invariants in every case except the $p=2$, $4\nmid n$ case, where $\Delta$ is not split but one can instead use that $\Sym^2 W^\vee$ and $\wedge^2 W^\vee$ both have $(V_1^{e_2})^\vee \otimes (V_2^{e_2})^\vee $ as a summand and the complementary summand has no $H$-invariants. 
We will identify $( \wedge_2 W^\vee)^H$ with pairs of matrices over $\kappa$.


We have  $0\neq \omega\in (V_1\tensor V_2)^H$, with $\omega^t \in V_2\tensor V_1$ the image under switching factors.
We have $(\kappa \tensor 1)\omega= (V_1^\vee\tensor V_2^\vee)^H$.   
The $\kappa$ action on $V_2$ is defined such that for $k\in \kappa$, we have $(k\tensor 1) \omega =(1\tensor k )\omega$.
Since $V_1,V_2$ are not self-dual,  we have $(V_i\tensor V_i)^H=0$ for $i=1,2$.
A pair $(M,N)$ of an $e_1\times e_2$ matrix over $\kappa$ and an $e_2\times e_1$  matrix over $\kappa$ with entries $m_{ij},n_{ij}\in\kappa$ corresponds to 
$(\underline{m},\underline{n})=\sum_{ij} (m_{ij}\tensor 1)\omega_{ij} + (n_{ij}\tensor 1)\omega^t_{ij}   \in (W^\vee \tensor W^\vee)^H$, where $\omega_{ij}$ is the copy of $\omega$ from the $i$th copy of $V_1$ on the left tensored with the $j$th copy of $V_2$ on the right, and $\omega^t_{ij}$ is the copy of $\omega^t$ from the $i$th copy of $V_2$ on the left tensored with the $j$th copy of $V_1$ on the right.

We can choose an $\F_p$ basis of  $V_1$, a basis of $V_2$, and a corresponding basis $w_i$ of $W^\vee$, and then 
$\wedge_2 W^\vee$ is the set of elements of the form $\sum_{i<j} a_{ij}(w_i\tensor w_j-w_j\tensor w_i)$, 
i.e. those whose corresponding $(a_{ij})_{ij}$ matrices are alternating.  
These matrices with $a_{ij}$ entries are block matrices, where each
$\dim_{\F_p} V_1\times \dim_{\F_p} V_1$ block corresponds to an element of $V_i\tensor V_j$,
where $i,j\in \{1,2\}$.
So an element of $(W^\vee\tensor W^\vee)^H$ has an $(a_{ij})_{ij}$ matrix, each of whose blocks 
corresponds to one entry of our matrices $M,N$, with the structure that the blocks of $(a_{ij})_{ij}$ the matrix entries are arranged
$$
\begin{bmatrix}
0 & M\\
N &0
\end{bmatrix}.
$$
Thus these $(a_{ij})_{ij}$ for elements of 
$(W^\vee\tensor W^\vee)^H$
matrices are alternating if and only if they are skew-symmetric.
Thus an element $(\underline{m},\underline{n})\in (W^\vee\tensor W^\vee)^H$ is in 
$(\wedge_2 W^\vee)^H$ if and only if it is the negative of its transpose, which is true if and only if $N=-M^t$.
Thus the elements of $(\wedge_2 W^\vee)^H$ correspond to $M\in M_{e_2\times e_1}(\kappa)$.

For $N=-M^t$, we can take $g= \sum_{ij} (m_{ij}\tensor 1)\omega_{ij}$ as an $H$-invariant bilinear form such that $\delta g=(\underline{m},\underline{n})$.  So for $[g]\in H^2(W,\Z/n)^H$, when computing $d_{2,2}^{0,2}([g])$ we can take $c$ to be $0$
and hence $d_{2,2}^{0,2}([g])=0.$
Then by Lemma~\ref{L:d202gen}, we have that
$d_2^{0,2}\circ \Delta^{-1}$ of the element of $\wedge^2 W^\vee$ corresponding to $(M,-M^t)$
is the image of $\alpha\in H^2(H,V_1)^{e_1}\times H^2(H,V_1)^{e_2}$ under the matrix 
$$
\begin{bmatrix}
0 & M\\
-M^t &0
\end{bmatrix}
$$
that sends the $j$th copy of $H^2(H,V_2)$ to the $i$th copy of $H^2(H,V_1)$ via the map $m_{ij}\omega$
and the $i$th copy of $H^2(H,V_1)$ to the $i$th copy of $H^2(H,V_2)$ via the map $-m_{ij}\omega^t$.
Thus $d_2^{0,2}$ of the element corresponding to $M$is $0$  if and only if $\im \alpha_1\sub \kappa^{e_1}$ is in the kernel of $M$ 
and $\im \alpha_2\sub \kappa^{e_2}$ is in the kernel of $-M^t$, i.e. the last $r_1$ columns and $r_2$ rows of $M$ are $0$.
Thus $|\ker d_2^{0,2}|=|\kappa|^{(e_1-r_1)(e_2-r_2)}$.

Now let $W'$ be the quotient of $W$ by the last $r_i$ copies of $V_i$ for $i=1,2$, and $G'$ be the quotient of $G$ by those copies of $V_1,V_2$.  
Let $\beta \in \ker d_2^{0,2}$. We
saw above that the corresponding $M$ has the last $r_1$ columns and $r_2$ rows all 0.  It follows that $\beta$ is the pullback of a class 
$\tilde{\beta}\in \wedge\Hom(W'\tensor W',\Z/n)^H \sub H^2(W',\Z/n)^H$.   Since we have taken the quotient of $W$ by the image of $\alpha$, we have that
$G'$ is a trivial extension class.  By  Lemma~\ref{L:d202gen}, we conclude $d_3^{0,2}(\tilde{\beta})=0$, and then by Lemma~\ref{pullback-differential} that $d_3^{0,2}(\beta)=0$.
\end{proof}

\section{Measures from moments}\label{S:momprob}

In this section,  we will use the results from \cite{Sawin2022} to determine the unique measure on profinite $n$-oriented $\Gamma$-groups that has moments as found in the previous section and prove  Theorem \ref{T:intro-measure-general}.

The results of  \cite{Sawin2022} require one to be working in a diamond category, as defined in  \cite[Section 1.1]{Sawin2022}.  
We now review some necessary terminology.

\subsection{Terminology of diamond categories}\label{SS:diamond}
Let $\category$ be a small category. 
A \emph{quotient} of $G$ is a pair $(F, \pi^G_F)$ of an object $F\in \category$ and an epimorphism $\pi^G_F:G\ra F$.  Quotients 
$(F, \pi^G_F)$ and $(F,' \pi^G_{F'})$ are isomorphic when there is an isomorphism $\rho:F\ra F'$ such that $\pi^G_{F'}=\rho \pi^G_F$.

For $G$ an object of $\category$, quotients of $G$, taken up to isomorphism, form a partially ordered set, where $ H\geq  F$ if there is $h \colon H \to F$ compatible with the morphisms from $G$. 

A \emph{lattice} is a partially ordered set where any two elements  $x,y$ have a least upper bound (join) denoted $x \vee y$, and
a greatest lower bound (meet), denoted  $x\wedge y$. A \emph{modular} lattice is one where $a \leq b $ implies $a \vee (x \wedge b) = (a \vee x) \wedge b$ for all $x$.

For $F$ a quotient of $G$, let $ [F,G]$ be the partially ordered set of (isomorphism classes of) quotients of $G$ that are $\geq F$.

An object $G\in \category$ is called \emph{minimal} if it has only one quotient (itself).

A set $\mathcal L$ of objects of $\category$, stable under isomorphism, is called \emph{downward-closed} if, whenever $H$ is a quotient of $G$, if $ G\in \mathcal L$ then $H \in\mathcal L$.  It is called \emph{join-closed} if, whenever $H$ and $F$ are quotients of $G$ and $H \vee F$ exists (in the partially ordered set of isomorphism classes of quotients of $G$), 
if $H \in \mathcal L $ and $F\in \mathcal L$ then $H \vee F\in \mathcal L$.

A \emph{formation} is a downward-closed join-closed set of isomorphism classes of $\category$ containing every minimal object of $\category$.
For sets $B,\C$ of isomorphism classes of objects of $\category$, 
we say $B$ \emph{generates} $\C$ if $\C$ is the 
smallest downward-closed join-closed set containing $B$   and every minimal object of $\category$.  
A \emph{level} $\C$ is a set of isomorphism classes of objects of  $\category$ generated by a finite set $B$.

For any formation $\W$,  we define the \emph{$\W$-quotient} $G^{\W}$ of $G$ to be the maximal quotient of $G$ that is in $\W$, i.e. the join of all quotients of $G$ that are in $\W$. 

An epimorphism $G \to F$ is \emph{simple} if it is not an isomorphism and, whenever it is written as a composition of two epimorphisms, one is an isomorphism (i.e. if the interval $[F,G]$ contains only $F$ and $G$).

\begin{definition}\label{D:diamond}
A \emph{diamond category} is a small category $\category$ such that
\begin{enumerate}
\item \label{D:modular}
For each object $G \in \category$, the set of (isomorphism classes of) quotients of $G$ form a finite modular lattice.
\item \label{D:Aut}
Each object in $\category$ has finite automorphism group.
\item \label{D:finite}
For each level $\C$ of $\category$ and each $G\in \C$, there are finitely many elements of $\C$ with a simple epimorphism to $G$.
\item The category $\category$ has at most countably many isomorphism classes of minimal objects.
\end{enumerate}
\end{definition}


In this paper, we will be concerned with the category $\Cat$ of finite $n$-oriented $\Gamma$-groups.  In this category, a morphism of finite $n$-oriented $\Gamma$-groups
$(G,s_G) \ra (H,s_H)$ is a \emph{surjective} $\Gamma$-equivariant group homomorphism $\pi: G\ra H$ such that $\pi_* s_G=s_H$.  In particular, all morphisms of $\Cat$ are epimorphisms.
We use bold font $\bG=(G,s_G)$ to denote the $n$-oriented $\Gamma$-group, which has an underlying $\Gamma$-group $G$ and $s_G\in H_3(G,\Z/n)^\Gamma$.
We write $\Sur(\bG,\bH)$ for the morphisms in $\Cat$ from $\bG$ to $\bH$.

\begin{lemma}\label{our-cat-is-diamond}
The category of finite $\Gamma$-groups is a diamond category.
The category $\Cat$ of finite $n$-oriented $\Gamma$-groups is a diamond category. \end{lemma}

\begin{proof} The category of finite groups is a diamond category \cite[Lemma 6.10]{Sawin2022}.  Since the definition of diamond category only depends on the objects and epimorphisms in the category,  the category of finite groups with morphisms taken to be surjective group homomorphisms is also a diamond category.  
Then from \cite[Lemma 6.25]{Sawin2022}, the category of finite groups of order prime to $|\Gamma|$, with  morphisms taken to be surjective group homomorphisms, is a diamond category.  
Then the category of such groups with a $\Gamma$-action
(as well as the category of finite groups with a $\Gamma$-action)
is a diamond category by \cite[Lemma 6.20]{Sawin2022}, and the category of such groups with a $\Gamma$-action and $n$-orientation is a diamond category by \cite[Lemma 6.22]{Sawin2022}.
\end{proof}

\begin{remark}
The category $\Cat$ may not have arbitrary fiber products, as there is not necessarily an $n$-orientation on the fiber product group that pushes forward to the given $n$-orientations.  However,  we can very concretely describe the joins and meets in this category.
The lattice of quotients of $\bG$ is the opposite of the lattice of normal $\Gamma$-subgroups of $G$ ordered by 
inclusion.   The join corresponds to intersecting subgroups and the meet corresponds to taking their products.
For any two quotients $H_1,H_2$ of $G$,  if $N_i$ is the kernel of the $\Gamma$-group homomorphism $G\ra H_i$, then 
$G/(N_1N_2)$,  with $n$-orientation pushed forward from $G$,  and natural structure as a quotient of $G$, is the join of $H_1$ and $H_2$.  Moreover,
$G/(N_1\cap N_2)$, with $n$-orientation pushed forward from $G$,  and natural structure as a quotient of $G$, is the meet of $H_1$ and $H_2$.  
\end{remark}

The (group theoretic) kernel of a morphism $\bG\ra \bH$ in $\Cat$ is a $\Gamma$-group, with a conjugation action by $G$ that makes it a $(G\rtimes \Gamma)$-group and also a $[H\rtimes\Gamma]$-group.  (While the kernel has an actual $\Gamma$ 
  action and not just an outer action,  we sometimes 
   forget this because we do not always have a good notation by which to refer to it and it will be unnecessary in some parts of our argument.)
We can check that a morphism $\bG\ra \bH$ in $\Cat$ is simple if and only if its kernel $N$ is a simple $(G\rtimes \Gamma)$-group, or equivalently its kernel $N$ is a minimal nontrivial normal $\Gamma$-subgroup of $G$. 
Similarly, $\bG\ra \bH$ in $\Cat$ is simple if and only if its kernel $N$ is a simple $[H\rtimes \Gamma]$-group.
We say a morphism $\bG\ra \bH$ in $\Cat$ is \emph{semisimple} if the kernel is a direct product of simple $(G\rtimes \Gamma)$-groups.
\details{To see this agrees with semisimple from that paper:
We use \cite[Lemma 2.4 (4)]{Sawin2022}.  
 The meet of simple $\bG\ra \bG_i$ corresponds to taking the product of pairwise distinct minimal nontrivial normal $\Gamma$-subgroups $N_i$ of $G$.
We proceed inductively.  Note that $N_i$ either has trivial intersection with $N_1\cdots N_{i-1}$ or is contained in it (by minimality of $N_i$).  By dropping the $N_i$ that are contained in the product of previous subgroups, we can assume that $N_1\cdots N_{i-1} \cap N_i=1$ for all $i$, and then we have 
that the meet of simple $\bG\ra \bG_i$ has kernel $N_1 \cdots N_k$ with $N_1\cdots N_k$ isomorphic to $N_1\times \cdots \times N_k$ as a $(G\rtimes \Gamma)$-group.  Conversely, it is easy to check that if the kernel of $\bG\ra \bH$ is a product of irreducible $(G\rtimes \Gamma)$-groups, then $\bH$ is the meet of simple $\bG\ra \bG_i$.
}

For a morphism $\pi: \bG\ra \bH$ in $\Cat$, we let $\mu(\pi)=\mu(F,G)$ be the M\"{o}bius function on the poset $[\bH,\bG]$.
If $\pi$ is not semisimple, then $\mu(\pi)=0$ \cite[Lemma 2.6]{Sawin2022} (using \cite[Lemma 2.4]{Sawin2022} to see that the definition of semisimple in \cite{Sawin2022} agrees with the one here).
\begin{lemma}\label{L:calcmu}
Let $\pi: \bG\ra \bH$ be a semisimple morphism in $\Cat$ with kernel isomorphic to the $(G\rtimes \Gamma)$-group $\prod_{i=1}^r V_i^{e_i} \times \prod_{i=1}^s N_i$,
where $V_i$ are pairwise non-isomorphic simple abelian $(G\rtimes \Gamma)$-groups and the $N_i$ are simple nonabelian $(G\rtimes \Gamma)$-groups.
Let $q_i$ be the size of the endomorphism field of $V_i$, as a $(G\rtimes \Gamma)$-module.  
Then 
the lattice $[\bH,\bG]$ is the product for $i$ from $1$ to $\ell$ of the lattice of subspaces of $\mathbb F_{q_i}^{e_i}$ together with $m$ copies of the two-element lattice, and 
$\mu(\pi)=(-1)^s \prod_{i=1}^{r} (-1)^{e_i} q_i^{\binom{e_i}{2}}$.

In the same setting, if $G$ and $H$ are in a level $\cL$, then the $V_i$ and $N_i$ are all admissible. 
\end{lemma}
\begin{proof}
The proof of the first claim is almost identical to the argument for \cite[Lemma 3.10]{Sawin2022}.

For the admissibility of the $V_i$, note that some extension of $H$ by $V_i$ is a quotient of $G$ and hence is in $\cL$, so by Lemma \ref{L:basicH2} the trivial extension of $H$ by $V_i$ is in $\cL$ (since every subspace contains zero), so by definition $V_i$ is admissible. For the admissibility of the $N_i$, note similarly that some extension $H'$ of $H$ by $V_i$ is in $\cL$.  By Lemma \ref{L:nonabrig} $H' \rtimes \Gamma$ has the form $H \rtimes \Gamma \times_{ \Out(N_i) } \Aut(N_i)$, so by definition $N_i$ is admissible.
\end{proof}

For $\pi$ as in Lemma~\ref{L:calcmu}, we
 define $Z(\pi)=2^{r+s}$ (the description of the lattice in
 Lemma~\ref{L:calcmu} along with
 \cite[Lemma 3.11]{Sawin2022} shows this agrees with the definition in 
\cite[Section 1.1]{Sawin2022}). 

\subsection{Distribution from the moments}\label{SS:distribution}
Let $\cL$ be a level of the category of finite $\Gamma$-groups.
Let $\cW$ be the set of isomorphism classes of finite $n$-oriented $\Gamma$-groups whose underlying $\Gamma$-group is in $\cL$.
\details{To be more like the 3-man paper we would switch $G$ and $H$.  To be more like the moments paper we could turn $H$ to $F$ and keep $G$.  In the 3-man paper we made $I_\C$ a set with one from each isom class, but in the moments paper with just went with $\C$ being a set of isomorphism classes, which is how we are so far doing it here.}
We fix a multiset $U$ of elements of $\Gamma$, and an associated function of $\Gamma$-groups
$$G^{\cdot U} :=\frac{\prod_{\gamma\in U} |G^\gamma|}{|G^\Gamma|}.$$
Note that $G^{\cdot U}$ does not change if we replace an element of $U$ by a conjugate. 
Let   
$$
 M_{\bG}:=  \frac{ |H^2(G \rtimes \Gamma, \mathbb Z/n) | }{ | H^3(G \rtimes \Gamma, \mathbb Z/n)|  G^{\cdot U}  }
$$
if $G_\Gamma=1$ and $0$ otherwise.
We are then interested in the quantities, for $\bH\in\cW$,
\begin{align*}
v_{\cW,\bH}:=\sum_{\bG \in \cW}  \sum_{ \pi :\bG\ra\bH}  \frac{\mu(\pi)}{|\Aut(\bH)||\Aut(\bG)|} M_\bG
\quad \textrm{and}\quad \tilde{v}_{\cW,\bH}:=\sum_{\bG \in \cW}  \sum_{ \pi :\bG\ra\bH}  \frac{|\mu(\pi)|Z(\pi)^3}{|\Aut(\bH)||\Aut(\bG)|} M_\bG
.
\end{align*}

The inner sum is over morphisms $\pi$ from $\bG$ to $\bH$ in $\cW$, or equivalently because of the factor of $\mu(\pi)$, semisimple morphisms.  These sums are the key quantities used for determining a measure from its moments $M_\bG$.
As in \cite{Sawin2022}, we say that the set $(M_\bG)_\bG$ of moments is \emph{well-behaved} if the $\tilde{v}_{\cW,\bH}$ are finite.

\begin{theorem}\label{T:wavefinite}
Let $\cL$ be a level of the category of finite $\Gamma$-groups.
Let $\cW$ be the set of isomorphism classes of finite $n$-oriented $\Gamma$-groups whose underlying $\Gamma$-group is in $\cL$. Let $U$ be a multiset of elements of $\Gamma$. If all nonzero
 representations of $\Gamma$ of characteristic dividing $n$ contain some nontrivial vector fixed by at least one element of $U$, then for all $\bH\in\cW$, we have that
  $\tilde{v}_{\cW,\bH}$ is finite
 and $v_{\cW,\bH}\geq 0$.
\end{theorem}

Proving Theorem~\ref{T:wavefinite} and determining expressions for $v_{\cW,\bH}$ will be the focus of this section.
This will allow us, using \cite{Sawin2022}, to conclude a theorem on distributions of profinite $n$-oriented $\Gamma$-groups.
Let $\operatorname{s-Pro}(\Cat)$ be the category of
small profinite $n$-oriented $\Gamma$-groups.
It is not hard to check that $\operatorname{s-Pro}(\Cat)$ is the category of small pro-objects of the category $\Cat$ of finite $n$-oriented $\Gamma$-groups, in the sense discussed in \cite[Section 5.2]{Sawin2022}.
Given an $\bX\in \operatorname{s-Pro}(\Cat)$ and a level $\mathcal{J}$ of $\Cat$,
$\bX$ has a maximal quotient $\bX^{\mathcal{J}}\in \mathcal{J}$ \cite[Lemma 5.6]{Sawin2022}.
For $\bG\in\Cat$, we write $\Sur(\bX,\bG)$ for the surjective (continuous) homomorphisms of $\Gamma$-groups $X\ra G$ taking $s_X$ to $s_G$.  Let $\mathcal{P}_{\Gamma,n}$ be the set of isomorphism classes of objects of $\operatorname{s-Pro}(\Cat)$.
We consider the $\sigma$-algebra on $\mathcal{P}_{\Gamma,n}$ generated by $\{X | X^\mathcal{J} \isom \bH\}$ as $\mathcal{J}$ varies
through levels of $\Cat$ and $\bH$ through elements of $\mathcal{J}$, and the Borel topology of $\mathcal{P}_{\Gamma,n}$ for this $\sigma$-algebra. 

\begin{theorem}\label{T:measure}
Let $U$ be a multiset of elements of $\Gamma$.
If all nonzero representations of $\Gamma$ of characteristic dividing $n$ contain some nontrivial vector fixed by at least one element of $U$, then there exists a unique measure $\nu$ on $\mathcal{P}_{\Gamma,n}$,
for the $\sigma$-algebra above,
such that for every $\bG\in \Cat$
$$
\int_{X\in \mathcal{P}_{\Gamma,n}} \Sur(X,\bG) d\nu =M_\bG.
$$
For any level $\cL$ of the category of finite $\Gamma$-groups, and $\cW$ the set of isomorphism classes of finite $n$-oriented $\Gamma$-groups whose underlying $\Gamma$-group is in $\cL$, we have
$$
\nu(\{X|X^\cL\isom  \bH\})=v_{\cW,\bH}.
$$

If there is a sequence $\nu^t$ of measures on $\mathcal{P}_{\Gamma,n}$ such that 
$$
\lim_{t\ra\infty} \int_{X\in \mathcal{P}_{\Gamma,n}} \Sur(X,\bG) d\nu^t =M_\bG,
$$
then $\nu$ exists and the $\nu^t$ weakly converge to $\nu$.
\end{theorem}

\begin{proof}[Proof of Theorem~\ref{T:measure}]
We use \cite[Lemma 5.7]{Sawin2022} to relate $\mathcal{P}$ to the pro-isomorphism classes of \cite{Sawin2022}.

By \cite[Corollary 6.24]{Sawin2022}, if $\tilde{v}_{\cW,\bH}$ is finite for every  
set of isomorphism classes $\cW$ of finite $n$-oriented $\Gamma$-groups whose underlying $\Gamma$-group is in a fixed level,
then $(M_\bG)_\bG$ is well-behaved.  Theorem~\ref{T:wavefinite} gives us the finiteness of such $\tilde{v}_{\cW,\bH}$.
The uniqueness of $\nu$ follows from \cite[Theorem 1.7]{Sawin2022} and the formulas for 
$\nu(\{X|X^\cW\isom  \bH\})$ follow from \cite[Corollary 5.8]{Sawin2022} (along with \cite[Corollary 6.24]{Sawin2022} to see that $\cW$ is a narrow formation).

Existence follows from \cite[Theorem 1.7]{Sawin2022} if we can show that $v_{\mathcal{M},\bH}\geq 0$ for every level $\mathcal{M}$ of $C$ and $\bH\in \mathcal{M}$, where these are defined as the $v_{\cW,\bH}$ but with $\mathcal{M}$ replacing $\cW$.  By \cite[Corollary 6.24]{Sawin2022},   
$v_{\mathcal{M},\bH}\geq 0$ for all $\mathcal{M},\bH$ follows from the inequalities $v_{\cW,\bH}\geq 0$
from Theorem~\ref{T:wavefinite}. 

If we have $\nu^t$ as in the theorem, then \cite[Theorem 1.8]{Sawin2022} says that $\nu$ exists and the $\nu^t$ weakly converge to $\nu$.
\end{proof}

Our proof of Theorem~\ref{T:wavefinite} has much of the same high level structure as the proof of Proposition 4.4 from \cite{Sawin2024}, on the distribution of 
the fundamental groups of $3$-manifolds (or more precisely, the distribution of the groups $G$ with an element $s\in H_3(G,\Z)$ coming from the fundamental class of the $3$-manifold).  One could view \cite{Sawin2024} as a simpler analog of what happens here, but it is not a special case.   In \cite{Sawin2024}, there is no $\Gamma$ action, but our requirement that $G_\Gamma=1$ means that there is no non-trivial special case in this paper for $\Gamma$ being the trivial group.
This condition that $G_\Gamma=1$ necessitates additional arguments in each piece of our argument here.  Our working with $\Z/n\Z$ coefficients necessitates several new arguments.   However, the biggest difference is that in the loose analogy, \cite{Sawin2024} works in the setting that $U$ is empty.  On the one hand, our proofs here do not go through in the $U$ empty case,  and \cite{Sawin2024} required us to prove parity theorems on fundamental groups of $3$-manifolds and solve the moment problem in a parity restricted setting.  On the other hand,  $U$ being empty philosophically relates to a result \cite[Lemma 7.1]{Sawin2024} of which the analog here would be always having $H^2(H\rtimes \Gamma, V_i^{e_i})^{\cL,s_H}=0$.  \details{We see in every Lemma below giving a $w_{V_i}$ value that for the probability to not be $0$, that we get an upper bound on the size of
$H^2(H\rtimes\Gamma,V_i)^{\cL}$ in terms of $U$, with bigger $U$ allowing for more $z:=\dim_{kappa_i} H^2(H\rtimes\Gamma,V_i)^{\cL}.$}
 The fact that we must consider non-zero extension classes $\alpha\in H^2(H\rtimes \Gamma, V_i^{e_i})^{\cL,s_H}$ requires an entirely new kind of analysis of the Lyndon-Hochschild-Serre spectral sequence, which was mainly the content of Section~\ref{S:LHS}.  The much more delicate behavior of the spectral sequence leads to several more cases here among the $V_i$,  and also requires more complicated $q$-series arguments than in \cite{Sawin2024}.
\details{Our Lemma~\L:sH is analogous to \cite[Lemma 7.2]{Sawin2024}, except we have to add in the $\Gamma$ action.  Our Lemma~\ref{L:Mefstep1}
setting up the spectral sequence is analogous to  \cite[Lemma 7.7]{Sawin2024}.
Our lemmas~\ref{automorphism-of-extension} and \ref{automorphism-of-F}, factoring automorphism counts, correspond to \cite[Lemmas 7.9, 7.10]{Sawin2024}.
In the first case, we need a new version with the $\Gamma$ factor, and in the second case we just cite.
Our Lemma~\ref{pullback-differential} is analogous to  \cite[Lemma 7.11]{Sawin2024}, ours is just a slight generalization, on the functoriality of the Lyndon-Hochschild-Serre spectral sequence, and we refer to the proof in \cite{Sawin2024}.  Our Lemma~\ref{inner-sum-product} is analogous to  \cite[Lemma 7.12]{Sawin2024},
factoring the inner sums over $\Ex(\bH,F)$, but our $\G_\Gamma=1$ condition and $\Z/n\Z$ coefficients requires a few new arguments.   
Our Lemma~\ref{L:Mfactor}, the final factoring of the $M(e,f)$, is analagous to  \cite[Lemma 7.16]{Sawin2024}, just different in the details.
Then \cite{Sawin2024} computes factors for each type of the minimal extensions as follows: non-abelian in 7.18,7.19, then no dual appearing in 7.21 and 7.23, self-dual in 7.25 and 7.26, and dual pairs in 7.25 and 7.27.   Here we have more types,  with the addition of representations of characteristic not dividing $n$ and anomalous representations. The notion of ``anomalous'' and ``intermediate'' are related, but not the same.

}

\subsection{Expressing the $v_{\cW,\bH}$ in terms of irreducibles}\label{S:getv}
Let $\cL$ be a level of the category of finite $\Gamma$-groups.
Let $\cW$ be the set of isomorphism classes of finite $n$-oriented $\Gamma$-groups whose underlying $\Gamma$-group is in $\cL$.
We fix $\cL$, $\cW$ and 
an $\bH\in\cW$ 
 through the end of Section~\ref{S:momprob}.
 In this subsection we will express $v_{\cW,\bH}$ and $\tilde{v}_{\cW,\bH}$ in terms of irreducible $[H\rtimes\Gamma]$-groups.

\begin{lemma}
Given $\cL$, $\cW$, and $\bH$ as at the start of Section~\ref{S:getv}, there are only finitely many isomorphism classes of admissible finite simple abelian $H \rtimes \Gamma$ groups and admissible finite simple nonabelian $[H\rtimes \Gamma]$-groups

\end{lemma}
\begin{proof}By definition, every admissible finite simple group is the kernel of a morphism
of a morphism $\bG\ra \bH$ for some $\bG\in \cW$.
By \cite[Corollary 6.24]{Sawin2022}, $\cW$ is \emph{narrow,} which means there are only finitely many simple morphisms to $\bH$ from objects of $\cW$.  The lemma follows.
\end{proof}

We enumerate all isomorphism classes of admissible abelian irreducible $(H\rtimes \Gamma)$-groups $V_1,\dots,V_r$.
 Let $q_i$ be the size of the endomorphism field $\kappa_i$ of $V_i$ as a 
$(H\rtimes \Gamma)$-module.
The $V_i$ are representations of $H\rtimes \Gamma$ over a finite field $\F_{p_i}$ for some prime $p_i$ that we call the \emph{characteristic} of $V_i$.
We say $V_i$ is \emph{trivial} if it is trivial as a representation and \emph{nontrivial} otherwise.
We write $V_i^\vee$ for the dual representation of $V_i$ over $\F_{p_i}$,  and say that $V_i$ is \emph{self-dual} if $V_i\isom V_i^\vee$, as representations over $\F_{p_i}$.
Similarly, we enumerate
all isomorphism classes of admissible non-abelian irreducible $[H\rtimes \Gamma]$-groups $N_1,\dots,N_s$.

For $\pi: \bG \to \bH $ a semisimple morphism in $C$ with $\bG\in \cW$, we can write $\ker \pi$ as a product $\prod_{i=1}^r V_i^{e_i} \times \prod_{i=1}^s N_i^{f_i}$ by Lemma \ref{L:calcmu}. In this case, we say that $\pi$ has \emph{type} $(e,f),$ where $e$ and $f$ denote the tuples $(e_1,\dots,e_r)$ and $(f_1,\dots,f_s)$ respectively.  


Let
  \[ M (e, f): = \sum_{ \bG \in {\cW}} \sum_{ \substack{\pi:\bG\ra \bH  \\ \textrm{type }(e,f)} } \frac{M_{\bG}}{|\Aut(\bG)|} .\] 

\begin{lemma}\label{L:vwithexplicitmu}
Given $\cL$, $\cW$, and $\bH$ as at the start of Section~\ref{S:getv}, 
\begin{align*}
v_{\cW,\bH}&=\frac{1}{|\Aut(\bH)|}  \sum_{e_1,\dots,e_r,f_1,\dots,f_s\geq 0} (-1)^{\sum_{i} f_i }\prod_{i=1}^{r} (-1)^{e_i} q_i^{\binom{e_i}{2}}  M (e, f)
\\ \tilde{v}_{\cW,\bH}&\leq  \frac{8^r}{|\Aut(\bH)|}  \sum_{e_1,\dots,e_r,f_1,\dots,f_s\geq 0}8^{\sum_{i} f_i }\prod_{i=1}^{r} q_i^{\binom{e_i}{2}}  M (e, f).
\end{align*}
\end{lemma}  
\begin{proof}
The lemma follows from the definitions of $v_{\cW,\bH},\tilde{v}_{\cW,\bH}, Z(\pi), M (e, f)$ and Lemma~\ref{L:calcmu}.
\end{proof}
  
We now turn to computing $M(e,f)$.
  Fix for now non-negative integers $e_1,\dots,e_r,f_1,\dots,f_s$, and let $F$ be the $[H\rtimes \Gamma]$-group $\prod_{i=1}^r V_i^{e_i} \times \prod_{i=1}^{s} N_i^{f_i}$. 
 We consider an exact sequence $$ 1 \to  F \to G \stackrel{\pi}{\to} H \to 1$$ of $\Gamma$-groups,  in which the $[H\rtimes \Gamma]$ structure on $F$ from the exact sequence agrees with the given $[H\rtimes \Gamma]$ structure.   We call such a thing an \emph{extension} of $\Gamma$-groups of $H$ by $F$. 
  Two such extensions 
$ 1 \to  F \to G \stackrel{\pi}{\to} H \to 1$ and 
$ 1 \to  F \to G' \stackrel{\pi'}{\to} H \to 1$ 
  are isomorphic if there is an isomorphism of $\Gamma$-groups from $G$ to $G'$ that restricts to the identity of $F$ and $H$.  We write $\Ext_{\Gamma}(H,F)$ for the set of
  isomorphism classes of such extensions.  We often write $(G,\pi)$ for such an extension (even though the map $F\ra G$ is also part of the data, we leave it implicit).
We write $\Aut_{F,H}(G)$ for the automorphisms of such an extension, i.e. $\Gamma$-automorphisms of $G$ that restrict to the identity  on $F$ and $H$.
 
Since $(|F|,|\Gamma|)=1$,  the isomorphism class of $F$ as a $[\Gamma]$-group determines a unique \emph{isomorphism class} of $F$ as a  $\Gamma$-group.  This follows because the map $\Aut(F)\ra \Out(F)$ has kernel coprime to $|\Gamma|$,  and thus by the Schur-Zassenhaus theorem the lifts of $\Gamma\ra\Out(F)$
to $\Gamma \ra\Aut(F)$ all differ by conjugation by an element of $\Aut(F)$.  So even if $F$ is only a $[\Gamma]$-group it makes sense to write $F^{\cdot U}$.
 
Consider the Lyndon-Hochschild-Serre spectral sequence calculating $H^{p+q} ( G \rtimes \Gamma, \mathbb Z/n)$, whose second page satisfies $E_2^{p,q} = H^p ( H\rtimes \Gamma, H^q (F, \mathbb Z/n))$.  
The next lemma explains how to calculate $M(e,f)$ using information about this spectral sequence.

For a map of $\Gamma$-groups $\pi:G\ra H$, there is a map $\pi_* :H_3(G,\Z/n)^\Gamma\ra H_3(H,\Z/n)^\Gamma$.

\begin{lemma} \label{L:Mefstep1}
Given $\cL$, $\cW$, and $\bH$ as at the start of Section~\ref{S:getv},  and
$F$  the $[H\rtimes \Gamma]$-group $\prod_{i=1}^r V_i^{e_i} \times \prod_{i=1}^{s} N_i^{f_i}$ for some $e_i$ and $f_i$,
\begin{align}\label{moment-spectral-formula} &|\Aut_{[H\rtimes \Gamma]}(F)| M(e,f) = \\ \notag
& \frac{ | H^2 (H \rtimes \Gamma, \mathbb Z/n) | }{|H^3 ( H \rtimes \Gamma, \mathbb Z/n) |H^{\cdot U}}
\sum_{\substack{
(G,\pi)\in \Ext_\Gamma(H,F)\\
G\in \cL\\G_\Gamma=1 \\  s_H\in \im \pi_*}} 
 \frac{1}{|\Aut_{F,H}(G)|}    
   \frac{ |H^1( H \rtimes \Gamma, H^1( F, \mathbb Z/n) )| | E_3^{0,2} |  }{  |H^1( F, \mathbb Z/n)^{H \rtimes \Gamma} |  F^{\cdot U} }
 .\end{align}
\end{lemma}

We take sums as in the above so frequently that we write $\Ex(\bH,F)$
 for the set of 
$(G,\pi)\in \Ext_\Gamma(H,F)$ such that 
$G\in \cL$, and $G_\Gamma=1$, and  $s_H\in \im \pi_*$.

\begin{proof} 
By definition,
$$M(e,f) =\sum_{ \bG \in {\cW}} \sum_{ \substack{\pi:\bG\ra \bH  \\ \textrm{type }(e,f)} } \frac{M_G}{|\Aut(\bG)|}.$$
Each  $\pi$ in the sum for $M(e,f)$ can be extended to an exact sequence $1\ra F \ra G \ra H\ra 1$ in $|\Aut_{[H\rtimes \Gamma]}(F)|$ ways
(compatible with the $[H\rtimes \Gamma]$ structure on $F$).
Given a $G\in\cL$, we have that $\Aut_\Gamma(G)$ 
acts on $H_3(G,\Z/n)^\Gamma$, 
with the orbits corresponding to isomorphism classes of $n$-oriented $\Gamma$-groups,
and the stabilizer of $s_G$ being $\Aut(\bG).$
Thus,
 \[ M(e,f) =
\sum_{\substack{G\in \cL\\G_\Gamma=1}} 
  \sum_{\substack{  1\to F \to G \stackrel{\pi}{\to} H \to 1 \\ \textrm{type }(e,f)} }
  \frac{1}{|\Aut_{[H\rtimes \Gamma]}(F)||\Aut_\Gamma(G)|}  \sum_{\substack{ s_G\in  H_3( G,  \mathbb Z/n)^\Gamma \\  \pi_* s_G=s_H }}     \frac{ |H^2(G \rtimes \Gamma, \mathbb Z/n) | }{ | H^3(G \rtimes \Gamma, \mathbb Z/n)|  G^{\cdot U}  }.\]  
The group $\Aut_\Gamma(G)$ acts on the exact sequences that the second sum is over,
with orbits being isomorphism classes of extensions, and the stabilizer being $\Aut_{F,H}(G)$.
Thus,
 \[ M(e,f) =
\sum_{(G,\pi)\in\Ex(\bH,F)}
\frac{1}{|\Aut_{[H\rtimes \Gamma]}(F)||\Aut_{F,H}(G)|}  \sum_{\substack{ s_G\in  H_3( G , \mathbb Z/n)^\Gamma \\  \pi_* s_G=s_H }}     \frac{ |H^2(G \rtimes \Gamma, \mathbb Z/n) | }{ | H^3(G \rtimes \Gamma, \mathbb Z/n)|  G^{\cdot U} }\]  where we may introduce the condition that $s_H \in \im \pi_*$ that is included in the definition of $\Ex(\bH,F)$ since the inner sum vanishes anyways if $s_H \notin \im \pi_*$.

Because none of the terms in the sum over $s_G$ depend on $s_G$, we can replace the third sum above with the count of $s_G\in H_3(G,\Z/n)^\Gamma $ that satisfy $\pi_* s_G=s_H$. This number is $|\ker \pi_*|$ because $s_H\in \im  \pi_*$. 
Using Lemmas~\ref{L:GammaH} and \ref{L:homdual},  it follows that $ | H^3(G \rtimes \Gamma, \mathbb Z/n)|= | H_3(G,\Z/n)^\Gamma|=|\ker \pi_*||\im \pi_*|$.
Thus
 \[ M(e,f) =
\sum_{(G,\pi)\in\Ex(\bH,F)}
\frac{1}{|\Aut_{[H\rtimes \Gamma]}(F)||\Aut_{F,H}(G)|}    \frac{ |H^2(G \rtimes \Gamma, \mathbb Z/n) | }{ |\im \pi_*|  G^{\cdot U}  }.\]

Now we can calculate some of the arrows of the Lyndon-Hochschild-Serre spectral sequence calculating $H^{p+q} ( G \rtimes \Gamma, \mathbb Z/n)$ from $E_2^{p,q} = H^p ( H\rtimes \Gamma, H^q (F, \mathbb Z/n))$.
When $G_\Gamma=1$, we have $0=\Hom_\Gamma(G,\Z/n)=H^1(G,\Z/n)^\Gamma$.  Then since $|\Gamma|$ is relatively prime to $n$, we have
\[ H^1 ( G \rtimes \Gamma, \mathbb Z/n)= H^1(G,\Z/n)^\Gamma =0 .\]
Thus
 \[ \ker d_2^{0,1} = E_3^{0,1} = E_{\infty}^{0,1} =0.  \] 

By the basic structure of the spectral sequence, 
\begin{align*}
| E_{\infty}^{2,0}|=| E_{3}^{2,0}|=\frac{|E_{2}^{2,0}||\ker d_2^{0,1}|}{|E_2^{0,1}|}, \quad
 |E_{\infty}^{1,1} |=\frac{ |E_2^{1,1}| }{| \Img d_2^{1,1}|},   \quad
 | E_{\infty}^{0,2}| =\frac{ |E_3^{0,2} | }{| \Img d_3^{0,2}| }.
\end{align*}
We therefore have \[ |H^2( G\rtimes \Gamma, \mathbb Z/n) | = | E_{\infty}^{2,0}| \cdot |E_{\infty}^{1,1} | \cdot | E_{\infty}^{0,2}| = \frac{ |E_2^{2,0}|}{ |E_2^{0,1}|} \cdot \frac{ |E_2^{1,1}| }{| \Img d_2^{1,1}|}   \cdot \frac{ |E_3^{0,2} | }{| \Img d_3^{0,2}| }. \]

Now $E_{\infty}^{3,0} = \Img \pi^* : H^3(H\rtimes \Gamma,\Z/n)\ra H^3(G\rtimes \Gamma,\Z/n) $
by the edge map of the spectral sequence.  Also, by Lemma~\ref{L:homdual} and the duality theory of finite abelian groups, we have $|\im \pi_*|=|\im \pi^*|.$
Thus, $| \Img d_2^{1,1}| | \Img d_3^{0,2}| | \Img \pi_* | = |E_2^{3,0} | = |H^3( H\rtimes \Gamma, \mathbb Z/n) |$. 

This gives

\[  \frac{ |H^2(G \rtimes \Gamma, \mathbb Z/n) | }{ | \Img \pi_* |  G^{\cdot U} } = \frac{ |E_2^{2,0}| |E_2^{1,1}|  |E_3^{0,2} |}{ |E_2^{0,1}|  |H^3 ( H \rtimes \Gamma, \mathbb Z/n) | G^{\cdot U}} \] \[= \frac{ | H^2 (H \rtimes \Gamma, \mathbb Z/n) |  |H^1( H \rtimes \Gamma, H^1( F, \mathbb Z/n) )| | E_3^{0,2} |  }{  |H^1( F, \mathbb Z/n)^{H \rtimes \Gamma} | |H^3 ( H \rtimes \Gamma, \mathbb Z/n) | G^{\cdot U} } \]

Using $G^{\cdot U}= F^{\cdot U} H^{\cdot U} $ (e.g. see \cite[Lemma 3.3 (1)]{Liu2024}),
we now obtain \eqref{moment-spectral-formula}. \end{proof}

\begin{lemma}\label{automorphism-of-extension} 
For $H,F,G$ as in Lemma~\ref{L:Mefstep1},  \[|\Aut_{F,H}(G)| = \prod_{i=1}^r |H^1 ( H \rtimes \Gamma, V_i)^{e_i}|  |V_i^{\Gamma}|^{e_i}  / | V_i^{H \rtimes \Gamma } |^{e_i} .\]
 \end{lemma}
 
 \begin{proof} 
For any simple $[H\rtimes \Gamma]$-group factor $S_j$ of $F$, we can consider the quotient $G_j$ of $G$ by all the other simple factors of $F$.  We have that $G_j\in \Ex(\bH,S_i)$,  and that $G$ is the fiber product of the $G_j$ over $H$.
\details{There is a map to the fiber product, which is easily checked to be an injection, and by counting it is a surjection.}
Any element of $ \Aut_{F,H}(G)$ fixes all the simple factors of $F$ and hence acts on each $G_i$ separately.
Hence $|\Aut_{F,H}(G)|=\prod_j |\Aut_{F,H}(G_j)|$.
 
 A factor $G_j$ corresponding to a non-abelian $S_j$ has no nontrivial automorphisms fixing $S_j$ and $H$ since the extension is canonically $H\times_{\Out(S_j)}\Aut(S_j)$ by Lemma~\ref{L:nonabrig} and the fact that such $S_j$, as a group, is a product of non-abelian simple groups and hence has trivial center.

Given an extension $G_j$ by some $V_i$,  for each $h\in H$ we pick a  lift $\tilde{h}$ in $G_j$. 
Then,  a group automorphism $\alpha$ of  $G_j$, fixing $H$ and $V_j$, 
 is determined by a map of sets $\phi: H\ra V$, where $\phi(h)=\alpha(\tilde{h})\tilde{h}^{-1}$.
 One can check such a $\phi$ gives a group automorphism if and only if it is a cocycle, and a $\Gamma$-automorphism if and only if it is a $\Gamma$-equivariant map.  
\details{Note $\alpha$ is determined on $V=V_i^{e_i}$, and $\phi(h)$ is in $V$.  Since every element of $E$ can be written uniquely as $v\tilde{h}$ for some $v\in V$ and $h\in H$,  this gives a map $E\ra E$, by $\alpha(v\tilde{h})=v\phi(h)\tilde{h}$.  This is a homomorphism iff
$v\phi(h)\tilde{h} v'\phi(h')\tilde{h}' =\alpha(v\tilde{h}v'\tilde{h'})$.  We have
$$\alpha(v\tilde{h}v'\tilde{h'})=v \alpha(h)(v') \alpha(\beta(h,h') \widetilde{hh'})=
v \alpha(h)(v') \beta(h,h') \phi(hh')\widetilde{hh'}
 $$
and
$$
v\phi(h)\tilde{h} v'\phi(h')\tilde{h}' =v\phi(h) \alpha(h)(v') \alpha(h)(\phi(h')) \beta(h,h')  \widetilde{hh'}.
$$
I.e. $\alpha$ is an automorphism iff
$$
v \alpha(h)(v') \beta(h,h') \phi(hh')\widetilde{hh'}=v\phi(h) \alpha(h)(v') \alpha(h)(\phi(h')) \beta(h,h')  \widetilde{hh'},
$$
i.e $  \phi(hh')=\phi(h) \alpha(h)(\phi(h')) ,$ i.e. $\phi$ is a co-cycle. 
What is the $\Gamma$ action?  Say $\gamma(v\tilde{h})=\gamma(v)\gamma_h\widetilde{\gamma(h)}$ for some set map $\gamma: H \ra V$.
Then $\phi$ gives a $\Gamma$ invariant map if the following two things are equal
$$
\alpha(\gamma(v\tilde{h}))=\alpha(\gamma(v)\gamma_h\widetilde{\gamma(h)} )=
\gamma(v)\gamma_h \phi(\gamma(h))\widetilde{\gamma(h)}
$$ 
and
$$
\gamma(\alpha(v\tilde{h}))=\gamma(v\phi(h)\tilde{h})=\gamma(v) \gamma(\phi(h)) \gamma_h\widetilde{\gamma(h)}
$$
i.e. iff 
$$
\gamma(v)\gamma_h \phi(\gamma(h))\widetilde{\gamma(h)}=\gamma(v) \gamma(\phi(h)) \gamma_h\widetilde{\gamma(h)}
$$
i.e.
$
\phi(\gamma(h))= \gamma(\phi(h)).
$} 
The number of $1$-cocycles is the number of $1$-coboundaries times the size of $H^1(H,V_i)$, and since $|\Gamma|$ and $|V_i|$ are relatively prime,
the number of $\Gamma$-equivariant cocycles is the number of $\Gamma$-equivariant coboundaries times $|H^1(H,V_i)^\Gamma|.$
Since $|\Gamma|$ and $|V_i|$ are relatively prime, we have $H^1(H,V_i)^\Gamma \isom H^1 ( H \rtimes \Gamma, V_i)$
by Lemma~\ref{L:GammaH}.
The group of $1$-coboundaries is the quotient of $V_i$ by $V_i^H$, and thus since $|\Gamma|$ and $|V_i|$ are relatively prime, 
the number of $\Gamma$-equivariant coboundaries is $|V_i^{\Gamma}|  / | V_i^{H \rtimes \Gamma } |.$
\details{A coboundary $h\mapsto hv-v$ is $\Gamma$-equivariant if and only if $\gamma(hv-v)=\gamma(h)v-v$
for all $h$ and $\gamma$
}
So the total number of $\Gamma$-automorphisms of 
$G_j$ is $|H^1 ( H \rtimes \Gamma, V_i)|  |V_i^{\Gamma}| / | V_i^{H \rtimes \Gamma } |,$ which proves the lemma.
\end{proof} 


 \begin{lemma}[{\cite[Lemma 7.10]{Sawin2024}}] \label{automorphism-of-F} 
 For  $F$  the $[H\rtimes \Gamma]$-group $\prod_{i=1}^r V_i^{e_i} \times \prod_{i=1}^{s} N_i^{f_i}$ for some $e_i$ and $f_i$,
\[\abs{\Aut_{[H\rtimes \Gamma]}(F)} = \prod_{i=1}^r \abs{GL_{e_i}(q_i) }\prod_{i=1}^s  \abs{N_i}^{f_i} \abs{ Z_{\Out(N_i)} (H\rtimes \Gamma)}^{f_i} f_i! .\]
\end{lemma}

Using Lemma \ref{automorphism-of-extension}, we see that every term in \eqref{moment-spectral-formula} 
in Lemma~\ref{L:Mefstep1}
is independent of the choice of extension of $H$ by $F$, except possibly for $| E_3^{0,2}|$.  Thus, the most difficult part of our sum to evaluate will be
\[ \sum_{(G,\pi)\in\Ex(\bH,F)}  | E_3^{0,2} |  .\]

The next few lemmas let us write this sum as a product of local factors by showing a multiplicativity property, which will culminate in writing $M(e,f)$ as a product of local factors.

\begin{lemma}\label{L:sH}
Let $\bH,F,G$ be as in Lemma~\ref{L:Mefstep1}
and let
$d_r^{p,q}$ be the differentials of the 
Lyndon-Hochschild-Serre spectral sequence calculating $H^{p+q} ( G \rtimes \Gamma, \mathbb Z/n)$ from $E_2^{p,q} = H^p ( H\rtimes \Gamma, H^q (F, \mathbb Z/n))$.
Viewing $s_H$ as a morphism $H^3(H\rtimes \Gamma,\Z/n\Z)\ra \Z/n\Z$ by 
composing the map $H^3(H,\Z/n\Z)\ra \Z/n\Z$ from
Lemma~\ref{L:homdual} with the pullback map $H^3(H\rtimes \Gamma,\Z/n\Z)\ra H^3(H,\Z/n\Z)$,
we have $s_H\in \im \pi_*$ if and only if $s_H \circ d_{2}^{1,1}=0$ and $s_H \circ d_{3}^{0,2}=0$. 
\end{lemma}

\begin{proof}
Since $n$ is relatively prime to $|\Gamma|$, 
we have that $s_H\in \im \pi_*$ if and only if it is in the image of the map $H_3(G,\Z/n)\ra H_3(H,\Z/n)$
(without $\Gamma$ invariants).  
By Lemma~\ref{L:homdual},  $s_H$ is in the image of $H_3(G,\Z/n)$ if and only if $s_H:H^3(H,\Z/n\Z)\ra \Z/n\Z$ is the 
pullback of a morphism $H^3(G,\Z/n\Z)\ra \Z/n\Z.$  By properties of finite abelian groups,  this happens if and only if
$s_H(\ker (H^3(H,\Z/n)\ra H^3(G,\Z/n)))=0$.   
Also by Lemma~\ref{L:homdual},  $s_H$ being in a $\Gamma$-invariant element of $H_3(H,\Z/n)$ implies that
$s_H:H^3(H,\Z/n\Z)\ra \Z/n\Z$ is a $\Gamma$-invariant map, i.e. factors through $H^3(H,\Z/n\Z)_\Gamma$.
Then since $n$ is relatively prime to $|\Gamma|$,  we have $s_H(\ker (H^3(H,\Z/n)\ra H^3(G,\Z/n)))=0$ if and only if
$s_H(\ker (H^3(H,\Z/n)^\Gamma \ra H^3(G,\Z/n)))=0.$

By Lemma~\ref{L:GammaH},
$H^3(G\rtimes \Gamma, \Z/n) \ra H^3(G,\Z/n)^\Gamma$ is an isomorphism, and similarly for $H$.
Thus $s_H(\ker (H^3(H,\Z/n)^\Gamma \ra H^3(G,\Z/n)))=0$ if and only if 
$s_H(\ker (H^3(H\rtimes\Gamma,\Z/n) \ra H^3(G\rtimes \Gamma,\Z/n)))=0.$
By considering the edge maps of the spectral sequence, we see that
$s_H(\ker (H^3(H\rtimes \Gamma,\Z/n)\ra H^3(G\rtimes\Gamma,\Z/n)))=0$ if and only if 
$s_H \circ d_{2}^{1,1}=0$ and $s_H \circ d_{3}^{0,2}=0$.  
\end{proof}

\begin{lemma}\label{inner-sum-product} 
Let $\cL$ be a level of the category of finite $\Gamma$-groups.
Let $H$ be a $\Gamma$-group and let $F_a$ and $F_b$ be two semisimple $[H \rtimes \Gamma]$-groups. Assume that, for each abelian simple $(H\rtimes\Gamma)$-group $V_i$ appearing in $F_a$ which either has characteristic dividing $n$ or is the trivial representation, 
the dual representation $V_i^\vee$  does not appear in $F_b$. 

Then 
\[\sum_{\substack{
(G,\pi)\in \Ex(\bH,F_a\times F_b)}}    | E_{3,G}^{0,2} |  = 
\Biggl( \sum_{\substack{
(G_a,\pi_a)\in \Ex(\bH,F_a)}}  | E_{3,G_a}^{0,2} | \Biggr) 
 \Biggl(\sum_{\substack{
(G_b,\pi_b)\in \Ex(\bH,F_b)}}   | E_{3,G_b}^{0,2} | \Biggr)   .\]
Here, $E_{3,G}^{0,2}$ refers to the spectral sequence for $1\ra F_a\times F_b \ra G \ra H\ra 1$, and 
$E_{3,G_a}^{0,2}$ refers to the spectral sequence for $1\ra F_a \ra G_a \ra H\ra 1$,
and 
$E_{3,G_b}^{0,2}$ refers to the spectral sequence for $1\ra F_b \ra G_b \ra H\ra 1$
\end{lemma}

\begin{proof} 
Every extension $G$ of $H$ by $F_a \times F_b$ is the fiber product of an extension $G_a=G/F_b$ of $H$ by $F_a$ and an extension $G_b=G/F_a$ of $H$ by $F_b$. 
This gives a bijection  $ \Ext_\Gamma(H,F_a\times F_b) \ra \Ext_\Gamma(H,F_a) \times \Ext_\Gamma(H,F_b)$.
Using that $\cL$ is join-closed and downward-closed, 
it follows that $G\in\cL$ if and only if $G_a,G_b\in\cL$.

Thus, matching terms on both sides, it suffices to show that 
$s_H\in \im \pi_*$ if and only if $s_H\in \im (\pi_a)_*$ and $s_H\in \im (\pi_b)_*$, 
that $G_\Gamma=1$ if and only if 
$(G_a)_\Gamma=1$ and $(G_b)_\Gamma=1$,
and the following isomorphism
\begin{equation}\label{E3-product} E_{3,G}^{0,2}\isom E_{3, G_a}^{0,2} \times E_{3, G_b}^{0,2}.\end{equation}

Let $F=F_a\times F_b$.
To evaluate $E_{3,G}^{0,2}$, we first check that the product of natural pullback maps
\begin{equation}\label{E02splitting} H^0 (H \rtimes \Gamma, H^2(F_{a}, \mathbb Z/n)) \times H^0 (H \rtimes \Gamma, H^2(F_{b}, \mathbb Z/n))\to H^0 (H \rtimes \Gamma, H^2(F, \mathbb Z/n))\end{equation} is an isomorphism. 

To do this,   first note that the maps $F_a,F_b\ra F$ and $F\ra F_a,F_b$, give $H\rtimes \Gamma$-equivariant maps
$$H^q (F_a, \mathbb Z) \times H^q (F_b, \mathbb Z) \stackrel{i_q}{\to} H^q (F, \mathbb Z) \stackrel{j_q}{\ra} H^q (F_a, \mathbb Z) \times H^q (F_b, \mathbb Z)$$
such that the composite is the identity for $q>0$.   So in particular, we have an
$H\rtimes \Gamma$-equivariant isomorphism
$H^q (F, \mathbb Z)\isom H^q (F_a, \mathbb Z) \times H^q (F_b, \mathbb Z)  \times \ker j_q$.  An analogous statement is true with $\Z/n\Z$ coefficients, with maps $i_q^n,j^n_q$.

\details{Because, in the composite, each component is the identity map in itself and, for $q>0$, the $0$ map into the other component, because it factors through the trivial group.}
For any finite group $A$,  we have functorial isomorphisms
$H^2(A, \mathbb Z)\isom H^1(A, \mathbb Q/\mathbb Z) \isom \Hom (A, \mathbb Q/\mathbb Z)$.
It follows that the map \[j_2: H^2 (F, \mathbb Z)\to  H^2 (F_a, \mathbb Z) \times H^2 (F_b, \mathbb Z)\] is an isomorphism.
\details{since we have an identity map on a group factoring through a group of the same size, the intermediate maps must be isomorphisms.}


For any finite group $A$,  we have $H^1(A,\Z)=0$.  Thus by the K\"unneth formula with principal ideal domain coefficients, 
we have a functorial exact sequence, 
 \[0\ra  H^3 (F_a, \mathbb Z) \times H^3 (F_b, \mathbb Z)  \stackrel{i_3}{\to}
H^3 (F, \mathbb Z) \ra   \operatorname{Tor}_1 ( H^2 (F_a, \mathbb Z), H^2(F_b, \mathbb Z))\ra 0\]
and thus an $H\rtimes \Gamma$-equivarant  isomorphism $\ker j_3 \isom  \operatorname{Tor}_1 ( H^2 (F_a, \mathbb Z), H^2(F_b, \mathbb Z))$.

 Applying the exact sequence $0 \to \mathbb Z \to \mathbb Z \to \mathbb Z/n\to 0$, we have exact sequences
\begin{align*}
0&\ra & H^2(F, \mathbb Z)/n &\ra & H^2(F, \mathbb Z/n) &\ra & H^3(F, \mathbb Z)[n] &\ra 0\\
0&\ra & H^2(F_a, \mathbb Z)/n \times H^2(F_b, \mathbb Z)/n &\ra & H^2(F_a, \mathbb Z/n) \times H^2(F_b, \mathbb Z/n)
&\ra & H^3(F_a, \mathbb Z)[n] \times H^3(F_b, \mathbb Z)[n] &\ra 0.
\end{align*} 
Using the map between these exact sequences given by $F_a\ra F$ and $F_b\ra F$,  
from the snake lemma and observations above it follows that we have an $H\rtimes \Gamma$-equivariant isomorphism $\ker j^n_2\isom \operatorname{Tor}_1 ( H^2 (F_a, \mathbb Z), H^2(F_b, \mathbb Z))[n].$  Thus, from the above observations,
we have an $H\rtimes \Gamma$-equivariant isomorphism
  \[H^2(F, \mathbb Z/n) \isom H^2(F_a, \mathbb Z/n)\times H^2(F_b, \mathbb Z/n) \times \operatorname{Tor}_1 ( H^2 (F_a, \mathbb Z), H^2(F_b, \mathbb Z))[n]. \]

Since $H^2(F_a, \mathbb Z)\isom\Hom(F_a,\Q/\Z)$
and $H^2(F_b, \mathbb Z)$ are
products of vector spaces over finite fields, 
 there is a functorial isomorphism from their $\operatorname{Tor}_1$  to their tensor product. 
\details{This can be done on the Sylow $p$ subgroups one prime at a time.  Here $A$ has a resolution $V\stackrel{p}{\ra} V \ra A$, where $V$ is a free $\Z$-module.
Then the kernel of $V\tensor B \stackrel{p}{\ra} V\tensor B$ is precisely $V\tensor B\isom A\tensor B$ since $pB=0$.}
We claim $H^2(F_a, \mathbb Z)\tensor H^2(F_b, \mathbb Z)$ contains no nontrivial element that is both $n$-torsion and $H \rtimes \Gamma$-invariant. Such an element would give a nontrivial $H \rtimes \Gamma$-invariant $\mathbb Z/n$-valued bilinear form on $F_a \times F_b$, which cannot exist because of our assumption on the irreducible factors of $F_a$ and $F_b$.   Thus we conclude that the map in Equation~\eqref{E02splitting} is an isomorphism.

More straightforwardly, the product of natural pullback maps \[ H^1 (F_a, \mathbb Z/n) \times H^1(F_b, \mathbb Z/n) \to H^1(F_a\times F_b, \mathbb Z/n)\] is an isomorphism because these cohomology groups are the same as sets of homomorphisms to $\mathbb Z/n$, hence \begin{equation}\label{pullback-product-1} H^p (H \rtimes \Gamma, H^1(F_{a}, \mathbb Z/n)) \times H^p (H \rtimes \Gamma, H^1(F_{b}, \mathbb Z/n))\to H^p (H \rtimes \Gamma, H^1(F_{ab}, \mathbb Z/n))\end{equation}  is an isomorphism for all $p$.

We adopt an analogous notation for the differentials as for the terms of the spectral sequences.  
Using \eqref{E02splitting}, the $p=2$ case of \eqref{pullback-product-1}, and Lemma \ref{pullback-differential}, it follows that
 \[ d_{2,G}^{0,2}: H^0 (H \rtimes \Gamma, H^2(F, \mathbb Z/n)) \to H^2 (H \rtimes \Gamma, H^1(F, \mathbb Z/n))\] is the product of $ d_{2,G_a}^{0,2}$ and $ d_{2,G_b}^{0,2}$. Hence $E_{3,G}^{0,2}= \ker d_{2,G}^{0,2}$ is the product of the kernel $E_{3, G_{a}}^{0,2}$ of $d_{2,G_a}^{0,2}$  and the kernel $E_{3, G_{b}}^{0,2}$ of $d_{2,G_b}^{0,2}$, verifying \eqref{E3-product}. 

Recall by Lemma~\ref{L:sH}, we have $s_H\in \im \pi_*$ if and only if $s_H \circ d_{2,G}^{1,1}=0$ and $s_H \circ d_{3,G}^{0,2}=0$, and the
 the analogous statements are true for $G_a$ and $G_b$.
Using Lemma \ref{pullback-differential} and the $p=1$ case of \eqref{pullback-product-1}, the map $d_{2,G}^{1,1}$ is the sum of $d_{2,G_a}^{1,1}$ and $d_{2,G_b}^{1,1}$, hence $s_H \circ d_{2,G}^{1,1}=0$ if and only if $s_H \circ d_{2,G_a}^{1,1}=0$ and $s_H \circ d_{2,G_b}^{1,1}=0$. Similarly, using Lemma \ref{pullback-differential} and \eqref{E3-product}, $d_{3,G}^{0,2}$ is the sum of $d_{3,G_a}^{0,2}$ and $d_{3,G_b}^{0,2}$, hence $\tau \circ d_{3,G}^{0,2}=0$ if and only if $\tau \circ d_{3,G_a}^{0,2}=0$ and $\tau \circ d_{3,G_b}^{0,2}=0$.   Thus we conclude that $s_H\in \im \pi_*$ if and only if $s_H\in \im (\pi_a)_*$ and $s_H\in \im (\pi_b)_*$.

Finally, if $G_{\Gamma}=1$ then $(G_a)_\Gamma=1$ and $(G_b)_{\Gamma}=1$ since $G_a$ and $G_b$ are quotients of $G$. 
Suppose $(G_a)_\Gamma=1$ and $(G_b)_{\Gamma}=1$, which in particular implies $H_\Gamma=1$.  The image of $F_a$ and $F_b$ in $G_\Gamma$ are both normal subgroups.
If we quotient $G_\Gamma$ by the image of $F_a$, we obtain $(G_b)_{\Gamma}=1$ (and similarly for $F_b$), and so
image of $F_a$ in $G_\Gamma$ is all of $G_\Gamma$ and similarly for $F_b$.
So given two elements of $G_\Gamma$,  we can write one as an image of an element of $F_a$ and the other as an image of an element of $F_b$ to see they commute, and hence $G_\Gamma$ is abelian.   So the image of $F_a$ in $G_\Gamma$ factors through $F_a^{\operatorname{ab}}$ (where $\operatorname{ab}$ denotes the abelianization), and the images of any element of $F$ and any element of $G$ commute in $G_\Gamma$, so the image of $F_a$ in $G_\Gamma$ factors through $(F_a^{\operatorname{ab}})_{H\rtimes \Gamma}$.  Each simple factor of $F_a$ that is not abelian and a trivial representation of $H\rtimes \Gamma$
is sent to $1$ in the map to $(F_a^{\operatorname{ab}})_{H\rtimes \Gamma}$.  All of this is also true for $F_b$.  So, since
for each prime $p$, one of $F_a$ or $F_b$ does not contain a characteristic $p$ trivial representation of $H\rtimes \Gamma$,  we conclude $G_\Gamma=1$.
We have shown that $G_\Gamma=1$ if and only if 
$(G_a)_\Gamma=1$ and $(G_b)_\Gamma=1$, which was the last statement needed to prove the lemma.
\end{proof}

Lemma~\ref{inner-sum-product} will allow us to factor $M(e,f)$, and we now will give notation for the factors.  

We let \[V_{i}^{\cdot U'} :=\frac{ |V_i^\Gamma|}{|V_i^{H\rtimes \Gamma}|  } V_i^{\cdot U}=\frac{1}{ |V_i^{H\rtimes \Gamma}|  } \prod_{\gamma\in U} |V_i^\gamma|.\]

\begin{mdframed}
For $V_i$ of characteristic not dividing $n$, let
\[ M_i (e_i) = \frac{ 1 }{  \abs{GL_{e_i}(q_i) }|H^1 ( H \rtimes \Gamma, V_i)|^{e_i}  (V_i^{\cdot U'})^{e_i} }    \sum_{\substack{
(G,\pi)\in \Ex(\bH,V_i^{e_i})}}    | E_{3}^{0,2} |  .\]

For $V_i$ self-dual of characteristic dividing $n$,  let
\[ M_i (e_i) =  \frac{1 }{ \abs{GL_{e_i}(q_i) }\ (V_i^{\cdot U'})^{e_i}  | V_i^{H \rtimes \Gamma } |^{e_i}}  \sum_{\substack{
(G,\pi)\in \Ex(\bH,V_i^{e_i})}}    | E_{3}^{0,2} |  .\]

For $V_i$ of characteristic dividing $n$ such that $V_i^\vee\not\isom V_j$ for any $j$, let
\[ M_i (e_i) =  \frac{|H^1 ( H \rtimes \Gamma, V_i^\vee)|^{e_i} }{ \abs{GL_{e_i}(q_i) }
|H^1 ( H \rtimes \Gamma, V_i)|^{e_i}
 (V_i^{\cdot U'})^{e_i}  }  \sum_{\substack{
(G,\pi)\in \Ex(\bH,V_i^{e_i})}}    | E_{3}^{0,2} |  .\]

If $V_i$ falls into one of the above three categories, i.e. $\operatorname{char}(V_i)\nmid n$, or $V_i$ is self-dual, or $V_i^\vee$ is not among the $V_j$, we say that $V_i$ is \emph{solo}.

For $V_i$ and $V_{i'}$ of characteristic dividing $n$, dual to each other, and non-isomorphic, define
\[ M_{i,i'} (e_i, e_{i'} ) = \frac{ |H^1 ( H \rtimes \Gamma, V_i)|^{e_{i'}-e_i}  |H^1 ( H \rtimes \Gamma, V_{i'})|^{e_{i} -e_{i'} }  }{  
   \abs{GL_{e_i}(q_i) } \abs{GL_{e_{i'}}(q_{i'}) }\ (V_i^{\cdot U'})^{e_i}(V_{i'}^{\cdot U'})^{e_{i'}} } 
  \sum_{
(G,\pi)\in \Ex(\bH,V_i^{e_i} \times V_{i'} ^{e_{i'}})} | E_3^{0,2} | .\] 

For any $N_i$, let
\[ \eta_i (f_i) = \frac{1}{  \abs{N_i}^{f_i} \abs{ Z_{\Out(N_i)} (H\rtimes \Gamma)}^{f_i} f_i! (N_i^{\cdot U})^{f_i} }   \sum_{(G,\pi)\in \Ex(\bH,N_i^{f_i} )}  | E_3^{0,2} |. \]
\end{mdframed}

\begin{lemma} \label{L:Mfactor}
Given $\cL$, $\cW$, and $\bH$ as at the start of Section~\ref{S:getv},  and non-negative integers
$e_1,\dots,e_r,f_1,\dots,f_s$,
\begin{equation*} \label{moment-product-formula} \begin{split}  M(e,f)   = M_{\bH} \prod_{ \substack{ i \in \{1,\dots ,r \} \\ V_i \textrm{ solo}}} M_i ( e_i) \prod_{\substack{ \{i, i'\} \subseteq \{1,\dots, r\}  \\  i \neq i' \\ V_i \cong V_{i'}^\vee \\ \chr(V_i) \mid n }} M_{i,i'} (e_i, e_{i'})  \prod_{i=1}^s \eta_i (f_i ) .\end{split}\end{equation*}  
\end{lemma}

\begin{proof}
Recall that
$$
 M_{\bH}:=  \frac{ |H^2(H \rtimes \Gamma, \mathbb Z/n) | }{ | H^3(H \rtimes \Gamma, \mathbb Z/n)| H^{\cdot U}  }.
$$
 By Lemma~\ref{L:Mefstep1}, it suffices to check that
\[ \sum_{(G,\pi)\in \Ex(\bH,F) }\frac{1}{|\Aut_{[H\times \Gamma]}(F)|
|\Aut_{F,H}(G)|} \frac{ | H^1( H \rtimes \Gamma, H^1( F, \mathbb Z/n)) | | E_3^{0,2} |  }{  |H^1( F, \mathbb Z/n)^{H \rtimes \Gamma} | F^{\cdot U} }\]
\[=  \prod_{ \substack{ i \in \{1,\dots ,r \} \\ V_i \textrm{ solo}}} M_i ( e_i) \prod_{\substack{ \{i, i'\} \subseteq \{1,\dots, r\}  \\  i \neq i' \\ V_i \cong V_{i'}^\vee \\ \chr(V_i) \mid n }} M_{i,i'} (e_i, e_{i'})  \prod_{i=1}^s \eta_i (f_i ) .\]

Lemmas \ref{automorphism-of-extension} and \ref{automorphism-of-F}  give factorizations of the 
$|\Aut_{[H\times \Gamma]}(F)|$ and $|\Aut_{F,H}(G)|$ terms:
$$
\frac{1}{|\Aut_{[H\times \Gamma]}(F)|
|\Aut_{F,H}(G)|} = \prod_{i=1}^r 
\frac{| V_i^{H \rtimes \Gamma } |^{e_i}}{ \abs{GL_{e_i}(q_i) } |H^1 ( H \rtimes \Gamma, V_i)^{e_i}|  |V_i^{\Gamma}|^{e_i}}
\prod_{i=1}^s  \frac{1}{\abs{N_i}^{f_i} \abs{ Z_{\Out(N_i)} (G)}^{f_i} f_i! }.
$$ Lemma~\ref{inner-sum-product} gives us a factorization of the sum over $|E_3^{0,2}|$.
We now split the remaining terms as products:
 \[ H^1( F, \mathbb Z/n)  = \prod_{ \substack{ i \in \{1,\dots, r\} \\ \chr(V_i) \mid n }} \left( V_i^\vee \right)^{e_i} \] so
 \begin{equation}\label{E11-as-product} | H^1( H \rtimes \Gamma, H^1( F, \mathbb Z/n)) |= \prod_{ \substack{ i \in \{1,\dots, r\} \\ \chr(V_i) \mid n }} H^1( H\rtimes \Gamma, V_i^\vee)^{e_i} \end{equation}
 and
 \begin{equation}\label{E01-as-product} |H^1( F, \mathbb Z/n)^{H \rtimes \Gamma}|= \prod_{ \substack{ i \in \{1,\dots, r\} \\ \chr(V_i) \mid n \\ V_i \textrm{ trivial $H\rtimes\Gamma$ rep} }}  |V_i|^{e_i} .\end{equation}
 Finally
 \begin{equation}\label{Y-as-product} F^{\cdot U} = \prod_{i=1}^r (V_i^{\cdot U})^{ e_i} \prod_{i=1}^s (G_i^{\cdot U})^{f_i}  .\end{equation}


We can check in each case that all the factors associated to $i$ in these products have been absorbed into the definition of $M_i (e_i)$, $M_{i,i'} (e_i, e_{i'})$, or $N_i(f_i)$, as appropriate.
\end{proof}

Now we will give notation for the factors of the $v_{\cW,\bH}$ and $\tilde{v}_{\cW,\bH}.$ 
\begin{mdframed}
We define, for a solo $V_i$, 
\begin{equation}
w_{V_i}:=\sum_{e_i\geq 0}(-1)^{e_i} q_i^{\binom{e_i}{2}} \quad
 M_i ( e_i)
\quad \textrm{and}\quad  \tilde{w}_{V_i}:=\sum_{e_i\geq 0} q_i^{\binom{e_i}{2}}
 M_i ( e_i)
 ,
\end{equation}
for a $V_i,V_{i'}$ for non-isomorphic dual representations of characteristic dividing $n$,
we define
\begin{align}
w_{V_i}w_{V_{i'}}&:=\sum_{e_i, e_{i'}\geq 0}  (-1)^{e_i+e_{i'}} q_i^{\binom{e_i}{2}+\binom{e_{i'}}{2}} 
 M_{i,i'} (e_i, e_{i'})
\quad \textrm{and} \\
\notag
\tilde{w}_{V_i}\tilde{w}_{V_{i'}}&:=\sum_{e_i, e_{i'}\geq 0} q_i^{\binom{e_i}{2}+\binom{e_{i'}}{2}} 
 M_{i,i'} (e_i, e_{i'}).
\end{align}
(This is a slight abuse of notation, as we never define $w_{V_i}$ and $w_{V_{i'}}$ separately, because they always appear together in all formulas, and similarly for  $\tilde{w}_{V_i}$ and $\tilde{w}_{V_{i'}}.$ We could define each of them to be the square-root of the relevant expression above.)
For each $N_i$, we define
\begin{equation}
w_{N_i}:=\sum_{f_i\geq 0} (-1)^{f_i}
 \eta_i (f_i )
 \quad \textrm{and}\quad  
\tilde{w}_{N_i}:=\sum_{f_i\geq 0} 8^{f_i}
 \eta_i (f_i )
 . 
\end{equation}
\end{mdframed}

\begin{corollary}\label{C:vfactor}
Given $\cL$, $\cW$, and $\bH$ as at the start of Section~\ref{S:getv},  
\begin{align*}
&v_{\cW,\bH}=\frac{M_\bH}{|\Aut(\bH)|}  
\prod_{ i=1}^{r} w_{V_i}
 \prod_{i=1}^s  w_{N_i}
  \quad \textrm{and}\quad  
 \tilde{v}_{\cW,\bH}\leq \frac{8^r M_\bH}{|\Aut(\bH)|}  
\prod_{ i=1}^{r} \tilde{w}_{V_i}
 \prod_{i=1}^s  \tilde{w}_{N_i}
 .
 \end{align*}
\end{corollary}

\begin{proof}
By Lemma~\ref{L:vwithexplicitmu}.
$$v_{\cW,\bH}=\frac{1}{|\Aut(\bH)|}  \sum_{e_1,\dots,e_r,f_1,\dots,f_s\geq 0} (-1)^{\sum_{i} f_i }\prod_{i=1}^{r} (-1)^{e_i} q_i^{\binom{e_i}{2}}  M (e, f).$$  Lemma~\ref{L:Mfactor} gives a factorization of the $M(e,f)$ term, which gives the corollary.
The argument is the same for  $\tilde{v}_{\cW,\bH}$.
\end{proof}

In the next several subsections we explicitly compute the terms $M_i(e_i),M_{i,i'}(e_i,e_{i'})$, and then $w_{V_i},w_{N_i}$.
From these calculations we will be able to use Corollary~\ref{C:vfactor} to prove Theorem~\ref{T:wavefinite} on the finiteness of 
 $\tilde{v}_{\cW,\bH}$ and the non-negativity of $v_{\cW,\bH}$.   From the definition of the $v_{\cW,\bH},\tilde{v}_{\cW,\bH}$ and the $M_\bH$, we can see that $v_{\cW,\bH}=\tilde{v}_{\cW,\bH}=0$ when $H_\Gamma\ne 1. $ So  we will restrict our work in the following subsections to the case that $H_\Gamma=1$.

 \subsection{Non-abelian groups}
 Recall from Lemma~\ref{L:nonabrig} that $\Aut(N_i)\times_{ \Out(N_i) } (H\rtimes \Gamma)$ is the unique extension of $H $ by the $[H\rtimes \Gamma]$-group $N_i$. 
Consider 
 the Lyndon-Hochschild-Serre spectral sequence computing $H^*(  \Aut(N_i)\times_{ \Out(N_i) } (H \rtimes \Gamma),\Z/n))$ with $E^{p,q}_2=H^p(H\rtimes \Gamma, H^q(N_i,\Z/n))$. 
Since $N_i$, as a group, is a product of non-abelian simple groups, we have $H^1(N_i,\Z/n)=0$ and hence $d_2^{0,2}=d_2^{1,1}=0$.  
  Let \[ \delta_{N_i}=d^{0,2}_3 \colon H^2 ( N_i, \mathbb Z/n)^{H\rtimes \Gamma} \to H^3 ( H \rtimes \Gamma, \mathbb Z/n). \] 

Given a $[H\rtimes \Gamma]$-group $N$,  let $L_{N}$ be the number of lifts of $H\rtimes \Gamma\ra \Out(N)$ to $H\rtimes\Gamma\ra \Aut(N)$.

\begin{lemma}\label{non-abelian-tilde} 
Let $\cL$, $\cW$, and $\bH$ be as at the start of Section~\ref{S:getv},  and further assume $H_\Gamma=1$.
 If $H \rtimes \Gamma\ra \Out(N_i)$ is trivial or $s_H \circ \delta_{N_i}$ is nontrivial,
 \[ \eta_i (f_i) = \begin{cases} 1 & f_i=0 \\ 0 & f_i >0 \end{cases} .\] For any other $i$ from $1$ to $s$,  \[   \eta_i (f_i) = \left(\frac{L_{N_i} | H^2 ( N_i, \mathbb Z/n)^{H \rtimes \Gamma}|}{  \abs{N_i} \abs{ Z_{\Out(N_i)} (H\rtimes \Gamma)}  N_i^{\cdot U}}  \right)^{f_i} \frac{1}{f_i!}.
 \]  \end{lemma}


 \begin{proof}  
Recall
$$ \eta_i (f_i) = \frac{1}{  \abs{N_i}^{f_i} \abs{ Z_{\Out(N_i)} (H\rtimes \Gamma)}^{f_i} f_i! (N_i^{\cdot U})^{f_i} }   \sum_{(G,\pi)\in \Ex(\bH,N_i^{f_i} )}  | E_3^{0,2} |$$
 and $\Ex(\bH,F)$ is the set of 
$(G,\pi)\in \Ext_\Gamma(H,F)$ such that 
$G\in \cL$, and $G_\Gamma=1$, and  $s_H\in \im \pi_*$.

 We have $H^1(N_i^{f_i}, \mathbb Z/n)=0$ so $H^p(H, H^1(N_i^{f_i}, \mathbb Z/n))=0$ for all $p$. Thus the differentials $d_2^{1,1}$ and $d_2^{0,2}$ vanish.
  By Lemma~\ref{L:sH}, we have $s_H\in\im \pi_*$ if and only if $s_H \circ d_3^{0,2} =0$.
  By Equation~\eqref{E02splitting},  $d_3^{0,2}$ is the sum of $f_i$ copies of $\delta_{N_i}$.
  So $s_H \circ \delta_{N_i}\neq 0$ then we have $s_H\in\im \pi_*$ if and only if $f_i=0$,  in which case $\eta_i (f_i) = \begin{cases} 1 & f_i=0 \\ 0 & f_i >0 \end{cases}$, while if $s_H \circ \delta_{N_i}= 0$ then we always have $s_H\in\im \pi_*$.

When $f_i>0$, we claim $G_\Gamma=1$ if and only if $H \rtimes \Gamma \ra \Out(N_i)$ is nontrivial.
If  $H \rtimes \Gamma \ra \Out(N_i)$ is trivial, then $G=N_i \times (H\rtimes \Gamma)$.  Since $(|N|,|\Gamma|)=1$, we must also have
$\Gamma\ra\Aut(N_i)$ is trivial, and so $G_\Gamma=N_i \times H_\Gamma$, which is nontrivial.

Now we assume $G_\Gamma$ is nontrivial.
Since $H_\Gamma=1$, we have that the image of $N=N_i^{f_i}$ must generate $G_\Gamma$, and 
 so some copy of $N_i$ must have nontrivial image in $G_\Gamma$.  The kernel of the map from $N_i$ to its image in 
$G_\Gamma$ is an intersection of two normal subgroups of $G$, so it is a normal subgroup of $G$, and hence fixed by the $H$-action on normal subgroups of $N_i$.  This kernel is also fixed by the $\Gamma$-action on normal subgroups of $N_i.$  Since $N_i$ is a simple $[H\rtimes \Gamma]$-group, it follows this kernel is either trivial or $N_i$, and since we assumed the image was nontrivial, the kernel must be trivial.  It follows that $\Gamma$ acts trivially on $N_i$, which implies that $H\rtimes \Gamma\ra \Out(N_i)$ factors through $H_\Gamma=1$, and hence $H\rtimes \Gamma\ra \Out(N_i)$ is trivial.  

Hence, we conclude the statement of the lemma when $H \rtimes \Gamma\ra \Out(N_i)$ is trivial.

By definition of the $N_i$,  the extension of $H$ by $N_i$ is in $\cL$ and thus the extension of $H$ by $N_i^{f_i}$ 
is the $f_i$-fold fiber product of the extension of $H$ by $N_i$, and hence also in $\cL$ since $\cL$ is join-closed.
 
 Next we count $|\Ext_\Gamma(H,N_i^{f_i})|.$  By Lemma~\ref{L:nonabrig}, there is a unique group extension $G$ of $H$ by $N_i^{f_i}$, and by Lemma~\ref{L:nonabGam}, the number of $\Gamma$-group structures on $G$ compatible with the structures on $H$ and $N_i^{f_i}$ is 
$L_{N_i}^{f_i}$.
 
 Since $d^{0,2}_2=0$,
  \[E_3^{0,2}= H^0 ( H \rtimes \Gamma, H^2 ( N_i^{f_i}, \mathbb Z/n)) =  H^0 ( H \rtimes \Gamma, H^2 ( N_i, \mathbb Z/n))^{f_i} .\]
 Putting this all into the definition of $\eta_i (f_i)$, we conclude the lemma.
  \end{proof}

\begin{lemma}\label{non-abelian-w}  
Let $\cL$, $\cW$, and $\bH$ be as at the start of Section~\ref{S:getv},  and further assume $H_\Gamma=1$.
For any $N_i$, we have
\[ w_{N_i}\oftau  = e^{ - \frac{  L_{N_i} \left| H^2( N_i, \mathbb Z/n)^{H \rtimes \Gamma} \right| }{ |N_i| | Z_{ \Out(N_i)} ( H \rtimes \Gamma) |N_i^{\cdot U}  }} .\] 
if $H \rtimes \Gamma\ra\Out(N_i)$ is nontrivial and $s_H \circ \delta_{N_i}=0$, and $1$ if either condition fails.  Further, $\tilde{w}_{N_i}$ is finite.
\end{lemma}

\begin{proof} 
This follows immediately from the definitions of
$w_{N_i}$, $\tilde{w}_{N_i}$ and  Lemma~\ref{non-abelian-tilde}.\end{proof}

\subsection{Preparation for extensions by abelian groups}

 In the sums for the $M_i(e_i)$, we will need to sum over $\Ex(\bH,V_i^{e_i})$,  the set of $(G,\pi)\in \Ext_\Gamma(H,V_i^{e_i})$ such that $G\in\cL$,  and $G_\Gamma=1$, and $s_H\in \im \pi_*$.   Even though $V_i$ was chosen because there is some extension $1\ra V_i\ra G \ra H \ra 1$ such that $G\in \cL$,
 it does not follow that every such extension of $H$ by $V_i$ is in $\cL$.  
 We will make a definition below that captures $G\in\cL$ and $s_H\in \im \pi_*$. 

 \begin{definition} \label{D:whichH2}
 For an abelian $H\rtimes \Gamma$-group $A$, we write  
$ H^2( H \rtimes \Gamma, A)^{\cL}$ for the subset of classes in $ H^2( H \rtimes \Gamma, A)$ whose corresponding $\Gamma$-group extension $1\ra A\ra G\ra H\ra 1$ via Lemma~\ref{L:diffext} has $G\in \cL$.
We write $H^2(H\rtimes \Gamma,A)^{\cL,s_H}$ for the subset of classes $\alpha\in H^2(H\rtimes \Gamma,A)^{\cL}$ whose corresponding $\Gamma$-extension
$G$ 
satisfies
\begin{itemize}
\item $s_H(\alpha\cup \beta)=0$ for all $\beta\in H^1(H\rtimes \Gamma, \Hom(A,\Z/n))$, and
\item for any $H\rtimes \Gamma$-equivariant homomorphism $A\ra \Z/n$,  we have that the image of $\alpha$ in
$H^2(H\rtimes \Gamma,\Z/n)$ is in the kernel of the
composite of
 Bockstein map $\mathcal{B}:H^2(H\rtimes \Gamma,\Z/n) \ra H^3(H\rtimes \Gamma,\Z/n)$ given by the connecting map in the long exact sequence for $0\ra \Z/n \ra  \Z/n^2 \ra \Z/n\ra 0$, and $s_H: H^3(H\rtimes \Gamma,\Z/n) \ra \Z/n$ (viewed as a homomorphism using Lemmas~\ref{L:GammaH} and \ref{L:homdual}).
\end{itemize}
\end{definition}

It is straightforward to check that this definition agrees with the definition in the introduction.  
After taking the composition with the map $H^2(H\rtimes \Gamma,\F_\ell) \stackrel{\times n/\ell}{\ra } H^2(H\rtimes \Gamma,\Z/n)$, the Bockstein msp
from $0\ra \Z/n \ra  \Z/n^2 \ra \Z/n\ra 0$ agrees with the Bockstein map from $0\ra \Z/\ell \ra  \Z/n\ell \ra \Z/\ell\ra 0$.

\begin{lemma}\label{L:shcriterion}
For $(G,\pi)\in \Ext_\Gamma(H,V_i^{e_i})$, we have
$G\in\cL$ and 
$s_H \in \im \pi$ 
 if and only if  
$(G,\pi)$
corresponds via Lemma~\ref{L:diffext} to a class in $H^2(H\rtimes \Gamma,V_i^{e_i})^{\cL,s_H}$.
\end{lemma}
\begin{proof}
By Lemma~\ref{L:sH}, we have $s_H \in \im \pi$  if and only if $s_H \circ d_{2}^{1,1}=0$ and $s_H \circ d_{3}^{0,2}=0$. 
The map $d_2^{1,1} :H^1(H\rtimes \Gamma, \Hom(F,\Z/n))\ra H^3(H\rtimes \Gamma,\Z/n\Z)$ is given by the cup product with the class of the extension in $H^2(H\rtimes \Gamma,F^{ab})$ \cite[Theorem 2.4.4]{Neukirch2000}.
So $s_H\circ d_2^{1,1}=0$ if and only if
$s_H(\alpha\cup \beta)=0$ for all $\beta\in H^1(H\rtimes \Gamma, \Hom(F,\Z/n))$. 

Lemma~\ref{L:H2struc} defines a homomorphism $B:\Hom(V_i^{e_i},\Z/n)^{H\rtimes\Gamma}  \ra H^2(V_i^{e_i},\Z/n)^{H\rtimes\Gamma}$.
Lemma \ref{L:d202gen} says that
 $d_3^{0,2}([B(\phi)])$ is the image of $[\alpha]\in H^2(H\rtimes \Gamma,V_i^{e_i})$ under $\phi_*:H^2(H\rtimes \Gamma,V_i^{e_i}) \ra H^2(H\rtimes \Gamma,\Z/n)$ and then $\mathcal{B}$. By Lemma~\ref{P:dforbilinear4},  $\im d_{3}^{0,2}=d_3^{0,2}(\im B)$.   Thus  $s_H\circ d_3^{0,2}=0$ if and only if the second bullet point above holds for the extension class.
\end{proof}

Lemma~\ref{L:basicH2} tells us that $ H^2( H \rtimes \Gamma, V_i^{e_i})^{\cL}$ is a $\kappa_i$-subspace of  $H^2( H \rtimes \Gamma, V_i^{e_i})$ and
that $ H^2( H \rtimes \Gamma, V_i^{e_i})^{\cL}=(H^2( H \rtimes \Gamma, V_i)^{\cL})^{e_i}.$  
Further, we can check that $ H^2( H \rtimes \Gamma, V_i^{e_i})^{\cL,s_H}$ is a $\kappa_i$-subspace of  $H^2( H \rtimes \Gamma, V_i^{e_i})$ and
that $ H^2( H \rtimes \Gamma, V_i^{e_i})^{\cL,s_H}=(H^2( H \rtimes \Gamma, V_i)^{\cL,s_H})^{e_i}.$  

Also, we have that,  the composite of $H^2(A,\F_p)\ra H^2(A,\Z/n)$ with $\mathcal{B}$  is the Bockstein map
$H^2(A,\F_p)\ra H^3(A,\Z/n)$  for $0\ra \Z/n \ra  \Z/np \ra \F_p\ra 0.$
 
We define an adjusted $q$-Pochhammer symbol $(q)_n:=\prod_{i=1}^n (1-q^{-i})$ for use in the next several sections.  By convention $(q)_0=1$.

 \subsection{Representations of characteristic prime to $n$} 
 
Note if $F=V_i^{e_i},$ where $V_i$ is of characteristic prime to $n$, then
$H^p(H\rtimes \Gamma,H^q(F,\Z/n))=0$ for $q>0$.  In particular, this implies that $|E_3^{0,2}|=1$ and
$\pi^*: H^3(H\rtimes \Gamma,\Z/n)\ra H^3(G\rtimes \Gamma,\Z/n)$ is an isomorphism and hence by Lemma~\ref{L:GammaH},
$\pi_*: H_3(G,\Z/n)^\Gamma\ra H_3(H,\Z/n)^\Gamma$ is an isomorphism. 

 \begin{lemma}\label{prime-n-nontrivial-factor} 
Let $\cL$, $\cW$, and $\bH$ be as at the start of Section~\ref{S:getv},  and further assume $H_\Gamma=1$.
 If $\chr(V_i) \nmid n$ and $V_i$ is nontrivial, \[ M_i (e_i)= 
\frac{1 }{ \abs{GL_{e_i}(q_i)}} 
 \left( \frac{ | H^2( H \rtimes \Gamma, V_i)^{\cL} |}{  |H^1( H \rtimes \Gamma, V_i )| V_i^{\cdot U'}  } \right) ^{e_i}  .\] \end{lemma}

 \begin{proof}  
 By definition, 
 \[ M_i (e_i) = \frac{ 1}{  \abs{GL_{e_i}(q_i) }|H^1 ( H \rtimes \Gamma, V_i)|^{e_i}  (V_i^{\cdot U'})^{ e_i } }  \sum_{\substack{
(G,\pi)\in \Ex(\bH,V_i^{e_i})}}    | E_{3}^{0,2} |  ,\] 
where $\Ex(\bH,V_i^{e_i})$ is the set of $(G,\pi)\in \Ext_\Gamma(H,V_i^{e_i})$ such that $G\in\cL$,  and $G_\Gamma=1$, and $s_H\in \im \pi_*$.  
By Lemma~\ref{L:GGamma}, we have $G_\Gamma=1$.
Since $\pi_*$ is an isomorphism, we always have $s_H\in\im\pi_*$. 
 
We then have
\begin{align*}
 M_i (e_i) &= \frac{  1}{  \abs{GL_{e_i}(q_i) }  |H^1 ( H \rtimes \Gamma, V_i)|^{e_i} (V_i^{\cdot U'})^{e_i}}  \sum_{\substack{
(G,\pi)\in \Ex(\bH,V_i^{e_i})}}    | E_{3}^{0,2} | \\
&= \frac{  1 }{  \abs{GL_{e_i}(q_i) } \ |H^1 ( H \rtimes \Gamma, V_i)|^{e_i} (V_i^{\cdot U'})^{e_i}}  \sum_{\substack{
(G,\pi)\in \Ex(\bH,V_i^{e_i})}}    1 \\
\\ 
& =\frac{1 }{ \abs{GL_{e_i}(q_i)}} \left( \frac{ | H^2( H \rtimes \Gamma, V_i)^{\cL} |}{  |H^1( H \rtimes \Gamma, V_i )| V_i^{\cdot U'} } \right) ^{e_i}  .
 \end{align*}
 \end{proof}

\begin{lemma}\label{prime-n-nontrivial-w} 
Let $\cL$, $\cW$, and $\bH$ be as at the start of Section~\ref{S:getv},  and further assume $H_\Gamma=1$.
If $\chr(V_i) \nmid n$ and $V_i$ is nontrivial, 
$$
w_{V_i} = \prod_{j=1}^\infty \left(1-q_i^{-j} \frac{ | H^2( H \rtimes \Gamma, V_i)^{\cL} |}{  |H^1( H \rtimes \Gamma, V_i )| V_i^{\cdot U'}  } 
 \right). 
$$
Further, $\tilde{w}_{V_i}$ is finite.
\end{lemma}

\begin{proof} 
We have by definition and Lemma \ref{prime-n-nontrivial-factor},
\begin{align*}
w_{V_i} =&\sum_{e_i\geq 0}(-1)^{e_i} q_i^{\binom{e_i}{2}} 
 M_i ( e_i)\\
&=\sum_{e_i\geq 0}(-1)^{e_i} q_i^{\binom{e_i}{2}} 
\frac{1 }{ \abs{GL_{e_i}(q_i)}} 
 \left( \frac{ | H^2( H \rtimes \Gamma, V_i)^{\cL} |}{  |H^1( H \rtimes \Gamma, V_i )| V_i^{\cdot U'}  } \right) ^{e_i}.
\end{align*}
Applying the infinite $q$-binomial theorem,
\begin{equation}\label{E:qbinomP}
\prod_{j=1}^\infty (1+q^{-j}v) 
=\sum_{k=0}^\infty \frac{q^{\binom{k}{2}-k^2} v^k}{(q)_k},
\end{equation}
gives the first statement of the lemma.  An analogous calculation without the sign gives that $\tilde{w}_{V_i}$ is finite.
\end{proof} 
 
 \begin{lemma}\label{prime-n-trivial-factor}
Let $\cL$, $\cW$, and $\bH$ be as at the start of Section~\ref{S:getv},  and further assume $H_\Gamma=1$.
 If {$\chr(V_i) \nmid n$} and $V_i$ is trivial, \[M_i (e_i) =
\frac{ |\Surj( H^2( H \rtimes \Gamma, V_i)^{\cL}, {V_i}^{e_i} )|  }{ \abs{GL_{e_i}(q_i)} (V_i^{\cdot U'})^{e_i}} 
.\] \end{lemma}
 
 \begin{proof}  
 The proof is similar to that of Lemma~\ref{prime-n-nontrivial-factor}. 
One can check that $G_\Gamma\ne 1$ if and only if there is no proper sub-representation  $W\sub V_i^{e_i} $,  such that 
$G/W$ is a trivial extension of $H$.  Such extensions $G$ are those that correspond to elements of $H^2(H\rtimes \Gamma, V_i^{e_i})=H^2(H\rtimes \Gamma, V_i)^{e_i}$ whose $e_i$ coordinates in $H^2(H\rtimes \Gamma, V_i)$ are linearly independent.  The number of such extensions in $\cL$ is hence the number of surjections $|\Surj( H^2( H \rtimes \Gamma, V_i)^{\cL}, {V_i}^{e_i} )|$ (as abelian groups). 
 
 Since $H_{\Gamma}=1$, we have $H^1(H\rtimes \Gamma, V_i)=\Hom(H\rtimes \Gamma, V_i)=0$.
 Putting this together,
 we conclude the lemma.
 \end{proof}

\begin{lemma}\label{L:wtrivial}
Let $\cL$, $\cW$, and $\bH$ be as at the start of Section~\ref{S:getv},  and further assume $H_\Gamma=1$.
Let $z=\dim_{\kappa_i} H^2( H \rtimes \Gamma, V_i)^{\cL}$ and $u$ be such that $q_i^u=V_i^{\cdot U'}$.
If $V_i$ is trivial, and $\chr(V_i)\nmid n$, then
\[ w_{V_i} =
\prod_{k=0}^{z-1} (1-q_i^{k-u}).
\]
Further, $\tilde{w}_{V_i}$ is finite.
If $U$ is nonempty, this expression is $0$ if $z\geq u+1$ and positive otherwise.  
\end{lemma}
\begin{proof}
Let $q=|V_i|$.
By definition and Lemma~\ref{prime-n-trivial-factor}, we have
\begin{align*}
w_{V_i}&=\sum_{e_i\geq 0}(-1)^{e_i} q_i^{\binom{e_i}{2}}
\frac{|\Surj( H^2( H \rtimes \Gamma, V_i)^{\cL,s_H}, {V_i}^{e_i} )| }{ \abs{GL_{e_i}(q_i) } (V_i^{\cdot U})^{e_i} \    }\\
&=\sum_{e=0}^{z} (-1)^{e} q^{\binom{e}{2}}
\frac{q^{ze}(q)_z/(q)_{z-e} }{ q^{e^2+ue}(q)_e   }\\
&=\sum_{e=0}^{z} q^{\binom{e}{2}+ze-e^2} (-q^{-u})^e
\frac{(q)_z}{ (q)_{z-e}(q)_e   }\\
&=\prod_{k=0}^{z-1} (1-q^{k-u}),
\end{align*}
where the last equality is by the $q$-binomial theorem
\begin{equation}\label{E:binomialfinP}
 \prod _{k=0}^{n-1 }{(1+q^{k}v)}
=\sum_{k=0}^{n}  q^{\binom{k}{2}+nk -k^2} v^k
\frac{ (q)_n  } { (q)_{n-k}(q)_k }.
\end{equation} 
An analogous calculation without the sign gives that $\tilde{w}_{V_i}$ is finite.

If $U$ is nonempty then $V_i^{\cdot U'} \geq 1 $ so $u\geq 0$. Furthermore $u$ is an integer. The expression $(1-q_i^{k-u})=0$ if $k= u$, which happens for some $k$ if  $z \geq u+1$, and positive if $k<u$, which happens for all $k$ if $z< u+1$.
\end{proof}

\subsection{Representations whose duals do not appear of characteristic dividing $n$}

\begin{lemma}\label{L:nodual-M}
Let $\cL$, $\cW$, and $\bH$ be as at the start of Section~\ref{S:getv},  and further assume $H_\Gamma=1$.
 If $V_i\not\isom V_j^\vee$ for any $j$, and $\chr(V_i)\mid n$, then
\[ M_i (e_i) =
\frac{1}{ \abs{GL_{e_i}(q_i) }} \left(
  \frac{|H^1 ( H \rtimes \Gamma, V_i^\vee)|  |H^2(H\rtimes \Gamma,V_i)^{\cL,s_H}| }{
|H^1 ( H \rtimes \Gamma, V_i)|
 V_i^{\cdot U'}  }\right)^{e_i}   .\]
\end{lemma}
\begin{proof}
By definition, 
\[ M_i (e_i) =  \frac{|H^1 ( H \rtimes \Gamma, V_i^\vee)|^{e_i} }{ \abs{GL_{e_i}(q_i) }
|H^1 ( H \rtimes \Gamma, V_i)|^{e_i}
  (V_i^{\cdot U'})^{ e_i }  }  \sum_{\substack{
(G,\pi)\in \Ex(\bH,V_i^{e_i})}}    | E_{3}^{0,2} |  ,\]
where $\Ex(\bH,V_i^{e_i})$ is the set of $(G,\pi)\in \Ext_\Gamma(H,V_i^{e_i})$ such that $G\in\cL$,  and $G_\Gamma=1$, and $s_H\in \im \pi_*$.  

Lemma~\ref{L:H2struc} gives an exact sequence $0\ra \im B \ra H^2(V_i^{e_i},\Z/n)\ra \wedge_2 (V_i^*)^{e_i}\ra 0$,
where $B$ is a homomorphism from $\Hom(V_i^{e_i},\Z/n)$ and $\wedge_2 A$ is defined as the subgroup of $A\tensor A$ generated by elements of the form
$a\tensor b-b\tensor a$ for $a,b\in A$.

Since $V_i$ is not-self dual,  we have $(\wedge_2 (V_i^*)^{e_i})^{H\rtimes \Gamma}
\subseteq ( ((V_i^*)^{e_i})^{\tensor 2})^{H \rtimes \Gamma}=0.$ 
This, with $V_i$ nontrivial and so $(\im B)^{H\rtimes \Gamma}=0$,
implies that $H^2(V_i^{e_i},\Z/n)^{H\rtimes \Gamma}=0$.
Thus $|E_3^{0,2}|=1$.  

For $(G,\pi)\in \Ext_\Gamma(H,V_i^{e_i})$, by Lemma~\ref{L:GGamma} and the fact that $V_i$ is not self-dual and thus not trivial, we have $G_\Gamma=1$.
The lemma then follows from Lemma~\ref{L:shcriterion}.
\end{proof}

\begin{lemma}\label{L:nodual-w} 
Let $\cL$, $\cW$, and $\bH$ be as at the start of Section~\ref{S:getv},  and further assume $H_\Gamma=1$. If $V_i\not\isom V_j^\vee$ for any $j$, and $\chr(V_i)\mid n$, then
$$
w_{V_i} = \prod_{j=1}^\infty \left(1-q_i^{-j}
\left(
  \frac{|H^1 ( H \rtimes \Gamma, V_i^\vee)|  |H^2(H\rtimes \Gamma,V_i)^{\cL,s_H}| }{
|H^1 ( H \rtimes \Gamma, V_i)|
  V_i^{\cdot U'}  }\right)
 \right). 
$$
Further, $\tilde{w}_{V_i}$ is finite.
\end{lemma}

\begin{proof} 
The lemma follows by definition, Lemma~\ref{L:nodual-M}, and the $q$-binomial theorem, as in the proof of Lemma~\ref{prime-n-nontrivial-w}.
\end{proof}

 \subsection{Non-anomalous self-dual representations of characteristic dividing $n$}

\begin{lemma}\label{L:Moddselfnontrivial}
Let $\cL$, $\cW$, and $\bH$ be as at the start of Section~\ref{S:getv},  and further assume $H_\Gamma=1$.
 If $V_i$ is self-dual, nontrivial,  non-anomalous, and $\chr(V_i)\mid n$, then
\[ M_i (e_i) =
\frac{1 }{ \abs{GL_{e_i}(q_i) }\  (V_i^{\cdot U'})^{ e_i }  }  \sum_{\substack{
\alpha \in (H^2(H\rtimes \Gamma,V_i)^{\cL,s_H})^{e_i} }}    q_i^{\frac{(e_i-r_\alpha)(e_i-r_\alpha-\epsilon_{V_i})}{2}}  
,\]
where $r_\alpha:=\dim_{\kappa_i} \im\alpha$ in the sense of Definition~\ref{D:im}.
\end{lemma}
\begin{proof}
By definition, 
\[ M_i (e_i) =  \frac{1 }{ \abs{GL_{e_i}(q_i) }\ (V_i^{\cdot U'})^{ e_i }  \abs{ V_i^{ H \rtimes \Gamma}}^{e_i}   }  \sum_{\substack{
(G,\pi)\in \Ex(\bH,V_i^{e_i})}}    | E_{3}^{0,2} |  ,\] where
 $\Ex(\bH,V_i^{e_i})$ is the set of
  $(G,\pi)\in \Ext_\Gamma(H,V_i^{e_i})$ such that $G\in\cL$,  and $G_\Gamma=1$, and $s_H\in \im \pi_*$.  
By Lemma~\ref{L:GGamma}, we have $G_\Gamma=1$ since $V_i$ is nontrivial. Also since $V_i$ is nontrivial we have $\abs{ V_i^{ H \rtimes \Gamma}}^{e_i} =1$.
By Proposition~\ref{P:dforbilinear4}, we have $|E_3^{0,2}|=|\kappa_i|^{\frac{(e_i-r_\alpha)(e_i-r_\alpha-\epsilon_{V_i})}{2}} $.
The lemma then follows from Lemma~\ref{L:shcriterion}.
\end{proof}

 \begin{lemma}\label{odd-selfdual-nontrivial-w} 
Let $\cL$, $\cW$, and $\bH$ be as at the start of Section~\ref{S:getv},  and further assume $H_\Gamma=1$.  Let $Z=H^2 ( H \rtimes \Gamma, V_i)^{\cL, s_H}$, and $z=\dim_{\kappa_i} Z$, and  
let $u$ be such that $q_i^u=V_i^{\cdot U'}.$
Assume that $V_i$ contains a nontrivial vector fixed by at least one element of $U$. 
 If $\chr(V_i) \mid n$, and $V_i$ is self-dual,  nontrivial, and non-anomalous, 
 then
 \[ w_{V_i} \oftau  =
 \prod _{k=0}^{\infty }(1+q_i^{-k-\frac{\epsilon_{V_i}+1}{2} -u})^{-1}
 \prod _{k=0}^{z-1 }{(1-q_i^{k -u})}.
 \] 
Further, $\tilde{w}_{V_i}$ is finite. 
If $z\geq u+1$, then $w_{V_i}=0$, and otherwise $w_{V_i}$ is positive.  
 \end{lemma}

\begin{proof} 
Let $\kappa=\kappa_{V_i}$, and $q=|\kappa|$, and $\epsilon=\epsilon_{V_i}$.
We are going to interchange an order of summation, which will be justified later by considering the same manipulations if no signs appeared to see that the sum is absolutely convergent.
Letting $e=e_i-r_\alpha$, by definition and Lemma~\ref{L:Moddselfnontrivial}
 \begin{align*} 
w_{V_i}=&
 \sum_{e_i=0}^{\infty} \frac{ (-1)^{e_i} q^{\binom{e_i}{2}}}{  \abs{GL_{e_i}(q) }\  (V_i^{\cdot U'})^{ e_i }  } \sum_{ \alpha \in (H^2 ( H \rtimes \Gamma, V_i)^{\cL, s_H} )^{e_i} }   q^{ \frac{ (e_i - r_\alpha) (e_{i} - r_\alpha - \epsilon) }{2}  } \\
 =& \sum_{e=0}^{\infty} \sum_{r =0}^{\infty}  
 \sum_{\substack{ \alpha \in (H^2 ( H \rtimes \Gamma, V_i)^{\cL, s_H} )^{e+r}\\ r_\alpha=r }  }
  \frac{ (-1)^{e+r} 
q^{ -\binom{e+r+1}{2}+ \frac{e(e- \epsilon) }{2} -u(e+r) }  
  }{  (q)_{e+r}   }. 
 \end{align*}
The number of $\alpha\in Z^{e+r}$ with $r_\alpha=r$,  is the number  
of $r$-dimensional subspaces of $\kappa^{e+r}$, which is 
$\frac{q^{er} (q)_{e+r}  } {  (q)_e (q)_r } $,
 times the number of surjections from the $\kappa$-dual of $Z$ to an $r$-dimensional $\kappa$ vector space, which is 
$\frac{q^{zr} (q)_z } {  (q)_{z-r}  }$
when $r\leq z$ and $0$ otherwise.
We then have
\begin{align*}
w_{V_i}
&=\sum_{e=0}^{\infty} \sum_{r =0}^{z} 
 \frac{q^{er} (q)_{e+r}  } {  (q)_e (q)_r } \frac{ q^{zr}(q)_z } {  (q)_{z-r}  }  \frac{ (-1)^{e+r} 
q^{ -\binom{e+r+1}{2}+ \frac{e(e- \epsilon) }{2} -u(e+r) }  
  }{  (q)_{e+r}   }.\\
 &=\sum_{e=0}^{\infty} \frac{(-1)^e q^{ -\frac{(\epsilon+1)e}{2}-ue}}{ (q)_e } \sum_{r =0}^{z}
(-1)^r q^{\binom{r}{2}-ur +zr-r^2}\frac{ (q)_z  } {  (q)_{z-r}  (q)_r }.
\end{align*}
Note that our assumption on $U$ implies that $V_i^{\cdot U'}>1$ so that $u>0$. By the same calculation without the signs, since $u>0$, we see that $\tilde{w}_{V_i}$ is finite and the sum above is absolutely convergent.
We 
 apply the infinite $q$-binomial theorem for negative powers
\begin{equation}\label{E:binomialnegP}
 \prod _{k=0}^{\infty }{\frac {1}{1-q^{-k}v}}
=\sum_{k=0}^{\infty }\frac {v^{k}}{(q)_k}
\end{equation} 
to the sum over $e$
and 
 the  $q$-binomial theorem \eqref{E:binomialfinP}
to the sum over $r$, and obtain
\begin{align*}
w_{V_i}=
 & \prod _{k=0}^{\infty }(1+q^{-k-\frac{\epsilon+1}{2} -u})^{-1}
 \prod _{k=0}^{z-1 }{(1-q^{k -u})}.
\end{align*}
Using that $u>0$, the product on the right is $0$ if $z-1-u\geq 0$, and is otherwise positive.
 \end{proof}

\begin{lemma}\label{L:Mntrivial}
Let $\cL$, $\cW$, and $\bH$ be as at the start of Section~\ref{S:getv},  and further assume $H_\Gamma=1$. If $V_i$ is trivial, and $\chr(V_i)\mid n$, then
\[ M_i (e_i) =
\frac{|\Surj( H^2( H \rtimes \Gamma, V_i)^{\cL,s_H}, {V_i}^{e_i} )| }{ \abs{GL_{e_i}(q_i) }(V_i^{\cdot U'})^{e_i} \   }.
\]
\end{lemma}
\begin{proof}
By definition, 
\[ M_i (e_i) =  \frac{1 }{ \abs{GL_{e_i}(q_i) }\  (V_i^{\cdot U'})^{ e_i }  \abs{V_i^{H \rtimes \Gamma}}^{e_i} }  \sum_{\substack{
(G,\pi)\in \Ex(\bH,V_i^{e_i})}}    | E_{3}^{0,2} |  ,\] where
 $\Ex(\bH,V_i^{e_i})$ is the set of
  $(G,\pi)\in \Ext_\Gamma(H,V_i^{e_i})$ such that $G\in\cL$,  and $G_\Gamma=1$, and $s_H\in \im \pi_*$.  
As in the proof of Lemma~\ref{prime-n-trivial-factor},  there are $|\Surj( H^2( H \rtimes \Gamma, V_i)^{\cL,s_H}, {V_i}^{e_i} )|$
choices of such $(G,\pi)$ that have $G_\Gamma=1$, and all of these have $r_\alpha=e_i$.
By Proposition~\ref{P:dforbilinear4}, we have $|E_3^{0,2}|=|V_i|^{e_i}|\kappa_i|^{\frac{(e_i-r_\alpha)(e_i-r_\alpha-\epsilon_{V_i})}{2}} =|V_i|^{e_i} = \abs{V_i^{H \rtimes \Gamma}}^{e_i}$.
\end{proof}

\begin{lemma}\label{L:wntrivial}
Let $\cL$, $\cW$, and $\bH$ be as at the start of Section~\ref{S:getv},  and further assume $H_\Gamma=1$. 
Let $z=\dim_{\kappa_i} H^2( H \rtimes \Gamma, V_i)^{\cL,s_H}$ and $u$ be such that $q_i^u=V_i^{\cdot U'}$.
If $V_i$ is trivial, and $\chr(V_i)\mid n$, then
\[ w_{V_i} =
\prod_{k=0}^{z-1} (1-q_i^{k-u}),
\]
which is $0$ if $z\geq u+1$ and positive otherwise.  Further, $\tilde{w}_{V_i}$ is finite.
If $U$ is nonempty, this expression is $0$ if $z\geq u+1$ and positive otherwise.  
\end{lemma}
\begin{proof}
The proof is the same as the proof of Lemma~\ref{L:wtrivial}.
\end{proof}

 \subsection{Non-self-dual representations whose duals appear of characteristic dividing $n$} 

\begin{lemma}\label{L:Modddualpairs}
Let $\cL$, $\cW$, and $\bH$ be as at the start of Section~\ref{S:getv},  and further assume $H_\Gamma=1$.
 If $V_i$ and $V_{i'}$ are dual, with $i\ne j$, 
and $\chr(V_i)\mid n$, then
\[ M_{i,i'} (e_i,e_{i'}) =
 \frac{ |H^1 ( H \rtimes \Gamma, V_i)|^{e_{i'}-e_i}  |H^1 ( H \rtimes \Gamma, V_{i'})|^{e_{i} -e_{i'} }  }{  
   \abs{GL_{e_i}(q_i) } \abs{GL_{e_{i'}}(q_{i'}) } (V_i^{\cdot U'})^{ e_i }   (V_{i'}^{\cdot U'})^{ e_{i'} }
}  
  \sum_{
(\alpha_i,\alpha_{i'})\in Z} q_i^{(e_i-r_{\alpha_i})(e_{i'}-r_{\alpha_{i'}})}
,\]
where 
$Z= (H^2(H\rtimes \Gamma,V_i)^{\cL,s_H})^{e_i}\times (H^2(H\rtimes \Gamma,V_{i'})^{\cL,s_H})^{e_{i'}}$ and 
$r_\alpha:=\dim_\kappa \im\alpha$ in the sense of Definition~\ref{D:im}. 
\end{lemma}
\begin{proof}
Let $F=V_i^{e_i} \times V_{i'} ^{e_{i'}}$.
Recall the definition,
\[ M_{i,i'} (e_i, e_{i'} ) = \frac{ |H^1 ( H \rtimes \Gamma, V_i)|^{e_{i'}-e_i}  |H^1 ( H \rtimes \Gamma, V_{i'})|^{e_{i} -e_{i'} }  }{  
   \abs{GL_{e_i}(q_i) } \abs{GL_{e_{i'}}(q_{i'}) } (V_i^{\cdot U'})^{ e_i }   (V_{i'}^{\cdot U'})^{ e_{i'} }
}  
  \sum_{
(G,\pi)\in \Ex(\bH,F)} | E_3^{0,2} |,\] 
where $ \Ex(\bH,F)$ is the set of
  $(G,\pi)\in \Ext_\Gamma(H,F)$ such that $G\in\cL$,  and $G_\Gamma=1$, and $s_H\in \im \pi_*$.  
By Lemma~\ref{L:GGamma}, we have $G_\Gamma=1$ since $V_i$ and $V_j$ are nontrivial.
By Lemma~\ref{L:sH}, $s_H\in \im \pi_*$ if and only if $s_H \circ d_{2}^{1,1}=0$ and $s_H \circ d_{3}^{0,2}=0$.
By Proposition~\ref{P:dfordualpairsgen}, $d_{3}^{0,2}=0$ and $|E_3^{0,2}|=|\kappa_i|^{(e_i-r_{\alpha_i})(e_{i'}-r_{\alpha_{i'}})} $.
The classes $(G,\pi)\in \Ext_\Gamma(H,F)$ such that $G\in\cL$ and $s_H \circ d_{2}^{1,1}=0$
correspond via Lemma~\ref{L:diffext} to classes in $H^2(H\rtimes \Gamma,V_i^{e_i})^{\cL,s_H}$.
\end{proof}

\begin{lemma}\label{odd-paired-w}
Let $\cL$, $\cW$, and $\bH$ be as at the start of Section~\ref{S:getv},  and further assume $H_\Gamma=1$.
 Let $\chr(V_i)\mid n$ and $V_i \cong V_{i'}^\vee $ but $i\neq i'$.
 Let 
 $u$ be such that $q^u=V_i^{\cdot U'}$
 and $z=\dim_{\kappa_{V_i}} H^2(H\rtimes \Gamma,V_i)^{\cL,s_H}$, and
 $z'=\dim_{\kappa_{V_i}} H^2(H\rtimes \Gamma,V_{i'})^{\cL,s_H}$.
 Let $h=\dim_\kappa H^1(H\rtimes\Gamma,V_{i'}) - \dim_\kappa H^1(H\rtimes\Gamma,V_i)$. 
  Then
\[ w_{V_i}\oftau  w_{V_{i'}}\oftau  =
\frac{(q_i)_\infty(q_i)_{2u}}{(q_i)_{u-h-z}(q_i)_{u+h-z'}} 
,\]
if $ -u+z'\leq h \leq u-z$ and $0$ otherwise.
Further, $\tilde{w}_{V_i}\tilde{w}_{V_{i'}}$ is finite. 
\end{lemma}

\begin{proof} 
Let $\kappa=\kappa_{V_i}$, and $q=|\kappa|$.
Let $Z= (H^2(H\rtimes \Gamma,V_i)^{\cL,s_H})^{e_i}\times (H^2(H\rtimes \Gamma,V_{i'})^{\cL,s_H})^{e_{i'}}$.
Since the characteristic of $V_i$ does not divide $|\Gamma|,$ for any subgroups $\Gamma'$ of $\Gamma$, we have $V_i$ and $V_{i'}$ 
are products of irreducible $\Gamma'$ representations, and $V_i^{\Gamma'}$ and $V_{i'}^{\Gamma'}$ are comprised of the trivial factors and hence the same size.  Thus $V_i^{\cdot U'}=(V_i^\vee) ^{\cdot U'}.$

We have by definition and Lemma~\ref{L:Modddualpairs},
\begin{align*}
w_{V_i}w_{V_{i'}}&=\sum_{e_i, e_{i'}\geq 0}  (-1)^{e_i+e_{i'}} 
  \frac{ q^{\binom{e_i}{2}+\binom{e_{i'}}{2} +h(e_{i}-e_{i'} )  -u(e_i+e_{i'})}   }{  
   \abs{GL_{e_i}(q_i) } \abs{GL_{e_{i'}}(q_{i'}) } } 
    \sum_{
(\alpha_i,\alpha_{i'})\in Z} q_i^{(e_i-r_{\alpha_i})(e_{i'}-r_{\alpha_{i'}})}
 .
\end{align*}
As in Lemma~\ref{odd-selfdual-nontrivial-w}, the number of $\alpha_i\in (H^2(H\rtimes \Gamma,V_i)^{\cL,s_H})^{e_i}$
with $r_{\alpha_i}=r$ is 
$$
\frac{q^{(e_i-r)r+zr}(q)_{e_i}(q)_{z}}{(q)_{e_i-r}(q)_{r}(q)_{z-r}}
$$
if $r\leq z$ and $r\leq e_i$,
and similarly for $i'$.
So
\begin{align*}
w_{V_i}w_{V_{i'}}=& \sum_{e_i, e_{i'}\geq 0}  (-1)^{e_i+e_{i'}} 
  \frac{ q^{\binom{e_i}{2}-e_i^2+\binom{e_{i'}}{2} -e_{i'}^2+h(e_{i} - e_{i'} )  -u(e_i+e_{i'})}   }{  
    (q)_{e_i}(q)_{e_{i'}}}   \times
\\
&     \sum_{r=0}^{\min(e_i,z)}\sum_{r'=0}^{\min(e_{i'},z')}  
\frac{q^{(e_i-r)r+zr}(q)_{e_i}(q)_{z}}{(q)_{e_i-r}(q)_{r}(q)_{z-r}}
    \frac{q^{(e_{i'}-r')r'+z'r'}(q)_{e_{i'}}(q)_{z'}}{(q)_{e_{i'}-r'}(q)_{r'}(q)_{z'-r'}}
    q^{(e_i-r_{i})(e_{i'}-r_{{i'}})}
 .
\end{align*}
We are going to interchange an order of summation, which will be justified later by considering the same manipulations if no signs appeared to see that the sum is absolutely convergent. 
Letting $e=e_i-r$ and $e'=e_{i'}-r'$,
\begin{align*}
&w_{V_i}w_{V_{i'}}=\\
& \sum_{e,e'\geq 0}  \sum_{r=0}^{z}\sum_{r'=0}^{z'}  (-1)^{e+e'+r+r'} 
  \frac{ q^{\binom{e}{2}-e^2+\binom{e'}{2} -(e')^2+h(e-e')  -u(e+e') +ee' -u(r+r')  +h(r-r')
 + \binom{r}{2}-r^2+\binom{r'}{2} -(r')^2
  +zr +z'r'  
  } (q)_{z} (q)_{z'} }{  
   (q)_{e}(q)_{e'}(q)_{r}(q)_{z-r}(q)_{r'}(q)_{z'-r'} }
 .
\end{align*}
We can factor out
$$
 \sum_{r=0}^{z}
  (-1)^{r} 
  \frac{ q^{hr-ur
 + \binom{r}{2}-r^2
  +zr
  } (q)_{z}  }{  
 (q)_{r}(q)_{z-r} }= \prod_{k=0}^{z-1} (1-q^{k+h-u}) ,
$$
with the equality by the $q$-binomial theorem \eqref{E:binomialfinP}, as well as the analogous sum over $r'$ (which has a $-$  in front of the $h$).  What remains is
\begin{equation}\label{E:justes}
\sum_{e,e'\geq 0} (-1)^{e+e'} 
  \frac{ q^{\binom{e}{2}-e^2+\binom{e'}{2} -(e')^2+h(e-e')  -u(e+e') +ee' 
  } }{  
   (q)_{e}(q)_{e'}}.
\end{equation}

We now argue that the sums we have been considering converge absolutely,  equivalently that $\tilde{w}_{V_i}\tilde{w}_{V_{i'}}$ is finite. 
If we had been considering the same sums without signs, the above argument would reduce the finiteness of $\tilde{w}_{V_i}\tilde{w}_{V_{i'}}$  to the absolute convergence of \eqref{E:justes}.  If $f$ is a quadratic polynomial with 
negative leading coefficient and $f(\Z)\sub \Z$, then $\sum_{e} q^{f(e)}\leq 4q^m$, where $m$ is the maximum value
taken by $f$.  This is because the sum on either side of the maximum is bounded by a geometric series summing to $q^m(1-q^{-1})$, and $q\geq 2$.  
\details{The maximum values of $ax^2+bx+c$ is $-b^2/(4a)+c$}
Thus
\begin{align*}
\sum_{e,e'\geq 0} 
  \frac{ q^{\binom{e}{2}-e^2+\binom{e'}{2} -(e')^2+h(e-e')  -u(e+e') +ee' 
  } }{  
   (q)_{e}(q)_{e'}}
   &\leq
   \frac{4}{(2)_\infty^2}  \sum_{e,e'\geq 0} 
q^{-\frac{e^2}{2} +(- \frac{1}{2}+h -u+e' )e +(-h-u- \frac{1}{2})e'     -\frac{(e')^2}{2}}    \\
   &\leq
   \frac{4}{(2)_\infty^2}  \sum_{e'\geq 0} 
q^{ \frac{(- \frac{1}{2}+h -u+e' )^2}{2} +(-h-u- \frac{1}{2})e'     -\frac{(e')^2}{2}}.    
\end{align*}
The last sum above is a geometric series with ratio $q^{-2u-1}$ and hence converges, showing that
the sums we have been considering converge absolutely,  and $\tilde{w}_{V_i}\tilde{w}_{V_{i'}}$ is finite. 

Returning to \eqref{E:justes}, we pull out the sum
\begin{equation}\label{E:onlye}
\sum_{e\geq 0} (-1)^{e} 
  \frac{ q^{\binom{e}{2}-e^2+he  -ue +ee' 
  } }{  
   (q)_{e}}=\prod_{j=1}^\infty (1-q^{-j+h-u+e'})=
   \begin{cases}
   0 & \textrm{ if } h-u+e'> 0 \\
   \frac{(q)_\infty}{(q)_{u-h-e'}} & \textrm{ if } h-u+e'\leq 0
   \end{cases},
\end{equation}
using the infinite $q$-binomial theorem \eqref{E:qbinomP} for the first equality.

Thus, by the $q$-binomial theorem \eqref{E:binomialfinP}, the sum in \eqref{E:justes} is equal to
$$
\sum_{e'= 0}^{u-h} (-1)^{e'} 
  \frac{ q^{\binom{e'}{2} -(e')^2-he'  -ue'
  } (q)_\infty}{  
   (q)_{e'}(q)_{u-h-e'}}=\frac{(q)_\infty}{(q)_{u-h}} \prod _{k=0}^{u-h-1 }{(1-q^{k-2u})}
$$
if $h\leq u$ and $0$ if $h>u$.  Note the product above is $0$ when $-h>u$.
Thus, we conclude
$$
w_{V_i}w_{V_{i'}}=\frac{(q)_\infty(q)_{2u}}{(q)_{u-h}(q)_{u+h}} \prod_{k=0}^{z-1} (1-q^{k-u+h})  \prod_{k=0}^{z'-1} (1-q^{k-u-h}). 
$$
if $-u\leq h \leq u$ and $0$ otherwise.  We note the final products are non-zero if and only if $h\leq u-z$ and $-h\leq u-z'$, and conclude
$$
w_{V_i}w_{V_{i'}}=\frac{(q)_\infty(q)_{2u}}{(q)_{u-h-z}(q)_{u+h-z'}} 
$$
if $ -u+z'\leq h \leq u-z$ and $0$ otherwise.
 \end{proof}

  \subsection{Anomalous self-dual representations}

 \begin{lemma}\label{L:Mintermediate}
Let $\cL$, $\cW$, and $\bH$ be as at the start of Section~\ref{S:getv},  and further assume $H_\Gamma=1$. 
Assume that $V_i$ contains a nontrivial vector fixed by at least one element of $U$. 
 Let $V_i$ be self-dual,  anomalous, and with $\chr(V_i)\mid n$.
Let $u$ be such that $q_i^u=V_i^{\cdot U'}$. 
and $z=\dim_{\kappa_i} H^2(H\rtimes \Gamma,V_i)^{\cL,s_H}$.
There is a particular class $\omega_*^{-1}(\Phi)\in H^2(H\rtimes \Gamma,V_i)$ defined in Proposition~\ref{P:dforbilinear4}.
If $\omega_*^{-1}(\Phi)\in H^2(H\rtimes \Gamma,V_i)^{\cL,s_H}$,
$$
w_{V_i}=
\prod _{k=0}^{\infty }(1+q_i^{-k-u})^{-1}
\prod_{k=0}^{z-1}(1-q_i^{k-u}) 
$$
If $\omega_*^{-1}(\Phi)\not\in H^2(H\rtimes \Gamma,V_i)^{\cL,s_H}$,
 \[ w_{V_i}   =
 \prod _{k=0}^{\infty }(1+q_i^{-k-1 -u})^{-1}
 \prod _{k=0}^{z-1 }{(1-q_i^{k -u})}.
 \] 
In either case, $\tilde{w}_{V_i}$ is finite, and if $z\geq u+1$,  then $w_{V_i}=0$ and otherwise $w_{V_i}$ is positive.   
\end{lemma}

\begin{proof}
Let $q=q_i$.
By definition, 
\[ M_i (e_i) =  \frac{1 }{ \abs{GL_{e_i}(q_i) }\  (V_i^{\cdot U'})^{ e_i }  }  \sum_{\substack{
(G,\pi)\in \Ex(\bH,V_i^{e_i})}}    | E_{3}^{0,2} |  
\quad \textrm{and} \quad w_{V_i}=\sum_{e_i\geq 0}(-1)^{e_i} q_i^{\binom{e_i}{2}}
 M_i ( e_i),\] where
 $\Ex(\bH,V_i^{e_i})$ is the set of
  $(G,\pi)\in \Ext_\Gamma(H,V_i^{e_i})$ such that $G\in\cL$,  and $G_\Gamma=1$, and $s_H\in \im \pi_*$.  
By Lemma~\ref{L:GGamma}, we have $G_\Gamma=1$ since $V_i$ is nontrivial.
So by Lemma~\ref{L:shcriterion}, the classes in $\Ex(\bH,V_i^{e_i})$ correspond to the classes in  $H^2(H\rtimes \Gamma,V_i^{e_i})^{\cL,s_H}=(H^2(H\rtimes \Gamma,V_i)^{\cL,s_H})^{e_i}$.

We can write the elements of  $(H^2(H\rtimes \Gamma,V_i)^{\cL,s_H})^{e_i}$ as $e_i\times z$ matrices over $\kappa_i$, where 
$\im \alpha$ is the column space of the associated matrix and
$\operatorname{span} \alpha$
is the row space of the associated matrix.  
Proposition~\ref{P:dforbilinear4} tells us that $| E_{3}^{0,2} | $ is $2 q^{\frac{(e_i-r)(e_i-r+1)}{2}}$ when $\omega_*^{-1}(\Phi)\in \operatorname{span} \alpha$ and $q^{\frac{(e_i-r)(e_i-r-1)}{2}}$ when $\omega_*^{-1}(\Phi)\not\in \operatorname{span} \alpha.$
Let $r_\alpha:=\dim_{\kappa_i} \im \alpha$.
We have
 \begin{align} \label{E:noeps1}
& \sum_{e_i=0}^{\infty} \frac{ (-1)^{e_i} q^{\binom{e_i}{2}}}{  \abs{\GL_{e_i}(q) }\  (V_i^{\cdot U'})^{ e_i }  } \sum_{ \alpha \in (H^2 ( H \rtimes \Gamma, V_i)^{\cL, s_H} )^{e_i} }   q^{ \frac{ (e_i - r_\alpha) (e_{i} - r_\alpha - 1) }{2}  } 
 =\prod _{k=0}^{\infty }\frac{1}{1+q^{-k-1 -u}}
 \prod _{k=0}^{z-1 }{(1-q^{k -u})}
 \end{align}
(and that the sum converges absolutely) exactly as in Lemma~\ref{odd-selfdual-nontrivial-w}.
If $\omega_*^{-1}(\Phi)\not\in H^2(H\rtimes \Gamma,V_i)^{\cL,s_H}$, then the above sum is $w_{V_i}$
and $\tilde{w}_{V_i}$ is finite by the same argument as in Lemma~\ref{odd-selfdual-nontrivial-w}.

Now we assume $\omega_*^{-1}(\Phi)\in H^2(H\rtimes \Gamma,V_i)$.
The argument for absolute convergence from Lemma~\ref{odd-selfdual-nontrivial-w} also works if
 the terms $q^{\frac{(e_i-r)(e_i-r-1)}{2}}$ above are replaced by $2 q^{\frac{(e_i-r)(e_i-r+1)}{2}}$.
 This implies that $\tilde{w}_{V_i}$ is finite and gives absolute convergence of the  sums we will consider below.

By Lemma~\ref{L:PhiASp}, we have that the $\omega_*^{-1}(\Phi)\in H^2(H\rtimes \Gamma,V_i)$ of Proposition~\ref{P:dforbilinear4} is non-zero.
To count $\alpha$ of rank $r$ where $\omega_*^{-1}(\Phi)\in \operatorname{span} \alpha$
 we first count  $e_i\times z$ matrices over $\kappa_i$
of rank $r$ whose row space contains 
a particular non-zero vector.  
There are $q^{(z-r)(r-1)}\frac{(q)_{z-1}}{(q)_{z-r}(q)_{r-1}}$
rank $r$ subspaces of $\kappa^z$ containing a particular non-zero vector for each $1\leq r \leq z$.  
 For each of these possible row spaces, there are
$q^{e_ir}\frac{(q)_{e_i}}{(q)_{e_i-r}}$  matrices of dimensions $e_i\times z$  with that row space if $e_i\geq r$,  and no such matrices if $e_i<r$.
So for $1\leq r \leq \min(e_i,z)$ there are 
$$
q^{(z-r)(r-1)+e_ir}\frac{(q)_{z-1}(q)_{e_i}}{(q)_{z-r}(q)_{r-1}(q)_{e_i-r}}
$$
$\alpha$ of rank $r$ where $\omega_*^{-1}(\Phi)\in \operatorname{span} \alpha$, and for other $r$ there are $0$ such $\alpha$.

We compute 
 \begin{align*} 
& \sum_{e_i=0}^{\infty} \frac{ (-1)^{e_i} q^{\binom{e_i}{2}}}{  \abs{\GL_{e_i}(q) }\  (V_i^{\cdot U'})^{ e_i }  } 
\sum_{ \substack{\alpha \in (H^2 ( H \rtimes \Gamma, V_i)^{\cL, s_H} )^e\\
\Phi\in \operatorname{span} \alpha} } 
 (2 q^{\frac{(e_i-r_\alpha)(e_i-r_\alpha+1)}{2}}-q^{\frac{(e_i-r_\alpha)(e_i-r_\alpha-1)}{2}}) \\
 &= \sum_{e_i=0}^{\infty} \frac{ (-1)^{e_i} q^{\binom{e_i}{2}}}{ q^{e_i^2}(q)_{e_i} q^{e_iu} } 
\sum_{r=1}^{\min(e_i,z)}
 q^{(z-r)(r-1)+e_ir}\frac{(q)_{z-1}(q)_{e_i}}{(q)_{z-r}(q)_{r-1}(q)_{e_i-r}}
 (2 q^{\frac{(e_i-r)(e_i-r+1)}{2}}-q^{\frac{(e_i-r)(e_i-r-1)}{2}}).
 \end{align*}
Since $H^2(H\rtimes \Gamma,V_i)^{\cL,s_H}$ contains a non-zero vector, $z\geq 1$.  We let $e=e_i-r$ and have that the above sum is 
 \begin{align*} 
 &= \sum_{e=0}^{\infty}
\sum_{r=1}^{z}
 (-1)^{e+r} q^{\binom{e+r}{2}-(e+r)^2-(e+r)u+(z-r)(r-1)+(e+r)r}\frac{(q)_{z-1}(q)_{e+r}}{(q)_{e+r}(q)_{z-r}(q)_{r-1}(q)_{e}}
 (2 q^{\frac{e(e+1)}{2}}-q^{\frac{e(e-1)}{2}}) \\ 
  &= \sum_{e=0}^{\infty}
\sum_{r=1}^{z}
 (-1)^{e+r} q^{\binom{r}{2}+er-(e+r)^2-(e+r)u+(z-r)(r-1)+(e+r)r+\binom{e}{2}}\frac{(q)_{z-1}}{(q)_{z-r}(q)_{r-1}(q)_{e}}
 (2 q^e-1) \\ 
   &= \sum_{e=0}^{\infty}
\frac{ (-1)^{e+1} q^{-e-eu-u}  (2 q^e-1)}{(q)_{e}}
\sum_{r=1}^{z}
(-1)^{r-1} 
q^{\binom{r-1}{2}-(r-1)(u-1)+(z-1)(r-1)-(r-1)^2} 
 \frac{(q)_{z-1}}{(q)_{z-r}(q)_{r-1}}
 \\ 
   &= \sum_{e=0}^{\infty}
\frac{ (-1)^{e} q^{-e(u+1)-u}  (1-2 q^e)}{(q)_{e}}
\prod_{k=0}^{z-2}(1-q^{k-u+1}) \quad \quad \textrm{by the $q$-binomial theorem \eqref{E:binomialfinP} }
 \\ 
    &= q^{-u}\left( \prod_{k=0}^{\infty} \frac{1}{1+q^{-k-u-1}} -2 \prod_{k=0}^{\infty} \frac{1}{1+q^{-k-u}}
     \right)
\prod_{k=0}^{z-2}(1-q^{k-u+1})   \quad \quad \textrm{by the infinite $q$-binomial theorem \eqref{E:binomialnegP} }\\
 \\ 
    &= q^{-u}\left( 1+q^{-u}-2 \right)\prod_{k=0}^{\infty} \frac{1}{1+q^{-k-u}}
\prod_{k=0}^{z-2}(1-q^{k-u+1})
 \\ 
    &= -q^{-u}\prod_{k=0}^{\infty} \frac{1}{1+q^{-k-u}}
\prod_{k=0}^{z-1}(1-q^{k-u}).
 \end{align*}

\details{
\begin{align*}
&\binom{e+r}{2}-(e+r)^2-(e+r)u+(z-r)(r-1)+(e+r)r+\binom{e}{2} \\
&=\binom{r}{2}+er-(e+r)^2-(e+r)u+(z-r)(r-1)+(e+r)r+2\binom{e}{2} \\
&=\binom{r}{2}+er-e^2-r^2-2er-(e+r)u+(z-r)(r-1)+er+r^2 +e^2-e\\
&=\binom{r}{2}-e-(e+r)u+(z-r)(r-1) \\
&=\binom{r-1}{2}+r-1-e-(e+r)u+z(r-1)-r(r-1) \\
&=\binom{r-1}{2}-e-(e+r)u+z(r-1)-(r-1)(r-1) \\
&=\binom{r-1}{2}-e-eu-ru+z(r-1)-(r-1)(r-1) \\
&=\binom{r-1}{2}-e-eu-ru+u-u+z(r-1)-(r-1)(r-1) \\
&=\binom{r-1}{2}-e-eu-(r-1)u-u+z(r-1)-(r-1)(r-1) \\
&=\binom{r-1}{2}-e-eu-(r-1)u-u+(z-1)(r-1)+r-1-(r-1)(r-1) \\
&=\binom{r-1}{2}-e-eu-(r-1)(u-1)-u+(z-1)(r-1)-(r-1)^2 \\
\end{align*}
}

We obtain $w_{V_i}$ by summing the above with \eqref{E:noeps1}, which gives
\begin{align*}
w_{V_i}=& \prod _{k=0}^{\infty }\frac{1}{1+q^{-k-1 -u}}
 \prod _{k=0}^{z-1 }{(1-q^{k -u})}-q^{-u}\prod_{k=0}^{\infty} \frac{1}{1+q^{-k-u}}
\prod_{k=0}^{z-1}(1-q^{k-u})\\
=& 
 \left(
(1+q^{-u})
 -q^{-u} \right)
\prod _{k=0}^{\infty }\frac{1}{1+q^{-k-u}}
\prod_{k=0}^{z-1}(1-q^{k-u}).
\end{align*}
Note that our assumption on $U$ implies that $V_i^{\cdot U'}>1$ so that $u>0$, which implies the final claim about when $w_{V_i}$ is $0$ or positive.
\end{proof}
 
 \subsection{Proof of Theorem~\ref{T:wavefinite}}\label{SS:conclude-moments}
 Corollary~\ref{C:vfactor} reduced the finiteness of $\tilde{v}_{\cW,\bH}$ to that of the $\tilde{w}_{V_i}$ and $\tilde{w}_{N_i}$.
 These later finiteness statements are proven in Lemmas~\ref{non-abelian-w}, \ref{prime-n-nontrivial-w}, \ref{L:wtrivial}, \ref{L:nodual-w}, \ref{odd-selfdual-nontrivial-w}, \ref{L:wntrivial}, \ref{odd-paired-w}, \ref{L:Mintermediate}. The assumption that each representation of $\Gamma$ of characteristic dividing $n$ has a nontrivial vector invariant under at least one element of $U$ is used in some of these lemmas, and immediately implies that $U$ is nonempty, which is used in Lemma \ref{L:wtrivial}.
 
Corollary~\ref{C:vfactor} reduced the non-negativity of ${v}_{\cW,\bH}$ to that of the ${w}_{V_i}$ and ${w}_{N_i}$. 
 These later non-negativity statements follow from the lemmas listed above, with a few additional observations.  
 When using Lemma~\ref{prime-n-nontrivial-w}, we have that
 $$
w_{V_i} = \prod_{j=1}^\infty \left(1-q_i^{-j} \frac{ | H^2( H \rtimes \Gamma, V_i)^{\cL} |}{  |H^1( H \rtimes \Gamma, V_i )| V_i^{\cdot U'}  } 
 \right). 
$$
We note that  $H^2( H \rtimes \Gamma, V_i)^{\cL}, H^1( H \rtimes \Gamma, V_i ),$ and $(V_i^{\cdot U'})$ are $\End_{H\rtimes \Gamma}(V_i)$-vector spaces, and thus every term in the product is $1-q_i^m$ for some integer $m$.
These exponents $m$  are either always negative, in which case the product is positive, or one factor has  $m=0$, which makes the product $0$.  Similar reasoning applies to the expression for $w_{V_i}$ given in Lemma~\ref{L:nodual-w}. 
\details{Lemma~\ref{odd-selfdual-nontrivial-w}: the first product is positive.  If $z=0$, the second product is $1$.  Otherwise the exponents on $q$ start at $-d$, which is non-positive and increase by on in each factor, so they only become positive if one is $0$ and the whole expression is $0$.
}

\subsection{Criteria for nonzero probability}

\begin{proposition}\label{nzp-criterion} Let $\cL$ be a level of the category of finite $\Gamma$-groups.
Let $\cW$ be the set of isomorphism classes of finite $n$-oriented $\Gamma$-groups whose underlying $\Gamma$-group is in $\cL$. Let $U$ be a multiset of elements of $\Gamma$. Assume all nonzero
 representations of $\Gamma$ of characteristic dividing $n$ contain some nontrivial vector fixed by at least one element of $U$. Fix $\bH\in\cW$ with $H_\Gamma$ trivial.
 
 Then we have $v_{\cW, \bH}>0$ if and only if, for each finite simple abelian $H \rtimes \Gamma$-group $V_i$ of characteristic prime to $n$ that can be the kernel of a simple morphism $\bG \to \bH$ with $\bG \in \cW$ we have
 \[ | H^2( H \rtimes \Gamma, V_i)^{\cL} | \leq  |H^1( H \rtimes \Gamma, V_i )| V_i^{\cdot U'}\] and for each finite simple abelian $H \rtimes \Gamma$-group $V_i$ of characteristic dividing $n$ that can be the kernel of a simple morphism $\bG \to \bH$ with $\bG \in \cW$ we have \[|H^1 ( H \rtimes \Gamma, V_i^\vee)|  |H^2(H\rtimes \Gamma,V_i)^{\cL,s_H}|  \leq  |H^1 ( H \rtimes \Gamma, V_i)|
  V_i^{\cdot U'}   .\]
  \end{proposition}
 
 \begin{proof} By Corollary~\ref{C:vfactor} we have ${v}_{\cW,\bH}\neq 0 $ if and only if $M_{\bH}$ is nonzero and all the ${w}_{V_i}$ and ${w}_{N_i}$ are nonzero. That $M_\bH$ is nonzero follows from the assumption that $H_\Gamma$ is trivial. It follows from Lemma \ref{non-abelian-w} that the $w_{N_i}$ are all nonzero, so it suffices to check for $V_i$ solo that $w_{V_i}$ is nonzero if and only if the stated criteria are satisfied for $V_i$ and for $V_i$ paired that $w_{V_i} w_{V_i'}$ is nonzero if and only if the stated criteria are satisifed for $V_i$ and $V_{i'}$.
 
 In the case of Lemma~\ref{prime-n-nontrivial-w}, each term $q_i^{-j} \frac{ | H^2( H \rtimes \Gamma, V_i)^{\cL} |}{  |H^1( H \rtimes \Gamma, V_i )| V_i^{\cdot U'}  }$  has the form $q_i^m$ for some $m$ and the product of the $1- q_i^m$ is nonzero if and only if no $m$ is zero, which happens and only if each $m$ is $\leq -1$, which occurs if and only if $| H^2( H \rtimes \Gamma, V_i)^{\cL} | \leq  |H^1( H \rtimes \Gamma, V_i )| V_i^{\cdot U'}$. 
 
Lemma~\ref{L:wtrivial} states that for $V_i$ trivial of characteristic prime to $n$, $w_{V_i} \neq 0$ if and only if $z \leq u$, in other words if $ | H^2( H \rtimes \Gamma, V_i)^{\cL} | \leq   V_i^{\cdot U'}$. Since $|H^1( H \rtimes \Gamma, V_i )|=1$ in this case by the assumption $H_\Gamma=1$, this is equivalent to the criterion  $ | H^2( H \rtimes \Gamma, V_i)^{\cL} | \leq  |H^1( H \rtimes \Gamma, V_i )| V_i^{\cdot U'}$.

In the case of Lemma~\ref{L:nodual-w}, $w_{V_i}$ is given by a product with all terms of the form $1-q_i^m$ where the product is nonzero if and only if all $m \leq -1$, which occurs if and only if $|H^1 ( H \rtimes \Gamma, V_i^\vee)|  |H^2(H\rtimes \Gamma,V_i)^{\cL,s_H}|  \leq  |H^1 ( H \rtimes \Gamma, V_i)| V_i^{\cdot U'}   $.
   
 Lemma~\ref{odd-paired-w} states that for $V_i, V_{i'}$ dual but nonisomorphic of characteristic dividing $n$, we have $w_{V_i} w_{V_{i'}}\neq 0$ if and only if $-u + z' \leq h \leq u-z$, which taking exponentials with base $q$  and applying the definitions of $u,z,z', h$ occurs if and only if
 \[ \frac{ |H^2(H \rtimes \Gamma, V_{i'})^{\cL, s_H} }{V_i^{\cdot U'}}  \leq \frac{ |H^1 (H\rtimes \Gamma, V_{i'} ) |}{ |H^1 (H\rtimes \Gamma, V_i) |}  \leq \frac{V_i^{\cdot U'}} { |H^2(H \rtimes \Gamma, V_{i})^{\cL, s_H} }.\] The first inequality is equivalent to the criterion $|H^1 ( H \rtimes \Gamma, V_{i'} ^\vee)|  |H^2(H\rtimes \Gamma,V_{i'})^{\cL,s_H}|  \leq  |H^1 ( H \rtimes \Gamma, V_{i'} )|
  V_{i'} ^{\cdot U'}   $ for $V_i'$ since $  V_{i'} ^{\cdot U'}  =   V_{i} ^{\cdot U'}  $ and $V_{i'}^\vee \cong V_i$, while the second inequality is equivalent to the criterion  $|H^1 ( H \rtimes \Gamma, V_i^\vee)|  |H^2(H\rtimes \Gamma,V_i)^{\cL,s_H}|  \leq  |H^1 ( H \rtimes \Gamma, V_i)| 
  V_i^{\cdot U'}  $ for $V_i$.
  
  In the case that $V_i$ is self-dual of characteristic dividing $n$, we have $H^1 ( H \rtimes \Gamma, V_i^\vee) \cong H^1 ( H \rtimes \Gamma, V_i)$ and thus the stated criterion simplifies to $|H^2(H\rtimes \Gamma,V_i)^{\cL,s_H}|  \leq   V_i^{\cdot U'}   $. This matches the condition for $w_{V_i} \neq 0$ given in Lemmas \ref{odd-selfdual-nontrivial-w}, \ref{L:wntrivial}, and \ref{L:Mintermediate}. \end{proof}
  
  \subsection{The main statements}\label{ss-statements}
  
  \begin{proof}[Proof of Theorem \ref{T:intro-measure-general}] 
  
  The statement of Theorem \ref{T:intro-measure-general} follows immediately from Theorem \ref{T:measure} and Corollary \ref{C:vfactor} once we check that the formulas for $w_{V_i}$ and $w_{N_i}$ given in \S\ref{ss-intro-formulas} match the definitions of $w_{V_i}$ and $w_{N_i}$. We check this by observing that the formulas given in \S\ref{ss-intro-formulas} match the formulas for $w_{V_i}$ and $w_{N_i}$ proven, in various cases, in Lemmas~\ref{non-abelian-w}, \ref{prime-n-nontrivial-w}, \ref{L:wtrivial}, \ref{L:nodual-w}, \ref{odd-selfdual-nontrivial-w}, \ref{L:wntrivial}, \ref{odd-paired-w}, \ref{L:Mintermediate}.\end{proof}
  
  \begin{proof}[Proof of Theorem \ref{T:intro-measure-nf}] Let $U$ be the multiset consisting of $\gamma_v$ for each $v \in \operatorname{Spec} k \otimes \mathbb R$. Let us first check that $U$ satisfies the hypothesis of Theorem \ref{T:intro-measure-general} that each nonzero representation of $\Gamma$ of characteristic dividing $n$ contains a nonzero vector fixed by some element of $U$. If $n=1$ then there are no nonzero representations of characteristic dividing $n$ and the hypothesis is vacuously satisfied.
  If some $\gamma_v$ is trivial, then $\gamma_v$ fixes every nonzero vector and the hypothesis is again satisfied.
    
  Theorem \ref{T:intro-measure-nf} now follows from specializing Theorem \ref{T:intro-measure-general} to this value of $U$, except for the final claim about convergence conditional on Conjecture \ref{intro-moments-conjecture}. 
   The counting of fields gives sequence of measures $\nu_B$. Conjecture \ref{intro-moments-conjecture} implies that the moments of the $\nu_B$ converge to the moments $b_{\mathbf H}$, which by Theorem \ref{T:measure} implies that the $\nu_B$ converges to $ \nu_{\Gamma, n, \overline{\gamma}}$. Because the set $\{X|X^\cW\isom  \bH\}$ is open and closed, it follows that $\nu_B ( \{X|X^\cW\isom  \bH\})$ converges to $\nu_{\Gamma, n, \overline{\gamma}}( \{X|X^\cW\isom  \bH\})$, which gives the desired statement.

  \end{proof}

    \begin{proof}[Proof of Theorem \ref{T:intro-measure-ff}] This follows from specializing Theorem \ref{T:intro-measure-general} to $U =\{1\}$, except for the final claim about convergence. Observe that for a function $F$ of parameters $b,q$ and limit $\ell$, we have $ \lim_{q \to\infty} \limsup_{b\to\infty} F(b,q) = \ell$ and $ \lim_{q \to\infty} \liminf_{b\to\infty} F(b,q) = \ell$ if and only if, for each sequence of pairs $b_i,q_i$, as long as $b_i$ grows sufficiently fast with respect to $q_i$, we have $\lim_{i \to\infty} F(b_i,q_i)=\ell$.
    
    We therefore fix a level $\cL$ and a group $\bH$, and a sequence of pairs $b_i, q_i$, where we may assume that $q_i$ grows arbitrarily fast with respect to $b_i$, and that $(q_i, \abs{\Gamma} M)=1$ and $(q_i-1,nM)=n$. It suffices to prove that 
    \begin{equation}\label{sequence-convergence-equation} \lim_{i\to\infty} \frac{\sum_{m\leq b_i}  \abs{ \{ K \in  E_\Gamma( q_i^m , \mathbb F_{q_i}(t)  \mid  \bGal ( K^{ \operatorname{un}, \abs{\Gamma}'} / K)^\cW \isom \bH   \}} }{ \sum_{m\leq b_i}\abs{ E_\Gamma( q_i^m , \mathbb F_{q_i}(t))}} =  \nu_{\Gamma, n, \{1\}}(\{X|X^\cW\isom  \bH\}).\end{equation}  We now define a sequence of measures $\nu_i$ by summing over $m \leq b_i$ and $K \in E_\Gamma( q_i^m , \mathbb F_{q_i}(t))  $ a delta measure supported on $ \bGal ( K^{ \operatorname{un}, \abs{\Gamma}'} / K)^\cW$ and dividing by $\sum_{m\leq b_i}\abs{ E_\Gamma( q_i^m , \mathbb F_{q_i}(t))}$. Equation \eqref{sequence-convergence-equation} follows from the claim that $\nu_i$ converges to the projection of  $\nu_{\Gamma, n, \{1\}}$ onto $\cW$. By \cite[Theorem 1.6]{Sawin2022}, this follows if we check that $\bG$-moment of $\nu_i$, for each $\bG \in \cW$, converges to the $\bG$-moment of $\nu_{\Gamma, n, \{1\}}$, which is $\begin{cases} \frac{ \abs{H^\Gamma} \abs{ H^2( H\rtimes \Gamma , \mathbb Z/n) } }{  \abs{H} \abs{ H^3( H\rtimes \Gamma , \mathbb Z/n) }}& \textrm{if } H_\Gamma = 1 \\ 0 & \textrm{if } H_\Gamma \neq 1 \end{cases} $.
    
    The $\bG$-moment of $\nu_i$ is given by 
   \[  \frac{\sum_{m\leq b_i} \sum_{K \in E_\Gamma( q_i^m , \mathbb F_{q_i}(t))}  \Sur ( \bGal ( K^{ \operatorname{un}, \abs{\Gamma}'} / K), \mathbf H)  }{ \sum_{m\leq b_i}\abs{ E_\Gamma( q_i^m , \mathbb F_{q_i}(t))}} \] so it suffices to prove that

   \begin{equation}\label{sequence-convergence-moments} \lim_{i\to\infty} \frac{\sum_{m\leq b_i} \sum_{K \in E_\Gamma( q_i^m , \mathbb F_{q_i}(t))}  \Sur ( \bGal ( K^{ \operatorname{un}, \abs{\Gamma}'} / K), \mathbf H)  }{ \sum_{m\leq b_i}\abs{ E_\Gamma( q_i^m , \mathbb F_{q_i}(t))}} = \begin{cases} \frac{ \abs{H^\Gamma} \abs{ H^2( H\rtimes \Gamma , \mathbb Z/n) } }{  \abs{H} \abs{ H^3( H\rtimes \Gamma , \mathbb Z/n) }}& \textrm{if } H_\Gamma = 1 \\ 0 & \textrm{if } H_\Gamma \neq 1 \end{cases}.\end{equation}

   Let us now check that we have $ (q_i, |\Gamma| |G| )=1 $ and $ q_i \equiv 1\bmod n $ and $(q_i-1, |G|) = (n,|G|) $. The first claim follows from $(q_1, |\Gamma| M)=1$ and the fact that $|G|$ divides a power of $M$ by definition of $M$. The second claim follows from $(q_i-1,nM)=n$. The third claim fails if there is a prime $p$ where the $p$-adic valuation of $q_i-1$ and the $p$-adic valuation of $|G|$ are both greater than the $p$-adic valuation of $n$. In particular, this requires $p$ to divide $M$, and so implies that the $p$-adic valuation of $(q_i-1, n M)$ is greater than the $p$-adic valuation of $n$, contardicting the assumption that $(q_i-1,nM)=n$.
  
  \eqref{sequence-convergence-moments} now follows, for $b_i$ growing sufficiently fast with respect to $q_i$, from Theorem \ref{ff-theorem} and the above observation on the relationship between double limits and single limits. (Theorem \ref{ff-theorem} states an identity that holds for $q$ sufficiently large depending on $H,\Gamma$, and thus always holds in the limit as $q$ tends to infinity.) \end{proof}
  
  \begin{proof}[Proof of Theorem \ref{intro-bb}] We specialize Theorem \ref{T:intro-measure-ff} to the setting where $\Gamma = \mathbb Z/3$, $n=2$, and $\cL$ is a level consisting of $2$-groups with an action of $\Gamma$. We have $M=2$ so the condition $(q,\abs{\Gamma} M)=1$ reads $(q,6)=1$ and is equivalent to $q$ being odd and not divisible by $3$, while the condition $(q-1, nM)=n$ reads $(q-1,4)=2$ and is equivalent to $q$ being congruent to $3$ modulo $4$. So these conditions match the conditions $q \equiv 3\bmod 4$ and $3\nmid q$ from the statement of Theorem \ref{intro-bb}.
  
   $\Gal(K^{\operatorname{ur},2}/K)$ is the maximal $2$-group quotient of $\Gal(K^{ \operatorname{ur}, \abs{\Gamma}'}/K)$. The set of profinite groups whose maximal $2$-group quotient is isomorphic to $H$ is the finite disjoint union of, for each $\Gamma$-group $H'$ whose underlying group is isomorphic to $H$, the set of profinite groups whose maximal $2$-group quotient is isomorphic as a $\Gamma$-group to $H'$.  By Lemma \ref{maximal-p-group-quotient} this set has the form $\{ X\mid X^{\cL} \cong H'\}$ for $\cL$ a level in the category of $\Gamma$-groups consisting only of $p$-groups with an action of $\Gamma$, including $H'$ and every extension of $H'$ by $V_i$ for $V_i$ any group of the form $\mathbb F_p^d$ with an irreducible action of $\Gamma$.  Thus the probability that $\Gal(K^{ \operatorname{ur}, \abs{\Gamma}'}/K)$ is the sum over $H'$ of the probability that $\Gal(K^{ \operatorname{ur}, \abs{\Gamma}'}/K)$ lies in $\{ X\mid X^{\cL} \cong H'\}$. By Theorem \ref{T:intro-measure-ff}, the double limit (with limsup or liminf) of the probability that $\Gal(K^{ \operatorname{ur}, \abs{\Gamma}'}/K)$ lies in $\{ X\mid X^{\cL} \cong H'\}$ is equal to $\nu_{\mathbb Z/3, 2, \{1\}} ( \{ X\mid X^{\cL} \cong H'\} )$ (because the set  $\{ X\mid X^{\cL} \cong H'\}$ is itself a disjoint union of  $\{ X\mid X^{\cL} \cong \bH'\}$ for oriented $\Gamma$-groups $\bH'$ whose underlying group is $H'$).  Hence if we take \[ p_H = \sum_ {H'} \nu_{\mathbb Z/3, 2, \{1\} }( \{ X\mid X^{\cL} \cong H'\} ) \] we obtain the first part of Theorem \ref{intro-bb}.  We have $ \nu_{\mathbb Z/3, 2, \{1\} }( \{ X\mid X^{\cL} \cong H'\} )=0$ if $(H')_\Gamma\ne 1$ by Theorem~\ref{T:intro-measure-general}.
   
  For $H$ the Klein four group or the eight-element quaternion group, there is a unique non-trivial $\Gamma$-action on $H$ up to isomorphism, and $ \nu_{\mathbb Z/3, 2, \{1\} }( \{ X\mid X^{\cL} \cong H'\} )$ were calculated respectively in Lemmas \ref{klein-four-probability} and \ref{quaternion-probability}, and these agree with the formula for $p_H$ claimed in Theorem \ref{intro-bb}.
  
Clearly $p_H$ is zero if $H$ is not a $2$-group. The remaining $2$-groups of order at most $8$ are the trivial group, $\mathbb Z/2$,  $\mathbb Z/4$,  $(\Z/2)^3$, $\Z/4\times\Z/2$,  $\Z/8$, and the dihedral group $D_4$. 
The 
groups $\mathbb Z/2$, $\mathbb Z/4$,  $\Z/4\times\Z/2$, $\Z/8$, and the dihedral group $D_4$  do not have a nontrivial automorphism of order $3$, so the only $\Gamma$-group structure is the trivial $\Gamma$-action.   The group $(\Z/2)^3$ has one non-trivial $\Gamma$-action up to automorphism (the regular representation of $\Gamma$), but that action still has non-trivial coinvariants.  
We have $p_H =0$ if $H _\Gamma $ is nontrivial since the moment is already trivial for these groups, so we get that $p_H=0$ for the remaining $H$. 

Finally, for the trivial group, we apply  Equation \eqref{cct-prob} to $H$ the trivial group and $\cL$ the level generated by, in the notation of \S\ref{ss-bb}, the groups $V_1$ and $V_2$. All the multiplicative factors outside the sum are $1$, the group $H_3( H \rtimes \Gamma, \mathbb F_2)$ being summed over has a single element $s_H=0$, the factor $w_{V_1}(s_H)$ inside the sum gives $1$ by Equation \eqref{wv1}, and the factor  $w_{V_2}(s_H)$ inside the sum gives $\prod_{k=0}^{\infty} (1+ 4^{ -k - \frac{1}{2} - 1})^{-1}$ by Equation \eqref{wv2}. Hence the probability in this case is $\prod_{k=0}^{\infty} (1+ 4^{ -k - \frac{1}{2} - 1})^{-1}$.  \end{proof}

\section{Non-existence results}\label{S:NE}

Let $k$ be a number field, $\Gamma$ a finite group, $K$ an extension of $k$ with Galois group $\Gamma$. Let $V$ be a finite abelian group of order prime to $|\Gamma|$ with an action of 
$\Gal ( K^{ \operatorname{un}, \abs{\Gamma}'} / k)$. 
 Let $H^i_c(\operatorname{Spec}  \mathcal O_K, V)$ be the compactly supported \'{e}tale cohomology groups in the sense of \cite[p. 165]{Milne2006}. Let $H^i_{\textrm{par}}(\operatorname{Spec}  \mathcal O_K, V)$ be the image of the natural map $H^i_c(\operatorname{Spec}  \mathcal O_K, V)\to H^i(\operatorname{Spec}  \mathcal O_K, V)$ (the \emph{parabolic cohomology groups}.)

\begin{lemma}\label{parabolic-Euler-characteristic} With notation as in the start of Section~\ref{S:NE},
we have 
\[ \frac{ | H^0(\operatorname{Spec} \mathcal O_K, V)^\Gamma| |H^2_c(\operatorname{Spec} \mathcal O_K, V)^\Gamma|}{ |H^1_{\textrm{par}}(\operatorname{Spec} \mathcal O_K, V)^\Gamma|   |H^3_c(\operatorname{Spec} \mathcal O_K, V)^\Gamma |} =  \prod_{ v \textrm{ archimedean place of } k} | V ^{ \Gal(k_v^{\operatorname{sep}}/k_v)} |,\]
where the action of the absolute Galois group $\Gal(k_v^{\operatorname{sep}}/k_v)$ on $V$ is through the map $\Gal(k_v^{\operatorname{sep}}/k_v)\ra \Gal ( K^{ \operatorname{un}, \abs{\Gamma}'} / k)$ (only defined up to conjugacy).
\end{lemma} 

\begin{proof} Let $U$ be an open subset of $\operatorname{Spec} \mathcal O_k$ avoiding all the primes that ramify in $K$ or divide the order of $V$. By \cite[II, Theorem 2.13(b)]{Milne2006} the right-hand side is equal to
\[ \frac{ |H^0_c( U, V) |  |H^2_c( U, V) | }{ |H^1_c( U, V) |   |H^3_c( U, V) |} \] where $H^i_c(U,V)$ arises from the mapping cone of the natural map from $H^i(U,V)$ to the sum over places $v$ of $F$ not in $U$ of  $H^i ( k_v,V)$, where $H^i(k_v,V)$ is understood as Galois cohomology for non-archimedean places and Tate cohomology for archimedean places.

Let $\tilde{U}$ be the inverse image of $U$ in $\operatorname{Spec} \mathcal O_K$. We have  $H^i(U, V) = H^i(\tilde{U}, V)^\Gamma$
by the Hochschild-Serre spectral sequence \cite[Ch. III, Theorem 2.20]{Milne1980} because $V$ has order prime to $|\Gamma|$ and hence so do all its cohomology groups.  We can upgrade this to $H^i_c(U, V) = H^i_c(\tilde{U}, V)^\Gamma$ using the long exact sequence for compactly supported cohomology 
\cite[Ch. 2, Prop. 2.3(a)]{Milne2006}, and the similar fact that
$H^i(\Gal(k_v^{\operatorname{sep}}/k_v)),V) \ra H^i(\Gal(K_w^{\operatorname{sep}}/K_w)),V)^{\Gamma_v}$ for any place $v$ of $k$, place $w$ of $K$ above $v$, and stabilizer $\Gamma_v$ of $w$ in $\Gamma$.

A long exact sequence relates $H^i_c(\tilde{U}, V) \to H^i_c( \mathcal O_K ,V)$ and the sum over non-archimedean places $w$ of $K$ not in $\tilde{U}$ of $H^i( R_w, V)$ where $R_w$ is the residue field \cite[Ch. 2, Prop. 2.3(d)]{Milne2006}.The cohomology of the residue field is concentrated in degrees $0,1$ so nonzero maps $H^i_c( \mathcal O_K, V) \to \prod_{w\not\in \tilde{U}} H^i (R_w, V)$ occur only for $i=0,1$ and nonzero maps $ \prod_{w\not\in \tilde{U}} H^{i-1} (R_w, V) \to H^i_c(\tilde{U},V)$ occur only for $i=1,2$. In particular these only occur in degrees between $0$ and $3$, giving

\[ \frac{ |H^0_c(\tilde{U},V)^\Gamma||H^2_c(\tilde{U},V)^\Gamma|}{|H^1_c(\tilde{U},V)^\Gamma||H^3_c(\tilde{U},V)^\Gamma|} = \frac{ |H^0_c(\operatorname{Spec} \mathcal O_K,V)^\Gamma||H^2_c(\operatorname{Spec} \mathcal O_K,V)^\Gamma|}{|H^1_c(\operatorname{Spec} \mathcal O_K,V)^\Gamma||H^3_c(\operatorname{Spec} \mathcal O_K,V)^\Gamma|}  \frac{ | ( \prod_{w\not\in \tilde{U}}  H^0(R_w, V))^\Gamma|}{| ( \prod_{w\not\in \tilde{U}}  H^1(R_w,V))^\Gamma|} .\]

The invariants $( \prod_{w\not\in \tilde{U}}  H^0(R_w, V))^\Gamma$ split as a product over one representative $w$ of each $\Gamma$-orbit of $H^0(R_w,V)^{\Gamma_w}$, where $\Gamma_w\sub \Gamma$ is the stabilizer of $w$, and the same is true for $( \prod_{w\not\in \tilde{U}}  H^1(R_w, V))^\Gamma$. Thus to show $ \frac{ | (\prod_{w\not\in \tilde{U}} H^0(R_w V))^\Gamma|}{| (\prod_{w\not\in \tilde{U}} H^1(R_w,V))^\Gamma|} =1$ it suffices to show $H^0(R_w,V)^{\Gamma_w}$ and $H^1(R_w,V)^{\Gamma_w}$ have the same cardinality for all $w$, which follows from the $\Gamma_w$-equivariant long exact sequence
\[0 \to H^0(R_w,V) \to V \stackrel{{\operatorname{Frob}_w-1}}{\ra} V \to H^1(R_w,V) \to 0.\]

This gives
\[ \frac{ | H^0_c(\operatorname{Spec} \mathcal O_K, V)^\Gamma| |H^2_c(\operatorname{Spec} \mathcal O_K, V)^\Gamma|}{ |H^1_c(\operatorname{Spec} \mathcal O_K, V)^\Gamma|   |H^3_c(\operatorname{Spec} \mathcal O_K, V)^\Gamma }|  =  \prod_{ v \textrm{ archimedean place of } k} | V ^{ \Gal(k_v^{\operatorname{sep}}/k_v)} |.\]

To obtain the statement, we use the long exact sequence \cite[Ch. 2, Prop. 2.3(a)]{Milne2006}, which together with the vanishing of $H^{-1} (\operatorname{Spec} \mathcal O_K, V)$ gives a long exact sequence
\begin{align*} 0 &\to  (\prod_{w \textrm{ arch}} H^{-1} ( K_w,V ) )^\Gamma \to H^0_c(\operatorname{Spec} \mathcal O_K, V)^\Gamma \to H^0(\operatorname{Spec} \mathcal O_K, V)^\Gamma \to  (\prod_{w \textrm{ arch}} H^{0} ( K_w,V ) )^\Gamma\\
 &\to H^1_c(\operatorname{Spec} \mathcal O_K, V)^\Gamma \to H^1(\operatorname{Spec} \mathcal O_K, V)^\Gamma
 \end{align*}
and thus
\begin{align*}
0 &\to  (\prod_{w \textrm{ arch}} H^{-1} ( K_w,V ) )^\Gamma \to H^0_c(\operatorname{Spec} \mathcal O_K, V)^\Gamma \to H^0(\operatorname{Spec} \mathcal O_K, V)^\Gamma \to  (\prod_{w \textrm{ arch}} H^{0} ( K_w,V ) )^\Gamma\\
 &\to H^1_c(\operatorname{Spec} \mathcal O_K, V)^\Gamma \to H^1_{\textrm{par}} (\operatorname{Spec} \mathcal O_K, V)^\Gamma\to 0. \end{align*}
Since the Galois group of every archimidean place is $\mathbb Z/2$ or $0$, the Tate cohomology groups are isomorphic in each degree and thus $  (\prod_w H^{-1} ( K_w,V ) )^\Gamma$ and  $(\prod_w H^{0} ( K_w,V ) )^\Gamma $ have the same cardinality. This gives
\[ \frac{ | H^0(\operatorname{Spec} \mathcal O_K, V)^\Gamma| }{ |H^1_{\textrm{par}}(\operatorname{Spec} \mathcal O_K, V)^\Gamma|   }=  \frac{ | H^0_c(\operatorname{Spec} \mathcal O_K, V)^\Gamma| }{ |H^1_c(\operatorname{Spec} \mathcal O_K, V)^\Gamma|  } \] 
and the statement.\end{proof} 

We now compare these cohomology groups to the Galois cohomology groups $H^i ( \Gal ( K^{ \operatorname{un}, \abs{\Gamma}'} / K), V)$. The next two lemmas are number field analogs of standard facts in topology (e.g. see \cite[Lemma 2.1]{Sawin2024}.)

\begin{lemma}\label{H1-comparison-prof} With notation as in the start of Section~\ref{S:NE}, there is a $\Gamma$-equivariant isomorphism \[ H^1( \Gal ( K^{ \operatorname{un}, \abs{\Gamma}'} / K), V) \to H^1_{\textrm{par}} (\Spec \mathcal O_K, V).\] \end{lemma}

\begin{proof} 
By the long exact sequence \cite[Ch. 2, Prop. 2.3(a)]{Milne2006}, we have that $H^1_{\textrm{par}} (\Spec \mathcal O_K, V)$ is precisely the subgroup of elements of $H^1 (\Spec\mathcal O_K, V)$ that are trivial in all maps to $H^1(K_v,V)$ for $v$ an archimedian place of $K$. 
Remark~\ref{R:defineor} gives a map $H^1( \Gal ( K^{ \operatorname{un}, \abs{\Gamma}'} / K), V) \to H^1 (\Spec\mathcal O_K, V)$. 
Any cohomology class arising from this map must restrict to a trivial class at each archimedean place, since $K^{ \operatorname{un}, \abs{\Gamma}'}/K$ is split completely at each archimedean place, and so  the image of the map is contained in $H^1_{\textrm{par}} (\Spec\mathcal O_K, V)$. 

Let $K^{ \operatorname{unf}}$ be the  maximal algebraic extension of $K$ that is unramified at all finite places.
By the definition, one can see that the map $H^1( \Gal ( K^{ \operatorname{un}, \abs{\Gamma}'} / K), V) \to H^1 (\Spec\mathcal O_K, V)$ then factors through $H^1( \Gal ( K^{ \operatorname{unf}} / K), V)$.
The resulting map $ H^1( \Gal ( K^{ \operatorname{unf}} / K), V) \ra H^1 (\Spec\mathcal O_K, V)$ is an isomorphism (see, e.g., \cite[Theorem 15]{Youcis2015}).
Note that by the inflation-restriction exact sequence 
$H^1( \Gal ( K^{ \operatorname{un}, \abs{\Gamma}'} / K), V)\ra H^1( \Gal ( K^{ \operatorname{unf}} / K), V)$ is injective.  
  Thus it remains to show that
$H^1( \Gal ( K^{ \operatorname{un}, \abs{\Gamma}'} / K), V)$ is the subgroup of
$H^1( \Gal ( K^{ \operatorname{unf}} / K), V)$ that is trivial in all maps to $H^1(K_v,V)$ for $v$ archimedian.

Let $G=\Gal ( K^{ \operatorname{unf}} / K)$ and $N=\Gal(K^{ \operatorname{unf}}/K^{ \operatorname{un}, \abs{\Gamma}'})$. 
Consider an element $\alpha$ of $H^1( \Gal ( K^{ \operatorname{unf}} / K), V)$ that is trivial in all maps to $H^1(K_v,V)$ for $v$ archimedian.
It is represented by a cocycle $G\ra V$ that gives a section $s$ of $V\rtimes G \ra G$. 
We compose $s$ with the quotient to obtain $G\ra V\rtimes G \ra V\rtimes G/N$.  The composite map has image whose order is prime to $|\Gamma|$.  Moreover if $G_v \sub G$ is the decomposition group for an archimedian place, since $G_v$ acts trivially on $V$, we have that $s_{G_v}$ is a homomorphism and is trivial by our assumption on $\alpha$.  So we have $G\ra V\rtimes G/N$ which is trivial on all $G_v$ for $v$ archimedian and has image of order prime to $|\Gamma|$.  Thus this map factors through $G/N$, and we see our cocycle is the image of a cocycle $G/N\ra V$, as desired.  
%
%
%
\end{proof}

\begin{lemma}\label{H2-comparison-prof} With notation as in the start of Section~\ref{S:NE}, there is a $\Gamma$-invariant injection \[ H^2( \Gal ( K^{ \operatorname{un}, \abs{\Gamma}'} / K), V) \to H^2_{c} (\mathcal O_K, V).\] \end{lemma}

\begin{proof} 
The homomorphism is given by Lemma~\ref{L:compactsupports}, and the $\Gamma$-invariance can be seen from the definition. It remains to check injectivity.

Fix $\alpha \in H^2( \Gal ( K^{ \operatorname{un}, \abs{\Gamma}'} / K), V) $ whose image in $H^2_{c} (\mathcal O_K, V)$ is zero. In particular, the image of $\alpha$ in $H^2 (\mathcal O_K, V)$ is zero. Fix a finite unramified $\abs{\Gamma}'$ Galois extension $K'$ of $K$ such that $\alpha$ arises from a class in $H^2(\Gal(K'/K), V)$. Then the induced class in $H^2 (\mathcal O_K, V)$ arises from a Cech cocycle of the covering $\operatorname{Spec} \mathcal O_{K'} \to \operatorname{Spec} \mathcal O_K$. The map from the Cech cohomology of this covering to the sheaf cohomology is the edge map of the Cech-to-sheaf spectral sequence of the covering $\operatorname{Spec} \mathcal O_{K'} \to \operatorname{Spec} \mathcal O_K$. Since $\alpha$ is sent to $0$, the associated cocyle must be in the image of a differential of this spectral sequence, which must be $d_2^{0,1}: E^{0,1}_2 \to E_2^{2,0}$. The source of this differential is a subgroup of $H^1( \operatorname{Spec} \mathcal O_{K'} , V)$.

The class $\beta$ in $H^1( \operatorname{Spec} \mathcal O_{K'} , V)$ whose image under $d_2^{0,1}$ is $\alpha$ corresponds to a $V$-torsor on $\mathcal O_{K'}$. Let $K^*$ be a Galois extension of $K$ over which the torsor $V$ splits. Since $V$ is unramified over $K$, we can take $K^*$ to be unramified at all finite places, but not necessarily at the infinite places. Since this torsor splits over $K^*$, the pullback of $\beta$ to $H^1( \operatorname{Spec} \mathcal O_{K^*} , V)$ vanishes. Since the differentials of this spectral sequence are compatible with pullbacks, the pullback of $\alpha$ to the Cech cohomology of $\operatorname{Spec} \mathcal O_{K^*} \to \operatorname{Spec} \mathcal O_K$ is the image of $d_{2}^{0,1}$ applied to the pullback of $\beta$, which is $0$, and thus must be $0$.

The Cech complex of  $\operatorname{Spec} \mathcal O_{K^*} \to \operatorname{Spec} \mathcal O_K$ is isomorphic to the group cohomology complex of $\Gal(K^*/K)$, and our map from group cohomology cocycles to Cech cohomology uses this isomorphism, so the pullback of $\alpha$ to $C^2 (\Gal(K^*/K), V)$ is a coboundary $d\gamma$ for some $\gamma \in C^1 (\Gal(K^*/K), V)$.

The image of $\alpha$ under the homomorphism of Remark~\ref{R:defineor} is represented by a pair consisting of the Cech cocycle on the covering $\operatorname{Spec} \mathcal O_{K^*} \to \operatorname{Spec} \mathcal O_K$ associated to $\alpha$ and, for each infinite place $v$ of $K$, a cochain $\delta_v$ for the Cech cover $\operatorname{Spec} K' \otimes_K K_v \to K_v$ of $K_v$ whose differential is the pullback of $\alpha$. Since $\alpha= d\gamma$, this cocycle is the coboundary $(d\gamma, (\gamma)_{v\mid\infty})$ plus the coycle $(0, (\delta_v-\gamma)_{v\mid\infty})$. The cocycle $(0, (\delta_v-\gamma)_{v\mid\infty})$ is the image of the cocycle $(\delta_v-\gamma)_{v\mid\infty} \in \prod_{v\mid\infty} H^1( K_v, V)$ under the connecting homomorphism $\prod_{v\mid\infty} H^1( K_v, V) \to H^2_c( \mathcal O_K, V)$.

It follows that $(\delta_v-\gamma)_{v\mid\infty}$ is in the image of $H^1(\mathcal O_K, V) \to \prod_{v\mid\infty} H^1( K_v, V)$.  A class in $H^1(\mathcal O_K, V)$ whose image is $(\delta_v-\gamma)_{v\mid\infty}$ corresponds to a $V$-torsor on $\Spec \mathcal O_K$. By a possible modification of $K^*$, we may assume that this torsor is also split over $K^*$, so that the class in $H^1(\mathcal O_K, V)$ is representated by a group cohomology cocycle $\epsilon \in C^1 ( \Gal(K^*)/K , V)$. For each $v$, we conclude that the pullback of $\epsilon$ to $H^1( K_v, V)$ agrees with $\delta_v-\gamma$.

Explicitly, the cocycle $\alpha$ endows $\Gal(K'/K) \times V$ with the structure of a group extension of $\Gal(K'/K)$ with $V$. Then $\gamma$ is a function $\Gal(K^*/K) \to V$ such that $g \mapsto (g \Gal(K^*/K'), \gamma(g))$ is a group homomorphism. Without loss of generality, we may assume $\alpha(1,1)=0$. For any infinite place $v$ that ramifies in $K^*$, we have $K_v = \mathbb R$ and $K^* \otimes_K K_v$ is a product of copies of $\mathbb C$. We can compute Cech cohomology classes in $H^*(K_v,V)$ by pulling back from the cover $\operatorname{Spec} K^* \otimes_K K_v \to \operatorname{Spec} K_v$ to a single copy of $\mathbb C$.  We can compute the pullbacks of $\alpha$ and $\delta_v$ by first pulling back to $\operatorname{Spec} K' \otimes_K K_v \to \operatorname{Spec} K_v$, which is a product of copies of $K_v$ since $K'$ is split at $v$, and the pullback of $\alpha$ to one of these copies is zero since $\alpha(1,1)=0$, so the pullback of $\delta_v$ to one of these copies is zero as well. 

Thus, for each infinite place $v$, the pullback of $\epsilon$ to $ H^1( \Gal(K_v),V)$ agrees with the pullback of $-\gamma$ to $H^1(\Gal(K_v),V)$. Since $\Gal(K_v)$ acts trivially on $V$, $H^1(\Gal(K_v),V)$ is simply the set of homomorphisms from $\Gal(K_v)$ to $V$, so $\gamma+\epsilon$ is the trivial homomorphism $\Gal(K_v)\to V$.  Hence $\gamma(g)+\epsilon(g)=0$ for $g\in \Gal(K_v)$.  Thus the function $f\colon \Gal(K^*/K) \to \Gal(K'/K) \times V$ sending $g \mapsto (g \Gal(K^*/K'), \gamma(g)+\epsilon(g))$ is trivial on $\Gal(K_v)$ for each infinite place $v$. The cocycle condition for $\epsilon$ implies that $f$ is a group homomorphism. Since $f$ is trivial on the Galois group of each infinite place, $f$ is unramified. Because the image of $f$ has order prime to $|\Gamma|$, $f$ factors through $\Gal ( K^{ \operatorname{un}, \abs{\Gamma}'} / K)$ so $\alpha$ splits over $\Gal ( K^{ \operatorname{un}, \abs{\Gamma}'} / K)$ and thus $\alpha=0$, as desired.\end{proof}

\begin{lemma}\label{group-scheme-comparison} 
With notation as in the start of Section~\ref{S:NE},
assume for some positive integer $n$ that $V$ has characteristic dividing $n$ and $k$ contains the $n$th roots of unity. Let $V^\vee$ be the dual representation of $V$ and let $DV$ be the Cartier dual of the group scheme associated to $V$. Then there is a natural map $V^\vee \to DV$ of group schemes over $\mathcal O_K$ and the induced $\Gamma$-equivariant map \[ H^1 ( \mathcal O_K, V^\vee) \to H^1(\mathcal O_K, DV)\] is an injection. 
\end{lemma}

\begin{proof} We have $V^\vee \cong \Hom (V, \mathbb Z/n)$ and $DV \cong \Hom(V, \mu_n)$ but fixing a generator of the $n$th roots of unity in $k$ gives a map $\mathbb Z/n\to\mu_n$ and thus $V^\vee \to DV$. The fact that the fixed generator lies in $k$ and not just in $K$ insures that the induced map is $\Gamma$-equivariant. This map is an isomorphism on the open set deleting the primes dividing $n$. A class in $H^1 ( \mathcal O_K, V^\vee)$
in the kernel of the map $H^1 ( \mathcal O_K, V^\vee) \to H^1(\mathcal O_K, D_V)$
 represents a torsor that becomes a trivial torsor when composed with $V^\vee \to DV$, hence trivial upon restriction to the open set, hence trivial everywhere since $\mathcal O_K$ is integrally closed, showing the map is injective. \end{proof}

\begin{lemma}\label{H0-H3} 
With notation as in the start of Section~\ref{S:NE}, assume
 $V$ is a simple $\Gal ( K^{ \operatorname{un}, \abs{\Gamma}'} / k)$-group. Let $p$ be the characteristic of $V$. If $k$ contains the $p$th roots of unity then
\[  |H^3_c( \operatorname{Spec} \mathcal O_K, V) ^\Gamma|=  |H^0 ( \operatorname{Spec} \mathcal O_K, V)^\Gamma|.\]

If $K^{ \operatorname{un}, \abs{\Gamma}'} $ does not contain the $p$th roots of unity,
\[  |H^3_c( \operatorname{Spec} \mathcal O_K, V) ^\Gamma| =1\] 
\end{lemma}
We will later see in Theorem~\ref{rou-descent} that if $k$ does not contain the $p$th roots of unity and has no nontrivial unramified extensions of order prime to $\abs{\Gamma}$ then $K^{ \operatorname{un}, \abs{\Gamma}'} $ contains the $p$th roots of unity only if $K$ contains the $p$th roots of unity.

\begin{proof}The group  $H^0 ( \operatorname{Spec} \mathcal O_K, V)^\Gamma$ is equal to the $\Gal ( K^{ \operatorname{un}, \abs{\Gamma}'} / k)$-invariants of $V$. Since $V$ is a smooth group scheme, $H^3_c( \operatorname{Spec} \mathcal O_K, V) ^\Gamma=H^3_{c, \operatorname{fppf}} ( \operatorname{Spec} \mathcal O_K, V) ^\Gamma $. By Artin-Verdier duality for flat cohomology, $H^3_{c,\operatorname{fppf}}( \operatorname{Spec} \mathcal O_K, V) ^\Gamma$ is dual to $H^0_{\operatorname{fppf}}( \operatorname{Spec} \mathcal O_K, DV) ^\Gamma$ where $DV$ is the Cartier dual of $V$. Then $H^0_{\operatorname{fppf}}( \operatorname{Spec} \mathcal O_K, DV) ^\Gamma$ is isomorphic to the Galois invariants of $DV$, hence dual to  the Galois coinvariants of the Tate twist $V(-1)$ of $V$. 

If $k$ contains the $p$th roots of unity then $V \cong V(-1)$ as $\Gal ( K^{ \operatorname{un}, \abs{\Gamma}'} / k)$-representations and, because $V$ is simple, its invariants and coinvariants are equal, so these are equal.

For the second statement, if $H^3_c( \operatorname{Spec} \mathcal O_K, V) ^\Gamma$ were nontrivial, the Galois action on $V(-1)$ would be trivial, so $V$ would be isomorphic as a Galois module to $\mathbb F_p(1)$, showing the Galois group of $k(\mu_p)$ is a quotient of $ \Gal(K^{ \operatorname{un}, \abs{\Gamma}'} / k)$ and hence that the $p$th roots of unity lie in $K^{ \operatorname{un}, \abs{\Gamma}'}$. \end{proof}

We are now ready to prove our results giving constraints on the Galois cohomology of $\Gal ( K^{ \operatorname{un}, \abs{\Gamma}'} / K) \rtimes \Gamma$.

\begin{theorem}\label{ne-nd-prof} With notation as in the start of Section~\ref{S:NE},  let $V$ be a simple $\Gal ( K^{ \operatorname{un}, \abs{\Gamma}'} / k)$-group. Assume that $K^{ \operatorname{un}, \abs{\Gamma}'} $ does not contain the $p$th roots of unity. Then
\[ |H^2( \Gal ( K^{ \operatorname{un}, \abs{\Gamma}'} / k) , V)| \leq |H^1(\Gal ( K^{ \operatorname{un}, \abs{\Gamma}'} / k) , V)|    \frac{ \prod_{ v \textrm{ archimedean place of } k} | V ^{\Gal(k_v^{\operatorname{sep}}/k_v)}|}{ |V^{ \Gal ( K^{ \operatorname{un}, \abs{\Gamma}'} / k)}|} .\]\end{theorem}

\begin{proof} This follows immediately from combining Lemmas \ref{parabolic-Euler-characteristic}, \ref{H1-comparison-prof}, \ref{H2-comparison-prof}, and \ref{H0-H3}, together with the observation that $|H^0 ( \operatorname{Spec} \mathcal O_K, V)^\Gamma |= |V^{ \Gal ( K^{ \operatorname{un}, \abs{\Gamma}'} / k)}|$, to show. \end{proof}

Let $n$ be the number of roots of unity in $k$ of order prime to $\abs{\Gamma}$.
For $V$ a simple  $\Gal ( K^{ \operatorname{un}, \abs{\Gamma}'} / k)$-group of characteristic $p$ where $p\mid n$ , let $H^2 (  \Gal ( K^{ \operatorname{un}, \abs{\Gamma}'} / k) , V)^{s_K}$ be the subgroup of classes in $H^2 (  \Gal ( K^{ \operatorname{un}, \abs{\Gamma}'} / k) , V)$  whose cup product with any element of $H^1 (   \Gal ( K^{ \operatorname{un}, \abs{\Gamma}'} / k) , V^\vee)$ has zero Artin-Verdier trace, and, if $V \cong \mathbb F_p$, whose image under the Bockstein homomorphism $H^2 (   \Gal ( K^{ \operatorname{un}, \abs{\Gamma}'} / k) , \mathbb F_p) \to H^3 (   \Gal ( K^{ \operatorname{un}, \abs{\Gamma}'} / k) , \mathbb Z/n)$  has zero Artin-Verdier trace.

\begin{theorem}\label{ne-d-prof} With notation as in the start of Section~\ref{S:NE},  let $V$ be a simple $\Gal ( K^{ \operatorname{un}, \abs{\Gamma}'} / k)$-group of characteristic $p$. Assume $k$ contains the $p$th roots of unity and $k(\mu_{pn})$ is a ramified extension of $k$. Then
\[  |H^1 (  \Gal ( K^{ \operatorname{un}, \abs{\Gamma}'} / k) , V^\vee)| |H^2 ( \Gal ( K^{ \operatorname{un}, \abs{\Gamma}'} / k) , V)^{s_K}| \] \[\leq |H^1 (  \Gal ( K^{ \operatorname{un}, \abs{\Gamma}'} / k) , V)|   \frac{  \prod_{ v \textrm{ archimedean place of } k} | V ^{ \Gal(k_v^{\operatorname{sep}}/k_v)}| }{|V^{\Gal ( K^{ \operatorname{un}, \abs{\Gamma}'} / k)}|}  .\]\end{theorem}   

\begin{proof} It follows from Lemmas \ref{parabolic-Euler-characteristic}, \ref{H1-comparison-prof}, and \ref{H0-H3} that 
\[ |H^2_{c} (\mathcal O_K, V)| =  |H^1 ( \Gal ( K^{ \operatorname{un}, \abs{\Gamma}'} / k) , V)|    \prod_{ v \textrm{ archimedean place of } k} | V ^{ \Gal(k_v^{\operatorname{sep}}/k_v)}|\] so it suffices to show
\[  |H^1 ( \Gal ( K^{ \operatorname{un}, \abs{\Gamma}'} / k) , V^\vee)| |H^2 (\Gal ( K^{ \operatorname{un}, \abs{\Gamma}'} / k) , V)^{s_K}| \leq \frac{  |H^2_{c} (\mathcal O_K, V)|}{|V^{ \Gal ( K^{ \operatorname{un}, \abs{\Gamma}'} / k)}|}  .\]
Lemma \ref{H2-comparison-prof} gives an injective map
\begin{equation}\label{first-injective-map} H^2 ( \Gal ( K^{ \operatorname{un}, \abs{\Gamma}'} / k) , V)^{s_K} \to H^2_{c} (\mathcal O_K, V) \end{equation}
while Lemmas \ref{H1-comparison-prof}, \ref{group-scheme-comparison}, and Artin-Verdier duality give a composition of injective maps
\begin{equation}\label{second-injective-map} H^1 (  \Gal ( K^{ \operatorname{un}, \abs{\Gamma}'} / k) , V^\vee) \to H^1_{\textrm{par}} (\mathcal O_K, V^\vee)  \to H^1(\mathcal O_K, V^\vee) \to H^1(\mathcal O_K, DV) \to (H^2_{c} (\mathcal O_K, V))^\vee \end{equation}
and by definition of $H^2 (  \Gal ( K^{ \operatorname{un}, \abs{\Gamma}'} / k) , V)^{s_K}$ the images of \eqref{first-injective-map} and \eqref{second-injective-map} are orthogonal and hence the product of their cardinalities is at most the cardinality of $H^2_{c} (\mathcal O_K, V) $. This handles the case where $V$ is nontrivial and thus  $|V^{ \Gal ( K^{ \operatorname{un}, \abs{\Gamma}'} / k)}|=1$.

In the complementary case where $V$ is trivial, i.e. $V \cong \mathbb F_p$, it suffices to show that the images of \eqref{first-injective-map} and \eqref{second-injective-map} are not orthognal complements of each other as then the product of their cardinalities will be at most $|H^2_{c} (\mathcal O_K, V) |/p$, as desired. To do this, observe that the image of \eqref{first-injective-map} consists of classes in $H^2_{c} (\mathcal O_K, \mathbb F_p)$ whose image under the Bockstein homomorphism $H^2_{c} (\mathcal O_K, \mathbb F_p) \to H^2_c( \mathcal O_K, \mathbb Z/n)$ has zero Artin-Verdier trace. By \cref{exact-sequence-duality}, the trace of the image of a class $\alpha$ under the Bockstein morphism is the Artin-Verdier pairing of $\alpha$ with the image under the dual Bockstein morphism $H^0(\mathcal O_K, \mu_n) \to H^1(\mathcal O_K, \mu_p)$ of the class in $H^0(\mathcal O_K, \mu_n)$ associated to a generator in $\mu_n$. 

This gives a class in $H^1(\mathcal O_K, \mu_p)= H^1(\mathcal O_K, DV)$ and induces a linear form in $(H^2_{c} (\mathcal O_K, V))^\vee $. To show the images of \eqref{first-injective-map} and \eqref{second-injective-map} are not orthogonal complements it suffices to show that this linear form is not in the image of \eqref{second-injective-map}. Since the last map in the construction of \eqref{second-injective-map} is injective, it suffices to show that the Bockstein class in $H^1(\mathcal O_K, \mu_p)$ does not lie in the image of $H^1_{\textrm{par}}(\mathcal O_K, \mathbb F_p)$. If it was, it would vanish over some everywhere unramified extension $K'$ of $K$ and hence the map $H^0 (\mathcal O_{K'},\mu_{np}) \to H^0( \mathcal O_{K'},\mu_n) $ would be surjective. (Note if $p=2$ then ``everywhere unramified" includes the infinite place, since the torsors parameterized by parabolic $H^1$ are those split at the infinite place.) Since $K'$ contains the $n$th roots of unity and $p\mid n$ this implies that $K'$ contains the $pn$'th roots of unity. Thus $K( \mu_{pn})$ is unramified over $K$. Thus the order of the inertia group of $k(\mu_{pn})/k$ divides the degree of $K/k$ and thus divides $\abs{\Gamma}$, but on the other hand divides the degree of $\mathbb Q(\mu_{pn})/\mathbb Q(\mu_n)$ which is $p$, and $p$ is the cardinality of $V$ which is prime to $|\Gamma|$, so this is $1$ and $k(\mu_{pn})$ is an unramified extension of $k$, contradicting our assumption.\end{proof}

\begin{lemma}\label{H1-comparison-finite} 
Let $\cL$ be a level in the category of finite $\Gamma$-groups, $H$ the maximal quotient of $\Gal ( K^{ \operatorname{un}, \abs{\Gamma}'} / K) $ in $\mathcal L$,  and $s_H$ the induced orientation on $H$.
For $V$ a $H\rtimes \Gamma$-group, the natural map $H^1 ( H\rtimes \Gamma, V) \to  H^1( \Gal ( K^{ \operatorname{un}, \abs{\Gamma}'} / K) \rtimes \Gamma, V)$ is injective.  If $ V \rtimes H \in \mathcal L$, this map is an isomorphism. \end{lemma}

\begin{proof} Classes in $H^1(G, V)$ represent lifts of $G \to \Aut(V) \to G \to V \rtimes (H \rtimes \Gamma) $, so the induced map on $H^1$ is always injective for a surjective homomorphism of groups. If $V \rtimes H \in \mathcal L$, for any lift $\Gal ( K^{ \operatorname{un}, \abs{\Gamma}'} / K) \rtimes \Gamma  \to V \rtimes (H \rtimes \Gamma)$  the image of $\Gal ( K^{ \operatorname{un}, \abs{\Gamma}'} / K)$ is a  $\Gamma$-invariant subgroup of $V \rtimes H$ whose projection to $H $ is surjective, hence either $V \rtimes H$ or $H$, so in either case in $\mathcal L$, and thus must be a quotient of $H$, showing the map is surjective in this case. \end{proof}

\begin{lemma}\label{H2-comparison-finite} 
Let $\cL$ be a level in the category of finite $\Gamma$-groups, $H$ the maximal quotient of $\Gal ( K^{ \operatorname{un}, \abs{\Gamma}'} / K) $ in $\mathcal L$,  and $s_H$ the induced orientation on $H$.
For $V$ an admissible finite simple $\Gamma$-group, the map \[H^2(H \rtimes \Gamma, V)^\mathcal L \to H^2( \Gal ( K^{ \operatorname{un}, \abs{\Gamma}'} / K) \rtimes \Gamma, V)\] is injective. \end{lemma}

\begin{proof} Each class in $H^2(H \rtimes \Gamma, V)$ represents an extension of $H \rtimes \Gamma$ by $V$, equivalently, a $\Gamma$-equivariant extension of $H$ by $V$, and its pullback to $H^2( \Gal ( K^{ \operatorname{un}, \abs{\Gamma}'} / K) \rtimes \Gamma, V)$ is trivial if and only if the extension splits $\Gamma$-equivariantly over $\Gal ( K^{ \operatorname{un}, \abs{\Gamma}'} / K) $, equivalently, the map $\Gal ( K^{ \operatorname{un}, \abs{\Gamma}'} / K) \to H$ lifts $\Gamma$-equivariantly to the extension. If the extension class lies in $H^2(H \rtimes \Gamma, V)^\mathcal L$ then the associated extension is in $\mathcal L$, so any lift must factor through $H$, showing the extension splits already over $H$ and its class is trivial.\end{proof}

\begin{lemma}\label{ne-nd-finite} 
Let $\cL$ be a level in the category of finite $\Gamma$-groups, $H$ the maximal quotient of $\Gal ( K^{ \operatorname{un}, \abs{\Gamma}'} / K) $ in $\mathcal L$,  and $s_H$ the induced orientation on $H$.
Let $V$ be an admissible simple $H \rtimes \Gamma$-group of characteristic $p$. Assume that $K^{ \operatorname{un}, \abs{\Gamma}'} $ does not contain the $p$th roots of unity. Then
\[ |H^2(H\rtimes \Gamma , V)^\mathcal L | \leq |H^1(  H \rtimes \Gamma , V)|   \frac{   \prod_{ v \textrm{ archimedean place of } k} | V ^{ \Gal(k_v)}|}{ |V^{H\rtimes \Gamma}| }.\]\end{lemma}

\begin{proof} This follows immediately from combining Lemmas \ref{ne-nd-prof}, \ref{H1-comparison-finite}, and \ref{H2-comparison-finite}. \end{proof}

\begin{lemma}\label{ne-d-finite} 
Let $\cL$ be a level in the category of finite $\Gamma$-groups, $H$ the maximal quotient of $\Gal ( K^{ \operatorname{un}, \abs{\Gamma}'} / K) $ in $\mathcal L$,  and $s_H$ the induced orientation on $H$.
Let $V$ be an admissible simple $H \rtimes \Gamma$-group of characteristic $p$. Assume $k$ contains the $p$th roots of unity.  Then
\[  |H^1 ( H \rtimes \Gamma , V^\vee)| |H^2 ( H \rtimes \Gamma , V)^{\mathcal L, s_H }| \leq |H^1 ( H\rtimes \Gamma , V)|   \frac{  \prod_{ v \textrm{ archimedean place of } k} | V ^{ \Gal(k_v)}|} { |V^{H\rtimes \Gamma}| }.\]\end{lemma}   

\begin{proof} In view of Theorem \ref{ne-d-prof} and Lemma \ref{H1-comparison-finite} it suffices to prove
\[  |H^1 ( H \rtimes \Gamma , V^\vee)| |H^2 ( H \rtimes \Gamma , V)^{\mathcal L, s_H } |  \leq  |H^1 (  \Gal ( K^{ \operatorname{un}, \abs{\Gamma}'} / K) \rtimes \Gamma , V^\vee)| |H^2 (  \Gal ( K^{ \operatorname{un}, \abs{\Gamma}'} / K) \rtimes \Gamma , V)^{s_K}| .\]

We first handle the case where $V$ is nontrivial.

For $p$ the characteristic of $V$, the right-hand side can be viewed as $p$ raised to the power of $\dim H^1 (  \Gal ( K^{ \operatorname{un}, \abs{\Gamma}'} / K) \rtimes \Gamma , V^\vee) $ plus $\dim H^2 (  \Gal ( K^{ \operatorname{un}, \abs{\Gamma}'} / K) \rtimes \Gamma , V)$ minus the rank of the bilinear form $ H^1 (  \Gal ( K^{ \operatorname{un}, \abs{\Gamma}'} / K) \rtimes \Gamma , V^\vee) \times H^2 (  \Gal ( K^{ \operatorname{un}, \abs{\Gamma}'} / K) \rtimes \Gamma , V) \to \mathbb F_p$ arising from cup product and Artin-Verdier duality.

Similarly, the left-hand side is  at most $p$ raised to the power of $\dim H^1 ( H \rtimes \Gamma , V^\vee)$ plus $\dim H^2 ( H \rtimes \Gamma , V)^{\mathcal L}$ minus the rank of the pullback bilinear form $H^1 ( H \rtimes \Gamma , V^\vee) \to \dim H^2 ( H \rtimes \Gamma , V)^{\mathcal L} \to \mathbb F_p$.  By Lemma \ref{H1-comparison-finite} this is the pullback along injective maps $H^1 ( H\rtimes \Gamma, V) \to  H^1( \Gal ( K^{ \operatorname{un}, \abs{\Gamma}'} / K) \rtimes \Gamma, V)$ and $H^2(H \rtimes \Gamma, V)^\mathcal L \to H^2( \Gal ( K^{ \operatorname{un}, \abs{\Gamma}'} / K) \rtimes \Gamma, V)$ of the previous bilinear form.

It is straightforward to check by linear algebra that the sum of the dimensions of two vector spaces minus the rank of a bilinear form on them can only decrease when we pull back the bilinear form along injective maps.

In the case when $V$ is nontrivial, the same claims are true except replacing $ H^2 (  \Gal ( K^{ \operatorname{un}, \abs{\Gamma}'} / K) \rtimes \Gamma , V)$ by the subspace of classes whose image under the Bockstein homomorphism has trivial Artin-Verdier trace, and doing the same for $H^2 ( H \rtimes \Gamma , V)^{\mathcal L}$. The same argument applies.\end{proof}

\begin{theorem}\label{rou-descent} Let $k$ be a number field, $\Gamma$ a finite group, and $p$ a prime. Assume that $k$ lacks nontrivial unramified extensions of degree prime to $\abs{\Gamma}$. Then  $K^{ \operatorname{un}, \abs{\Gamma}'} $ contains the $p$th roots of unity if and only if $K$ contains the $p$th roots of unity.
\end{theorem}

\begin{proof}The ``if'' direction is clear so we establish ``only if'' and thus assume that $K^{ \operatorname{un}, \abs{\Gamma}'} $ contains the $p$th roots of unity.

Let $\ell$ be a prime dividing the order of $\Gal(k(\mu_p)/k)$.  We show that $\ell\mid \abs{\Gamma}$. Since $\Gal(k(\mu_p)/k)$ is abelian, it has some quotient of order $\ell$, corresponding to a degree $\ell$ extension of $k$. If this extension is unramified then by assumption $\ell \mid \abs{\Gamma}$ so we may assume this extension is ramified at some place $v$ and thus the inertia group of $k(\mu_p)/k$ at $v$ has order a multiple of $\ell$.

If $K^{ \operatorname{un}, \abs{\Gamma}'} $ contains $k(\mu_p)$ then the inertia group of $K^{ \operatorname{un}, \abs{\Gamma}'} /k$ at $v$ surjects onto the inertia group of $k(\mu_p)/k$ at $v$. Since $K^{ \operatorname{un}, \abs{\Gamma}'} $ is an everywhere unramified extension of $K$, the inertia group of $K^{ \operatorname{un}, \abs{\Gamma}'} /k$ at $v$ is isomorphic to the inertia group of $K/k$ at $v$. Hence $\ell$ divides the order of the inertia group of $k(\mu_p)/k$ at $v$ which divides the order of the inertia group of $K/k$ at $v$ which divides $\abs{\Gamma}$, and in particular $\ell\mid \abs{\Gamma}$, as desired.

 If $K^{ \operatorname{un}, \abs{\Gamma}'} $ contains $k(\mu_p)$ then the extension $ K(\mu_p)/K$ must have degree relatively prime to $\abs{\Gamma}$ but on the other hand its degree divides the order of $\Gal(k(\mu_p) /k)$. Since each prime divding the order of $\Gal(k(\mu_p) /k)$ divides $\abs{\Gamma}$, it follows that $K(\mu_p)/K$ has degree $1$ and thus $\mu_p\in K$, as desired.\end{proof}

 We are now ready to verify the criteria for nonzero probability in Proposition \ref{nzp-criterion} holds for number fields under a mild assumption on roots of unity.

 \begin{proposition}\label{nzp-converse} Let $k$ be a number field and $\Gamma$ a finite group. Let $n$ be the number of roots of unity in $k$ of order prime to $\abs{\Gamma}$. Assume that $k$ lacks nontrivial unramified extensions of degree prime to $\abs{\Gamma}$.
 
 Let $K$ be an extension of $k$ with Galois group $\Gamma$. Let $U$ be a multiset of elements of $\Gamma$, consisting of, for each infinite place $v$ of $k$, a generator of $\Gal(k_v)$.  Assume that $K$ does not contain any roots of unity that are not contained in $k$.
 
 Let $\cL$ be a level of the category of finite $\Gamma$-groups.  Let $\cW$ be the set of isomorphism classes of finite $n$-oriented $\Gamma$-groups whose underlying $\Gamma$-group is in $\cL$. Let $H$ be the maximal quotient of $\Gal ( K^{ \operatorname{un}, \abs{\Gamma}'} / K) $ in $\mathcal L$,  and $s_H$ the orientation on $H$.
 
 Then for each finite simple abelian $H \rtimes \Gamma$-group $V_i$ of characteristic prime to $n$ that can be the kernel of a simple morphism $\bG \to \bH$ with $\bG \in \cW$ we have
 \begin{equation}\label{nzp-converse-nd} | H^2( H \rtimes \Gamma, V_i)^{\cL} | \leq  |H^1( H \rtimes \Gamma, V_i )| V_i^{\cdot U'}\end{equation} and for each finite simple abelian $H \rtimes \Gamma$-group $V_i$ of characteristic dividing $n$ that can be the kernel of a simple morphism $\bG \to \bH$ with $\bG \in \cW$ we have \begin{equation}\label{nzp-converse-d} |H^1 ( H \rtimes \Gamma, V_i^\vee)|  |H^2(H\rtimes \Gamma,V_i)^{\cL,s_H}|  \leq  |H^1 ( H \rtimes \Gamma, V_i)|  V_i^{\cdot U'}   .\end{equation}
  \end{proposition}
  The assumption that $k$ lacks any nontrivial everywhere unramified field extension of degree prime to $\abs{\Gamma}$ is (weaker than) one of the assumptions in the main conjecture Conjecture \ref{intro-moments-conjecture}, so it is not surprising that we need it to verify a prediction of that conjecture.

The assumption that $K$ does not contain any roots of unity that are not contained in $k$ is more subtle, but should always hold outside a set of density zero, under the methods of counting allowed in Conjecture~\ref{intro-moments-conjecture}.
 The field extensions with Galois group $\Gamma$, containing additional roots of unity all contain one of finitely many nontrivial subextensions, that being the cyclotomic extensions of degree at most $\abs{\Gamma}$, so indeed should have density zero under these counting methods.

\begin{proof} Take $V = V_i$ and let $p$ be the characteristic of $V_i$. Since every element of $\cW$ has order prime to $\abs{\Gamma}$, the assumption that $V_i$ can be the kernel of a simple morphism implies that $p$ is prime to $\abs{\Gamma}$. 

In the case that the characteristic $p$ of $V_i$ is prime to $n$, it follows by definition of $n$ that $k$ does not contain the $p$th roots of unity. By assumption, $K$ also does not contain the $p$th roots of unity. Theorem \ref{rou-descent}, $K^{ \operatorname{un}, \abs{\Gamma}'} $ also does not contain the $p$th roots of unity.  Thus we can apply Lemma \ref{ne-nd-finite}, which gives
\[ |H^2(H\rtimes \Gamma , V_i)^\mathcal L | \leq |H^1(  H \rtimes \Gamma , V_i)|   \frac{   \prod_{ v \textrm{ archimedean place of } k} | V_i ^{ \Gal(k_v)}|}{ |V_i^{H\rtimes \Gamma}| }.\]
We have
 \begin{equation}\label{uprime-matching}V_i^{\cdot U'}=\frac{\prod_{ v \textrm{ archimedean place of } k} | V _i^{ \Gal(k_v)}|}{|V_i^{H\rtimes \Gamma}|} \end{equation} by definition of $U$ and $V_i^{\cdot U'}$, and thus we obtain \eqref{nzp-converse-nd}.
 
 In the case that the characteristic $p$ of $V_i$ divides $n$, by definition of $n$, $k$ contains the $p$th roots of unity. Hence we can apply Lemma \ref{ne-d-finite} which gives
 \[  |H^1 ( H \rtimes \Gamma , V_i^\vee)| |H^2 ( H \rtimes \Gamma , V_i)^{\mathcal L, s_H }| \leq |H^1 ( H\rtimes \Gamma , V_i)|   \frac{  \prod_{ v \textrm{ archimedean place of } k} | V_i ^{ \Gal(k_v)}|} { |V_i^{H\rtimes \Gamma}| }\]  which by \eqref{uprime-matching} implies \eqref{nzp-converse-d}.\end{proof}

\bibliographystyle{alpha}
\bibliography{../MyLibrary.bib,../references.bib}

\end{document}